\def\arxivcsm{}

\documentclass[10pt,twoside]{article}

\usepackage{ifthen}
\usepackage{pgffor}
\usepackage{etoolbox}
\usepackage{parselines}
\usepackage{environ}
\usepackage{keyval}
\usepackage{listofitems}
\usepackage{forloop}
\usepackage{calc}
\usepackage{xstring}
\usepackage{changepage}
\usepackage{scrextend}
\usepackage{expl3}
\usepackage{stackengine}

\makeatletter
\ifthenelse{\isundefined{\lefteqs}}{
\usepackage{amsmath}
\usepackage{amsfonts}
\usepackage{mathtools}
}{
\usepackage[fleqn]{amsmath}
\usepackage[fleqn]{mathtools}
}
\usepackage{amssymb}
\@ifclassloaded{autart}{}{\usepackage{amsthm}}
\makeatother
\usepackage{thmtools}
\usepackage{ar}
\usepackage{empheq}
\usepackage{mathrsfs}
\usepackage{bm}
\usepackage{sansmath}

\makeatletter
\let\ltxxlabel\ltx@label
\makeatother

\usepackage{color}
\makeatletter
\@ifclassloaded{beamer}{}{
\usepackage[table]{xcolor}}
\makeatother

\usepackage{graphicx}
\makeatletter
\@ifclassloaded{uwthesis/uwthesis}{
  \usepackage{subfig}
}{
  \@ifclassloaded{IEEEtran}{
  }{
    \usepackage{subcaption}
  }
}
\@ifclassloaded{IEEEtran}{}{\usepackage{caption}}
\makeatother
\usepackage{dblfloatfix}

\usepackage{tabularx}
\usepackage{adjustbox}
\usepackage{multicol}
\usepackage{diagbox} 
\usepackage{varwidth} 
\usepackage{makecell} 
\usepackage{arydshln}

\usepackage{tikz}
\usepackage{tikz-3dplot}
\usetikzlibrary{backgrounds}
\usetikzlibrary{arrows}
\usetikzlibrary{arrows.meta}
\usetikzlibrary{intersections}
\usetikzlibrary{math}
\usetikzlibrary{fpu}
\usetikzlibrary{shapes}
\usetikzlibrary{positioning}
\usetikzlibrary{decorations.markings}
\usetikzlibrary{decorations.pathreplacing}
\usetikzlibrary{decorations.text}
\usetikzlibrary{tikzmark}
\usetikzlibrary{fadings}

\makeatletter
\@ifclassloaded{standalone}{}{\usepackage[tikz]{mdframed}}
\@ifclassloaded{beamer}{}{
  \@ifclassloaded{uwthesis/uwthesis}{}{
    \@ifclassloaded{IEEEtran}{}{
      \usepackage{tocloft}
    }}}
\usepackage{enumitem}
\makeatother
\usepackage{listings}
\usepackage{lettrine}
\usepackage{algorithm, algorithmicx}
\ifthenelse{\isundefined{\algouseends}}{
\usepackage[noend]{algpseudocode}}{
\usepackage{algpseudocode}}
\usepackage[many,magazine]{tcolorbox}
\usepackage[binary-units=true]{siunitx}
\usepackage{chemformula}
\usepackage[normalem]{ulem}
\usepackage{cancel}
\usepackage{accents}
\makeatletter
\@ifclassloaded{beamer}{\usepackage{multimedia}}{}
\makeatother
\usepackage{lipsum}
\usepackage{mfirstuc}
\usepackage{import}
\usepackage{nicematrix}

\makeatletter
\ifthenelse{\isundefined{\nobiblatex}}{
  \usepackage{bibentry} 
  \newcommand{\usebiblatex}[1]{%
    \usepackage[backend=biber,
    style=#1,
    sorting=none,
    url=false,
    sortcites=true,
    defernumbers=true]{biblatex}}
  \@ifclassloaded{beamer}{
    \usebiblatex{authoryear}}{
    \ifthenelse{\isundefined{\preprintcsm}}{%
      \usebiblatex{numeric}}{%
      \usebiblatex{ieee}}}
}{}
\makeatother

\usepackage[outline]{contour}
\usepackage{extdash}
\usepackage{textcomp}

\usepackage{float}
\usepackage{afterpage}
\ifthenelse{\isundefined{\landscapemode}}{
\def\landscapemode{pdflscape}
}{}
\usepackage{\landscapemode} 
\usepackage{varwidth}
\usepackage{nomencl}
\usepackage{longtable}
\usepackage{marginnote}
\makeatletter
\@ifclassloaded{beamer}{}{
\@ifclassloaded{autart}{}{
\@ifclassloaded{standalone}{}{
\@ifclassloaded{ieeeconf}{}{
\usepackage[explicit]{titlesec}
}}}}
\makeatother
\usepackage{soul}
\usepackage{stackengine}
\usepackage{scalerel}
\usepackage{lastpage}
\usepackage{fancyhdr}
\usepackage{nameref}

\makeatletter

\ifthenelse{\isundefined{\preprintcsm}}{}{
  \gdef\nohref{}}
\@ifclassloaded{beamer}{%
  \gdef\nohref{}}{}

\ifthenelse{\isundefined{\nohref}}{
  \usepackage[hypertexnames=false]{hyperref}}{}

\makeatother

\usepackage{xr}

\usepackage{xspace}
\usepackage{setspace}

\newcommand{\myappend}[1]{\ifthenelse{\equal{#1}{}}{}{,#1}}
\newcommand{\myprepend}[1]{\ifthenelse{\equal{#1}{}}{}{#1,}}

\ExplSyntaxOn
\makeatletter
\NewDocumentCommand{\parengen}{m m m m}
{
  \group_begin:
  \tl_set:Nn \l_left_paren_tl {#1}
  \tl_set:Nn \l_right_paren_tl {#2}
  \tl_set:Nn \l_paren_sz_tl {#3}
  \ifthenelse{\equal{\l_paren_sz_tl}{}}{
    \l_left_paren_tl #4 \l_right_paren_tl
  }{
    \ifthenelse{\equal{\l_paren_sz_tl}{auto}}{
      \left\l_left_paren_tl #4 \right\l_right_paren_tl
    }{
      \csname \l_paren_sz_tl\endcsname\l_left_paren_tl
      #4
      \csname \l_paren_sz_tl\endcsname\l_right_paren_tl
    }
  }
  \group_end:
}

\newcommand{\pare}[2][]{\parengen{(}{)}{#1}{#2}}
\newcommand{\brak}[2][]{\parengen{[}{]}{#1}{#2}}
\newcommand{\brac}[2][]{\parengen{\{}{\}}{#1}{#2}}

\newcommand{\abso}[2][]{\parengen{|}{|}{#1}{#2}}

\newcommand{\@@normgen}[2][]{\parengen{\|}{\|}{#1}{#2}}

\makeatother
\ExplSyntaxOff

\newcommand*{\bigs}[2]{\vcenter{\hbox{\scalebox{#1}{\ensuremath#2}}}}
\let\oldbig\big
\let\oldbigg\bigg
\let\oldBig\Big
\let\oldBigg\Bigg
\renewcommand{\big}{\oldbig}
\renewcommand{\bigg}{\bigs{1.3}}
\newcommand{\biggg}{\oldBig}
\renewcommand{\Big}{\oldbigg}
\renewcommand{\Bigg}{\oldBigg}

\makeatletter
\newcommand{\fade}[2][1]{%
  \edef\fade@amount{\the\numexpr 100-#1*30 \relax}%
  {\color{black!\fade@amount!white}#2}}
\makeatother

\newcommand{\ifelsemath}[3]{%
  \pgfmathparse{#1 ? 1:0}%
  \ifthenelse{\pgfmathresult>0}{#2}{#3}}

\newcommand{\where}{:}

\newcommand{\inv}[1][]{\ifthenelse{\equal{#1}{}}{^{-1}}{^{-#1}}}

\ExplSyntaxOn
\makeatletter
\NewDocumentCommand{\norm}{o m}
{
  \group_begin:
  \IfNoValueTF{#1}{
    \def\@nargin{0}
  }{
    \seq_set_split:Nnn \l_var_seq { , } { #1 }
    \edef\@nargin{\seq_count:N \l_var_seq}
  }
  \seq_pop_left:NN \l_var_seq \l_norm_type_tl
  \seq_pop_left:NN \l_var_seq \l_delimiter_size_tl
  \pgfmathparse{\@nargin<=1 ? 1:0}
  \ifthenelse{\pgfmathresult>0}{
    \tl_set:Nn \l_delimiter_size_tl {}
  }{}
  \pgfmathparse{\@nargin==0 ? 1:0}
  \ifthenelse{\pgfmathresult>0}{
    \tl_set:Nn \l_norm_type_tl {}
  }{}
  \ifthenelse{\equal{\l_delimiter_size_tl}{auto}}{%
    \@@normgen[auto]{#2}}{%
    \ifthenelse{\equal{\l_delimiter_size_tl}{}}{
      \@@normgen[]{#2}}{\@@normgen[\l_delimiter_size_tl]{#2}}}%
  \c_math_subscript_token{\l_norm_type_tl}
  \group_end:
}
\seq_new:N \l_var_seq
\tl_new:N \l_norm_type_tl
\tl_new:N \l_delimiter_size_tl
\makeatother
\ExplSyntaxOff

\newcommand{\vertiii}[1]{{\left\vert\kern-0.25ex\left\vert\kern-0.25ex\left\vert #1
        \right\vert\kern-0.25ex\right\vert\kern-0.25ex\right\vert}}

\newcommand{\seq}[1]{\{#1\}}

\newlength{\dhatheight}

\makeatletter
\define@key{projKeys}{set}{\def\set{#1}}
\define@key{proxKeys}{t}{\def\t{#1}}
\define@key{proxKeys}{f}{\def\f{#1}}
\makeatother

\DeclareMathOperator{\proximal}{prox}

\newcommand{\Prox}[2][]{%
  \setkeys{proxKeys}{t=t,#1}%
  \setkeys{proxKeys}{f=f,#1}%
  \mathchoice{\underset{\t\f}{\proximal}}%
  {\proximal_{\t\f}}{\proximal_{\t\f}}{\proximal_{\t\f}}%
  \ifthenelse{\equal{#2}{}}{}{\left(#2\right)}%
}
\DeclareMathOperator{\distance}{dist}
\newcommand{\dist}[2][]{%
  \distance_{#1}\ifthenelse{\equal{#2}{}}{}{(#2)}
}
\newcommand{\gauge}[2][]{%
  \gamma_{#1}\ifthenelse{\equal{#2}{}}{}{(#2)}
}

\newcommand{\transp}{{\scriptscriptstyle\mathsf{T}}}
\newcommand{\T}{^\transp}

\newcommand{\dd}[1][]{#1\mathrm{d}}
\newcommand{\dt}[1][]{\dd[#1]t}
\newcommand{\sdd}{\dd[\,]}
\newcommand{\sdt}{\dt[\,]}

\newcommand{\fun}[2][1]{%
  #2(%
  \foreach \index in {1, ..., #1} {%
    \ifthenelse{\equal{\index}{#1}}{%
      \cdot%
    }{%
      \cdot,%
    }%
  })}
\newcommand{\definedas}[1][tri]{%
  \ifthenelse{\equal{#1}{tri}}{\triangleq}{\coloneqq}}

\renewcommand{\implies}{\Rightarrow}
\renewcommand{\iff}{\Leftrightarrow}

\DeclareMathOperator*{\exptx}{exp}
\renewcommand{\exp}[2][exponent]{\ifthenelse{\equal{#1}{exponent}}{e^{#2}}{\exptx\left(#2\right)}}
\DeclareMathOperator*{\expsf}{\mathsf{exp}}
\newcommand{\expm}[2][exponent]{\ifthenelse{\equal{#1}{exponent}}{\mathsf{e}^{#2}}{\expsf\left(#2\right)}}

\makeatletter
\DeclareFontFamily{U}{tipa}{}
\DeclareFontShape{U}{tipa}{m}{n}{<->tipa10}{}
\newcommand{\arc@char}{{\usefont{U}{tipa}{m}{n}\symbol{62}}}%
\renewcommand{\arc}[1]{\mathpalette\arc@arc{#1}}
\newcommand{\arc@arc}[2]{%
  \sbox0{$\m@th#1#2$}%
  \vbox{
    \hbox{\resizebox{\wd0}{\height}{\arc@char}}
    \nointerlineskip
    \box0
  }%
}
\makeatother



\newcommand{\dash}{\Hyphdash}

\makeatletter
\newcommand{\pushright}[1]{%
  \ifmeasuring@#1\else\omit\hfill$\displaystyle#1$\fi\ignorespaces}
\newcommand{\pushleft}[1]{%
  \ifmeasuring@#1\else\omit$\displaystyle#1$\hfill\fi\ignorespaces}
\makeatother

\newcommand{\Union}{\bigcup}

\newcommand{\set}[1]{\mathcal #1}

\newcommand{\eucledian}[1][]{\ifthenelse{\equal{#1}{}}{\mathbb E}{\mathbb #1}}
\DeclareMathOperator{\@convhull}{conv}
\newcommand{\reals}[1][]{%
  \ifthenelse{\equal{#1}{extended}}{\overline}{}{\mathbb R}%
}
\newcommand{\real}{\mathbb R}
\newcommand{\pos}[1][\reals]{#1_{++}}
\newcommand{\nonneg}[1][\reals]{#1_{+}}

\newcommand{\symm}{\mathbb S}  

\newcommand{\pd}{\symm_{++}}
\DeclareMathOperator{\dom}{dom}
\DeclareMathOperator{\epi}{epi}
\DeclareMathOperator{\interior}{int}

\DeclareMathOperator{\boundary}{\partial}

\newcommand{\algvar}[1]{\text{\IfSubStr{#1}{_}{%
    \StrSubstitute{#1}{_}{\textunderscore}}{#1}}}

\makeatletter
\definecolor{listinggray}{gray}{0.9}
\definecolor{lbcolor}{rgb}{0.9,0.9,0.9}
\colorlet{Darkgreen}{green!60!black}
\lstdefinelanguage{Julia}%
{morekeywords={abstract,break,case,catch,const,continue,do,else,elseif,%
    end,export,false,for,function,immutable,import,importall,if,in,%
    macro,module,otherwise,quote,return,switch,true,try,type,typealias,%
    using,while},%
  sensitive=true,%
  alsoother={\$},
  morecomment=[l]\#,%
  morecomment=[n]{\#=}{=\#},%
  morestring=[s]{"}{"},%
  morestring=[m]{'}{'},%
}[keywords,comments,strings]
\lstdefinestyle{listing@C++}{
  language=[GNU]C++,
}
\lstdefinestyle{listing@Julia}{%
  language = Julia,
}
\lstdefinestyle{listing@Python}{%
  language = Python,
}
\newcommand\langname@bash{}
\def\langname@bash{bash}
\newcommand\prompt@bash{\texttt{\$}\ } 
\newcommand\addedToEveryPar@bash{}
\lst@AddToHook{EveryPar}{\addedToEveryPar@bash}
\lst@AddToHook{PreInit}{%
  \ifx\lst@language\langname@bash%
  \let\addedToEveryPar@bash\prompt@bash%
  \fi
}
\lstdefinestyle{listingterminal}{%
  language=bash,
  backgroundcolor=\color{black!90!blue},
  basicstyle=\color{white},
  commentstyle=\color{yellow!50},
  keywordstyle=\bfseries\color{orange},
  breaklines,
  escapechar=@
}

\lstnewenvironment{terminalcode}[1][]{\lstset{style=listingterminal,#1}}{}
\makeatother

\newcommand{\grad}{\nabla}

\newcommand{\diff}[2]{\nabla_{#1}#2}

\newcommand{\subdiff}[1][none]{%
  \ifthenelse{\equal{#1}{none}}{%
    \partial%
  }{%
    \partial_{#1}%
  }%
}

\makeatletter
\DeclareMathOperator{\@convhull}{conv}
\newcommand{\convhull}[1][]{%
  \ifthenelse{\equal{#1}{closed}}{\overline}{}\@convhull%
}
\DeclareMathOperator{\@closed@convex@envelope}{\overline{co}}
\newcommand{\cce}{\@closed@convex@envelope}
\makeatother

\newcommand{\support}[2][delta]{\ifthenelse{\equal{#1}{delta}}{%
    \delta^*_{#2}%
  }{\sigma_{#2}}}

\DeclareMathOperator{\nul}{null}
\DeclareMathOperator{\rank}{rank}

\DeclareMathOperator{\diag}{diag}

\DeclareMathOperator{\tr}{\mathbf{tr}}

\newcommand{\Matrix}[2][]{%
  \ifthenelse{\equal{#1}{}}{}{\setlength\arraycolsep{#1}}%
  \begin{bmatrix}#2\end{bmatrix}}
\renewcommand{\Array}[2][]{%
  \ifthenelse{\equal{#1}{}}{}{\setlength\arraycolsep{#1}}%
  \begin{matrix}#2\end{matrix}}

\makeatletter
\newcommand{\toset}{%
  \def\arr@offset{0.15em}
  \def\arr@len{0.7em}
  \def\arr@height{0.3em}
  \tikz[minimum height=0ex,outer sep=0,inner sep=0]
  \path[-{Latex[length=0.8mm]}]
  node (a) at (0,0) {}
  node (b) at (0,\arr@height) {}
  (a) edge ++(\arr@len,0)
  (b) edge ++(\arr@len,0)
  (a) edge[draw=none] ++(0,-\arr@offset);%
}
\makeatother

\renewcommand{\skew}[1][]{^{\myprepend{#1}\times}}

\DeclareMathOperator*{\Expquat}{\mathsf{Exp}}
\DeclareMathOperator*{\Logquat}{\mathsf{Log}}
\newcommand{\qconj}{^{*}}

\newcommand{\ev}[3][c]{%
  \ifthenelse{\equal{#1}{c}}{%
    \ifthenelse{\equal{#2}{}}{\forall[0,#3]}{\forall[#2,#3]}%
  }{%
    \ifthenelse{\equal{#1}{o}}{%
      \ifthenelse{\equal{#2}{}}{\forall(0,#3)}{\forall(#2,#3)}%
    }{%
      \ifthenelse{\equal{#1}{oc}}{%
        \ifthenelse{\equal{#2}{}}{\forall(0,#3]}{\forall(#2,#3]}%
      }{%
        \ifthenelse{\equal{#2}{}}{\forall[0,#3)}{\forall[#2,#3)}%
      }%
    }%
  }%
}

\newcommand{\optimal}[1]{#1^*}
\newcommand{\lse}[1][]{\mathsf{L}_{#1}}
\newcommand*\softmax{\lse}

\newcommand{\GetLabel}[1]{\expandafter\csname #1Label\endcsname}
\makeatletter
\define@key{printRefKeys}{otherLabel}{\def\otherLabel{#1}}
\define@key{printRefKeys}{concatenate}{\def\concatenate{#1}}
\makeatother
\newcommand{\PrintRefs}[3][]{%
  \setkeys{printRefKeys}{otherLabel=,#1}%
  \setkeys{printRefKeys}{concatenate=false,#1}%
  \xdef\MyModifiedLabel{#2}%
  \ifthenelse{\equal{\otherLabel}{}}{}{\xdef\MyModifiedLabel{\otherLabel}}%
  \xdef\MyRefCount{0}%
  \foreach \i in {#3} {%
    \tikzmath{\MyRefCount=int(\MyRefCount+1);}%
    \xdef\MyRefCount{\MyRefCount}%
  }%
  \edef\OrigRefCount{\MyRefCount}%
  \ifthenelse{\equal{\MyRefCount}{1}}{%
    #2~\ref{\GetLabel{\MyModifiedLabel}:#3}%
  }{%
    \xdef\MyCounter{0}%
    \ifthenelse{\equal{#2}{Corollary}}{%
      Corollaries%
    }{%
      #2s%
    }~%
    \foreach \AlgRef in {#3} {%
      \tikzmath{\MyCounter=int(\MyCounter+1);}%
      \xdef\MyCounter{\MyCounter}%
      \ifthenelse{\equal{\concatenate}{true}}{%
        \pgfmathparse{\MyCounter==1 ? 1 : 0}%
        \ifthenelse{\pgfmathresult>0}{%
          \ref{\GetLabel{\MyModifiedLabel}:\AlgRef}-%
        }{%
          \pgfmathparse{\MyCounter==\MyRefCount ? 1 : 0}%
          \ifthenelse{\pgfmathresult>0}{%
            \ref{\GetLabel{\MyModifiedLabel}:\AlgRef}%
          }{}%
        }%
      }{%
        \pgfmathparse{\MyCounter<\MyRefCount-1 ? 1 : 0}%
        \ifthenelse{\pgfmathresult>0}{%
          \ref{\GetLabel{\MyModifiedLabel}:\AlgRef}, %
        }{%
          \pgfmathparse{\MyCounter<\MyRefCount ? 1 : 0}%
          \ifthenelse{\pgfmathresult>0}{%
            \ref{\GetLabel{\MyModifiedLabel}:\AlgRef}%
            \pgfmathparse{\OrigRefCount==2 ? 1:0}%
            \ifthenelse{\pgfmathresult>0}{ and }{, and }%
          }{%
            \ref{\GetLabel{\MyModifiedLabel}:\AlgRef}%
          }%
        }%
      }%
    }%
  }%
}

\newcommand{\cref}[1]{\PrintRefs{Chapter}{#1}}

\newcommand{\figref}[1]{\PrintRefs{Figure}{#1}}
\newcommand{\tabref}[2][]{\PrintRefs[#1]{Table}{#2}}

\makeatletter
\define@key{algKeys}{start}{\def\startline{#1}}
\define@key{algKeys}{end}{\def\endline{#1}}
\define@key{algKeys}{show}{\def\showalg{#1}}
\makeatother
\renewcommand{\algref}[2][]{%
  \setkeys{algKeys}{start=,#1}%
  \setkeys{algKeys}{end=,#1}%
  \setkeys{algKeys}{show=true,#1}%
  \ifthenelse{\equal{\startline}{}}{}{line%
    \ifthenelse{\equal{\endline}{}}{}{s}}%
  \ifthenelse{\equal{\startline}{}}{}{~\ref{alg:#2:line:\startline}%
    \ifthenelse{\equal{\endline}{}}{}{-\ref{alg:#2:line:\endline}} %
    of }%
  \ifthenelse{\equal{\showalg}{true}}{\PrintRefs{Algorithm}{#2}}{}%
}
\newcommand{\pref}[2][]{\PrintRefs[#1]{Problem}{#2}}

\newcommand{\tref}[2][]{\PrintRefs[#1]{Theorem}{#2}}
\newcommand{\corref}[2][]{\PrintRefs[#1]{Corollary}{#2}}
\newcommand{\conref}[2][]{\PrintRefs[#1]{Condition}{#2}}
\newcommand{\aref}[2][]{\PrintRefs[#1]{Assumption}{#2}}

\newcommand{\dref}[2][]{\PrintRefs[#1]{Definition}{#2}}

\newcommand{\figlabel}[1]{\label{fig:#1}}
\newcommand{\tablabel}[1]{\label{table:#1}}

\makeatletter
\newcommand*{\addFileDependency}[1]{
  \typeout{(#1)}
  \@addtofilelist{#1}
  \IfFileExists{#1}{}{\typeout{No file #1.}}
}
\makeatother

\DeclareSIUnit{\radian}{rad}
\DeclareSIUnit\century{century}
\DeclareSIUnit\year{yr}
\DeclareSIUnit{\ton}{t}

\makeatletter

\newcommand{\add@list@item}[2]{
  \ifthenelse{\equal{#1}{}}{\xdef#1{{#2}}}{\xdef#1{#1,{#2}}}
}

\def\providecounter#1{%
  \@ifundefined{c@#1}%
  {\newcounter{#1}}{\setcounter{#1}{0}}}

\newcommand{\makelist}[2]{
  \xdef\@current@list@name{#1}
  \xdef\@current@list@counter{#1@listcounter}
  \providecounter{\@current@list@counter}
  \edef\@tmp@list{{#2}}
  \expandafter\forcsvlist\expandafter\list@saveitem\@tmp@list
}

\newcommand{\list@saveitem}[1]{%
  \stepcounter{\@current@list@counter}%
  \expandafter\def\csname\@current@list@name%
  \arabic{\@current@list@counter}\endcsname{#1}
}

\newcommand{\list@nth}[2]{\csname #1#2\endcsname}

\makeatother

\newcommand{\mrm}[2][m]{\ifthenelse{\equal{#1}{m}}{\mathrm{#2}}
  {\textnormal{#2}}}

\newcommand{\Behcet}{Beh\c{c}et}
\newcommand{\Acikmese}{A\c{c}{\i}kme\c{s}e}

\makeatletter
\@ifclassloaded{beamer}{}{\newcommand{\alert}[1]{\textbf{#1}}}
\makeatother
\newcommand{\lcvx}{LCvx\xspace}
\newcommand{\emphasize}[1]{\textit{#1}}
\newcommand{\scvx}{SCvx\xspace}

\newcommand{\defintext}[1]{\textbf{#1}}

\makeatletter
\newcommand{\linkdest}[1]{\Hy@raisedlink{\hypertarget{#1}{}}}
\makeatother

\newcommand{\makeflag}[2]{%
  \expandafter\newif\csname ifmake#1\endcsname
  \csname make#1#2\endcsname
}

\NewEnviron{makeif}[2][]{%
  \csname ifmake#2\endcsname
  \BODY
  \else
  #1
  \fi
}

\ifthenelse{\isundefined{\noindentadjust}}{
\setlength\parindent{0pt} 
\setlength{\parskip}{4pt} 
}{}
\allowdisplaybreaks       

\makeatletter
\ifthenelse{\isundefined{\nohref}}{
\@ifclassloaded{beamer}{}{
\ifthenelse{\isundefined{\preprintcsm}}{
  \usepackage[hypertexnames=false]{hyperref}}{}}{}}{}
\makeatother

\ifthenelse{\isundefined{\nohref}}{%
  \hypersetup{colorlinks,linkcolor={red},citecolor={green},urlcolor={blue}}
  }{}

\newtheoremstyle{faded}
{\topsep}%
{\topsep}%
{\color{black!75!white}}%
{}%
{\color{black}\bfseries}
{.}
{.5em}%
{\thmname{#1}~\thmnumber{#2}\thmnote{ (#3)}}%

\ifthenelse{\isundefined{\theoremlook}}{
\def\theoremlook{theorem}}{}

\ifthenelse{\isundefined{\definitionlook}}{
\def\definitionlook{definition}}{}

\theoremstyle{\theoremlook}
\newtheorem{theorem}{Theorem}

\newtheorem{corollary}{Corollary}

\theoremstyle{\definitionlook}

\newtheorem{condition}{Condition}

\newtheorem{definition}{Definition}
\newtheorem{assumption}{Assumption}

\makeatletter
\def\thmenvs{theorem,lemma,corollary,proposition,remark,property,%
  condition,method,example,problem,definition,assumption}
\foreach \@te in \thmenvs {
  \bgroup
  \globaldefs=1
  \expandafter\let\csname o\@te\expandafter\endcsname\csname\@te\endcsname
  \expandafter\let\csname eo\@te\expandafter\endcsname\csname end\@te\endcsname
  \renewenvironment{\@te}[1][]{%
    \csname o\@currenvir\endcsname%
    \ifthenelse{\equal{##1}{}}{}{\label{\@currenvir:##1}}
  }{%
    \expandafter\csname eo\@currenvir\endcsname%
  }
  \egroup
}
\makeatother

\newcommand{\defvar}[3][show]{\nomenclature{#2}{#3}%
  \ifthenelse{\equal{#1}{show}}{#2}{}}

\input{latex_common/optimization.tex}

\makeflag{prettycsm}{true}
\makeflag{twocolcsm}{true}
\ifthenelse{\isundefined{\preprintcsm}}{}{
  \makeflag{prettycsm}{false}
  \makeflag{twocolcsm}{false}
}
\ifthenelse{\isundefined{\arxivcsm}}{}{
  \makeflag{prettycsm}{false}
  \makeflag{arxivcsm}{true}
}

\ifmaketwocolcsm
\else

\usepackage{csm16}
\usepackage[margin=1in]{geometry}
\usepackage[nolists,nomarkers,tablesfirst]{endfloat} 

\usepackage{url}
\usepackage{graphicx,xcolor}
\usepackage{verbatim}
\makeatletter
\newcommand{\verbatimfont}[1]{\def\verbatim@font{#1}}%
\makeatother
\verbatimfont{\ttfamily\small}

\newcommand{\bi}{\begin{itemize}}\newcommand{\ei}{\end{itemize}}
\newcommand{\be}{\begin{equation}}\newcommand{\ee}{\end{equation}}
\newcommand{\bee}{\begin{enumerate}}\newcommand{\eee}{\end{enumerate}}
\newcommand{\bea}{\begin{eqnarray}}\newcommand{\eea}{\end{eqnarray}}
\newcommand{\beas}{\begin{eqnarray*}}\newcommand{\eeas}{\end{eqnarray*}}
\newcommand{\bc}{\begin{center}}\newcommand{\ec}{\end{center}}
\usepackage[left,pagewise]{lineno}
\usepackage[english]{babel}
\usepackage{blindtext}

\newif\ifPDF \ifx\pdfoutput\undefined\PDFfalse \else \ifnum\pdfoutput > 0\PDFtrue \else\PDFfalse \fi
\ifPDF
\usepackage[pdftex, plainpages = false, colorlinks=true, linkcolor=black, citecolor = green!50!blue, urlcolor = blue, filecolor=black, pagebackref=false, hypertexnames=false,  pdfpagelabels ]{hyperref}
\fi

\gdef\csmpreprintfigscale{1.1}

\fi

\makeatletter

\definecolor{csm@intro@first@char}{HTML}{df9434}
\definecolor{csm@sidebar@bg}{HTML}{f2e9c4}
\colorlet{csm@sidebar@arxiv@bg}{black!2!white}
\definecolor{csm@sidebar@first@char}{HTML}{a08c2f}
\definecolor{csm@fig@bg}{HTML}{e6f7fe}
\definecolor{@light@blue}{HTML}{26baf2}
\colorlet{csm@fig@label@body}{@light@blue}
\colorlet{csm@blurb@color}{@light@blue}
\definecolor{@golden}{HTML}{c8b03b}
\colorlet{csm@fig@label@sidebar}{@golden}
\colorlet{csm@abstract@color}{@golden}
\definecolor{csm@table@caption@bg}{HTML}{abc2e9}
\definecolor{csm@table@border}{HTML}{5e81c9}
\colorlet{csm@table@arxiv@border}{black}
\colorlet{csm@table@bg}{csm@fig@bg}
\colorlet{csm@table@arxiv@bg}{black!4!white}
\ifmakearxivcsm
\definecolor{TitleRed}{HTML}{db6245}
\fi

\makeatother

\ifmaketwocolcsm
\makeatletter
\def\@csm@parindent{1em}
\def\@csm@parskip{0pt plus1pt}

\makeatother
\fi

\ifmakeprettycsm
\let\temp\rmdefault
\usepackage{mathpazo}
\let\rmdefault\temp

\DeclareMathSizes{10}{9}{7}{5}

\DeclareSymbolFont{Greekletters}{OT1}{iwona}{m}{n}
\DeclareSymbolFont{greekletters}{OML}{iwona}{m}{it}
\DeclareMathSymbol{\Delta}{\mathord}{Greekletters}{"01}
\DeclareMathSymbol{\Theta}{\mathord}{Greekletters}{"02}
\DeclareMathSymbol{\Lambda}{\mathord}{Greekletters}{"03}
\DeclareMathSymbol{\Xi}{\mathord}{Greekletters}{"04}
\DeclareMathSymbol{\Pi}{\mathord}{Greekletters}{"05}
\DeclareMathSymbol{\Sigma}{\mathord}{Greekletters}{"06}
\DeclareMathSymbol{\Upsilon}{\mathord}{Greekletters}{"07}
\DeclareMathSymbol{\Phi}{\mathord}{Greekletters}{"08}
\DeclareMathSymbol{\Psi}{\mathord}{Greekletters}{"09}
\DeclareMathSymbol{\Omega}{\mathord}{Greekletters}{"0A}
\DeclareMathSymbol{\alpha}{\mathord}{greekletters}{"0B}
\DeclareMathSymbol{\beta}{\mathord}{greekletters}{"0C}
\DeclareMathSymbol{\gamma}{\mathord}{greekletters}{"0D}
\DeclareMathSymbol{\delta}{\mathord}{greekletters}{"0E}
\DeclareMathSymbol{\epsilon}{\mathord}{greekletters}{"0F}
\DeclareMathSymbol{\zeta}{\mathord}{greekletters}{"10}
\DeclareMathSymbol{\eta}{\mathord}{greekletters}{"11}
\DeclareMathSymbol{\theta}{\mathord}{greekletters}{"12}
\DeclareMathSymbol{\iota}{\mathord}{greekletters}{"13}
\DeclareMathSymbol{\kappa}{\mathord}{greekletters}{"14}
\DeclareMathSymbol{\lambda}{\mathord}{greekletters}{"15}
\DeclareMathSymbol{\mu}{\mathord}{greekletters}{"16}
\DeclareMathSymbol{\nu}{\mathord}{greekletters}{"17}
\DeclareMathSymbol{\xi}{\mathord}{greekletters}{"18}
\DeclareMathSymbol{\pi}{\mathord}{greekletters}{"19}
\DeclareMathSymbol{\rho}{\mathord}{greekletters}{"1A}
\DeclareMathSymbol{\sigma}{\mathord}{greekletters}{"1B}
\DeclareMathSymbol{\tau}{\mathord}{greekletters}{"1C}
\DeclareMathSymbol{\upsilon}{\mathord}{greekletters}{"1D}
\DeclareMathSymbol{\phi}{\mathord}{greekletters}{"1E}
\DeclareMathSymbol{\chi}{\mathord}{greekletters}{"1F}
\DeclareMathSymbol{\psi}{\mathord}{greekletters}{"20}
\DeclareMathSymbol{\omega}{\mathord}{greekletters}{"21}
\DeclareMathSymbol{\varepsilon}{\mathord}{greekletters}{"22}
\DeclareMathSymbol{\vartheta}{\mathord}{greekletters}{"23}
\DeclareMathSymbol{\varpi}{\mathord}{greekletters}{"24}
\DeclareMathSymbol{\varrho}{\mathord}{greekletters}{"25}
\DeclareMathSymbol{\varsigma}{\mathord}{greekletters}{"26}
\DeclareMathSymbol{\varphi}{\mathord}{greekletters}{"27}

\let\oldprime\prime
\renewcommand{\prime}{^\oldprime}

\newcommand\der{^{\mkern+2.5mu\raise-3pt\hbox{$\scriptstyle\prime$}}}
\else
\newcommand\der{'}
\fi

\ifmaketwocolcsm
\tikzset{
  list bullet/.style={
    scale=0.5
  }
}

\setlist[itemize,1]{label={%
    \tikz[outer sep=0,inner sep=0,baseline=-0.3em]{%
      \node[list bullet] at (0,0) {\faChevronRight};%
      \node[list bullet] at (0.25em,0) {\faChevronRight};}%
  },leftmargin=1.3em}
\fi

\ifmaketwocolcsm
\usepackage{geometry}

\ExplSyntaxOn
\ifmakeprettycsm
\tl_set:Nn \l_sidemargin_dim {1.5cm}
\tl_set:Nn \l_topmargin_dim {2.2cm}
\tl_set:Nn \l_botmargin_dim {1.6cm}
\tl_set:Nn \l_footmargin_dim {6mm}
\tl_set:Nn \l_columnsep_dim {3mm}
\else
\ifmakearxivcsm
\tl_set:Nn \l_sidemargin_dim {1.2cm}
\tl_set:Nn \l_topmargin_dim {2cm}
\tl_set:Nn \l_botmargin_dim {1.4cm}
\tl_set:Nn \l_footmargin_dim {6mm}
\tl_set:Nn \l_columnsep_dim {4mm}

\xdef\csmtopmargin{\tl_use:N \l_topmargin_dim}

\fi
\fi
\ExplSyntaxOff
\fi

\ifmaketwocolcsm
\usepackage{helvet}
\titleformat{\section}{%
  \normalfont\fontsize{10}{12}%
  \sffamily\bfseries\MakeUppercase}{}{0cm}{#1}
\titlespacing*{\section}{0cm}{4mm}{0mm}

\titleformat{\subsection}{%
  \normalfont\fontsize{10}{12}%
  \sffamily\bfseries\slshape}{}{0cm}{#1}
\titlespacing*{\subsection}{0cm}{4mm}{0mm}

\titleformat{\subsubsection}{%
  \normalfont\fontsize{10}{12}%
  \sffamily}{}{0cm}{#1}
\titlespacing*{\subsubsection}{0cm}{4mm}{0mm}
\fi

\makeatletter

\ifmaketwocolcsm
\newcommand{\firstword}[2][]{%
  \StrChar{#2}{1}[\first@char]%
  \StrBehind{#2}{\first@char}[\rest@of@word]%
  \renewcommand{\LettrineFontHook}{\sffamily\bfseries}%
  \lettrine[lines=5,#1]{%
    \ifmakeprettycsm
    \color{csm@intro@first@char}%
    \else
    \ifmakearxivcsm
    \color{TitleRed}%
    \fi
    \fi
    \first@char}{\rest@of@word}%
}
\else
\newcommand{\firstword}[2][]{#2}
\fi

\makeatother

\ifmakeprettycsm
\usepackage{fontawesome}

\ExplSyntaxOn
\tl_new:N \l_csm_date
\tl_set:Nn \l_csm_date {}

\newcommand{\csmdate}[1]{
  \tl_set:Nn \l_csm_date {#1}
}

\fancypagestyle{csmfancy}{%
  \fancyhf{}
  
  \fancyfoot[LO]{%
    \sffamily%
    \fontsize{9}{11}%
    \textbf{\thepage}~~%
    \fontsize{6}{8}%
    \textbf{IEEE~CONTROL~SYSTEMS}~~%
    \tikz[outer~sep=0,inner~sep=0]{%
      \node[scale=0.8] at (0,0) {\faChevronRight};%
      \node[scale=0.8] at (0.4em,0) {\faChevronRight};}%
    {\normalfont\sffamily{}~\tl_use:N \l_csm_date}%
  }%
  \fancyfoot[RE]{%
    \fontsize{6}{8}%
    \sffamily%
    {\tl_use:N \l_csm_date}~%
    \tikz[outer~sep=0,inner~sep=0]{%
      \node[scale=0.8] at (0,0) {\faChevronLeft};%
      \node[scale=0.8] at (0.4em,0) {\faChevronLeft};}%
    \textbf{~IEEE~CONTROL~SYSTEMS}%
    \fontsize{9}{11}%
    \textbf{~~\thepage}%
  }
}
\ExplSyntaxOff
\else
\ifmakearxivcsm
\usepackage{fontawesome}

\ExplSyntaxOn
\tl_new:N \l_csm_date
\tl_set:Nn \l_csm_date {}

\newcommand{\csmdate}[1]{
  \tl_set:Nn \l_csm_date {#1}
}

\fancypagestyle{csmfancy}{%
  \fancyhf{}
  
  \fancyfoot[LO]{%
    \sffamily%
    \fontsize{9}{11}%
    \textbf{\thepage}%
  }%
  \fancyfoot[RE]{%
    \fontsize{6}{8}%
    \sffamily%
    \fontsize{9}{11}%
    \textbf{\thepage}%
  }%
  \fancyhead[C]{%
    \footnotesize\sffamily\color{black!20!white}%
    Preprint%
  }%
}
\ExplSyntaxOff
\else
\newcommand{\csmdate}[1]{}
\fi
\fi

\allowdisplaybreaks

\ifmaketwocolcsm
\makeatletter

\newtheoremstyle{csm@theorem@style}
{4mm}
{1pt}
{}
{}
{\fontsize{9}{11}\sffamily}
{}
{\newline}
{}

\theoremstyle{csm@theorem@style}

\newtheorem{csm@theorem}{Theorem}
\newtheorem{csm@lemma}{Lemma}
\newtheorem{csm@corollary}{Corollary}
\newtheorem{csm@proposition}{Proposition}
\newtheorem{csm@remark}{Remark}
\newtheorem{csm@property}{Property}
\newtheorem{csm@condition}{Condition}
\newtheorem{csm@method}{Method}
\newtheorem{csm@example}{Example}
\newtheorem{csm@problem}{Problem}
\newtheorem{csm@definition}{Definition}
\newtheorem{csm@assumption}{Assumption}

\def\@thmenvs{theorem,lemma,corollary,proposition,remark,property,%
  condition,method,example,problem,definition,assumption}

\foreach \@te in \@thmenvs {
  \bgroup
  \globaldefs=1
  \expandafter\let\csname ocsm@\@te\expandafter\endcsname
  \csname csm@\@te\endcsname
  \expandafter\let\csname eocsm@\@te\expandafter\endcsname
  \csname endcsm@\@te\endcsname
  \renewenvironment{\@te}[1][]{%
    \csname ocsm@\@currenvir\endcsname%
    \ifthenelse{\equal{##1}{}}{}{\label{\@currenvir:##1}}
  }{%
    \hfill$\blacksquare$
    \expandafter\csname eocsm@\@currenvir\endcsname%
  }
  \egroup
}

\makeatother
\fi

\ExplSyntaxOn
\makeatletter

\tl_new:N \l_csm_figure_label_tl
\tl_new:N \l_csm_figure_caption_tl
\tl_new:N \l_csm_figure_pos_tl
\tl_new:N \l_csm_figure_cols_tl

\define@key{csm@figure@keys}{label}{\tl_set:Nn \l_csm_figure_label_tl {#1}}
\define@key{csm@figure@keys}{caption}{\tl_set:Nn \l_csm_figure_caption_tl {#1}}
\define@key{csm@figure@keys}{position}{\tl_set:Nn \l_csm_figure_pos_tl {#1}}
\define@key{csm@figure@keys}{columns}{\tl_set:Nn \l_csm_figure_cols_tl {#1}}

\ifmakeprettycsm
\newtcolorbox{figurebox}{
  enhanced~jigsaw,
  width=\linewidth,
  arc=2mm,
  boxrule=0pt,
  titlerule=0pt,
  left=3mm,
  right=3mm,
  bottom=3mm,
  top=3mm,
  boxsep=\tcb@boxsep,
  colback=csm@fig@bg,
  colbacktitle=csm@fig@bg
}

\else
\ifmakearxivcsm
\newtcolorbox{figurebox}{
  enhanced~jigsaw,
  width=\linewidth,
  arc=0pt,
  boxrule=0pt,
  titlerule=0pt,
  left=3mm,
  right=3mm,
  bottom=3mm,
  top=3mm,
  boxsep=\tcb@boxsep,
  colback=white,
  colbacktitle=white,
  opacityfill=0
}

\fi
\fi

\ifmaketwocolcsm

\newenvironment{figure@internal}[2][tbp]{%
  \def\fig@width{#2}%
  \ifthenelse{\equal{\fig@width}{}}{%
    \begin{figure}[#1]}{\begin{figure*}[#1]}}{%
      \ifthenelse{\equal{\fig@width}{}}{%
      \end{figure}}{\end{figure*}}}

\newenvironment{csmfigure}[1][]{%
  \setkeys{csm@figure@keys}{label=,#1}%
  \setkeys{csm@figure@keys}{caption=,#1}%
  \setkeys{csm@figure@keys}{position=tbp,#1}%
  \setkeys{csm@figure@keys}{columns=1,#1}%
  \edef\csm@fig@pos{\tl_use:N \l_csm_figure_pos_tl}%
  \pgfmathparse{\tl_use:N \l_csm_figure_cols_tl ==1 ? 1:0}%
  \ifthenelse{\pgfmathresult>0}{%
    \gdef\csm@figure@width{}%
  }{%
    \gdef\csm@figure@width{dbl}%
  }%
  \edef\begin@figure@cmd{{figure@internal}[\csm@fig@pos}%
  \expandafter\begin\begin@figure@cmd]{\csm@figure@width}%
    \begin{figurebox}
    }{%
    \end{figurebox}
    \caption{\tl_use:N \l_csm_figure_caption_tl}
    \figlabel{\tl_use:N \l_csm_figure_label_tl}
  \end{figure@internal}
}
\else
\gdef\csm@singlecol@emulate@width{0.7}
\newenvironment{csmfigure}[1][]{%
  \setkeys{csm@figure@keys}{label=,#1}%
  \setkeys{csm@figure@keys}{caption=,#1}%
  \setkeys{csm@figure@keys}{position=tbp,#1}%
  \setkeys{csm@figure@keys}{columns=1,#1}%
  \edef\csm@num@cols{\tl_use:N \l_csm_figure_cols_tl}
  \pgfmathparse{\csm@num@cols>1 ? 1:0}
  \ifthenelse{\pgfmathresult>0}{\def\csm@figarea@width{1}}{%
    \def\csm@figarea@width{\csm@singlecol@emulate@width}}
  \begin{figure}
    \centering
    \begin{minipage}{\csm@figarea@width\textwidth}
    }{%
    \end{minipage}
    \caption{\tl_use:N \l_csm_figure_caption_tl}
    \figlabel{\tl_use:N \l_csm_figure_label_tl}
  \end{figure}
}
\DeclareDelayedFloatFlavor{csmfigure}{figure}
\fi

\makeatother
\ExplSyntaxOff

\newcommand{\fakesubfigref}[3][]{%
  \figref{#2}(#3)\ifthenelse{\equal{#1}{}}{}{-(#1)}}

\ExplSyntaxOn
\makeatletter

\tl_new:N \l_csm_table_label_tl
\tl_new:N \l_csm_table_caption_tl
\tl_new:N \l_csm_table_pos_tl
\tl_new:N \l_csm_table_cols_tl

\define@key{csm@table@keys}{label}{\tl_set:Nn \l_csm_table_label_tl {#1}}
\define@key{csm@table@keys}{caption}{\tl_set:Nn \l_csm_table_caption_tl {#1}}
\define@key{csm@table@keys}{position}{\tl_set:Nn \l_csm_table_pos_tl {#1}}
\define@key{csm@table@keys}{columns}{\tl_set:Nn \l_csm_table_cols_tl {#1}}

\ifmakeprettycsm
\def\@table@inner@pad{3mm}
\def\@table@border@width{1pt}

\newtcolorbox{tablebox}[1][]{
  enhanced~jigsaw,
  fonttitle=\sffamily\bfseries\small,
  fontupper=\sffamily,
  width=\linewidth,
  arc=2mm,
  boxrule=\@table@border@width,
  titlerule=\@table@border@width,
  left=\@table@inner@pad,
  right=\@table@inner@pad,
  bottom=\@table@inner@pad,
  top=\@table@inner@pad,
  toptitle=\@table@inner@pad,
  bottomtitle=\@table@inner@pad,
  boxsep=\tcb@boxsep,
  colback=csm@table@bg,
  coltitle=black,
  colbacktitle=csm@table@caption@bg,
  colframe=csm@table@border,
  title={#1},
  attach~boxed~title~to~top,
  boxed~title~style={colframe=csm@table@border,rounded~corners}
}
\else
\ifmakearxivcsm
\def\@table@inner@pad{3mm}
\def\@table@border@width{1pt}

\newtcolorbox{tablebox}[1][]{
  enhanced~jigsaw,
  fonttitle=\sffamily\bfseries\small,
  fontupper=\sffamily,
  width=\linewidth,
  arc=2mm,
  boxrule=\@table@border@width,
  titlerule=\@table@border@width,
  left=\@table@inner@pad,
  right=\@table@inner@pad,
  bottom=\@table@inner@pad,
  top=\@table@inner@pad,
  toptitle=\@table@inner@pad,
  bottomtitle=\@table@inner@pad,
  boxsep=\tcb@boxsep,
  colback=csm@table@arxiv@bg,
  coltitle=black,
  colbacktitle=csm@table@caption@bg,
  colframe=csm@table@arxiv@border,
  minipage~boxed~title,
  title={#1},
  attach~boxed~title~to~top~left,
  boxed~title~style={%
    rounded~corners,
    boxrule=0pt,
    colback=white,
    colframe=white,
    frame~hidden,
    left=0pt,
    right=0pt
  }
}
\fi
\fi

\ifmaketwocolcsm

\newenvironment{table@internal}[2][tbp]{%
  \def\table@width{#2}%
  \ifthenelse{\equal{\table@width}{}}{%
    \begin{table}[#1]}{\begin{table*}[#1]}}{%
      \ifthenelse{\equal{\table@width}{}}{%
      \end{table}}{\end{table*}}}

\newenvironment{csmtable}[1][]{%
  \setkeys{csm@table@keys}{label=,#1}%
  \setkeys{csm@table@keys}{caption=,#1}%
  \setkeys{csm@table@keys}{position=tbp,#1}%
  \setkeys{csm@table@keys}{columns=1,#1}%
  \edef\csm@table@pos{\tl_use:N \l_csm_table_pos_tl}%
  \pgfmathparse{\tl_use:N \l_csm_table_cols_tl ==1 ? 1:0}%
  \ifthenelse{\pgfmathresult>0}{%
    \gdef\csm@table@width{}%
  }{%
    \gdef\csm@table@width{dbl}%
  }%
  \edef\begin@table@cmd{{table@internal}[\csm@table@pos}%
    \expandafter\begin\begin@table@cmd]{\csm@table@width}%
    \refstepcounter{table}
    \begin{tablebox}[%
      {\bfseries TABLE~\arabic{table}}\hspace{1em}%
      \tl_use:N \l_csm_table_caption_tl]%
      \fontsize{9}{13}%
    }{%
    \end{tablebox}
    \tablabel{\tl_use:N \l_csm_table_label_tl}
  \end{table@internal}
}
\else
\newenvironment{csmtable}[1][]{%
  \setkeys{csm@table@keys}{label=,#1}%
  \setkeys{csm@table@keys}{caption=,#1}%
  \setkeys{csm@table@keys}{position=tbp,#1}%
  \setkeys{csm@table@keys}{columns=1,#1}%
  \edef\csm@num@cols{\tl_use:N \l_csm_table_cols_tl}
  \pgfmathparse{\csm@num@cols>1 ? 1:0}
  \ifthenelse{\pgfmathresult>0}{\def\csm@tabarea@width{1}}{%
    \def\csm@tabarea@width{\csm@singlecol@emulate@width}}
  \begin{table}[t]
    \caption{\tl_use:N \l_csm_table_caption_tl}
    \vspace*{.1in}
    \centering
    \begin{minipage}{\csm@tabarea@width\textwidth}
    }{%
    \end{minipage}
    \tablabel{\tl_use:N \l_csm_table_label_tl}
  \end{table}
}
\DeclareDelayedFloatFlavor{csmtable}{table}
\fi

\makeatother
\ExplSyntaxOff

\makeatletter

\let\oldhypertarget\hypertarget
\renewcommand{\hypertarget}[2]{%
  \oldhypertarget{#1}{#2}%
  \protected@write\@mainaux{}{%
    \string\expandafter\string\gdef
    \string\csname\string\detokenize{#1}\string\endcsname{#2}%
  }%
}

\newcommand{\csmhyperlink}[1]{%
  \hyperlink{#1}{\csname #1\endcsname}%
}

\makeatother

\ExplSyntaxOn
\makeatletter

\tl_new:N \l_csm_sidebar_title_tl
\tl_new:N \l_csm_sidebar_label_tl
\tl_new:N \l_csm_sidebar_side_tl
\tl_new:N \l_csm_sidebar_pos_tl
\tl_new:N \l_csm_sidebar_cols_tl

\define@key{csm@sidebar@keys}{title}{\tl_set:Nn \l_csm_sidebar_title_tl {#1}}
\define@key{csm@sidebar@keys}{label}{\tl_set:Nn \l_csm_sidebar_label_tl {#1}}
\define@key{csm@sidebar@keys}{side}{\tl_set:Nn \l_csm_sidebar_side_tl {#1}}
\define@key{csm@sidebar@keys}{position}{\tl_set:Nn \l_csm_sidebar_pos_tl {#1}}
\define@key{csm@sidebar@keys}{columns}{\tl_set:Nn \l_csm_sidebar_cols_tl {#1}}

\newcounter{sidebar@equation}%
\newcounter{sidebar@table}%
\newcounter{sidebar@figure}%
\newcounter{subsidebar@counter}

\newcounter{sbrefs}
\renewbibmacro*{finentry}{\stepcounter{sbrefs}\finentry}

\newcommand{\csm@sb@printbib}[2]{%
  \printbibliography[title=References,keyword=#1,resetnumbers=#2]}

\newcommand{\make@sidebar@bib}{%
  \edef\sidebar@bib@label{sidebar:\tl_use:N \l_csm_sidebar_label_tl}%
  \ifmaketwocolcsm
  \def\@sb@bibnum@reset{false}%
  \else
  \pgfmathparse{\thesbrefs==0 ? 1:0}
  \ifthenelse{\pgfmathresult>0}{%
    \def\@sb@bibnum@reset{true}
  }{%
    \def\@sb@bibnum@reset{false}
  }%
  \fi
  \begin{refcontext}[sorting=none,labelprefix=S]%
    \ifmaketwocolcsm
    \renewcommand*{\bibfont}{\fontsize{7}{9}\sffamily}%
    \else
    \renewcommand*{\bibfont}{\footnotesize}%
    \fi
    \DeclareFieldFormat{labelnumberwidth}{[##1]}%
    \begingroup\edef\sidebar@print@bib{\endgroup\noexpand
      \printbibliography[title=References,
      keyword=\sidebar@bib@label,
      resetnumbers=\@sb@bibnum@reset]
    }\sidebar@print@bib
  \end{refcontext}%
}

\newcommand{\sidebar@labels@on}{
  \global\chardef\dc@currentequation=\value{equation}%
  \global\chardef\dc@currentfigure=\value{figure}%
  \global\chardef\dc@currenttable=\value{table}%
  \let\c@equation\c@sidebar@equation
  \let\c@figure\c@sidebar@figure
  \let\c@table\c@sidebar@table
  \renewcommand{\theequation}{S\arabic{equation}}%
  \renewcommand{\thetable}{S\arabic{table}}%
  \renewcommand{\thefigure}{S\arabic{figure}}%
  \renewcommand{\theHequation}{S\arabic{equation}}%
  \renewcommand{\theHtable}{S\arabic{table}}%
  \renewcommand{\theHfigure}{S\arabic{figure}}%
}

\newcommand{\sidebar@labels@off}{
  \edef\equation@mem{\arabic{sidebar@equation}}%
  \edef\figure@mem{\arabic{sidebar@figure}}%
  \edef\table@mem{\arabic{sidebar@table}}%
  \setcounter{equation}{\dc@currentequation}%
  \setcounter{figure}{\dc@currentfigure}%
  \setcounter{table}{\dc@currenttable}%
  \setcounter{sidebar@equation}{\equation@mem}%
  \setcounter{sidebar@figure}{\figure@mem}%
  \setcounter{sidebar@table}{\table@mem}%
}

\ifmakeprettycsm

\def\tcb@left{1em}
\def\tcb@right{1em}
\def\tcb@bottom{1em}
\def\tcb@boxsep{0pt}
\def\tcb@leftfill{0cm}
\def\tcb@rightfill{0cm}

\def\continued@message{\textit{(continued...)}}

\newtcolorbox{sidebarbox}[1][]{
  enhanced~jigsaw,
  fonttitle=\sffamily\bfseries\large,
  width=#1,
  arc=2mm,
  boxrule=0pt,
  titlerule=0pt,
  toptitle=2mm,
  left=\tcb@left,
  right=\tcb@right,
  bottom=\tcb@bottom,
  boxsep=\tcb@boxsep,
  leftrule=\tcb@leftfill,
  grow~to~left~by=\tcb@leftfill,
  rightrule=\tcb@rightfill,
  grow~to~right~by=\tcb@rightfill,
  colback=csm@sidebar@bg,
  colframe=csm@sidebar@bg,
  colbacktitle=csm@sidebar@bg,
  coltitle=black
}

\else
\ifmakearxivcsm

\def\tcb@left{1em}
\def\tcb@right{1em}
\def\tcb@bottom{1em}
\def\tcb@boxsep{0pt}
\def\tcb@leftfill{0cm}
\def\tcb@rightfill{0cm}

\def\continued@message{\textit{(continued...)}}

\newtcolorbox{sidebarbox}[1][]{
  enhanced~jigsaw,
  fonttitle=\sffamily\bfseries\large,
  width=#1,
  boxrule=0pt,
  arc=0pt,
  titlerule=0pt,
  toptitle=2mm,
  left=\tcb@left,
  right=\tcb@right,
  bottom=\tcb@bottom,
  boxsep=\tcb@boxsep,
  bottomrule=1.5pt,
  colback=csm@sidebar@arxiv@bg,
  frame~style={left~color=white,
    right~color=white,
    top~color=white,
    bottom~color=black},
  colbacktitle=csm@sidebar@arxiv@bg,
  coltitle=black
}

\fi

\ifmaketwocolcsm

\newsavebox{\sidebar@box@full}
\newsavebox{\sidebar@box@left}
\newsavebox{\sidebar@box@right}

\newenvironment{sidebar}[1][]{%
  \setkeys{csm@sidebar@keys}{title=,#1}%
  \setkeys{csm@sidebar@keys}{label=,#1}%
  \setkeys{csm@sidebar@keys}{side=odd,#1}%
  \setkeys{csm@sidebar@keys}{position=t,#1}%
  \setkeys{csm@sidebar@keys}{columns=2,#1}%
  \edef\sidebar@side{\tl_use:N \l_csm_sidebar_side_tl}%
  \edef\fig@pos{\tl_use:N \l_csm_sidebar_pos_tl}%
  \pgfmathparse{\tl_use:N \l_csm_sidebar_cols_tl ==1 ? 1:0}%
  \ifthenelse{\pgfmathresult>0}{%
    \gdef\figure@width{}%
    \tikzmath{\sidebar@box@width=\columnwidth;
      \sidebar@width=\sidebar@box@width-\tcb@left-\tcb@right-\tcb@boxsep;}
    \xdef\sidebar@box@width{\sidebar@box@width pt}%
    \xdef\sidebar@width{\sidebar@width pt}%
  }{%
    \gdef\figure@width{dbl}%
    \tikzmath{\sidebar@box@width=\textwidth;
      \sidebar@width=\sidebar@box@width-\tcb@left-\tcb@right-\tcb@boxsep;}
    \xdef\sidebar@box@width{\sidebar@box@width pt}%
    \xdef\sidebar@width{\sidebar@width pt}%
  }%
  \begin{subsidebar}[reset]%
  }{%
  \end{subsidebar}%
  \pgfmathparse{\value{subsidebar@counter}==1 ? 1:0}%
  \ifthenelse{\pgfmathresult>0}{%
    \sidebar@draw{none}{\sidebar@box@full}%
  }{%
    \ifthenelse{\equal{\sidebar@side}{odd}}{%
      \def\sidebar@side@other{even}%
    }{%
      \def\sidebar@side@other{odd}%
    }%
    \sidebar@draw{\sidebar@side}{\sidebar@box@left}%
    \sidebar@draw[notitle]{\sidebar@side@other}{\sidebar@box@right}%
  }%
  \setcounter{subsidebar@counter}{0}
}

\newcommand{\sidebar@draw}[3][title]{
  \xdef\sidebar@title{\tl_use:N \l_csm_sidebar_title_tl}%
  \xdef\sidebar@label{sidebar:\tl_use:N \l_csm_sidebar_label_tl}%
  \edef\begin@sidebar@cmd{{figure@internal}[\fig@pos}%
  \expandafter\begin\begin@sidebar@cmd]{\figure@width}%
    \ifthenelse{\equal{#1}{title}}{%
      \tcbset{%
        title={\hypertarget{\sidebar@label}{\sidebar@title}}
      }%
    }{}%
    \ifthenelse{\equal{#2}{none}}{%
      \def\tcb@leftfill{0cm}%
      \def\tcb@rightfill{0cm}%
    }{%
      \ifthenelse{\equal{#2}{odd}}{%
        \def\tcb@leftfill{0cm}%
        \def\tcb@rightfill{10cm}%
      }{%
        \def\tcb@leftfill{10cm}%
        \def\tcb@rightfill{0cm}%
      }%
    }%
    \begin{sidebarbox}[\sidebar@box@width]%
      \usebox{#3}%
    \end{sidebarbox}%
  \end{figure@internal}%
}

\newenvironment{subsidebar}[1][noreset]{
  \def\sidebar@reset{#1}
  \pgfmathparse{\value{subsidebar@counter}>=3 ? 1:0}%
  \ifthenelse{\pgfmathresult>0}{%
    \PackageError{CSM~style}{Only~two~subsidebars~are~allowed}{}%
  }{}%
  \pgfmathparse{\value{subsidebar@counter}==0 ? 1:0}%
  \ifthenelse{\pgfmathresult>0}{%
    \def\sidebar@box{sidebar@box@full}%
  }{%
    \pgfmathparse{\value{subsidebar@counter}==1 ? 1:0}%
    \ifthenelse{\pgfmathresult>0}{%
      \def\sidebar@box{sidebar@box@left}%
    }{%
      \def\sidebar@box{sidebar@box@right}%
    }%
  }%
  \stepcounter{subsidebar@counter}%
  \begin{lrbox}{\csname\sidebar@box\endcsname}%
    \begin{minipage}{\sidebar@width}%
      \ifthenelse{\equal{\sidebar@reset}{noreset}}{}{\sidebar@labels@on}%
      \ifthenelse{\equal{\figure@width}{}}{}{%
        \begin{multicols}{2}%
        }%
        \fontsize{9}{13}
        \setlength\parindent{\@csm@parindent} 
        \setlength{\parskip}{\@csm@parskip} 
        \sffamily\sansmath
        \ifmakeprettycsm
        \captionsetup[figure]{name={%
            \bfseries\color{csm@fig@label@sidebar}FIGURE~}}
        \fi
        \pgfmathparse{\value{subsidebar@counter}==3 ? 1:0}%
        \ifthenelse{\pgfmathresult>0}{%
          \continued@message%
        }{}%
      }{%
        \pgfmathparse{\value{subsidebar@counter}!=2 ? 1:0}%
        \ifthenelse{\pgfmathresult>0}{%
          \make@sidebar@bib%
        }{}%
        \ifthenelse{\pgfmathresult>0}{}{%
          \continued@message%
        }%
        \ifthenelse{\equal{\figure@width}{}}{}{%
        \end{multicols}%
      }%
      \ifthenelse{\equal{\sidebar@reset}{noreset}}{}{\sidebar@labels@off}%
    \end{minipage}%
  \end{lrbox}%
  \global%
  \expandafter\setbox\csname\sidebar@box\endcsname%
  \expandafter\box\csname\sidebar@box\endcsname%
}

\else

\newenvironment{sidebar}[1][]{%
  \setkeys{csm@sidebar@keys}{title=,#1}%
  \setkeys{csm@sidebar@keys}{label=,#1}%
  \setkeys{csm@sidebar@keys}{side=odd,#1}%
  \setkeys{csm@sidebar@keys}{position=t,#1}%
  \setkeys{csm@sidebar@keys}{columns=2,#1}%
  \edef\@csm@sb@title{\tl_use:N \l_csm_sidebar_title_tl}
  \edef\@csm@sb@label{sidebar:\tl_use:N \l_csm_sidebar_label_tl}%
  \clearpage
  \newpage
  \section[\@csm@sb@title]{Sidebar:~%
    \hypertarget{\@csm@sb@label}{\@csm@sb@title}}
  \sidebar@labels@on 
}{
  \make@sidebar@bib
  \newpage
  \processdelayedfloats 
  \clearpage
  \sidebar@labels@off 
}

\newenvironment{subsidebar}[1][noreset]{}{}

\fi

\newcounter{sidebar@queue@size}%
\newcommand{\reset@sidebar@queue}{
  \gdef\csm@sb@queue@title{}
  \gdef\csm@sb@queue@label{}
  \gdef\csm@sb@queue@side{}
  \gdef\csm@sb@queue@position{}
  \gdef\csm@sb@queue@columns{}
  \gdef\csm@sb@queue@content{}
  \setcounter{sidebar@queue@size}{0}
}
\reset@sidebar@queue

\newcommand{\declaresidebar}[2][]{
  \setkeys{csm@sidebar@keys}{title=,#1}%
  \setkeys{csm@sidebar@keys}{label=,#1}%
  \setkeys{csm@sidebar@keys}{side=odd,#1}%
  \setkeys{csm@sidebar@keys}{position=t,#1}%
  \setkeys{csm@sidebar@keys}{columns=2,#1}%
  \add@list@item{\csm@sb@queue@title}{\tl_use:N \l_csm_sidebar_title_tl}
  \add@list@item{\csm@sb@queue@label}{\tl_use:N \l_csm_sidebar_label_tl}
  \add@list@item{\csm@sb@queue@side}{\tl_use:N \l_csm_sidebar_side_tl}
  \add@list@item{\csm@sb@queue@position}{\tl_use:N \l_csm_sidebar_pos_tl}
  \add@list@item{\csm@sb@queue@columns}{\tl_use:N \l_csm_sidebar_cols_tl}
  \add@list@item{\csm@sb@queue@content}{#2}
  \stepcounter{sidebar@queue@size}
  \ifmaketwocolcsm
  \flushsidebars
  \fi
}

\makeflag{printsb}{false}
\newcommand{\flushsidebars}{
  \pgfmathparse{\thesidebar@queue@size>0 ? 1:0}
  \ifthenelse{\pgfmathresult>0}{\makeprintsbtrue}{\makeprintsbfalse}
  \ifmakeprintsb
  \setcounter{sbrefs}{0}
  \makelist{csm@sb@titles}{\csm@sb@queue@title}
  \makelist{csm@sb@labels}{\csm@sb@queue@label}
  \makelist{csm@sb@sides}{\csm@sb@queue@side}
  \makelist{csm@sb@positions}{\csm@sb@queue@position}
  \makelist{csm@sb@columns}{\csm@sb@queue@columns}
  \makelist{csm@sb@contents}{\csm@sb@queue@content}
  \foreach \@sidebar@number in {1,...,\thesidebar@queue@size} {
    \begin{sidebar}[%
      title={\list@nth{csm@sb@titles}{\@sidebar@number}},
      label={\list@nth{csm@sb@labels}{\@sidebar@number}},
      side={\list@nth{csm@sb@sides}{\@sidebar@number}},
      position={\list@nth{csm@sb@positions}{\@sidebar@number}},
      columns={\list@nth{csm@sb@columns}{\@sidebar@number}}]
      \input{\list@nth{csm@sb@contents}{\@sidebar@number}}
    \end{sidebar}
  }
  \reset@sidebar@queue 
  \fi
}

\ifmaketwocolcsm

\newcommand{\sbfirstword}[1]{%
  \StrChar{#1}{1}[\first@char]%
  \StrBehind{#1}{\first@char}[\rest@of@word]%
  \lettrine[nindent=0pt,slope=0pt,findent=3pt]{%
    \ifmakeprettycsm
    \color{csm@sidebar@first@char}\bfseries%
    \else
    \ifmakearxivcsm
    \bfseries%
    \fi
    \fi
    \first@char}{\rest@of@word}%
}
\else
\newcommand{\sbfirstword}[1]{#1}
\fi

\newcommand{\sbref}[1]{%
  \@ifundefined{sidebar:#1}{%
    ``\textbf{??}''%
  }{%
    ``\csmhyperlink{sidebar:#1}''%
  }%
}

\makeatother
\ExplSyntaxOff

\ifmaketwocolcsm
\makeatletter

\ifmakeprettycsm
\def\blurb@pad{4mm}

\newtcolorbox{blurbbox}{
  enhanced jigsaw,
  fontupper=\sffamily\bfseries\large,
  width=\textwidth,
  arc=2mm,
  boxrule=0.5pt,
  titlerule=0pt,
  toptitle=2mm,
  left=\blurb@pad,
  right=\blurb@pad,
  bottom=\blurb@pad,
  top=\blurb@pad,
  boxsep=0pt,
  opacityback=0,
  colframe=csm@blurb@color,
  coltext=csm@blurb@color
}
\else
\def\blurb@pad{4mm}

\newtcolorbox{blurbbox}{
  enhanced jigsaw,
  fontupper=\sffamily\itshape\large,
  width=\textwidth,
  sharp corners=all,
  boxrule=0pt,
  titlerule=0pt,
  toptitle=2mm,
  left=\blurb@pad,
  right=\blurb@pad,
  bottom=\blurb@pad,
  top=\blurb@pad,
  boxsep=0pt,
  colback=black!5!white,
  colframe=white,
  coltext=black
}
\fi

\newenvironment{blurb}{
  \begin{figure*}[!t]
    \centering
    \begin{blurbbox}
      \centering
      \setlength{\baselineskip}{1.4em}
    }{
    \end{blurbbox}
  \end{figure*}
}

\makeatother
\else
\NewEnviron{blurb}{}
\fi

\ifmakeprettycsm
\makeatletter

\def\abstract@pad{2mm}

\newtcolorbox{abstractbox}{
  enhanced jigsaw,
  fontupper=\sffamily,
  width=\textwidth,
  arc=2mm,
  boxrule=0.5pt,
  titlerule=0pt,
  toptitle=2mm,
  left=\abstract@pad,
  right=\abstract@pad,
  bottom=\abstract@pad,
  top=\abstract@pad,
  boxsep=0pt,
  opacityback=0,
  colframe=csm@abstract@color,
  coltext=black!70
}

\newenvironment{csmabstract}{
  \begin{figure*}[t]
    \centering
    \begin{abstractbox}%
      \paragraph{\sffamily Abstract (hidden in print).}%
    }{%
    \end{abstractbox}
  \end{figure*}
}

\makeatother
\else
\NewEnviron{csmabstract}{}
\fi

\ExplSyntaxOn
\makeatletter

\tl_new:N \l_csm_nomenclature_title_tl
\tl_new:N \l_csm_nomenclature_label_tl
\tl_new:N \l_csm_nomenclature_first_col_width_tl
\tl_new:N \l_csm_nomenclature_position_tl

\define@key{csm@nomenclature@keys}{title}{\tl_set:Nn
  \l_csm_nomenclature_title_tl {#1}}
\define@key{csm@nomenclature@keys}{label}{\tl_set:Nn
  \l_csm_nomenclature_label_tl {#1}}
\define@key{csm@nomenclature@keys}{width}{\tl_set:Nn
  \l_csm_nomenclature_first_col_width_tl {#1}}
\define@key{csm@nomenclature@keys}{position}{\tl_set:Nn
  \l_csm_nomenclature_position_tl {#1}}

\newcommand{\csmnomenclature}[2][]{
  \setkeys{csm@nomenclature@keys}{title=,#1}%
  \setkeys{csm@nomenclature@keys}{label=,#1}%
  \setkeys{csm@nomenclature@keys}{width=0.2,#1}%
  \setkeys{csm@nomenclature@keys}{position=t,#1}%
  \edef\nomenc@title{\tl_use:N
    \l_csm_nomenclature_title_tl}%
  \edef\nomenc@label{\tl_use:N
    \l_csm_nomenclature_label_tl}%
  \edef\nomenc@width{\tl_use:N
    \l_csm_nomenclature_first_col_width_tl}%
  \edef\nomenc@position{\tl_use:N
    \l_csm_nomenclature_position_tl}%
  \def\@nomenclature@content{\BODY}
  \newwrite\file
  \immediate\openout\file=csm_\nomenc@label.tex
  \immediate\write\file{\noexpand\footnotesize}
  \immediate\write\file{\string\def\string\arraystretch{1.2}}
  \immediate\write\file{\string\newcolumntype{V}{%
      >{\string\hsize=\nomenc@width\string\linewidth}X}}
  \immediate\write\file{\string\begin{tabularx}{\string\columnwidth}{VX}}
  \immediate\write\file{  \string\input{#2}}
  \immediate\write\file{\string\end{tabularx}}
  \immediate\closeout\file
  \declaresidebar[%
  title={\nomenc@title},
  columns=1,
  label={\nomenc@label},
  position={\nomenc@position}]{csm_\nomenc@label.tex}
}

\makeatother
\ExplSyntaxOff

\ifmaketwocolcsm
\newenvironment{authorbio}[1][]{
  {\bfseries\itshape #1}~}{}
\else
\newenvironment{authorbio}[1][]{
  #1~}{}
\fi

\defbibenvironment{bibliography}
{\list
  {\printtext[labelnumberwidth]{%
      \printfield{labelprefix}%
      \printfield{labelnumber}}}
  {\setlength{\labelwidth}{\labelnumberwidth}%
    \setlength{\leftmargin}{0pt}
    \setlength{\labelsep}{\biblabelsep}%
    \setlength{\itemsep}{\bibitemsep}%
    \setlength{\parsep}{\bibparsep}}%
  }
{\endlist}
{\item}

\newcommand{\csmprintmainbib}{%
  \begin{refcontext}[sorting=none]%
    \renewcommand*{\bibfont}{\footnotesize}%
    \DeclareFieldFormat{labelnumberwidth}{[##1]}%
    \printbibliography[notkeyword=sidebarbib, resetnumbers=true]%
  \end{refcontext}%
}

\definecolor{tay_col_b}{RGB}{0,8,126}       
\definecolor{tay_col_g}{RGB}{0,136,82}      
\definecolor{tay_col_r}{RGB}{235,61,37}     

\definecolor{beamerYellow}{HTML}{f1d46a}
\definecolor{beamerDarkYellow}{HTML}{bbb870}
\definecolor{beamerRed}{HTML}{db6245}
\definecolor{beamerBlue}{HTML}{356397}
\definecolor{beamerDarkBlue}{HTML}{26415d}
\definecolor{beamerGrey}{HTML}{acacac}
\definecolor{beamerGreen}{HTML}{5da9a1}
\definecolor{beamerBlack}{HTML}{000000}

\colorlet{lcvxColor}{beamerRed}

\tikzset{
  every node/.style={
    inner sep=0,
    outer sep=0,
    minimum size=0
  },
  inline box/.style={
    rectangle,
    inner xsep=1mm,
    font=\normalfont\sffamily\small,
    minimum height=1em
  },
  delim line/.style={
    line width=1pt
  }
}

\newcommand{\alglocation}[1]{%
  \tikz[baseline=-0.6ex]{\node[
    circle,
    draw=black,
    minimum size=0.6em,
    inner sep=0.2mm,
    scale=0.8
    ] at (0,0) {#1};}}

\def\iterstartloc{1}
\def\solveloc{2}
\def\testloc{3}

\newcommand{\alglayer}[2]{%
  \tikz[baseline=-0.6ex]{\node[
    circle,
    draw=#1,
    minimum size=0.6em,
    inner sep=0.2mm,
    scale=0.8
    ] at (0,0) {#2};}}

\usepackage{enumitem}

\newcommand{\runningk}{\zeta}

\newcommand{\lcvxresultref}[2][singular]{\lcvx~\ifthenelse{\equal{#1}{singular}}{Result}{Results}~\ref{#2}}
\newcommand{\csmnameref}[1]{\hyperlink{sidebar_#1}{\csname #1\endcsname}}

\makeatletter
\@ifclassloaded{standalone}{}{%

	\newmdenv[
	outerlinewidth = 1,%
	linewidth = 0pt,%
	backgroundcolor = yellow!40,%
	outerlinecolor = black!50,%
	innertopmargin = \topskip,%
	splittopskip = \topskip,%
	]{csm_sidebar}
}
\makeatother

\DeclareMathOperator{\cyclicshift}{\Psi_c}

\newcommand{\gusto}{GuSTO\xspace}

\newcommand{\runn}{\Gamma}
\newcommand{\term}{\phi}

\newcommand{\flow}{\psi}
\newcommand{\ti}{0}
\newcommand{\tio}{t_0}
\newcommand{\tf}[1][]{t_{f#1}}
\newcommand{\gic}{g_{\mathrm{ic}}}
\newcommand{\gtf}{g_{\mathrm{tc}}}
\newcommand{\indi}{w}
\newcommand{\Hf}{H_{\mathrm{f}}}
\newcommand{\Kf}{K_{\mathrm{f}}}
\newcommand{\ellf}{\ell_{\mathrm{f}}}
\newcommand{\dimgic}{n_{ic}}
\newcommand{\dimgtf}{n_{tc}}
\newcommand{\dimss}{n_{s}}
\newcommand{\dimindi}{n_{\indi}}

\newcommand{\dimvc}{n_{\vc}}

\newcommand{\trajohone}{\{\xb(t),\,\ub(t),\,\pb\}_0^1}

\newcommand{\ocost}{J}

\newcommand{\ocostL}{L}

\newcommand{\diffLx}{A^{\scriptscriptstyle \runn}}
\newcommand{\diffLu}{B^{\scriptscriptstyle \runn}}
\newcommand{\diffLp}{F^{\scriptscriptstyle \runn}}

\newcommand{\timeintvl}{\Delta t}

\newcommand{\xk}[1][]{x_{k#1}}
\newcommand{\xkp}{x_{k+1}}
\newcommand{\uk}[1][]{u_{k#1}}
\newcommand{\ukp}{u_{k+1}}
\newcommand{\pk}{p}
\newcommand{\xb}{\bar{x}}
\newcommand{\ub}{\bar{u}}
\newcommand{\pb}{\bar{p}}
\newcommand{\xbk}{\xb_k}
\newcommand{\ubk}{\ub_k}
\newcommand{\pbk}{\pb}
\newcommand{\vc}{\nu}

\newcommand{\vck}[1][]{\vc_{k#1}}

\newcommand{\vcc}[1][]{\nu_{s#1}}
\newcommand{\vcany}{\hat\nu}
\newcommand{\vccic}{\nu_{\mrm{ic}}}
\newcommand{\vcctf}{\nu_{\mrm{tc}}}
\newcommand{\vcck}[1][]{\vcc[,k#1]}

\newcommand{\defect}{\delta}
\newcommand{\defectk}[1][]{\defect_{k#1}}

\newcommand{\xxic}{x_{\mathrm{ic}}}
\newcommand{\xxfc}{x_{\mathrm{tc}}}
\newcommand{\uuic}{u_{\mathrm{ic}}}
\newcommand{\uufc}{u_{\mathrm{tc}}}

\newcommand{\tilt}{\theta}
\newcommand{\tiltmax}{\tilt_{\max}}

\newcommand{\Nobs}{n_{\text{obs}}}

\newcommand{\body}{\mathcal{B}}
\newcommand{\Fbody}{\mathcal{F}_{\body}}

\newcommand{\inertial}{\mathcal{I}}
\newcommand{\Finertial}{\mathcal{F}_{\inertial}}

\newcommand{\rI}[1][]{r_{\inertial\myappend{#1}}}
\newcommand{\drI}{\dot{r}_{\inertial}}
\newcommand{\ric}{r_0}
\newcommand{\rfc}{r_f}

\newcommand{\vI}{v_{\inertial}}
\newcommand{\dvI}{\dot{v}_{\inertial}}
\newcommand{\vic}{v_0}
\newcommand{\vfc}{v_f}

\newcommand{\q}{q}
\newcommand{\qIB}{\q_{\scriptscriptstyle \body\gets\inertial}}
\newcommand{\dqIB}{\dot{\q}_{\scriptscriptstyle \body\gets\inertial}}
\newcommand{\qic}{\q_0}
\newcommand{\qfc}{\q_f}

\newcommand{\wB}{\omega_{\body}}
\newcommand{\dwB}{\dot{\omega}_{\body}}

\newcommand{\gI}{g_{\inertial}}

\renewcommand{\tr}{\eta}
\newcommand{\trx}{\alpha_x}
\newcommand{\tru}{\alpha_u}
\newcommand{\trp}{\alpha_p}
\newcommand{\trgrow}{\beta_{\text{gr}}}
\newcommand{\trshrink}{\beta_{\text{sh}}}
\newcommand{\Jwgrow}{\gamma_{\textit{fail}}}

\newcommand{\rat}{\rho}

\newcommand{\mfnl}{\mathcal{J}}
\newcommand{\mfl}{\mathcal{L}}

\newcommand{\Jq}{S}
\newcommand{\Jl}{\ell}
\newcommand{\Jc}{g}
\newcommand{\Jpen}{\Jc_{\mrm{sp}}}
\newcommand{\Jtr}{\Jc_{\mrm{tr}}}

\newcommand{\Jw}{\lambda}

\newcommand{\runnw}[1][]{\runn_{\Jw#1}}
\newcommand{\termw}{\term_{\Jw}}
\newcommand{\ocostw}[1][]{\csname #1\endcsname{\ocost}_{\Jw}}
\newcommand{\ocostwCvx}{\ocostw[breve]}
\newcommand{\ocostwNCvx}{\ocostw[tilde]}

\newcommand{\ocostLw}[1][]{\csname #1\endcsname{\ocostL}_{\Jw}}
\newcommand{\ocostLwNCvx}{\ocostLw[tilde]}

\newcommand{\ocostLwrunn}[1][]{\ocostL_{\Jw#1}^\runn}
\newcommand{\ocostwrunn}[1][]{\ocost_{\Jw#1}^\runn}
\newcommand{\ocostLwrunnN}[1][]{\ocostL_{\Jw#1}^{\runn,N}}
\newcommand{\ocostwrunnN}[1][]{\ocost_{\Jw#1}^{\runn,N}}

\newcommand{\hpen}{h_{\Jw}}

\newcommand{\mflw}{\mfl_{\Jw}}
\newcommand{\mfnlw}{\mfnl_{\Jw}}

\newcommand{\shrinkrate}{\mu}
\newcommand{\shrinkk}{k_{*}}

\newcommand{\Obsi}[1][i]{\mathcal{O}_{#1}}
\newcommand{\ISS}{\mathrm{ss}}
\newcommand{\ObsISS}{\mathcal{O}_{\ISS}}
\newcommand{\sObsISS}{\tilde{\mathcal{O}}_{\ISS}}
\newcommand{\niss}{n_{\ISS}}

\newcommand{\diss}[1][]{d_{\ISS\myappend{#1}}}
\newcommand{\sdiss}[1][]{\tilde d_{\ISS\myappend{#1}}}

\newcommand{\sdissweight}{\varepsilon_{\ISS}}
\newcommand{\disslb}[1][]{\delta_{\ISS\myappend{#1}}}
\newcommand{\disslbvec}[1][]{\Delta_{\ISS\myappend{#1}}}

\newcommand{\tabs}{\mathsf{t}}

\newcommand{\CTNLconvexpath}{
  \optieqref{scp_gen_cont}{convex_path_constraints_X}\dash
  \optieqref{scp_gen_cont}{convex_path_constraints_U}\xspace}

\newcommand{\tdil}{\alpha_{\tabs}}

\graphicspath{{./figures/}}
\DeclareSourcemap{
  \maps[datatype=bibtex, overwrite]{
    \map{
      \perdatasource{sidebars/references.bib}
      \step[fieldset=KEYWORDS, fieldvalue={, }, appendstrict]
      \step[fieldset=KEYWORDS, fieldvalue=sidebarbib, append]
    }
  }
}

\makeatletter
\let\blx@rerun@biber\relax
\makeatother

\begin{document}

\input{sections/front_cover.tex}

\newcommand{\csmprintnomenclature}{%
  \csmnomenclature[%
  title={Abbreviations},
  width=0.4,
  label=nomencl]{sections/abbreviations.tex}
}

\newcommand{\csmprintsymbollist}{%
  \csmnomenclature[%
  title={Notation},
  width=0.4,
  label=symbols]{sections/symbols.tex}
}


\begin{csmabstract}
  Reliable and efficient trajectory generation methods are a fundamental need for
autonomous dynamical systems of tomorrow. The goal of this article is to
provide a comprehensive tutorial of three major convex optimization-based
trajectory generation methods: lossless convexification (\lcvx), and two
sequential convex programming algorithms known as \scvx and \gusto. In this
article, trajectory generation is the computation of a dynamically feasible
state and control signal that satisfies a set of constraints while optimizing
key mission objectives. The trajectory generation problem is almost always
nonconvex, which typically means that it is not readily amenable to efficient
and reliable solution onboard an autonomous vehicle. The three algorithms that
we discuss use problem reformulation and a systematic algorithmic strategy to
nonetheless solve nonconvex trajectory generation tasks through the use of a
convex optimizer. The theoretical guarantees and computational speed offered by
convex optimization have made the algorithms popular in both research and
industry circles. To date, the list of applications include rocket landing,
spacecraft hypersonic reentry, spacecraft rendezvous and docking, aerial motion
planning for fixed-wing and quadrotor vehicles, robot motion planning, and
more. Among these applications are high-profile rocket flights conducted by
organizations like NASA, Masten Space Systems, SpaceX, and Blue Origin. This
article aims to give the reader the tools and understanding necessary to work
with each algorithm, and to know what each method can and cannot do. A publicly
available source code repository supports the numerical examples provided at
the end of this article. By the end of the article, the reader should be ready
to use each method, to extend them, and to contribute to their many exciting
modern applications.


\end{csmabstract}

\firstword[nindent=0pt,slope=4pt,findent=-8pt]{Autonomous} vehicles and robots
promise many exciting new applications that will transform our society. For
example, autonomous aerial vehicles (AAVs) operating in urban environments
could deliver commercial goods, emergency medical supplies, monitor traffic,
and provide threat alerts for national security \cite{DAndrea2014}. At the same
time, these applications present significant engineering challenges for
performance, trustworthiness, and safety. For instance, AAVs can be a
catastrophic safety hazard should they lose control or situational awareness
over a populated area. Space missions, self\dash driving cars, and applications
of autonomous underwater vehicles (AUVs) share similar concerns.

Generating a trajectory autonomously onboard the vehicle is not only desirable
for many of these applications but also indeed a necessity when considering the
deployment of autonomous systems either in remote areas with little to no
communication, or at scale in a dynamic and uncertain world. For example, it is
not possible to remotely control a spacecraft during a Mars landing scenario
\cite{sanmartin2013development,Steltzner2014,way2007mars}, nor is it practical
to coordinate the motion of tens of thousands of delivery drones from a
centralized location \cite{DAndrea2014}.
In these and many other scenarios, individual vehicles must be endowed with
their own high quality decision making capability. Failure to generate a safe
trajectory can result in losing the vehicle, the payload, and even human life.
Reliable methods for trajectory generation are a fundamental need if we are to
maintain public trust and the high safety standard that we have come to expect
from autonomous or automatic systems \cite{stein2003respect}.

Computational resource requirements are the second major consideration for
onboard trajectory generation. Historically, this has been the driving factor
for much of practical algorithm development
\cite{MalyutaARC,MindellBook}. Although the modern consumer desktop is an
incredibly powerful machine, industrial central processing units (CPUs) can
still be relatively modest \cite{carson2019splice}. This is especially true for
spaceflight, where the harsh radiation environment of outer space prevents the
rapid adoption of new computing technologies. For example, NASA's flagship Mars
rover ``Perseverance'' landed in 2021 using a BAE RAD750 PowerPC flight
computer, which is a 20 year\dash old technology \cite{Dueri2017,
  Mars2020RoverBrains}. When we factor in the energy requirements of powerful
CPUs and the fact that trajectory generation is a small part of all the tasks
that an autonomous vehicle must perform, it becomes clear that modern
trajectory generation is still confined to a small computational
footprint. Consequently, real\dash time onboard trajectory generation
algorithms must be computationally efficient.


We define trajectory generation to be the computation of a multi-dimensional
temporal state and control signal that satisfies a set of specifications, while
optimizing key mission objectives. This article is concerned exclusively with
dynamically feasible trajectories, which are those that respect the equations
of motion of the vehicle under consideration. Although it is commonplace to
track dynamically \textit{in}feasible trajectories using feedback control, at
the end of the day a system can only evolve along dynamically feasible paths,
whether those are computed upfront during trajectory generation or are the
result of feedback tracking. Performing dynamically feasible trajectory
generation carries two important advantages. First, it provides a method to
systematically satisfy constraints (i.e., specifications) that are hard, if not
impossible, to meet through feedback control. This includes, for example,
translation-attitude coupled sensor pointing constraints. Secondly, dynamically
feasible trajectories leave much less tracking error for feedback controllers
to ``clean up'', which usually means that tracking performance can be vastly
improved. These two advantages shall become more apparent throughout the
article.



Numerical optimization provides a systematic mathematical framework to specify
mission objectives as \alert{costs or rewards} to be optimized, and to enforce
state and control specifications as well as the equations of motion via
\alert{constraints}. As a result, we can express trajectory generation problems
as optimal control problems, which are infinite\dash dimensional optimization
problems over function spaces \cite{PontryaginBook,BerkovitzBook}. Since the
early 1960s, optimal control theory has proven to be extremely powerful
\cite{BettsBook,MalyutaARC}. Early developments were driven by aerospace
applications, where every gram of mass matters and trajectories are typically
sought to minimize fuel or some other mass\dash reducing metric (such as
aerodynamic load and thereby structural mass). This led to work on trajectory
algorithms for climbing aircraft, ascending and landing rockets, and spacecraft
orbit transfer, to name just a few
\cite{LawdenBook,MarecBook,BrysonBook,KirkBook,LonguskiBook,Meditch1964}. Following
these early applications, trajectory optimization problems have now been
formulated in many practical areas, including aerial, underwater, and space
vehicles, as well as for chemical processes
\cite{morari1999model,eaton1992model,campo1987robust}, building climate control
\cite{oldewurtel2012use,ma2011model,bengea2014implementation}, robotics, and
medicine \cite{BettsBook}, to mention a few.

At its core, optimization-based trajectory generation requires solving an
optimal control problem of the following form (at this point, we keep the
problem very general):
\begin{optimization}[
  variables={u,p,t_f},
  objective={J(x,u,p,t_f)},
  label={general_ocp}]%
  \optilabel{dynamics}
  & \dot x(t) = f(x,u,p,t), \\
  \optilabel{xu_constraints}
  & \pare[big]{x(t),u(t),p,t_f} \in \set C(t),~\forall t\in [0,t_f], \\
  \optilabel{bcs} & \pare[big]{x(0), p}\in\set X_0,~\pare[big]{x(t_f), p}\in\set X_f.
\end{optimization}

The cost \optiobjref{general_ocp} encodes the mission goal, the system dynamics are
modeled by the differential equation constraint \optieqref{general_ocp}{dynamics},
the state and control specifications are enforced through
\optieqref{general_ocp}{xu_constraints}, and the boundary conditions are fixed by
\optieqref{general_ocp}{bcs}. Note that \pref{general_ocp} makes a distinction between
a control vector $u(t)$, which is a temporal signal, and a so-called parameter
vector $p$, which is a static variable that encodes other decision variables
like temporal scaling.

In most cases, the solution of such an infinite\dash dimensional optimization
problem is neither available in closed form nor is it computationally tractable
to numerically compute (especially in real-time). Instead, different solution
methods have been proposed that typically share the following three main
components:
\begin{itemize}
\item Formulation: specification of how the functions $J$ and $f$, and the sets
  $\set C$, $\set X_0$, and $\set X_f$ are expressed mathematically;
\item Discretization: approximation of the infinite\Hyphdash dimensional state
  and control signal by a finite-dimensional set of basis functions;
\item Numerical optimization: iterative computation of an optimal solution of
  the discretized problem.
\end{itemize}

Choosing the most suitable combination of these three components is truly a
mathematical art form that is highly problem dependent, and not an established
plug-and-play process like least\dash squares regression
\cite{BoydConvexBook}. No single recipe or method is the best, and methods that
work well for some problems can fare much worse for others. Trajectory
algorithm design is replete with problem dependent tradeoffs in performance,
optimality, and robustness, among others. Still, optimization literature (to
which this article belongs) attempts to provide some formal guidance and
intuition about the process. For the rest of this article, we will be
interested in methods that solve the \alert{exact} form of \pref{general_ocp},
subject only to the initial approximation made by discretizing the
continuous-time dynamics.

Many excellent references discuss the discretization component
\cite{Betts1998,BettsBook,Kelly2017}. Once the problem is discretized, we can
employ numerical optimization methods to obtain a solution. This is where the
real technical challenges for reliable trajectory generation arise. Depending
on the problem formulation and the discretization method, the finite\dash
dimensional optimization problem that must be solved can end up being a
nonlinear (i.e., \alert{nonconvex}) optimization problem (NLP). However, NLP
optimization has high computational complexity, and there are no general
guarantees of either obtaining a solution or even certifying that a solution
does not exist
\cite{BoydConvexBook,NocedalBook,rockafellar1993watershed}. Hence, general NLP
methods may not be appropriate for control applications, since we need
deterministic guarantees for reliable trajectory generation in autonomous
control.


In contrast, if the discretized problem is convex, then it can be solved
reliably and with an efficiency that exceeds all other areas of numerical
optimization except for linear programming and least squares
\cite{BoydConvexBook,NocedalBook,WrightIPMBook,
  GillPracticalBook,RoosLPIPMBook,wright2005interior}. This is the key
motivation behind this article's focus on a convex optimization-based problem
formulations for trajectory generation. By building methods on top of convex
optimization, we are able to leverage fast iterative algorithms with
polynomial-time complexity \cite[Theorem~4.7]{Peng2002MaxItersGuarantee}. In
plain language, given any desired solution accuracy, a convex optimization
problem can be solved to within this accuracy in a predetermined number of
arithmetic operations that is a polynomial function of the problem size. Hence,
there is a deterministic bound on how much computation is needed to solve a
given convex optimization problem, and the number of iterations cannot grow
indefinitely like for NLP. 
These special properties, together with a mature theoretical understanding and
an ever\dash growing list of successful real-world use\dash cases, leave little
doubt that convex optimization-based solution methods are uniquely well-suited
to be used in the solution of optimal control problems.

Among the successful real-world use-cases of convex optimization-based
trajectory generation are several inspiring examples from the aerospace
domain. These include autonomous drones, spacecraft rendezvous and docking, and
most notably planetary landing. The latter came into high profile 
as an enabling technology for reusability of the SpaceX Falcon 9 and Heavy
rockets \cite{lars2016autonomous}. Even earlier, the NASA Jet Propulsion
Laboratory (JPL) demonstrated the use of a similar method for Mars pinpoint
landing aboard the Masten Xombie sounding rocket
\cite{Scharf2017,flights2012a,flights2012b}. Today, these methods are being
studied and adopted for several Mars, Moon, and Earth landing applications
\cite{lars2016autonomous,carson2019splice}. Although each application has its
own set of unique challenges, they all share the need to use the full
spacecraft motion envelope with limited sensing, actuation, and fuel/power
\cite{carson2019splice}. These considerations are not unique to space
applications, and can be found in almost all autonomous vehicles such as cars
\cite{Paden2016,Buehler2009DARPA}, walking robots
\cite{Buchanan2020,Sleiman2021}, and quadrotors \cite{Kaufmann2018Drone}.

Having motivated the use of convex optimization, we note that many trajectory
generation problems have common sources of nonconvexity, among which are:
nonconvex control constraints, nonconvex or coupled state-control constraints,
and nonlinear dynamics. The goal of a convex optimization-based trajectory
generation algorithm is to provide a systematic way of handling these
nonconvexities, and to generate a trajectory using a convex solver at its core.

Two methods stand out to achieve this goal. In special cases, it is possible to
reformulate the problem into a convex one through a variable substitution and a
``lifting'' (i.e., augmentation) of the control input into a higher-dimensional
space. In this case, Pontryagin's maximum principle \cite{PontryaginBook} can be
used to show that solving the new problem recovers a globally optimal solution
of the original problem. This gives the approach the name \alert{lossless
  convexification} (\lcvx), and the resulting problem can often be solved with
a single call to a convex solver. As one can imagine, however, \lcvx tends to
apply only to very specific problems. Fortunately, this includes some important
and practically useful forms of rocket landing and other trajectory generation
problems for spacecraft and AAV vehicles
\cite{Acikmese2007,Blackmore2012,MalyutaARC}. When \lcvx cannot be used, convex
optimization can be applied via \alert{sequential convex programming}
(SCP). This natural extension linearizes all nonconvex elements of
\pref{general_ocp}, and solves the convex problem in a local neighborhood where the
linearization is accurate. Roughly speaking, the problem is then re-linearized
at the new solution, and the whole process is repeated until a stopping
criterion is met. In the overall classification of optimization algorithms, SCP
is a trust region method
\cite{ConnTrustRegionBook,NocedalBook,Kochenderfer2019}. While SCP is a whole
class of algorithms, our primary focus is on two particular and closely related
methods called \scvx (also known as successive convexification) and \gusto
\cite{Mao2018,Bonalli2019a}.

Let us go through the numerous applications where the \lcvx, \scvx and \gusto
methods have been used. The \lcvx method was originally developed for rocket
landing \cite{Acikmese2007,harris2013maximum}. This was the method at the center
of the aforementioned NASA JPL multi-year flight test campaign for Mars
pinpoint landing \cite{Scharf2017}. We hypothesize that \lcvx is also a close
relative of the SpaceX Falcon 9 terminal landing algorithm
\cite{lars2016autonomous}. The method also appears in applications for
fixed-wing and quadrotor AAV trajectory generation
\cite{Blackmore2012,Szmuk2017}, spacecraft hypersonic reentry
\cite{Liu2016,Wang2017,Liu2019}, and spacecraft rendezvous and docking
\cite{Harris2014JGCD,MalyutaLCvxRejected}. The \scvx method, which applies to
far more general problems albeit with fewer runtime guarantees, has been used
extensively for the general rocket landing problem
\cite{Szmuk2016,SzmukReynolds2018,szmuk2018successive-conf,
  szmuk2019successive,Reynolds2019b,Reynolds2020}, quadrotor trajectory
generation \cite{Szmuk2017,Mao2017,szmuk2019real}, spacecraft rendezvous and
docking \cite{malyuta2020fast}, and cubesat attitude control
\cite{Reynolds2021}. Recently, as part of the NASA SPLICE project to develop a
next-generation planetary landing computer \cite{carson2019splice}, the \scvx
algorithm is being tested as an experimental payload aboard the Blue Origin New
Shepard rocket \cite{TippingPointFlight}. The \gusto method has been applied to
free-flyer robots such as those used aboard the international space station
(ISS) \cite{Bonalli2019a,Bonalli2019b,BanerjeeEtAl2020,LewBonalliECC2020}, car
motion planning \cite{BonalliLewTAC2021,BonalliLewSICON2021}, aircraft motion
planning \cite{Bonalli2019a}, and robot manipulator arms \cite{Bonalli2019b}.

\begin{csmfigure}[%
  caption={%
    A typical control architecture consists of trajectory generation and
    feedback control elements. This article discusses algorithms for trajectory
    generation, which traditionally provides reference and feedforward control
    signals. By repeatedly generating new trajectories, a feedback action is
    created that can itself be used for control. Repeated trajectory generation
    for feedback control underlies the theory of model predictive control.
  },%
  label={control_architecture}]%
  \centering%
  \includegraphics[width=1\columnwidth]{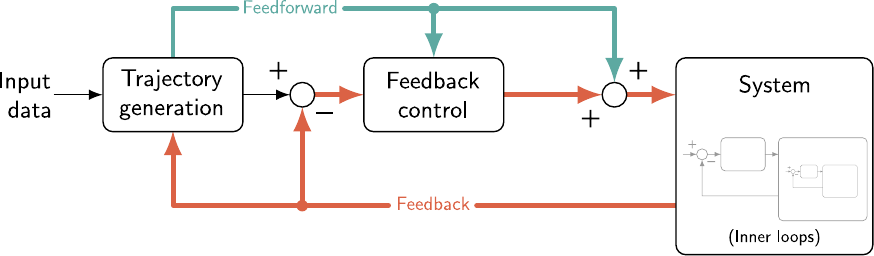}
\end{csmfigure}

This article provides a first\dash ever comprehensive tutorial of the \lcvx,
\scvx and \gusto algorithms. Placing these related methods under the umbrella
of a single article allows to provide a unified description that highlights
common ideas and helps the practitioner to know how, where, and when to deploy
each method. Previous tutorials on \lcvx and \scvx provide a complementary
technical discussion \cite{MaoTutorialACC, MaoTutorialMPC}.

There are two reasons for focusing on \lcvx, \scvx, and \gusto
specifically. First, the authors are the developers of the three
algorithms. Hence, we feel best positioned to provide a thorough description
for these particular methods, given our experience with their
implementation. Second, these three methods have a significant history of
real-world application. This should provide confidence that the methods
withstood the test of time, and have proven themselves to be useful when the
stakes were high. 
By the end of the article, our hope is to have provided the understanding and
the tools necessary in order to adapt each method to the reader's particular
engineering application.

Although our discussion for SCP is restricted to \scvx and \gusto, both methods
are closely related to other existing SCP algorithms. We hope that after
reading this tutorial, the reader will be well-positioned to understand most if
not all other SCP methods for trajectory generation. Applications of these SCP
alternatives are discussed in the recent survey paper \cite{MalyutaARC}.

Finally, we note that this article is focused on solving a trajectory
optimization problem like \pref{general_ocp} \textit{once} in real\dash time. As
illustrated in \figref{control_architecture}, this results in a single optimal
trajectory that can be robustly tracked by a downstream control system. The
ability to solve for the trajectory in real\dash time, however, can allow for
updating the trajectory as the mission evolves and more information is revealed
to the autonomous vehicle. 
Repetitive trajectory generation provides a feedback action that can itself be
used for control purposes. This approach is the driving force behind model
predictive control (MPC), which has been applied to many application domains
over the last three decades
\cite{mayne2014model,garcia1989model,eren2017model}. This article does not cover
MPC, and we refer the interested reader to existing literature
\cite{RawlingsMPCBook,eren2017model}.

\begin{csmfigure}[%
  caption={%
    The complete implementation source code for the numerical examples at the
    end of this article can be found in the \texttt{csm} branch of our
    open-source GitHub repository. The \texttt{master} branch provides even
    more algorithms and examples that are not covered in this article.
  },%
  label={github_qr},
  position={!b}]%
  \centering%
  \ifmaketwocolcsm




\makeatletter

\def\@sz{3cm}

\tikzset{
  qr/.style={
    rectangle,
    minimum size=\@sz
  },
  url/.style={
    anchor=north,
    scale=1.0,
    outer sep=2mm,
    font=\tt
  }
}

\begin{tikzpicture}[line join=round, line cap=round]

  \node[qr] (code) at (0,0) {\includegraphics[width=\@sz]{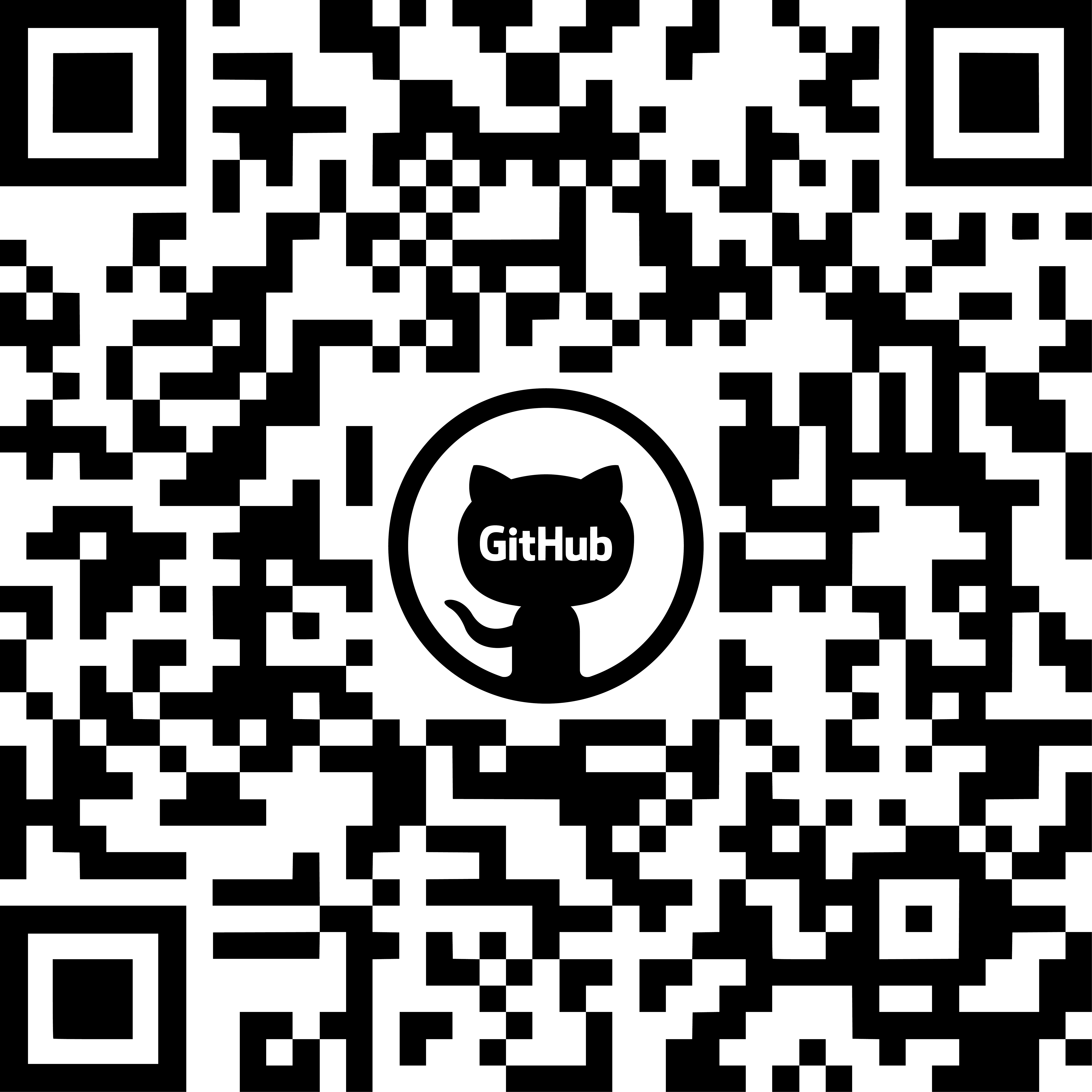}};

  \node[url] at (code.south)
  {\href{https://www.github.com/dmalyuta/scp\_traj\_opt/tree/csm}%
    {\color{beamerBlue} github.com/dmalyuta/scp\_traj\_opt/tree/csm}};

\end{tikzpicture}

\makeatother



  \else
  \includegraphics[scale=\csmpreprintfigscale]{qr_code}
  \fi
\end{csmfigure}

The rest of this article is organized as follows. We begin with a short section
on convex optimization, with the primary objective of highlighting why it is so
useful for trajectory generation. The article is then split into three main
parts. Part I surveys the major results of lossless convexification (\lcvx) to
solve nonconvex trajectory problems in one shot. Part II discusses sequential
convex programming (SCP) which can handle very general and highly nonconvex
trajectory generation tasks by iteratively solving a number of convex
optimization problems. In particular, Part II provides a detailed tutorial on
two modern SCP methods called \scvx and \gusto
\cite{Mao2018,Bonalli2019a}. Lastly, Part III applies \lcvx, \scvx, and \gusto
to three complex trajectory generation problems: a rocket\dash powered
planetary lander, a quadrotor, and a microgravity free\dash flying robotic
assistant. Some important naming conventions and notation that we use
throughout the article are defined in the \sbref{nomencl} and \sbref{symbols}.

To complement the tutorial style of this article, the numerical examples in
Part III are accompanied by open-source implementation code linked in
\figref{github_qr}. We use the Julia programming language because it is simple to
read like Python, yet it can be as fast as C/C++ \cite{bezanson2017julia}. By
downloading and running the code, the reader can recreate the exact plots seen
in this article.

\csmprintnomenclature

\csmprintsymbollist


\subsection{Convex Optimization Background}

Convex optimization seeks to minimize a convex objective function while
satisfying a set of convex constraints. The technique is expressive enough to
capture many trajectory generation and control applications, and is appealing
due to the availability of solution algorithms with the following properties
\cite{NocedalBook,BoydConvexBook}:
\begin{itemize}
\item A {globally optimal} solution is found if a feasible solution exists;
\item A {certificate of infeasibility} is provided when a feasible solution
  does not exist;
\item The runtime complexity is {polynomial} in the problem size;
\item The algorithms can {self-initialize}, eliminating the need for an expert
  initial guess.
\end{itemize}

The above properties are fairly general and apply to most, if not all,
trajectory generation and control applications. This makes convex programming
safer than other optimization methods for autonomous applications.

To appreciate what makes an optimization problem convex, we introduce some
basic definitions here and defer to \cite{BoydConvexBook,RockafellarConvexBook}
for further details. Two fundamental objects must be considered: a convex
function and a convex set. For reference, \figref{convex_example} illustrates a
notional convex set and function. By definition, $\mathcal C\subseteq\real^n$
is a convex set if and only if it contains the line segment connecting any two
of its points:
\begin{equation}
  \label{eq:generic_convex_set}
  x,y\in\mathcal C\Rightarrow [x,y]_\theta\in\mathcal C
\end{equation}
for all $\theta\in[0,1]$, where
$[x,y]_\theta\definedas\theta x+(1-\theta)y$. An important property is that
convexity is preserved under set intersection. This allows us to build
complicated convex sets by intersecting simpler sets. By replacing the word
``sets'' with ``constraints'', we can readily appreciate how this property
plays into modeling trajectory generation problems using convex optimization.

\begin{csmfigure}[%
  position=t,
  caption={%
    Illustration of a notional convex set (\protect\subref{fig:convex_example_set})
    and convex function (\protect\subref{fig:convex_example_function}). In both
    cases, the variable $\theta\in [0,1]$ generates a line segment between two
    points. The epigraph $\epi f\subseteq\real^n\times\real$ is the set of
    points which lie above the function, and itself defines a convex set.
  },
  label={convex_example},
  position={!t}]%
  \centering
  \begin{subfigure}[b]{1.0\linewidth}
    \centering
    \includegraphics[width=0.9\columnwidth]{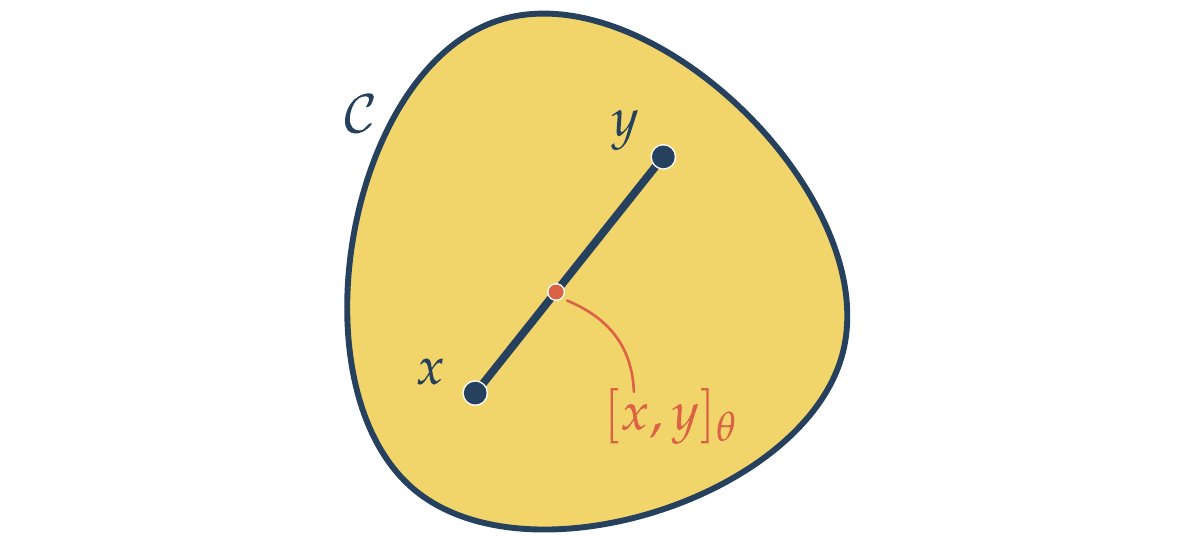}
    \caption{A convex set contains all line segments connecting its points.}
    \label{fig:convex_example_set}
  \end{subfigure}
  \vspace{3mm}

  \begin{subfigure}[b]{1.0\linewidth}
    \centering
    \includegraphics[width=0.9\columnwidth]{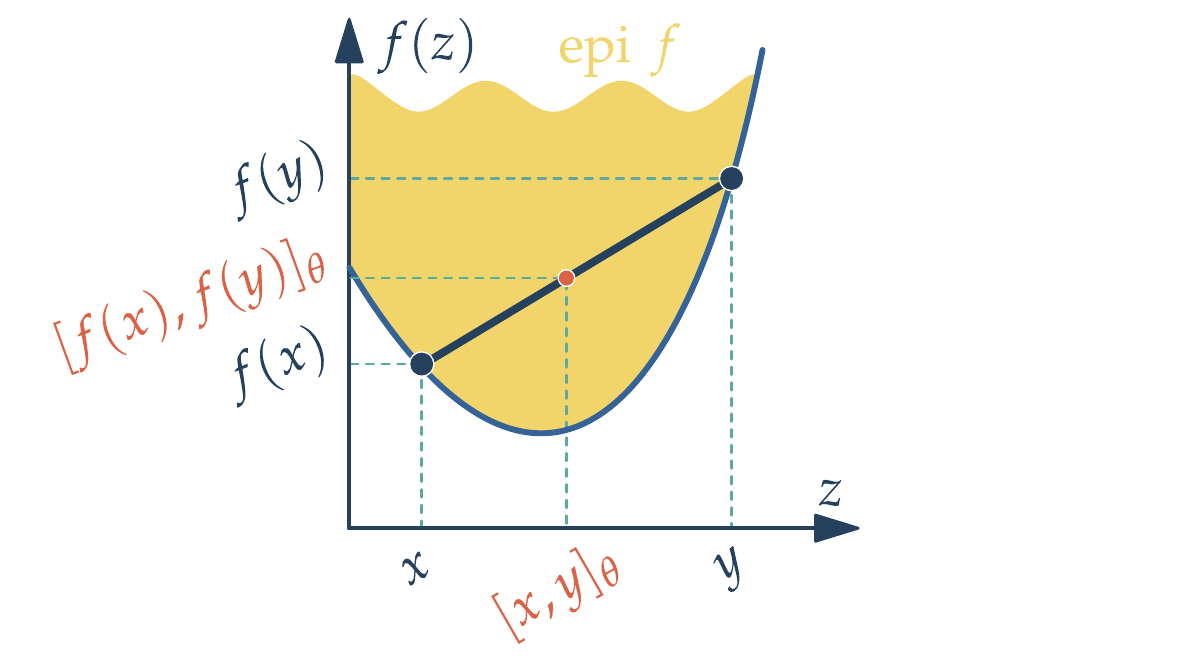}
    \caption{A convex function lies below all line segments connecting its
      points.}
    \label{fig:convex_example_function}
  \end{subfigure}
\end{csmfigure}

A function $f:\real^n\to\real$ is convex if and only if $\dom f$ is a convex
set and $f$ lies below the line segment connecting any two of its points:
\begin{equation}
  x,y\in\dom f\Rightarrow f([x,y]_\theta)\le [f(x),f(y)]_\theta
\end{equation}
for all $\theta\in[0,1]$. A convex optimization problem is simply the
minimization of a convex function subject to a number of convex constraints
that act to restrict the search space:
\begin{optimization}[
  label={cvx},
  variables={x\in\reals^n},
  objective={f_0(x)}]
  \optilabel{inequalities}
  & f_i(x)\le 0,~i=1,\ldots,n_{\mrm{ineq}}, \\
  \optilabel{equalities}
  & g_i(x) = 0,~i=1,\ldots,n_{\mrm{eq}},
\end{optimization}
where~$f_0:\real^n\to\real$ is a convex \alert{cost} function,
$f_i:\real^n\to\real$ are convex inequality constraints, and
$g_i:\real^n\to\real$ are affine equality constraints. The problem contains
$n_{\mrm{ineq}}$ inequality and $n_{\mrm{eq}}$ equality constraints. We stress
that the equality constraints must be affine, which means that each function
$g_i$ is a linear expression in $x$ plus a constant offset. The equations of
motion are equality constraints, therefore basic convex optimization restricts
the dynamics to be affine (i.e., linear time-varying at most). Handling
nonlinear dynamics will be a major topic of discussion throughout this article.

Each constraint defines a convex set so that, together,
\optieqref{cvx}{inequalities} and \optieqref{cvx}{equalities} form a convex feasible
set of values that the decision variable $x$ may take. To explicitly connect
this discussion back to the generic convex set introduced in
\eqref{eq:generic_convex_set}, we can write the feasible set as:
\begin{align}
  \notag
  \set C=\bigg\{
  x\in\reals^n~\where~
  &f_i(x)\le 0,~i=1,\ldots,n_{\mrm{ineq}}, \\
  \label{eq:convex_feasible_set}
  &g_i(x) = 0,~i=1,\ldots,n_{\mrm{eq}}\bigg\}.
\end{align}

A fundamental consequence of convexity is that any local minimum of a convex
function is a global minimum \cite[Section~4.2.2]{BoydConvexBook}. More
generally, convex functions come with a plethora of properties that allow
algorithm designers to obtain global information about function behavior from
local measurements. For example, a differentiable convex function is globally
lower bounded by its local first\dash order approximation
\cite{BoydConvexBook}. Thus, we may look at convexity as a highly beneficial
assumption on function behavior that enables efficient algorithm
design. Indeed, a landmark discovery of the twentieth century was that it is
convexity, not linearity, that separates ``hard'' and ``easy'' problems
\cite{rockafellar1993watershed}.

For practitioners, the utility of convex optimization stems not so much from
the ability to find the global minimum, but rather from the ability to find it
(or indeed any other feasible solution) \textit{quickly}. The field of
numerical convex optimization was invigorated by the interior-point method
(IPM) family of optimization algorithms, first introduced in 1984 by Karmarkar
\cite{Karmarkar1984}. Today, convex optimization problems can be solved by
primal\dash dual IPMs in a few tens of iterations
\cite{BoydConvexBook,WrightIPMBook}. Roughly speaking, we can say that
substantially large trajectory generation problems can usually be solved in
under one second \cite{Reynolds2020,Dueri2017,Scharf2017}. In technical
parlance, IPMs have a favorable polynomial problem complexity: the number of
iterations required to solve the problem to a given tolerance grows
polynomially in the number of constraints $n_{\mrm{ineq}}+n_{\mrm{eq}}$. With
some further assumptions, it is even possible to provide an upper bound on the
number of iterations \cite{PengSelfRegularityBook,Peng2002MaxItersGuarantee}.
We defer to \cite[Chapters~14 and 19]{NocedalBook} for futher details on convex
optimization algorithms. Throughout this article, our goal will be to leverage
existing convex problem solvers to create higher-level frameworks for the
solution of trajectory generation problems.


\section{Part I: Lossless Convexification}

\alert{Lossless convexification} (\lcvx) is a modeling technique that solves
nonconvex optimal control problems through a convex relaxation. In this method,
Pontryagin's maximum principle \cite{PontryaginBook} is used to show that a
convex relaxation of a nonconvex problem finds the globally optimal solution to
the original problem, hence the method's name. To date, the method has been
extended as far as relaxing certain classes of nonconvex control constraints,
such as an input norm lower bound (see \sbref{inputrelax}) and a nonconvex pointing
constraint (see \sbref{pointingrelax}).

The \lcvx method has been shown to work for a large class of state-constrained
optimal control problems, however a working assumption is that state
constraints are convex. Lossless relaxation of nonconvex state constraints
remains under active research, and some related results are available
\cite{Blackmore2012} (which we will cover in this section).

As the reader goes through Part I, it is suggested to keep in mind that the
main concerns of lossless convexification are:
\begin{itemize}
\item To find a convex lifting of the feasible input set;
\item To show that the optimal input of the lifted problem projects back to a
  feasible input of the original non-lifted problem.
\end{itemize}

\figref{lcvx_history} chronicles the development history of \lcvx. The aim of this
part of the article is to provide a tutorial overview of the key results, so
theoretical proofs are omitted in favor of a more practical and action\dash
oriented description. Ultimately, our aim is for the reader to come away with a
clear understanding of how \lcvx can be applied to their own problems.

In the following sections, we begin by introducing \lcvx for the input norm
lower bound and pointing nonconvexities using a baseline problem with no state
constraints. Then, the method is extended to handle affine and quadratic state
constraints, followed by general convex state constraints. We then describe how
the \lcvx method can also handle a class of dynamical systems with nonlinear
dynamics. For more general applications, embedding \lcvx into nonlinear
optimization algorithms is also discussed. At the very end of this part of the
article, we cover some of the newest \lcvx results from the past year, and
provide a toy example that gives a first taste of how \lcvx can be used in
practice.

\begin{csmfigure}[
  caption={%
    Chronology of lossless convexification theory development. Note the
    progression from state\dash unconstrained problems to those with
    progressively more general state constraints and, finally, to problems that
    contain integer variables.
  },
  label={lcvx_history},
  position=t]%
  \centering
  \ifmaketwocolcsm
  \sffamily\sansmath
  \fi
  \begin{tikzpicture}

  \ifmaketwocolcsm
  \tikzmath{
    \cx=-2; \cy=5;
    \dy=0.5; \hy=18*\dy;
    \dx=0.625;
    \ccx=0.3;
    \DY = 0.8;
    \Y0 = \cy-0.1;
  }
  \else
  \tikzmath{
    \cx=-2; \cy=5;
    \dy=0.6; \hy=18*\dy;
    \dx=0.625;
    \ccx=0.3;
    \DY = 0.95;
    \Y0 = \cy-0.1;
  }
  \fi

  \def\whiskdxl{1mm}
  \def\whiskdxr{1mm}

  \tikzset{
    whisker tip/.append style={
      inner sep=0,
      outer sep=0,
      circle,
      fill=black,
      minimum size=1mm
    },
    connection/.append style={
      line width=0.3mm,
      postaction={
        decorate,
        decoration={
          markings,
          mark=at position 0 with {\node[whisker tip] {};},
          mark=at position 0.999 with {\node[whisker tip] {};}
        }
      }
    },
    left edge/.append style={
      black,
      line width=0.3mm
    },
    description/.append style={
      align=left,
      anchor=west,
      inner ysep=1pt,
      scale=0.8
    }
  }

  \draw[-latex] (\cx,\cy) -- ++(0,-\hy) node{\phantom{H}};
  \foreach \x in {1,...,17} \draw[] (\cx,\cy-\x*\dy) -- (\cx+0.15,\cy-\x*\dy);

  \node[align=center,name=L2005] at (\cx+\dx,\cy-\dy) {$2005$};
  \node[align=center,name=L2007] at (\cx+\dx,\cy-3*\dy) {$2007$};
  \node[align=center,name=L2010] at (\cx+\dx,\cy-6*\dy) {$2010$};
  \node[align=center,name=L2011] at (\cx+\dx,\cy-7*\dy) {$2011$};
  \node[align=center,name=L2012] at (\cx+\dx,\cy-8*\dy) {$2012$};
  \node[align=center,name=L2013] at (\cx+\dx,\cy-9*\dy) {$2013$};
  \node[align=center,name=L2014] at (\cx+\dx,\cy-10*\dy) {$2014$};
  \node[align=center,name=L2019] at (\cx+\dx,\cy-15*\dy) {$2019$};
  \node[align=center,name=L2020] at (\cx+\dx,\cy-16*\dy) {$2020$};
  \node[align=center,name=L2021] at (\cx+\dx,\cy-17*\dy) {$2021$};

  \node[name=T2005,description] at (\ccx,\Y0) {Input lower\\
    bound~\cite{Acikmese2005,Acikmese2007}};
  \newcommand{\leftedge}[1]{
    \draw[left edge] ($(#1.south west)+(-\whiskdxr,0)$) to
    ($(#1.north west)+(-\whiskdxr,0)$);
  }
  \leftedge{T2005}

  \node[name=T2010,description] at (\ccx,\Y0-\DY) {Minimum landing\\
    error~\cite{Blackmore2010}};
  \leftedge{T2010}

  \node[name=T2011a,description] at (\ccx,\Y0-2*\DY) {Generalized
    input\\ nonconvexity~\cite{Acikmese2011}};
  \leftedge{T2011a}

  \node[name=T2011b,description] at (\ccx,\Y0-3*\DY) {Active state
    constraints\\at discrete times~\cite{Acikmese2011}};
  \leftedge{T2011b}

  \node[name=T2011c,description] at (\ccx,\Y0-4*\DY) {Input pointing\\
    constraint~\cite{Carson2011}};
  \leftedge{T2011c}

  \node[name=T2012,description] at (\ccx,\Y0-5*\DY-0.25) {Fusion with\\
    mixed-integer\\ programming~\cite{Blackmore2012}};
  \leftedge{T2012}

  \node[name=T2013,description] at (\ccx,\Y0-6.7*\DY) {Persistently
    active\\ quadratic state\\ constraints~\cite{Harris2013b}};
  \leftedge{T2013}

  \node[name=T2014,description] at (\ccx,\Y0-8.13*\DY) {Persistently
    active\\linear state\\ constraints~\cite{Harris2013b,Harris2014}};
  \leftedge{T2014}

  \node[name=T2019,description] at (\ccx,\Y0-9.35*\DY) {Semi-continuous\\
    input norms~\cite{MalyutaIFAC}};
  \leftedge{T2019}

  \node[name=T2021a,description] at (\ccx,\Y0-10.15*\DY) {Disconnected sets~\cite{Harris2021}};
  \leftedge{T2021a}

  \node[name=T2021b,description] at (\ccx,\Y0-10.75*\DY) {Fixed\dash final time~\cite{Kunhippurayil2021FixedTime}};
  \leftedge{T2021b}

  \node[name=T2021c,description] at (\ccx,\Y0-11.35*\DY) {Mixed constraints~\cite{Kunhippurayil2021Observability}};
  \leftedge{T2021c}

  \draw[connection] ($(L2005.east)+(\whiskdxl,0)$) to[out=0,in=180,looseness=0.3] ($(T2005.west)+(-\whiskdxr,0)$);
  \draw[connection] ($(L2007.east)+(\whiskdxl,0)$) to[out=0,in=180,looseness=0.3] ($(T2005.west)+(-\whiskdxr,0)$);
  \draw[connection] ($(L2010.east)+(\whiskdxl,0)$) to[out=0,in=180,looseness=0.3] ($(T2010.west)+(-\whiskdxr,0)$);
  \draw[connection] ($(L2011.east)+(\whiskdxl,0)$) to[out=0,in=180,looseness=0.3] ($(T2011a.west)+(-\whiskdxr,0)$);
  \draw[connection] ($(L2011.east)+(\whiskdxl,0)$) to[out=0,in=180,looseness=0.3] ($(T2011b.west)+(-\whiskdxr,0)$);
  \draw[connection] ($(L2011.east)+(\whiskdxl,0)$) to[out=0,in=180,looseness=0.3] ($(T2011c.west)+(-\whiskdxr,0)$);
  \draw[connection] ($(L2012.east)+(\whiskdxl,0)$) to[out=0,in=180,looseness=0.3] ($(T2012.west)+(-\whiskdxr,0)$);
  \draw[connection] ($(L2013.east)+(\whiskdxl,0)$) to[out=0,in=180,looseness=0.3] ($(T2013.west)+(-\whiskdxr,0)$);
  \draw[connection] ($(L2013.east)+(\whiskdxl,0)$) to[out=0,in=180,looseness=0.3] ($(T2011c.west)+(-\whiskdxr,0)$);
  \draw[connection] ($(L2013.east)+(\whiskdxl,0)$) to[out=0,in=180,looseness=0.3] ($(T2014.west)+(-\whiskdxr,0)$);
  \draw[connection] ($(L2014.east)+(\whiskdxl,0)$) to[out=0,in=180,looseness=0.3] ($(T2014.west)+(-\whiskdxr,0)$);
  \draw[connection] ($(L2019.east)+(\whiskdxl,0)$) to[out=0,in=180,looseness=0.3] ($(T2019.west)+(-\whiskdxr,0)$);
  \draw[connection] ($(L2020.east)+(\whiskdxl,0)$) to[out=0,in=180,looseness=0.3] ($(T2019.west)+(-\whiskdxr,0)$);
  \draw[connection] ($(L2021.east)+(\whiskdxl,0)$) to[out=0,in=180,looseness=0.3] ($(T2021a.west)+(-\whiskdxr,0)$);
  \draw[connection] ($(L2021.east)+(\whiskdxl,0)$) to[out=0,in=180,looseness=0.3] ($(T2021b.west)+(-\whiskdxr,0)$);
  \draw[connection] ($(L2021.east)+(\whiskdxl,0)$) to[out=0,in=180,looseness=0.3] ($(T2021c.west)+(-\whiskdxr,0)$);

\end{tikzpicture}


\end{csmfigure}

\subsection{No State Constraints}

We begin by stating perhaps the simplest optimal control problem for which an
\lcvx result is available. Its salient features are a distinct absence of state
constraints (except for the boundary conditions), and its only source of
nonconvexity is a lower-bound on the input given by
\optieqref{lcvx_o_nostate}{bounds}. A detailed description of \lcvx for this
problem may be found in \cite{Acikmese2005,Acikmese2007,Acikmese2011}.

\begin{optimization}[
  label={lcvx_o_nostate},
  variables={u,t_f},
  objective={m(t_f,x(t_f))+\runningk \int_0^{t_f} \ell(g_1(u(t)))\,\dt}]
  \optilabel{dynamics}
  & \dot x(t) = A(t)x(t)+B(t)u(t)+E(t)w(t), \\
  \optilabel{bounds}
  & \rho_{\min}\le g_1(u(t)),~g_0(u(t))\le\rho_{\max}, \\
  \optilabel{boundary}
  & x(0) = x_0,~b(t_f,x(t_f))=0.
\end{optimization}

In \pref{lcvx_o_nostate}, $t_f>0$ is the terminal time, $x(\cdot)\in\real^n$ is the
state trajectory, $u(\cdot)\in\real^m$ is the input trajectory,
$w(\cdot)\in\real^p$ is an exogenous additive disturbance,
$m:\real\times\real^n\to\real$ is a convex terminal cost, $\ell:\real\to\real$
is a convex and non-decreasing running cost modifier, $\runningk\in\{0,1\}$ is
a fixed user-chosen parameter to toggle the running cost,
$g_0,g_1:\real^m\to\real_+$ are convex functions, $\rho_{\min}>0$ and
$\rho_{\max}>\rho_{\min}$ are user-chosen bounds, and
$b:\real\times\real^n\to\real^{n_b}$ is an affine terminal constraint
function. Note that the dynamics in \pref{lcvx_o_nostate} define a linear time
varying (LTV) system. For all the \lcvx discussion that follows, we will also
make the following two assumptions on the problem data.

\begin{assumption}[terminal_state_not_overconstrained]
  If any part of the terminal state is constrained then the Jacobian
  $\grad_x b[t_f]\in\reals^{n_b\times n}$ is full row rank, i.e.,
  $\rank{\grad_x b[t_f]}=n_b$. This implies that the terminal state is not
  over-constrained.
\end{assumption}

\begin{assumption}[position_running_cost]
  The running cost is positive definite, in other words $\ell(z)>0$ for all
  $z\ne 0$.
\end{assumption}

When faced with a nonconvex problem like \pref{lcvx_o_nostate}, an engineer has two
choices. Either devise a nonlinear optimization algorithm, or solve a simpler
problem that is convex. The mantra of \lcvx is to take the latter approach by
``relaxing'' the problem until it is convex. In the case of \pref{lcvx_o_nostate},
let us propose the following relaxation, which introduces a new variable
$\sigma(\cdot)\in\real$ called the \alert{slack input}.

\begin{optimization}[
  label={lcvx_r_nostate},
  variables={\sigma,u,t_f},
  objective={m(t_f,x(t_f))+\runningk \int_0^{t_f} \ell(\sigma(t))\,\dt}]
  \optilabel{dynamics}
  & \dot x(t) = A(t)x(t)+B(t)u(t)+E(t)w(t), \\
  \optilabel{bounds}
  & \rho_{\min}\le \sigma(t),~g_0(u(t))\le\rho_{\max}, \\
  \optilabel{lcvx_equality}
  & {\color{lcvxColor}g_1(u(t))\le\sigma(t)}, \\
  \optilabel{boundary}
  & x(0) = x_0,~b(t_f,x(t_f))=0.
\end{optimization}

\declaresidebar[%
title={Convex Relaxation of an Input Lower Bound},
label={inputrelax}]{%
  sidebars/lcvx_lower_bound.tex}

The relaxation of the nonconvex input constraint \optieqref{lcvx_o_nostate}{bounds}
to the convex constraints
\optieqref{lcvx_r_nostate}{bounds}-\optieqref{lcvx_r_nostate}{lcvx_equality} is
illustrated in \sbref{inputrelax} for the case of a throttleable rocket
engine. Note that if \optieqref{lcvx_r_nostate}{lcvx_equality} is replaced with
equality, then \pref{lcvx_r_nostate} is equivalent to \pref{lcvx_o_nostate}. Indeed,
the entire goal of \lcvx is to \textit{prove} that
\optieqref{lcvx_r_nostate}{lcvx_equality} holds with equality at the globally
optimal solution of \pref{lcvx_r_nostate}. Because of the clear importance of
constraint \optieqref{lcvx_r_nostate}{lcvx_equality} to the \lcvx method, we shall
call it the \alert{\lcvx equality constraint}. This special constraint will be
highlighted in red in all subsequent \lcvx optimization problems.

Consider now the following set of conditions, which arise naturally when using
the maximum principle to prove lossless convexification. The theoretical details
are provided in \cite[Theorem~2]{Acikmese2011}.

\begin{condition}[lcvx_nostate_controllability]
  The pair $\{A(\cdot),B(\cdot)\}$ must be totally controllable. This means
  that any initial state can be transferred to any final state by a bounded
  control trajectory in any finite time interval $[0,t_f]$
  \cite{Blackmore2012,DAngelo1970}. If the system is time invariant, this is
  equivalent to $\{A,B\}$ being controllable, and can be verified by checking
  that the controllability matrix is full rank or, more robustly, via the PBH
  test \cite{Antsaklis2006}.
\end{condition}

\begin{condition}[lcvx_nostate_linindep]
  Define the quantities:
  \begin{subequations}
    \label{eq:lcvx_nostate_linindep_m_B}
    \begin{align}
      m_{\text{\lcvx}}
      &= \Matrix{\grad_x m[t_f] \\ \grad_t m[t_f]+\runningk\ell(\sigma(t_f))}%
      \in\real^{n+1}, \\
      \label{eq:lcvx_nostate_linindep_B}
      B_{\text{\lcvx}}
      &= \Matrix{\grad_x b[t_f]\T \\ \grad_t b[t_f]\T}%
      \in\real^{(n+1)\times n_b}.
    \end{align}
  \end{subequations}

  The vector $m_{\text{\lcvx}}$ and the columns of $B_{\text{\lcvx}}$ must be
  linearly independent.
\end{condition}

We can now state \tref{nostate}, which is the main lossless convexification
result for \pref{lcvx_o_nostate}. The practical implication of \tref{nostate} is
that the solution of \pref{lcvx_o_nostate} can be found in polynomial time by
solving \pref{lcvx_r_nostate} instead.

\begin{theorem}[nostate]
  The solution of \pref{lcvx_r_nostate} is globally optimal for
  \pref{lcvx_o_nostate} if
  \conref{lcvx_nostate_controllability,lcvx_nostate_linindep} hold.
\end{theorem}

There is a partial generalization of \tref{nostate}. First, restrict
\pref{lcvx_o_nostate} to the choices $\runningk=1$ and $g_0=g_1$. Next, introduce a
new pointing-like input constraint
\optieqref{lcvx_o_nostate_pointing}{pointing}. The quantities
$\hat n_u\in\real^m$ and $\hat n_g\in\real$ are user-chosen parameters. The new
problem takes the following form:
\begin{optimization}[
  label={lcvx_o_nostate_pointing},
  variables={u,t_f},
  objective={m(t_f,x(t_f))+\int_0^{t_f} \ell(g_0(u(t)))\,\dt}]
  & \dot x(t) = A(t)x(t)+B(t)u(t)+E(t)w(t), \\
  \optilabel{bounds}
  & \rho_{\min}\le g_0(u(t))\le\rho_{\max}, \\
  \optilabel{pointing}
  & \hat n_u\T u(t)\ge \hat n_g g_0(u(t)), \\
  & x(0) = x_0,~b(t_f,x(t_f))=0.
\end{optimization}

Using \optieqref{lcvx_o_nostate_pointing}{pointing}, one can, for example,
constrain an airborne vehicle's tilt angle. This constraint, however, is
nonconvex for $\hat n_g<0$. We take care of this nonconvexity, along with the
typical nonconvexity of the lower bound in
\optieqref{lcvx_o_nostate_pointing}{bounds}, by solving the following relaxation of
the original problem:
\begin{optimization}[
  label={lcvx_r_nostate_pointing},
  variables={\sigma,u,t_f},
  objective={m(t_f,x(t_f))+\int_0^{t_f} \ell(\sigma(t))\,\dt}]
  & \dot x(t) = A(t)x(t)+B(t)u(t)+E(t)w(t), \\
  & \rho_{\min}\le \sigma(t)\le\rho_{\max}, \\
  \optilabel{lcvx_equality}
  & {\color{lcvxColor}g_0(u(t))\le\sigma(t)}, \\
  \optilabel{pointing}
  & \hat n_u\T u(t)\ge \sigma(t) \hat n_g, \\
  & x(0) = x_0,~b(t_f,x(t_f))=0.
\end{optimization}

Just like in \pref{lcvx_r_nostate}, we introduced a slack input
$\sigma(\cdot)\in\reals$ to strategically remove nonconvexity. Note once again
the appearance of the \lcvx equality constraint
\optieqref{lcvx_r_nostate_pointing}{lcvx_equality}. Meanwhile, the relaxation of
\optieqref{lcvx_o_nostate_pointing}{pointing} to
\optieqref{lcvx_r_nostate_pointing}{pointing} corresponds to a halfspace input
constraint in the $(u,\sigma)\in\reals^{m+1}$ space. A geometric intuition
about the relaxation is illustrated in \sbref{pointingrelax} for a typical vehicle
tilt constraint. Lossless convexification can again be shown under an extra
\conref{lcvx_nostate_controllability_pointing}, yielding
\tref{nostate_pointing}. Theoretical details are provided in
\cite{Carson2011,acikmese2013lossless}.

\declaresidebar[%
title={Convex Relaxation of an Input Pointing Constraint},
label={pointingrelax}]{%
  sidebars/lcvx_pointing.tex}

\begin{condition}[lcvx_nostate_controllability_pointing]
  Let $N\in\real^{m\times (m-1)}$ be a matrix whose columns span the nullspace
  of $\hat n_u$ in \optieqref{lcvx_o_nostate_pointing}{pointing}. The pair
  $\{A(\cdot),B(\cdot)N\}$ must be totally controllable.
\end{condition}

\begin{theorem}[nostate_pointing]
  The solution of \pref{lcvx_r_nostate_pointing} is globally optimal for
  \pref{lcvx_o_nostate_pointing} if
  \conref{lcvx_nostate_controllability,lcvx_nostate_linindep,%
    lcvx_nostate_controllability_pointing} hold.
\end{theorem}

\subsection{Affine State Constraints}

The logical next step after \tref{nostate,nostate_pointing} is to ask whether
\pref{lcvx_o_nostate} can incorporate state constraints. It turns out that this
is possible under a fairly mild set of extra conditions. The results presented
in this section originate from \cite{Harris2014,Harris2013b,HarrisThesis}.

Affine inequality constraints are the simplest class of state constraints that
can be handled in \lcvx. The nonconvex statement of the original problem is a
close relative of \pref{lcvx_o_nostate}:

\begin{optimization}[
  label={lcvx_o_linstate},
  variables={u,t_f},
  objective={m(x(t_f))+\runningk \int_0^{t_f} \ell(g_1(u(t)))\,\dt}]
  \optilabel{dynamics}
  & \dot x(t) = Ax(t)+Bu(t)+Ew, \\
  \optilabel{bounds}
  & \rho_{\min}\le g_1(u(t)),~g_0(u(t))\le\rho_{\max}, \\
  \optilabel{affine_input}
  & Cu(t)\le c, \\
  \optilabel{affine_state}
  & Hx(t)\le h, \\
  \optilabel{boundary}
  & x(0) = x_0,~b(x(t_f))=0.
\end{optimization}

First, we note that \pref{lcvx_o_linstate} is \alert{autonomous}, in other words
the terminal cost in \optieqref{lcvx_o_linstate}{objective}, the dynamics
\optieqref{lcvx_o_linstate}{dynamics}, and the boundary constraint
\optieqref{lcvx_o_linstate}{boundary} are all independent of time. Note that a
limited form of time variance can still be included through an additional time
integrator state whose dynamics are $\dot z(t)=1$. The limitation here is that
time variance must not introduce nonconvexity in the cost, in the dynamics, and
in the terminal constraint \optieqref{lcvx_o_linstate}{boundary}. The matrix of
facet normals $H\in\real^{n_h\times n}$ and the vector of facet offsets
$h\in\real^{n_h}$ define a new polytopic (affine) state constraint set. A
practical use-case for this constraint is described in
\sbref{sidebar_affinestate}. Similarly, $C\in\real^{n_c\times m}$ and
$c\in\real^{n_c}$ define a new polytopic subset of the input constraint set, as
illustrated in \sbref{sidebar_inputcut}.


\declaresidebar[%
title={Landing Glideslope as an Affine State Constraint},
label={sidebar_affinestate}]{%
  sidebars/lcvx_glideslope.tex}

\declaresidebar[%
title={Using Halfspaces to Further Constrain the Input Set},
label={sidebar_inputcut}]{%
  sidebars/lcvx_input_cut.tex}

Let us propose the following convex relaxation of \pref{lcvx_o_linstate}, which
takes the familiar form of \pref{lcvx_r_nostate}:
\begin{optimization}[
  label={lcvx_r_linstate},
  variables={\sigma,u,t_f},
  objective={m(x(t_f))+\runningk \int_0^{t_f} \ell(\sigma(t))\,\dt}]
  & \dot x(t) = Ax(t)+Bu(t)+Ew, \\
  & \rho_{\min}\le \sigma(t),~g_0(u(t))\le\rho_{\max}, \\
  & {\color{lcvxColor}g_1(u(t))\le\sigma(t)}, \\
  & Cu(t)\le c, \\
  \optilabel{affine_state}
  & Hx(t)\le h, \\
  & x(0) = x_0,~b(x(t_f))=0.
\end{optimization}

To guarantee lossless convexification for this convex relaxation,
\conref{lcvx_nostate_controllability} can be modified to handle the new state and
input constraints \optieqref{lcvx_o_linstate}{affine_input} and
\optieqref{lcvx_o_linstate}{affine_state}. To this end, we use the following notion
of cyclic coordinates from mechanics \cite{Landau1969}.

\begin{definition}[cyclic]
  For a dynamical system $\dot x = f(x)$, with state $x\in\reals^n$, any
  components of $x$ that do not appear explicitly in $f(\cdot)$ are said to be
  \alert{cyclic coordinates}. Without loss of generality, we can decompose the
  state as follows:
  \begin{equation}
    \label{eq:state_cyclic}
    x = \Matrix{x_c \\ x_{nc}},
  \end{equation}
  where $x_c\in\reals^{n_c}$ are the cyclic and $x_{nc}\in\reals^{n_{nc}}$ are
  the non-cyclic coordinates, such that $n_c+n_{nc}=n$. We can then write
  $\dot x = f(x_{nc})$.
\end{definition}

Many mechanical systems have cyclic coordinates. For example, quadrotor drone
and fixed-wing aircraft dynamics do not depend on the position or yaw angle
\cite{FlightMechanicsBook}. Satellite dynamics in a circular low Earth orbit,
frequently approximated with the Clohessy\dash Wiltshire\dash Hill equations
\cite{DeRuiterBook}, do not depend on the true anomaly angle which locates the
spacecraft along the orbit.

With \dref{cyclic} in hand, we call a \alert{cyclic transformation}
$\cyclicshift:\real^n\to\real^n$ any mapping from the state space to itself
which translates the state vector along the cyclic coordinates. In other words,
assuming without loss of generality that the state is given by
\eqref{eq:state_cyclic}, we can write:
\begin{equation}
  \cyclicshift(x) = \Matrix{x_c+\Delta x_c \\ x_{nc}}
\end{equation}
for some translation $\Delta x_c\in\reals^n$. Let us now consider the polytopic
state constraint \optieqref{lcvx_r_linstate}{affine_state} and, in particular,
let $\set F_i = \{x\in\real^n: H_i\T x= h_i\}$ be its $i$-th facet ($H_i\T$ is
the $i$-th row of $H$). 
The following condition must then hold.

\begin{condition}[lcvx_linstate_cyclic_shift]
  Let $\set F_i = \{x\in\real^n: H_i\T x= h_i\}$ denote the $i$-th facet of the
  polytopic state constraint \optieqref{lcvx_r_linstate}{affine_state}, for any
  $i\in\{1,\dots,n_h\}$. If $h_i\ne 0$, then there must exist a cyclic
  transformation $\cyclicshift$ such that
  \begin{equation}
    H_i\T\cyclicshift(x) = 0.
  \end{equation}
\end{condition}

\begin{csmfigure}[%
  caption={%
    Illustration of a landing glideslope constraint ({\color{beamerRed}red},
    sideview of \figref{sidebar_affine_glideslope}) that undergoes a cyclic shift
    in position along the positive $x_2$ axis, to arrive at a new landing
    location ({\color{beamerBlue}blue}). Thanks to
    \conref{lcvx_linstate_cyclic_shift}, the \lcvx guarantee continues to hold for
    the new glideslope constraint, even though the new constraint facets are
    not subspaces.
  },
  label={landing_glideslope_cyclic_shift}]%
  \centering \includegraphics[width=0.6\columnwidth]{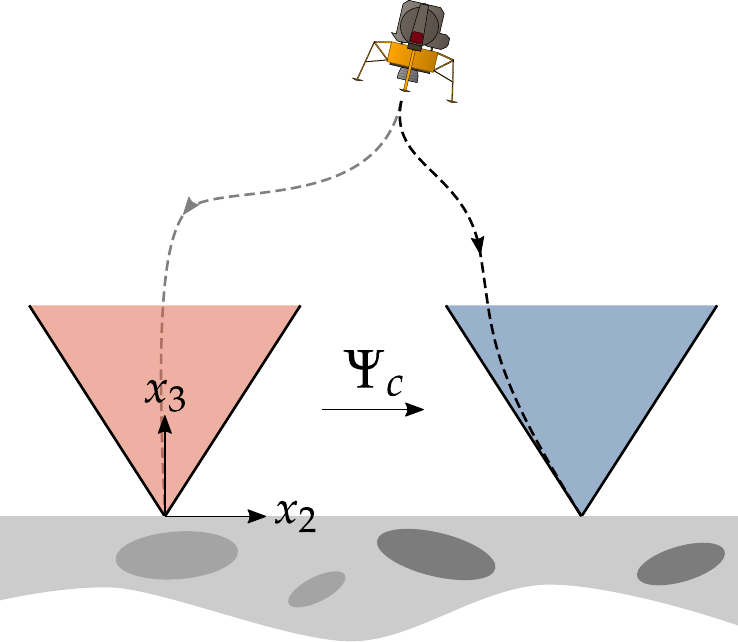}
\end{csmfigure}

To visualize the implication of \conref{lcvx_linstate_cyclic_shift}, simply
consider the case of the landing glideslope constraint from
\sbref{sidebar_affinestate}. Because the position of a spacecraft in a constant
gravity field is a cyclic coordinate, \conref{lcvx_linstate_cyclic_shift} confirms
our intuitive understanding that we can impose landing at a coordinate other
than the origin without compromising lossless
convexification. \figref{landing_glideslope_cyclic_shift} provides an illustration.

From linear systems theory \cite{TrentelmanBook,Harris2013a,Harris2014}, it
turns out that $x(t)\in\set F_i$ can hold for a non\dash zero time interval
(i.e. the state ``sticks'' to a facet) if and only if there exists a triplet of
matrices
$\{F_i\in\reals^{m\times n},G_i\in\reals^{m\times n_v},H_i\in\reals^{m\times p}\}$
such that:
\begin{equation}
  \label{eq:input_friends_mapping}
  u(t) = F_ix(t)+G_iv(t)+H_iw.
\end{equation}

\begin{csmfigure}[%
  caption={Given $x(0)\in\set F_i$, the dynamical system
    \optieqref{lcvx_o_linstate}{dynamics} evolves on $\set F_i$ if and only if $u(t)$ is
    of the form \eqref{eq:input_friends_mapping}.},
  label={friends_block_diagram}
  ]%
  \centering
  \ifmaketwocolcsm
  \includegraphics[width=0.75\columnwidth]{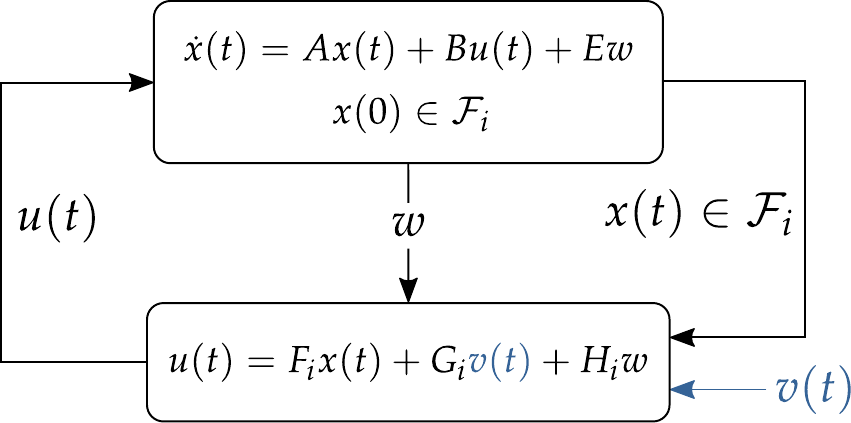}
  \else
  \includegraphics[width=0.7\columnwidth]{friends_block_diagram}
  \fi
\end{csmfigure}

The ``new'' control input $v(t)\in\real^{n_v}$ effectively gets filtered
through \eqref{eq:input_friends_mapping} to produce a control that maintains
$x(\cdot)$ on the hyperplane $\set F_i$. The situation is illustrated in
\figref{friends_block_diagram} in a familiar block diagram form. While the matrix
triplet is not unique, a valid triplet may be computed via standard linear
algebra operations. The reader may consult these operations directly in the
source code of the \lcvx examples provided at the end of this article.


The proof of \lcvx for \pref{lcvx_o_linstate} was originally developed in
\cite{Harris2013b,Harris2014}. The theory behind equation
\eqref{eq:input_friends_mapping} is among the most abstract in all of \lcvx, and we
shall not attempt a proof here. The ultimate outcome of the proof is that the
following condition must hold.

\begin{condition}[lcvx_linstate_controllability]
  For each facet $\set F_i\subseteq\reals^n$ of the polytopic state constraint
  \optieqref{lcvx_r_linstate}{affine_state}, the following ``dual'' linear system
  has no transmission zeros:
  \begin{align*}
    \dot\lambda(t) &= -(A+BF_i)\T\lambda(t)-(CF_i)\T \mu(t), \\
    y(t) &= (BG_i)\T\lambda(t) +(CG_i)\T\mu(t).
  \end{align*}
\end{condition}

Transmission zeros are defined in \cite[Section~4.5.1]{MIMOBook}. Roughly
speaking, if there are no transmission zeros then there cannot exist an initial
condition $\lambda(0)\in\reals^n$ and an input trajectory $\mu\in\reals^{n_c}$
such that $y(t)=0$ for a non\dash zero time interval.

We can now state when lossless convexification holds for
\pref{lcvx_r_linstate}. Note that the statement is very similar to
\tref{nostate}. Indeed, the primary contribution of \cite{Harris2013b,Harris2014}
was to introduce \conref{lcvx_linstate_controllability} and to show that \lcvx
holds by using a version of the maximum principle that includes state
constraints \cite{hartl1995survey,Milyutin1998}. The actual lossless
convexification procedure, meanwhile, does not change.

\begin{theorem}[linstate]
  The solution of \pref{lcvx_r_linstate} is globally optimal for
  \pref{lcvx_o_linstate} if
  \conref{lcvx_nostate_linindep,lcvx_linstate_controllability} hold.
\end{theorem}

Most recently, a similar \lcvx result was proved for problems with affine
equality state constraints that can furthermore depend on the input
\cite{Kunhippurayil2021Observability}. These so\dash called mixed constraints
are of the following form:
\begin{equation}
  \label{eq:lcvx_mixed_constraints}
  y(t) = C(t)x(t)+D(t)u(t),
\end{equation}
where $y$, $C$, and $D$ are problem data.

\subsection{Quadratic State Constraints}

In the last section, we showed an \lcvx result for state constraints that can
be represented by the affine description
\optieqref{lcvx_o_linstate}{affine_state}. The natural next question is whether
\lcvx extends to more complicated state constraints. It turns out that a
generalization of \lcvx exists for quadratic state constraints, if one can
accept a slight restriction to the system dynamics. This result was originally
presented in \cite{Harris2013a,HarrisThesis}. The nonconvex problem statement is
as follows:

\begin{optimization}[
  label={lcvx_o_quadstate},
  variables={u,t_f},
  objective={m(x(t_f))+\runningk \int_0^{t_f} \ell(g_1(u(t)))\,\dt}]%
  \optilabel{dynamics_1}
  & \dot x_1(t) = Ax_2(t), \\
  \optilabel{dynamics_2}
  & \dot x_2(t) = Bu(t)+w, \\
  \optilabel{bounds}
  & \rho_{\min}\le g_1(u(t)),~g_0(u(t))\le\rho_{\max}, \\
  \optilabel{max_speed}
  & x_2(t)\T H x_2(t)\le 1, \\
  & x_1(0)=x_{1,0},~x_2(0)=x_{2,0},~b(x(t_f))=0,
\end{optimization}
where $(x_1,x_2)\in\real^n\times\real^n$ is the state that has been partitioned
into two distinct parts, $u\in\reals^n$ is an input of the same dimension, and
the exogenous disturbance $w\in\reals^n$ is some fixed constant. The matrix
$H\in\real^{n\times n}$ is symmetric positive definite such that
\optieqref{lcvx_o_quadstate}{max_speed} maintains the state in an ellipsoid. An
example of such a constraint is illustrated in \sbref{sidebar_quadstate}.

\declaresidebar[%
title={Maximum Velocity as a Quadratic State Constraint},
label={sidebar_quadstate}]{%
  sidebars/lcvx_max_speed.tex}

Although the dynamics
\optieqref{lcvx_o_quadstate}{dynamics_1}-\optieqref{lcvx_o_quadstate}{dynamics_2} are
less general than \optieqref{lcvx_o_linstate}{dynamics}, they can still accommodate
problems related to vehicle trajectory generation. In such problems, the
vehicle is usually closely related to a double integrator system, for which
$A=I_n$ and $B=I_n$ such that $x_1$ is the position and $x_2$ is the velocity
of the vehicle. The control $u$ in this case is the acceleration.

The following assumption further restricts the problem setup, and is a
consequence of the lossless convexification proof \cite{Harris2013b}.

\begin{assumption}[lcvx_quadstate_ass]
  The matrices $A$, $B$, and $H$ in \pref{lcvx_o_linstate} are invertible. The
  functions $g_0(\cdot)$ and $g_1(\cdot)$ satisfy
  $\rho_{\min}\le g_1(-B\inv w)$ and $g_0(-B\inv w)\le\rho_{\max}$.
\end{assumption}

\aref{lcvx_quadstate_ass} has the direct interpretation of requiring that the
disturbance $w$ can be counteracted by an input that is feasible with respect
to \optieqref{lcvx_o_quadstate}{bounds}. The relaxed problem once again simply
convexifies the nonconvex input lower bound by introducing a slack input:

\begin{optimization}[
  label={lcvx_r_quadstate},
  variables={\sigma,u,t_f},
  objective={m(x(t_f))+\runningk \int_0^{t_f} \ell(\sigma(t))\,\dt}]
  \optilabel{dynamics_1}
  & \dot x_1(t) = Ax_2(t), \\
  \optilabel{dynamics_2}
  & \dot x_2(t) = Bu(t)+w, \\
  & \rho_{\min}\le \sigma(t),~g_0(u(t))\le\rho_{\max}, \\
  & {\color{lcvxColor}g_1(u(t))\le\sigma(t)}, \\
  & x_2(t)\T H x_2(t)\le 1, \\
  & x_1(0) = x_{1,0},~x_2(0) = x_{2,0},~b(x(t_f))=0.
\end{optimization}

Thanks to the structure of the dynamics
\optieqref{lcvx_r_quadstate}{dynamics_1}-\optieqref{lcvx_r_quadstate}{dynamics_2},
it can be shown that \conref{lcvx_nostate_controllability} is automatically
satisfied. On the other hand, \conref{lcvx_nostate_linindep} must be modified to
account for the quadratic state constraint.

\begin{condition}[lcvx_quadstate_linindep]
  If $\runningk=0$, the vector $\grad_x m[t_f]\in\real^{2n}$ and the columns of
  the following matrix must be linearly independent:
  \begin{equation}
    \tilde B_{\lcvx} =%
    \Matrix{\grad_{x_1} b[t_f]\T & 0 \\ \grad_{x_2} b[t_f]\T & 2Hx_2(t_f)}%
    \in\real^{2n\times (n_b+1)}.
  \end{equation}
\end{condition}

Note that \conref{lcvx_quadstate_linindep} carries a subtle but important
implication. Recall that due to \aref{terminal_state_not_overconstrained},
$\grad_x b[t_f]\T$ must be full column rank. Hence, if $\runningk=0$ then the
vector $\pare[big]{0,~2Hx_2(t_f)}\in\reals^{2n}$ and the columns of
$\grad_x b[t_f]\T$ must be linearly dependent. Otherwise, $\tilde B_{\lcvx}$ is
full column rank and \conref{lcvx_quadstate_linindep} cannot be satisfied. With
this in mind, the following \lcvx result was proved in
\cite[Theorem~2]{Harris2013a}.


\begin{theorem}[quadstate]
  The solution of \pref{lcvx_r_quadstate} is globally optimal for
  \pref{lcvx_o_quadstate} if \conref{lcvx_quadstate_linindep} holds.
\end{theorem}

\subsection{General Convex State Constraints}

The preceding two sections discussed problem classes where an \lcvx guarantee
is available even in the presence of affine or quadratic state constraints. For
obvious reasons, an engineer may want to impose more exotic constraints than
afforded by \pref{lcvx_o_linstate,lcvx_o_quadstate}. Luckily, an \lcvx guarantee is
available for general convex state constraints.

As may be expected, generality comes at the price of a somewhat weaker
result. In the preceding sections, the \lcvx guarantee was independent from the
way in which the affine and quadratic state constraints get activated:
instantaneously, for periods of time, or even for the entire optimal trajectory
duration. In contrast, for the case of general convex state constraints, an
\lcvx guarantee will only hold as long as the state constraints are active
\alert{pointwise} in time. In other words, they get activated at isolated time
instances and never persistently over a time interval. This result was
originally provided in \cite{Acikmese2011}. The nonconvex problem statement is:
\begin{optimization}[
  label={lcvx_o_genstate},
  variables={u,t_f},
  objective={m(x(t_f))+\runningk \int_0^{t_f} \ell(g_1(u(t)))\,\dt}]%
  \optilabel{dynamics}
  & \dot x(t) = Ax(t)+Bu(t)+Ew, \\
  \optilabel{bounds}
  & \rho_{\min}\le g_1(u(t)),~g_0(u(t))\le\rho_{\max}, \\
  \optilabel{state_constraint}
  & x(t)\in\mathcal X, \\
  & x(0) = x_0,~b(x(t_f))=0,
\end{optimization}
where $\mathcal X\subseteq\real^n$ is a convex set that defines the state
constraints. Without the state constraint, \pref{lcvx_o_genstate} is nothing but
the autonomous version of \pref{lcvx_o_nostate}. As for \pref{lcvx_o_linstate}, time
variance can be introduced in a limited way by using a time integrator state,
as long as this does not introduce nonconvexity. The relaxed problem uses the
by-now familiar slack variable relaxation technique for
\optieqref{lcvx_o_genstate}{bounds}:

\begin{optimization}[
  label={lcvx_r_genstate},
  variables={\sigma,u,t_f},
  objective={m(x(t_f))+\runningk \int_0^{t_f} \ell(\sigma(t))\,\dt}]
  & \dot x(t) = Ax(t)+Bu(t)+Ew, \\
  & \rho_{\min}\le \sigma(t),~g_0(u(t))\le\rho_{\max}, \\
  & {\color{lcvxColor}g_1(u(t))\le\sigma(t)}, \\
  & x(t)\in\mathcal X, \\
  & x(0) = x_0,~b(x(t_f))=0.
\end{optimization}

The \lcvx proof is provided in \cite[Corollary~3]{Acikmese2011}, and relies on
recognizing two key facts:
\begin{enumerate}
\item When $x(t)\in\interior{\set X}$ for any time interval $t\in[t_1,t_2]$, the
  state of the optimal control problem is unconstrained along that time
  interval;
\item For autonomous problems (recall the description after
  \pref{lcvx_o_linstate}), every segment of the trajectory is itself optimal
  \cite{Harris2013a}.
\end{enumerate}

\begin{csmfigure}[%
  caption={%
    The dashed {\color{beamerBlue}blue} curve represents any segment of the
    optimal state trajectory for \pref{lcvx_o_genstate} that evolves in the
    interior of the state constraint set
    \optieqref{lcvx_o_genstate}{state_constraint}. Because the optimal control
    problem is autonomous, any such segment is the solution to the
    state-unconstrained \pref{lcvx_o_nostate}. When $\runningk=1$ and in the limit
    as $a\to\infty$, \lcvx applies to the entire (open) segment inside
    $\interior{\set X}$ \cite{Acikmese2011}.
  },
  label={lcvx_genstate_interior_arc}
  ]%
  \centering
  \ifmaketwocolcsm
  \includegraphics[width=0.7\columnwidth]{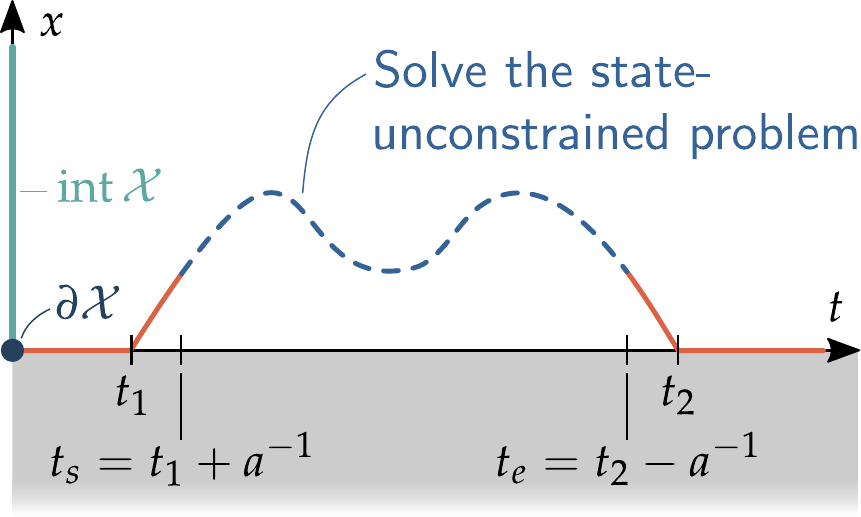}
  \else
  \includegraphics[width=0.7\columnwidth]{lcvx_genstate_interior_arc}
  \fi
\end{csmfigure}

As a result, whenever $x(t)\in\interior{\set X}$, the solution of
\pref{lcvx_r_genstate} is equivalent to the solution of
\pref{lcvx_r_nostate}. Consider an \alert{interior} trajectory segment, as
illustrated in \figref{lcvx_genstate_interior_arc}. The optimal trajectory for
the dashed portion in \figref{lcvx_genstate_interior_arc} is the solution of the
following fixed final state, free final time problem:
\begin{optimization}[
  label={lcvx_genstate_subarc},
  variables={u,t_f},
  objective={m(z(t_e))+\runningk \int_{t_s}^{t_e} \ell(g_1(u(t)))\,\dt}]
  & \dot z(t) = Az(t)+Bu(t)+Ew, \\
  & \rho_{\min}\le g_1(u(t)),~g_0(u(t))\le\rho_{\max}, \\
  & z(t_s) = x(t_s),~z(t_e)=x(t_e).
\end{optimization}

We recognize that \pref{lcvx_genstate_subarc} is an instance of \pref{lcvx_o_nostate}
and, as long as $\runningk=1$ (in order for \conref{lcvx_nostate_linindep} to
hold), \tref{nostate} applies. Because $a>0$ can be arbitrarily large in
\figref{lcvx_genstate_interior_arc}, lossless convexification applies over the open
time interval
$(t_1,t_2)$. 
Thus, the solution segments of the relaxed problem that lie in the interior of
the state constraint set are feasible and globally optimal for the original
\pref{lcvx_o_genstate}.

\declaresidebar[%
title={Permissible State Constraint Activation for General Convex State
  Constraints},
label={sidebar_genstate}]{%
  sidebars/lcvx_pointwise_state.tex}

The same cannot be said when $x(t)\in\boundary{\set X}$. During these segments,
the solution can become infeasible for \pref{lcvx_o_genstate}. However, as long as
$x(t)\in\boundary{\set X}$ at isolated time instances, \lcvx can be guaranteed
to hold. This idea is further illustrated in \sbref{sidebar_genstate}.


When $\runningk=0$, the situation becomes more complicated because
\conref{lcvx_nostate_linindep} does not hold for \pref{lcvx_genstate_subarc}. This is
clear from the fact that the terms defined in \conref{lcvx_nostate_linindep}
become:
\begin{equation*}
  m_{\text{\lcvx}} = \Matrix{\grad_z m[t_e] \\ 0}, \quad
  B_{\text{\lcvx}} = \Matrix{I_n \\ 0},
\end{equation*}
which are clearly not linearly independent since $B_{\text{\lcvx}}$ is full
column rank. Thus, even for interior segments the solution may be infeasible
for \pref{lcvx_o_genstate}. To remedy this, \cite[Corollary~4]{Acikmese2011}
suggests \algref{lcvx_genstate}. At its core, the algorithm relies on the following
simple idea. By solving \pref{lcvx_r_genstate} with the suggested modifications on
\algref[start=modifications]{lcvx_genstate}, every interior segment once again
becomes an instance of \pref{lcvx_o_nostate} for which \tref{nostate}
holds. Furthermore, due to the constraint $x(t_f)=\optimal{x}(\optimal{t_f})$,
any solution to the modified problem will be optimal for the original
formulation where $\runningk=0$ (since
$m(x(t_f))=m(\optimal{x}(\optimal{t_f}))$).

This modification can be viewed as a search for an equivalent solution for
which \lcvx holds. As a concrete example, \pref{lcvx_r_genstate} may be searching
for a minimum miss distance solution for a planetary rocket landing trajectory
\cite{Blackmore2010}. The ancillary problem in \algref{lcvx_genstate} can search for
a minimum fuel solution that achieves the same miss distance. Clearly, other
running cost choices are possible. Thus, the ancillary problem's running cost
becomes an extra tuning parameter.

\begin{algorithm}[t]
  \centering
  \begin{algorithmic}[1]
    \State Solve \pref{lcvx_r_genstate} to obtain $\optimal{x}(\optimal{t_f})$
    \label{alg:lcvx_genstate:line:step1}
    \If{$\runningk=0$}
    \State Solve \pref{lcvx_r_genstate} again, with the modifications:
    \begin{itemize}
    \item Use the cost $\int_0^{t_f}\ell(\sigma(t))\dd t$
    \item Set $b(x(t_f))=x(t_f)-\optimal{x}(\optimal{t_f})$ 
    \end{itemize}
    \label{alg:lcvx_genstate:line:modifications}
    \EndIf
  \end{algorithmic}
  \caption{Solution algorithm for \pref{lcvx_o_genstate}. When
    $\runningk=0$, a two-step procedure is used where an auxiliary problem with
    $\runningk=1$ searches over the optimal solutions to the original problem.}
  \label{alg:lcvx_genstate}
\end{algorithm}

We are now able to summarize the lossless convexification result for problems
with general convex state constraints.

\begin{theorem}[genstate]
  \algref{lcvx_genstate} returns the globally optimal solution of
  \pref{lcvx_o_genstate} if the state constraint
  \optieqref{lcvx_o_genstate}{state_constraint} is activated at isolated time
  instances, and \conref{lcvx_nostate_controllability,lcvx_nostate_linindep} hold.
\end{theorem}

\subsection{Nonlinear Dynamics}

A unifying theme of the previous sections is the assumption that the system
dynamics are linear. In fact, across all \lcvx results that we have mentioned
so far, the dynamics did not vary much from the first formulation in
\optieqref{lcvx_o_nostate}{dynamics}. Many engineering applications, however,
involve non-negligible nonlinearities. A natural question is then whether the
theory of lossless convexification can be extended to systems with general
nonlinear dynamics.

An \lcvx result is available for a class of nonlinear dynamical systems. The
groundwork for this extension was presented in \cite{Blackmore2012}. The goal
here is to show that the standard input set relaxation based on the \lcvx
equality constraint is also lossless when the dynamics are
nonlinear. Importantly, note that the dynamics themselves are not convexified,
so the relaxed optimization problem is still nonlinear, and it is up to the
user to solve the problem to global optimality. This is possible in special
cases, for example if the nonlinearities are approximated by piecewise affine
functions. This yields a mixed-integer convex problem whose globally optimal
solution can be found via mixed-integer programming \cite{achterberg2013mixed}.

With this introduction, let us begin by introducing the generalization of
\pref{lcvx_o_nostate} that we shall solve using \lcvx:
\begin{optimization}[
  label={lcvx_o_nonlinear},
  variables={u,t_f},
  objective={m(t_f,x(t_f))+\runningk \int_0^{t_f} \ell(g(u(t)))\,\dt}]%
  \optilabel{dynamics}
  & \dot x(t) = f(t,x(t),u(t),g(u(t))), \\
  \optilabel{bounds}
  & \rho_{\min}\le g(u(t))\le\rho_{\max}, \\
  & x(0) = x_0,~b(t_f,x(t_f))=0,
\end{optimization}
where $f:\real\times\real^n\times\real^m\times\real\to\real^n$ defines the
nonlinear dynamics. Just as for \pref{lcvx_o_nostate_pointing}, it is required that
$g_0=g_1\definedas g$. Consider the following convex relaxation of the input
constraint by using a slack input:
\begin{optimization}[
  label={lcvx_r_nonlinear},
  variables={\sigma,u,t_f},
  objective={m(t_f,x(t_f))+\runningk \int_0^{t_f} \ell(\sigma(t))\,\dt}]%
  \optilabel{dynamics}
  & \dot x(t) = f(t,x(t),u(t),\sigma(t)), \\
  & \rho_{\min}\le\sigma(t)\le\rho_{\max}, \\
  & {\color{lcvxColor}g(u(t))\le\sigma(t)}, \\
  & x(0) = x_0,~b(t_f,x(t_f))=0.
\end{optimization}

Note that the slack input $\sigma$ makes a new appearance in the dynamics
\optieqref{lcvx_r_nonlinear}{dynamics}. The more complicated dynamics require an
updated version of \conref{lcvx_nostate_controllability} in order to guarantee that
\lcvx holds.

\begin{condition}[lcvx_nonlinear_controllability]
  The pair $\{\grad_{x}f[t],\grad_{u}f[t]\}$ must be totally controllable on
  $[0,t_f]$ for all feasible sequences of $x(\cdot)$ and $u(\cdot)$ for
  \pref{lcvx_r_nonlinear} \cite{Blackmore2012,DAngelo1970}.
\end{condition}

Using the above condition, we can state the following quite general \lcvx
guarantee for problems that fit the \pref{lcvx_o_nonlinear} template.

\begin{theorem}[nonlinear]
  The solution of \pref{lcvx_r_nonlinear} is globally optimal for
  \pref{lcvx_o_nonlinear} if
  \conref{lcvx_nostate_linindep,lcvx_nonlinear_controllability} hold.
\end{theorem}

Alas, \conref{lcvx_nonlinear_controllability} is generally quite difficult to
check. Nevertheless, two general classes of systems have been shown to
automatically satisfy this condition thanks to the structure of their dynamics
\cite{Blackmore2012}. These classes accommodate vehicle trajectory generation
problems with double integrator dynamics and nonlinearities like mass
depletion, aerodynamic drag, and nonlinear gravity. The following discussion of
these system classes can appear hard to parse at first sight. For this reason,
we provide two practical examples of systems that belong to each class in
\sbref{lcvx_nonlinear}.

\declaresidebar[%
title={Examples of Losslessly Convexifiable Nonlinear Systems},
label={lcvx_nonlinear},
position=h]{%
  sidebars/lcvx_nonlinear.tex}

The first corollary of \tref{nonlinear} introduces the first class of systems. A
key insight is that the nullspace conditions of the corollary require that
$2m\ge n$, in other words there are at least twice as many control variables as
there are state variables. This is satisfied by some vehicle trajectory
generation problems where $2m=n$, for example when the state consists of
position and velocity while the control is an acceleration that acts on all the
velocity states. This is a common approximation for flying vehicles. We shall
see an example for rocket landing in Part III of the article.

\begin{corollary}[nonlinear_1]
  Suppose that the dynamics \optieqref{lcvx_o_nonlinear}{dynamics} are of the
  form:
  \begin{equation}
    \label{eq:lcvx_o_nonlinear_class1_dynamics}
    x=\Matrix{
      x_1 \\
      x_2
    },~
    f(t,x,u,g(u)) = \Matrix{
      f_1(t,x) \\
      f_2(t,x,u)
    },
  \end{equation}
  where $\nul(\grad_u f_2)=\{0\}$ and $\nul(\grad_{x_2} f_1)=\{0\}$. Then
  \tref{nonlinear} applies if \conref{lcvx_nostate_linindep} holds.
\end{corollary}

The next corollary to \tref{nonlinear} introduces the second class of systems, for
which $2m<n$ is allowed. This class is once again useful for vehicle trajectory
generation problems where the dynamics are given by
\eqref{eq:lcvx_o_nonlinear_class2_dynamics} and $g(u)$ is a function that measures
control effort. A practical example is when the state $x_2$ is mass, which is
depleted as a function of the control effort (such as thrust for a rocket).

\begin{corollary}[nonlinear_2]
  Suppose that the dynamics \optieqref{lcvx_o_nonlinear}{dynamics} are of the
  form:
  \begin{equation}
    \label{eq:lcvx_o_nonlinear_class2_dynamics}
    x=\Matrix{
      x_1 \\
      x_2
    },~
    f(t,x,u,g(u)) = \Matrix{
      f_1(t,x,u) \\
      f_2(t,g(u))
    }.
  \end{equation}

  Define the matrix:
  \begin{equation}
    \label{eq:lcvx_o_nonlinear_class2_M3}
    M\definedas\Matrix{
      (\grad_u f_1)\T \\
      \displaystyle\frac{\dd(\grad_u f_1)\T}{\dt}-(\grad_u f)\T (\grad_x f_1)\T
    }.
  \end{equation}

  Furthermore, suppose that the terminal constraint function $b$ is affine and
  $x_2(t_f)$ is unconstrained, such that $\grad_{x_2} b = 0$. Then \tref{nonlinear}
  applies if $\nul(M)=\{0\}$ and \conref{lcvx_nostate_linindep} holds.
\end{corollary}

It must be emphasized that \pref{lcvx_r_nonlinear} is still a nonlinear program
and that for \tref{nonlinear} to hold, a globally optimal solution of
\pref{lcvx_r_nonlinear} must be found. Although this cannot be done for general
nonlinear programming, if the dynamics $f$ are piecewise affine then the problem
can be solved to global optimality via mixed-integer programming
\cite{achterberg2007constrained,achterberg2013mixed,Schouwenaars2001}. In this
case, convexification of the nonconvex input lower bound reduces the number of
disjunctions in the branch-and-bound tree, and hence lowers the problem
complexity \cite{Blackmore2012}. Several examples of nonlinear systems that can
be modeled in this way, and which comply with \corref{nonlinear_1,nonlinear_2},
are illustrated in \sbref{lcvx_pwa}.

\declaresidebar[%
title={Approximating Nonlinear Systems with Piecewise Affine Functions},
label={lcvx_pwa},
position=t]{%
  sidebars/lcvx_pwa.tex}

\subsection{Embedded Lossless Convexification}

The reader will notice that the \lcvx theory of the previous sections deals
with special cases of problems whose nonconvexity is ``just right'' for an
\lcvx guarantee to be provable using the maximum principle. Although such
problems have found their practical use in problems like spaceflight
\cite{Scharf2017} and quadrotor path planning \cite{Szmuk2017}, it leaves out
many trajectory generation applications that do not fit the tight mold of
original problems and conditions of the previous sections.

\begin{csmfigure}[%
  caption={Illustration of how embedded \lcvx can be used to solve an optimal
    control problem that does not fit into any of the templates presented in Part I
    of this article.},
  label={embedded_lcvx_procedure}]%
  \centering
  \begin{subfigure}{1.0\linewidth}
    \centering
    \includegraphics[width=0.9\linewidth]{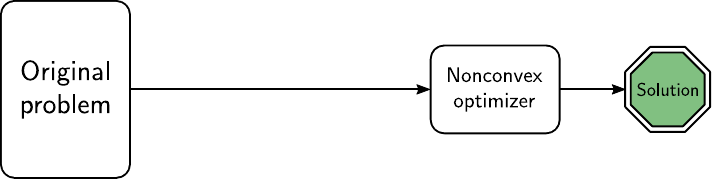}
    \caption{The solution process for a nonconvex optimal control problem,
      without using \lcvx.}
    \figlabel{embedded_lcvx_none}
  \end{subfigure}%
  \vspace{3mm}

  \begin{subfigure}{1.0\linewidth}
    \centering
    \includegraphics[width=0.9\linewidth]{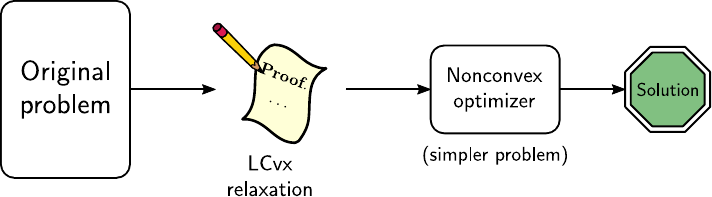}
    \caption{The solution process for a nonconvex optimal control problem, where
      \lcvx is embedded to convexify part of the original problem.}
    \figlabel{embedded_lcvx_with}
  \end{subfigure}
\end{csmfigure}

Despite this apparent limitation, \lcvx is still highly relevant for problems
that simply do not conform to one of the forms given in the previous
sections. For such problems, we assume that the reader is facing the challenge
of solving a nonconvex optimal control problem that fits the mold of
\pref{scp_gen_cont} (the subject of Part II of the article), and is considering
whether \lcvx can help. There is evidence that the answer is affirmative, by
using \lcvx theory only on the constraints that are losslessly
convexifiable. We call this \alert{embedded \lcvx}, because it is used to
convexify only part of the problem, while the rest is handled by another
nonconvex optimization method such as presented in Part II of this
article. Because \lcvx reduces the amount of nonconvexity present in the
problem, it can significantly improve the convergence properties and reduce the
computational cost to solve the resulting problem. An example of this approach
for quadrotor trajectory generation is demonstrated in Part III.

The basic procedure for applying embedded \lcvx is illustrated in
\figref{embedded_lcvx_procedure}. As shown in \figref{embedded_lcvx_with}, we
reiterate that \lcvx is not a computation scheme, but rather it is a convex
relaxation with an accompanying proof of equivalence to the original problem. As
such, it happens prior to the solution and simply changes the problem
description seen by the subsequent numerical optimization algorithm.

There are a number of examples of embedded \lcvx that we can mention. First,
the previous section on nonlinear dynamics can be intepreted as embedded
\lcvx. For example, \cite{Blackmore2012} solves a rocket landing problem where
only the nonconvex input constraint \optieqref{lcvx_o_nonlinear}{bounds} is
convexified. This leaves behind a nonconvex problem due to nonlinear dynamics,
and mixed-integer programming is used to solve it. Another example is
\cite{Zhang2017}, where \lcvx is embedded in a mixed-integer autonomous aerial
vehicle trajectory generation problem in order to convexify a stall speed
constraint of the form:
\begin{equation}
  \label{eq:stall_speed_constraint}
  0<v_{\min}\le\norm[2]{v_{cmd}(t)}\le v_{\max},
\end{equation}
where the input $v_{cmd}(\cdot)\in\reals^3$ is the commanded velocity, while
$v_{\min}$ and $v_{\max}$ are lower and upper bounds that guarantee a stable
flight envelope. The same constraint is also considered in
\cite{Blackmore2012}. In \cite{Liu2016}, the authors develop a highly nonlinear
planetary entry trajectory optimization problem where the control input is the
bank angle $\beta\in\reals$, parameterized via two inputs
$u_1\definedas\cos(\beta)$ and $u_2\definedas\sin(\beta)$. The associated
constraint $\norm[2]{u}^2=1$ is convexified to $\norm[2]{u}^2\le 1$, and
equality at the optimal solution is shown in an \lcvx-like fashion (the authors
call it ``assurance of active control constraint''). Similar methods are used
in \cite{Liu2019} in the context of rocket landing with aerodynamic controls. A
survey of related methods is available in \cite{Liu2017}. Finally, we will
mention \cite{Szmuk2016,Szmuk2017} where embedded \lcvx is used to convexify an
input lower bound and an attitude pointing constraint for rocket landing and
for agile quadrotor flight. Sequential convex programming from Part II is then
using to solve the remaining nonlinear optimal control problems. The quadrotor
application in particular is demonstrated as a numerical example in Part III.

As a result of the success of these applications, we foresee there being
further opportunities to use \lcvx as a strategy to derive simpler problem
formulations. The result would be a speedup in computation for
optimization-based trajectory generation.

\subsection{The Future of Lossless Convexification}


Lossless convexification is a method that solves nonconvex trajectory
generation problems with one or a small number of calls to a convex
solver. This places it among the most reliable and robust methods for nonconvex
trajectory generation. The future of \lcvx therefore has an obvious motivation:
to expand the class of problems that can be losslessly convexified. The most
recent result discussed in the previous sections is for problems with affine
state constraints \cite{Harris2014}, and this dates back to 2014. In the past
two years, \lcvx research has been rejuvenated by several fundamental
discoveries and practical methods that expand the method to new and interesting
problem types. This section briefly surveys these new results.

\subsubsection{Fixed\dash final Time Problems}

The first new \lcvx result applies to a fixed\dash final time and fixed\dash
final state version of \pref{lcvx_o_nostate} with no state constraints. To begin,
recognize that the classical \lcvx result from \tref{nostate} does not apply when
both $t_f$ and $x(t_f)$ are fixed. In this case, $B_{\text{\lcvx}}=I_{n+1}$ in
\eqref{eq:lcvx_nostate_linindep_B} and therefore its columns, which span all of
$\reals^{n+1}$, cannot be linearly independent from $m_{\text{\lcvx}}$. Thus,
traditionally one could not fix the final time and the final state
simultaneously. Very recently, Kunhippurayil et
al. \cite{Kunhippurayil2021FixedTime} showed that \conref{lcvx_nostate_linindep} is
in fact not necessary for the following version of \pref{lcvx_o_nostate}:
\begin{optimization}[
  label={lcvx_o_fixedtime},
  variables={u,t_f},
  objective={\int_0^{t_f} \ell(g(u(t)))\,\dt}]
  \optilabel{dynamics}
  & \dot x(t) = Ax(t)+Bu(t), \\
  \optilabel{bounds}
  & \rho_{\min}\le g(u(t))\le\rho_{\max}, \\
  \optilabel{boundary}
  & x(0) = x_0,~x(t_f)=x_f,
\end{optimization}
where $t_f$ is fixed and $x_f\in\reals^n$ specifies the final state. The
lossless relaxation is the usual one, and is just a specialization of
\pref{lcvx_r_nostate} for \pref{lcvx_o_fixedtime}:
\begin{optimization}[
  label={lcvx_r_fixedtime},
  variables={\sigma,u,t_f},
  objective={\int_0^{t_f} \ell(\sigma(t))\,\dt}]
  \optilabel{dynamics}
  & \dot x(t) = Ax(t)+Bu(t)w(t), \\
  \optilabel{bounds}
  & \rho_{\min}\le \sigma(t)\le\rho_{\max}, \\
  \optilabel{lcvx_equality}
  & {\color{lcvxColor}g(u(t))\le\sigma(t)}, \\
  \optilabel{boundary}
  & x(0) = x_0,~x(t_f)=x_f.
\end{optimization}

The following result is then proved in \cite{Kunhippurayil2021FixedTime}. By
dropping \conref{lcvx_nostate_linindep}, the result generalizes \tref{nostate} and
significantly expands the reach of \lcvx to problems without state constraints.

\begin{theorem}[fixed_time]
  The solution of \pref{lcvx_r_fixedtime} is globally optimal for
  \pref{lcvx_o_fixedtime} if \conref{lcvx_nostate_controllability} holds and
  $t_f$ is between the minimum feasible time and the time that minimizes
  \optiobjref{lcvx_o_fixedtime}. For longer trajectory durations, there exists a
  solution to \pref{lcvx_r_fixedtime} that is globally optimal for
  \pref{lcvx_o_fixedtime}.
\end{theorem}

Perhaps the most important part of \tref{fixed_time}, and a significant future
direction for \lcvx, is in its final sentence. Although a lossless solution
``exists'', how does one find it? An \textit{algorithm} is provided in
\cite{Kunhippurayil2021FixedTime} to find the lossless solution, that is, one
solution among many others which may not be lossless. This is similar to
\tref{genstate} and \algref{lcvx_genstate}: we know that slackness in
\optieqref{lcvx_r_fixedtime}{lcvx_equality} may occur, so we devise an algorithm
that works around the issue and is able to recover an input for which
\optieqref{lcvx_r_fixedtime}{lcvx_equality} holds with equality. Most traditional
\lcvx results place further restrictions on the original problem in order to
``avoid'' slackness, but this by definition limits the applicability of
\lcvx. By instead providing algorithms which recover lossless inputs from
problems that do not admit \lcvx naturally, we can tackle lossless
convexification ``head on'' and expand the class of losslessly convexifiable
problems. A similar approach is used for spacecraft rendezvous in
\cite{Harris2014JGCD}, where an iterative algorithm modifies the dynamics in
order to extract bang\dash bang controls from a solution that exhibits
slackness.

\subsubsection{Hybrid System Problems}

Many physical systems contain on\dash off elements such as valves, relays, and
switches
\cite{Bemporad1999a,MalyutaJGCD,Blackmore2012,Schouwenaars2001,Schouwenaars2006}. Discrete
behavior can also appear through interactions between the autonomous agent and
its environment, such as through foot contact for walking robots
\cite{Posa2013,Hwangbo2018}. Modeling discrete behavior is the province of
hybrid systems theory, and the resulting trajectory problems typically combine
continuous variables and discrete logic elements (i.e., ``and'' and ``or''
gates) \cite{Bemporad1999a,Schouwenaars2006,Goebel2009}. Because these is no
concept like local perturbation for values that, for example, can only be equal
to zero or one, problems with discrete logic are fundamentally more
difficult. Traditional solution methods use mixed\dash integer programming
\cite{achterberg2013mixed}. The underlying branch\dash and\dash bound method,
however, has poor (combinatorial) worst\dash case complexity. Historically,
this made it very difficult to put optimization with discrete logic onboard
computationally constrained and safety\dash critical systems throughout aerospace,
automotive, and even state\dash of\dash the\dash art robotics
\cite{TedrakeDARPA}.

Two recent results showed that \lcvx can be applied to certain classes of
hybrid optimal control problems that are useful for trajectory generation
\cite{MalyutaIFAC,Harris2021}. While the results are more general, the following
basic problem will help ground our discussion:
\begin{optimization}[
  label={lcvx_o_hybrid},
  variables={u,\gamma,t_f},
  objective={\int_0^{t_f} \sum_{i=1}^M\norm[2]{u_i(t)}\,\dt}]
  \optilabel{dynamics}
  & \dot x(t) = Ax(t)+B\sum_{i=1}^M u_i(t), \\
  \optilabel{bounds}
  & \gamma_i(t)\rho_{\min,i}\le \norm[2]{u_i(t)}\le\gamma_i(t)\rho_{\max,i}, \\
  \optilabel{binary}
  & \gamma_i(t)\in\{0,1\}, \\
  \optilabel{affine_input}
  & C_iu_i(t)\le 0, \\
  \optilabel{boundary}
  & x(0) = x_0,~b\pare[big]{t_f, x(t_f)}=x_f,
\end{optimization}
where $M$ is the number of individual input vectors and the binary variables
$\gamma_i$ are used to model the on\dash off nature of each input. Compared to
the traditional \pref{lcvx_o_nostate}, this new problem can be seen as a system
controlled by $M$ actuators that can be either ``off'' or ``on'' and norm\dash
bounded in the $[\rho_{\min,i}, \rho_{\max,i}]$ interval. The affine input
constraint \optieqref{lcvx_o_hybrid}{affine_input} represents an affine cone, and
is a specialized version of the earlier constraint
\optieqref{lcvx_o_linstate}{affine_input}. \figref{lcvx_mixed_example} illustrates the
kind of input set that can be modeled. Imitating the previous results, the
convex relaxation uses a slack input for each control vector:
\begin{optimization}[
  label={lcvx_r_hybrid},
  variables={\sigma,u,\gamma,t_f},
  objective={\int_0^{t_f} \sum_{i=1}^M\ell(\sigma_i(t))\,\dt}]
  \optilabel{dynamics}
  & \dot x(t) = Ax(t)+B\sum_{i=1}^M u_i(t), \\
  \optilabel{bounds}
  & \gamma_i(t)\rho_{\min,i}\le\sigma_i(t)\le\gamma_i(t)\rho_{\max,i}, \\
  \optilabel{lcvx_equality}
  & {\color{lcvxColor}\norm[2]{u_i(t)}\le\sigma_i(t)}, \\
  \optilabel{binary}
  & 0\le \gamma_i(t)\le 1, \\
  \optilabel{affine_input}
  & C_iu_i(t)\le 0, \\
  \optilabel{boundary}
  & x(0) = x_0,~b\pare[big]{t_f, x(t_f)}=x_f,
\end{optimization}
where the only real novelty is that the $\gamma_i$ variables have also been
relaxed to the continuous $[0, 1]$ interval.

\begin{csmfigure}[%
  caption={%
    Example of a feasible input set that can be modeled in \pref{lcvx_o_hybrid}. It
    is a nonconvex disconnected set composed of the origin, a
    {\color{beamerRed}point}, an {\color{beamerGreen}arc}, and a
    {\color{beamerBlue}nonconvex set} with an interior. For example, this setup
    can represent a satellite equipped with thrusters and drag plates, or a
    rocket with a thrust-gimbal coupled engine \cite{Harris2021,MalyutaIFAC}.
  },
  label={lcvx_mixed_example}]%
  \centering
  \ifmaketwocolcsm
  \includegraphics[width=0.65\linewidth]{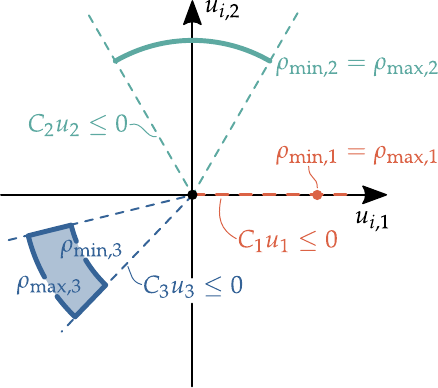}
  \else
  \includegraphics[width=0.65\linewidth]{disconnected_input_set}
  \fi
\end{csmfigure}

Taking \pref{lcvx_r_hybrid} as an example, the works of
\cite{MalyutaIFAC,Harris2021} prove lossless convexification from slightly
different angles. In \cite{Harris2021} it is recognized that \pref{lcvx_r_hybrid}
will be a lossless convexification of \pref{lcvx_o_hybrid} if the dynamical system
is ``normal'' due to the so\dash called bang\dash bang principle
\cite{BerkovitzBook,Hermes1969}. Normality is related to, but much stronger
than, the notion of controllability from
\conref{lcvx_nostate_controllability}. Nevertheless, it is shown that the dynamical
system can be perturbed by an arbitrarily small amount to induce
normality. This phenomenon was previously observed in a practical context for
rocket landing \lcvx with a pointing constraint, which we discussed for
\pref{lcvx_o_nostate_pointing} \cite{Carson2011}. Practical examples are shown for
spacecraft orbit reshaping, minimum\dash energy transfer, and CubeSat
differential drag and thrust maneuvering. It is noted that while mixed\dash
integer programming fails to solve the latter problem, the convex relaxation is
solved in $<0.1$ seconds.

The results in \cite{MalyutaIFAC,MalyutaLCvxRejected} also prove lossless
convexification for \pref{lcvx_r_hybrid}. However, instead of leveraging normality
and perturbing the dynamics, the nonsmooth maximum principle
\cite{Vinter2000,Clarke2010,hartl1995survey} is used directly to develop a set
of conditions for which \lcvx holds. These conditions are an interesting mix of
problem geometry (i.e., the shapes and orientations of the constraint cones
\optieqref{lcvx_r_hybrid}{affine_input}) and \conref{lcvx_nostate_controllability,
  lcvx_nostate_linindep}. Notably, they are more general than normality, so
they can be satisfied by systems that are not normal. Practical examples are
given for spacecraft rendezvous and rocket landing with a coupled thrust\dash
gimbal constraint. The solution is observed to take on the order of a few
seconds and to be more than 100 times faster than mixed\dash integer
programming.

We see the works \cite{MalyutaIFAC,Harris2021} as complementary:
\cite{MalyutaIFAC} shows that for some systems, the perturbation proposed by
\cite{Harris2021} is not necessary. On the other hand, \cite{Harris2021} provides
a method to recover \lcvx when the conditions of \cite{MalyutaIFAC}
fail. Altogether, the fact that an arbitrarily small perturbation of the
dynamics can recover \lcvx suggests a deeper underlying theory for how and why
problems can be losslessly convexified. We feel that the search for this theory
will be a running theme of future \lcvx research, and its eventual discovery
will lead to more general lossless convexification algorithms.

\subsection{Toy Example}

The following example provides a simple illustration of how \lcvx can be used
to solve a nonconvex problem. This example is meant to be a ``preview'' of the
practical application of \lcvx. More challenging and realistic examples are
given in Part III.

The problem that we will solve is minimum-effort control of a double integrator
system (such as a car) with a constant ``friction'' term $g$. This can be
written as a linear time-invariant instance of \pref{lcvx_o_nostate}:
\begin{optimization}[
  label={lcvx_o_toy},
  variables={u},
  objective={\int_0^{t_f}u(t)^2\,\dt}]
  \optilabel{dynamics_1}
  & \dot x_1(t) = x_2(t), \\
  \optilabel{dynamics_2}
  & \dot x_2(t)=u(t)-g, \\
  \optilabel{bounds}
  & 1\le |u(t)|\le 2, \\
  & x_1(0)=x_2(0)=0, \\
  & x_1(t_f)=s,~x_2(t_f)=0,~t_f=10.
\end{optimization}

The input $u(\cdot)\in\reals$ is the acceleration of the car. The constraint
\optieqref{lcvx_o_toy}{bounds} is a nonconvex one-dimensional version of the
constraint \eqref{eq:inputrelax:lb_convexify_1}. Assuming that the car has unit
mass, the integrand in \optiobjref{lcvx_o_toy} has units of Watts. The objective of
\pref{lcvx_o_toy} is therefore to move a car by a distance $s$ in $t_f=10$ seconds
while minimizing the average power. Following the relaxation template provided
by \pref{lcvx_r_nostate}, we propose the following convex relaxation to solve
\pref{lcvx_o_toy}:
\begin{optimization}[
  label={lcvx_r_toy},
  variables={\sigma,u},
  objective={\int_0^{t_f}\sigma(t)^2\,\dt}]%
  & \dot x_1(t) = x_2(t), \\
  & \dot x_2(t)=u(t)-g, \\
  \optilabel{bounds}
  & 1\le \sigma(t)\le 2, \\
  \optilabel{lcvx_equality}
  & {\color{lcvxColor} |u(t)|\le\sigma(t)}, \\
  \optilabel{boundary_start}
  & x_1(0)=x_2(0)=0, \\
  \optilabel{boundary_end}
  & x_1(t_f)=s,~x_2(t_f)=0,~t_f=10.
\end{optimization}

To guarantee that \lcvx holds, in other words that \pref{lcvx_r_toy} finds the
globally optimal solution of \pref{lcvx_o_toy}, let us first attempt to verify the
conditions of \tref{nostate}. In particular, we need to show that
\conref{lcvx_nostate_controllability,lcvx_nostate_linindep} hold. First, from
\optieqref{lcvx_o_toy}{dynamics_1}-\optieqref{lcvx_o_toy}{dynamics_2}, we can extract
the following state\dash space matrices:
\begin{equation}
  \label{eq:toy_AB}
  A =
  \Matrix{
    0 & 1 \\
    0 & 0
  },~B =
  \Matrix{
    0 \\
    1
  }.
\end{equation}

We can verify that \conref{lcvx_nostate_controllability} holds by either showing
that the controllability matrix is full rank, or by using the PBH test
\cite{Antsaklis2006}. Next, from \optiobjref{lcvx_r_toy} and
\optieqref{lcvx_r_toy}{boundary_start}-\optieqref{lcvx_r_toy}{boundary_end}, we can
extract the following terminal cost and terminal constraint functions:
\begin{equation}
  \label{eq:lcvx_r_toy_mb}
  m(t_f,x(t_f)) = 0,~b(t_f,x(t_f))=\Matrix{
    t_f-10 \\
    x_1(t_f)-s \\
    x_2(t_f)
  }.
\end{equation}

We can now substitute \eqref{eq:lcvx_r_toy_mb} into
\eqref{eq:lcvx_nostate_linindep_m_B} to obtain:
\begin{equation}
  \label{eq:m_B_lcvx_toy}
  m_{\text{\lcvx}} =
  \Matrix{
    0 \\ 0 \\ \sigma(t_f)^2
  },~
  B_{\text{\lcvx}} = I_3.
\end{equation}

Thus, $B_{\text{\lcvx}}$ is full column rank and its columns cannot be linearly
independent from $m_{\text{\lcvx}}$. We conclude that
\conref{lcvx_nostate_linindep} does not hold, so \tref{nostate} cannot be applied. In
fact, \pref{lcvx_o_toy} has both a fixed final time and a fixed final state. This
is exactly the edge case for which traditional \lcvx does not apply, as was
mentioned in the previous section on future \lcvx. Instead, we fall back on
\tref{fixed_time} which says that \conref{lcvx_nostate_linindep} is not needed as long
as $t_f$ is between the minimum and optimal times for \pref{lcvx_o_toy}. It turns
out that this holds for the problem parameters used in \figref{lcvx_toy}. The
minimum time is just slightly below $10~\si{\second}$ and the optimal time is
$\approx 13.8~\si{\second}$ for \figref{lcvx_toy_1} and
$\approx 13.3~\si{\second}$ for \figref{lcvx_toy_2}. Most interestingly, lossless
convexification fails (i.e., \optieqref{lcvx_r_toy}{lcvx_equality} does not hold
with equality) for $t_f$ values almost exactly past the optimal time for
\figref{lcvx_toy_1}, and just slightly past it for \figref{lcvx_toy_2}.

Although \pref{lcvx_r_toy} is convex, it has an infinite number of solution
variables because time is continuous. To be able to find an approximation of
the optimal solution using a numerical convex optimization algorithm, the
problem must be temporally discretized. To this end, we apply a first-order
hold (FOH) discretization with $N=50$ temporal nodes, as explained in
\sbref{discretization}.

\ifmaketwocolcsm
\gdef\toycols{1}
\else
\gdef\toycols{2}
\fi
\begin{csmfigure}[%
  caption={\lcvx solutions of \pref{lcvx_r_toy} for two scenarios. The close
    match of the analytic solution using the maximum principle (drawn as a
    continuous line) and the discretized solution using \lcvx (drawn as discrete
    dots) confirms that \lcvx finds the globally optimal solution of the
    problem.},
  label={lcvx_toy},
  position={!t},
  columns=\toycols]%
  \centering
  \ifmaketwocolcsm
  \def\lcvxtoyplotsz{\columnwidth}
  \def\lcvxtoyinnerszleft{0.8\textwidth}
  \def\lcvxtoyinnerszright{0.8\textwidth}
  \else
  \def\lcvxtoyplotsz{0.48\textwidth}
  \def\lcvxtoyinnerszleft{0.975\textwidth}
  \def\lcvxtoyinnerszright{\textwidth}
  \fi
  \begin{subfigure}[b]{\lcvxtoyplotsz}
    \centering
    \includegraphics[width=\lcvxtoyinnerszleft]{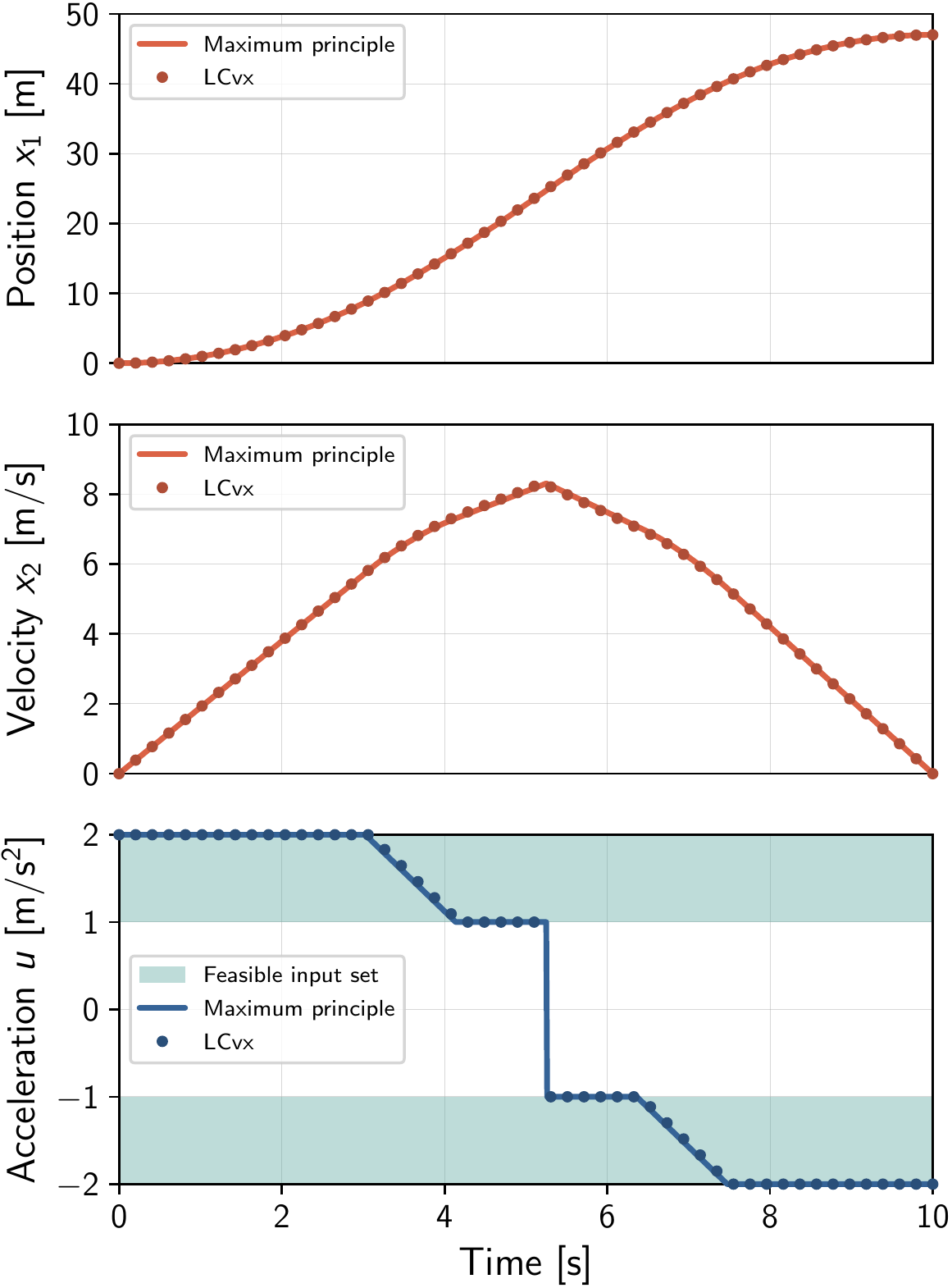}
    \caption{Solution of \pref{lcvx_r_toy} for $g=0.1$~m/s$^2$ and $s=47$~m.}
    \label{fig:lcvx_toy_1}
  \end{subfigure}%
  \ifmaketwocolcsm
  \ifmakearxivcsm
  \vspace{3mm}
  \fi
  \else
  \hfill%
  \fi
  \begin{subfigure}[b]{\lcvxtoyplotsz}
    \centering
    \includegraphics[width=\lcvxtoyinnerszright]{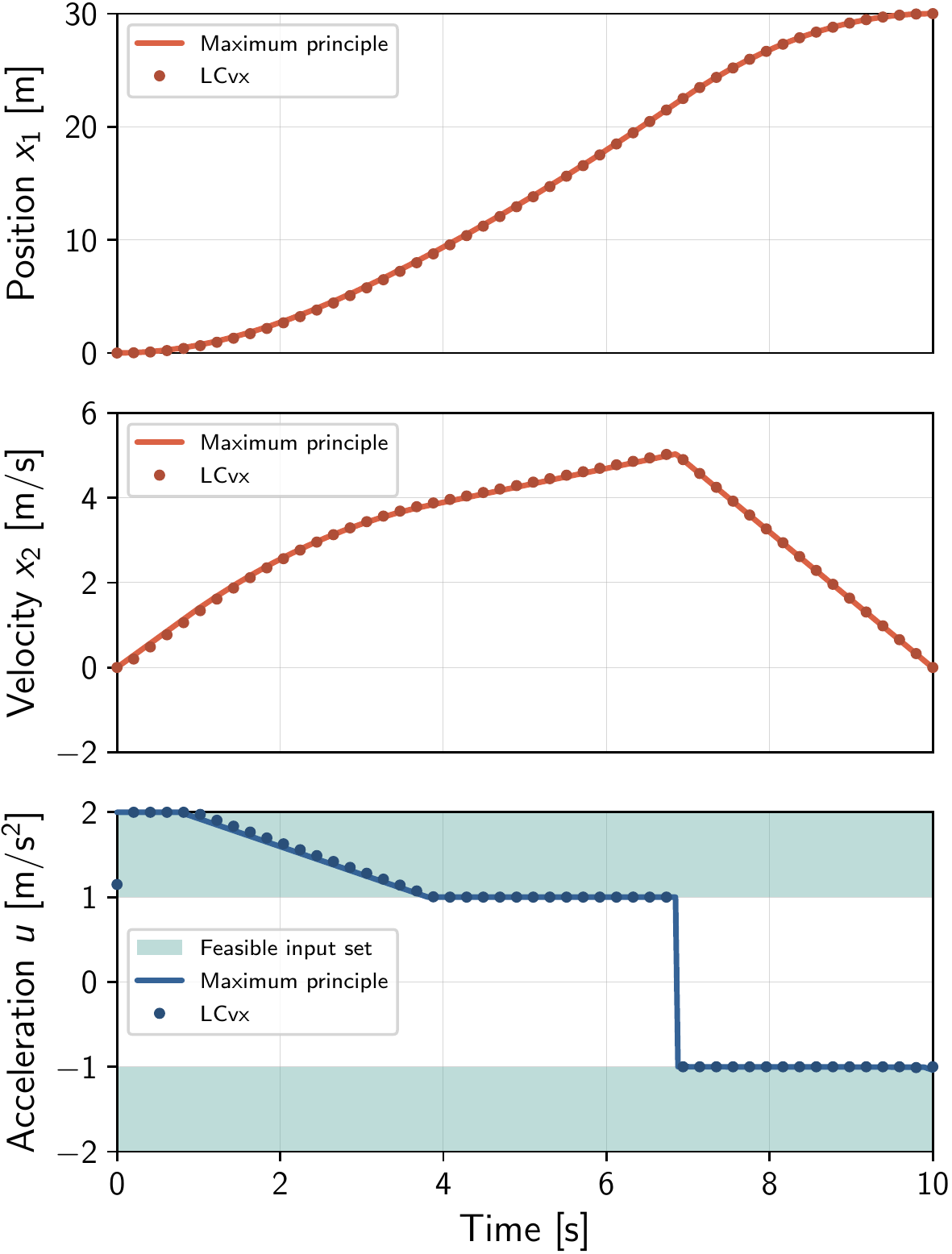}
    \caption{Solution of \pref{lcvx_r_toy} for $g=0.6$~m/s$^2$ and $s=30$~m.}
    \label{fig:lcvx_toy_2}
  \end{subfigure}%
\end{csmfigure}

Looking at the solutions in \figref{lcvx_toy} for two values of the friction
parameter $g$, we can see that the nonconvex constraint
\optieqref{lcvx_o_toy}{bounds} holds in both cases. We emphasize that this is
despite the trajectories in \figref{lcvx_toy} coming from the solution of
\pref{lcvx_r_toy}, where $|u(t)|<1$ is feasible. The fact that this does not occur
is the salient feature of \lcvx theory, and for this problem it is guaranteed
by \tref{fixed_time}. Finally, we note that \figref{lcvx_toy} also plots the analytical
globally optimal solution obtained via the maximum principle, where no
relaxation nor discretization is made. The close match between this solution
and the numerical \lcvx solution further confirms the theory, as well as the
accuracy of the FOH discretization method. Note that the mismatch at $t=0$ in
the acceleration plot in \figref{lcvx_toy_2} is a benign single-time-step
discretization artifact that is commonly observed in \lcvx numerical solutions.


\section{Part II: Sequential Convex Programming}

We now move on to a different kind of convex optimization\dash based trajectory
generation algorithm, known as sequential convex programming (SCP). The reader
will see that this opens up a whole world of possibilities beyond the
restricted capabilities of lossless convexification. One could say that if
\lcvx is a surgical knife to remove acute nonconvexity, then SCP is a catch-all
sledgehammer for nonconvex trajectory design \cite{MalyutaARC}.

A wealth of industrial and research applications, including high-profile
experiments, support this statement. Examples can be found in many engineering
domains, ranging from
aerospace~\cite{SzmukReynolds2018,Lee2017,Liu2014,Bonalli2019} and mechanical
design~\cite{Kocvara2002,Beck2010} to power grid technology~\cite{Wei2017},
chemical processes~\cite{Biegler2014}, and computer
vision~\cite{Jiang2007,Chan2006}. Just last year, the Tipping Point Partnership
between NASA and Blue Origin started testing an SCP algorithm (that we will
discuss in this section) aboard the New Shepard rocket
\cite{TippingPointFlight,carson2019splice}. Another application of SCP methods
is for the SpaceX Starship landing flip maneuver
\cite{SN10LandingRecap}. Although SpaceX's methods are undiscolsed, we know that
convex optimization is used by the Falcon 9 rocket and that SCP algorithms are
highly capable of solving such challenging trajectories
\cite{lars2016autonomous, CSMCodeMasterBranch}.

Further afield, examples of SCP can be found in
medicine~\cite{Rocha2018,Silva2017}, economics~\cite{Dorfman1969,Caputo2005},
biology~\cite{Banga2008,Goelzer2011}, and fisheries~\cite{Liski2001}. Of course,
in any of these applications, SCP is not the only methodology that can be used
to obtain good solutions. Others might include interior point methods
\cite{BettsBook,Nesterov1994}, dynamic programming \cite{BertsekasDPBook},
augmented Lagrangian techniques \cite{BertsekasConvexAlgoBook}, genetic or
evolutionary algorithms \cite{MitchellGABook}, and machine learning and neural
networks \cite{SuttonRLBook}, to name only a few. However, it is our view that
SCP methods are fast, flexible and efficient local optimization algorithms for
trajectory generation. They are a powerful tool to have in a trajectory
engineer's toolbox, and they will be the focus of this part of the article.

\begin{csmfigure}[%
  caption={%
    Block diagram illustration of a typical SCP algorithm. Every SCP-based
    trajectory generation method is comprised of three major components: a way
    to guess the initial trajectory ({\color{beamerBlue}Starting}), an
    iteration scheme which refines the trajectory until it is feasible and
    locally optimal ({\color{beamerRed}Iteration}), and an exit criterion to
    stop once the trajectory has been computed
    ({\color{beamerGreen}Stopping}). In a well-designed SCP scheme, the test
    (convergence) criterion is guaranteed to trigger, but the solution may be
    infeasible for the original problem.
  },%
  label={scp_loop},position=t,columns=2]%
  \centering%
  \ifmaketwocolcsm
  \includegraphics[scale=1.1]{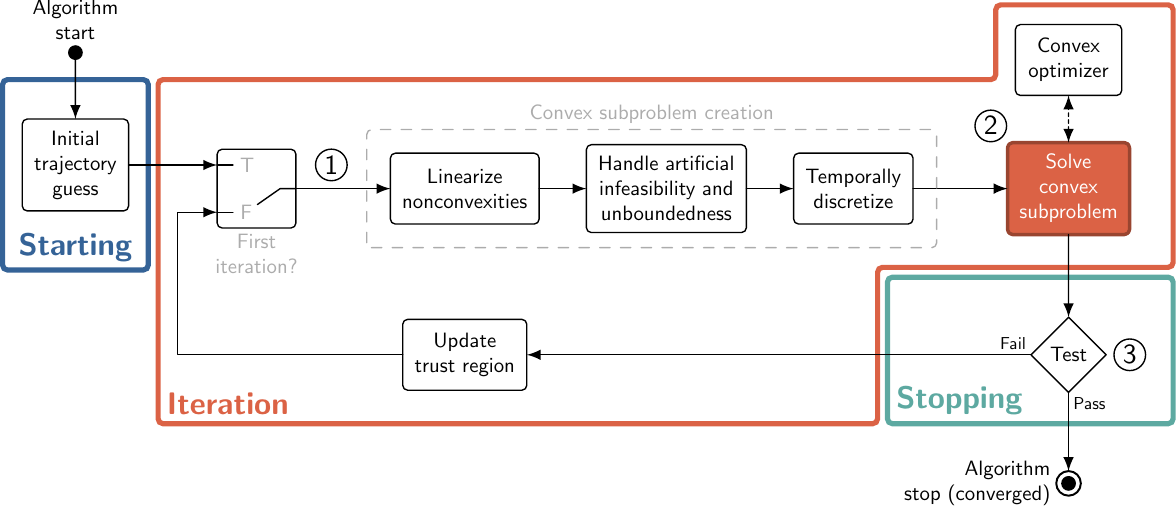}
  \else
  \includegraphics[scale=1.2]{scvx_loop}
  \fi
\end{csmfigure}

As the name suggests, at the core of SCP is the idea of iterative convex
approximation. Most, if not all, SCP algorithms for trajectory generation can
be cast in the form illustrated by \figref{scp_loop}. Strictly speaking, SCP
methods are nonlinear local optimization algorithms. In particular, the reader
will recognize that SCP algorithms are specialized trust region methods for
continuous-time optimal control problems
\cite{ConnTrustRegionBook,NocedalBook,Kochenderfer2019}. 

All SCP methods solve a sequence of convex approximations, called
\alert{subproblems}, to the original nonconvex problem and update the
approximation as new solutions are obtained. Going around the loop of
\figref{scp_loop}, all algorithms start with a user-supplied initial guess, which
can be very coarse (more on this later). At \alglocation{\iterstartloc}, the
SCP algorithm has available a so-called reference trajectory, which may be
infeasible with respect to the problem dynamics and constraints. The
nonconvexities of the problem are removed by a local linearization around the
reference trajectory, while convex elements are kept
unchanged. 
Well-designed SCP algorithms add extra features to the problem in order to
maintain subproblem feasibility after linearization. The resulting convex
continuous-time subproblem is then temporally discretized to yield a
finite-dimensional convex optimization problem. The optimal solution to the
discretized subproblem is computed at \alglocation{\solveloc}, where the SCP
algorithm makes a call to any appropriate 
convex optimization solver. The solution is tested at \alglocation{\testloc}
against stopping criteria. If the test passes, we say that the algorithm has
\alert{converged}, and the most recent 
solution from \alglocation{\solveloc} is returned. Otherwise, the solution is
used to update the trust region (and possibly other parameters) that are
internal to the SCP algorithm. The solution then becomes the new reference
trajectory for the next iteration of the algorithm.

The SCP approach offers two main advantages. First, a wide range of algorithms
exist to reliably solve each convex subproblem at
\alglocation{\solveloc}. Because SCP is agnostic to the particular choice of
subproblem optimizer, well-tested 
algorithms can be used. This makes SCP very attractive for safety\dash critical
applications, which are ubiquitous throughout disciplines like aerospace and
automotive engineering. Second, one can derive meaningful theoretical
guarantees on algorithm performance and computational complexity, as opposed to
general NLP optimization where the convergence guarantees are much
weaker. Taken together, these advantages have led to the development of very
efficient SCP algorithms with runtimes low enough to enable real-time
deployment for some applications \cite{Reynolds2020}.

A fundamental dilemma of NLP optimization is that one can either compute
locally optimal solutions quickly, or globally optimal solutions slowly. SCP
techniques are not immune to this trade-off, despite the fact that certain
subclasses of convex optimization can be viewed as ``easy'' from a
computational perspective due to the availability of interior point
methods. Some of the aforementioned applications may favor the ability to
compute solutions quickly (i.e., in near real-time), such as aerospace and
power grid technologies. Others, such as economics and structural truss design,
may favor global optimality and put less emphasis on solution time (although
early trade studies may still benefit from a fast optimization method). Given
the motivation from the beginning of this article, our focus is on the former
class of algorithms that provide locally optimal solutions in near real-time.

This part of the article will provide an overview of the algorithmic design
choices and assumptions that lead to effective SCP implementations. The core
tradeoffs include how the convex approximations are formulated, what structure
is devised for updating the solutions, how progress towards a solution is
measured, and how all of the above enables theoretical convergence and
performance guarantees.

We hope that the reader comes away with the following view of SCP: it is an
effective and flexible way to do trajectory optimization, and inherits some but
not all of the theoretical properties of convex
optimization. 
SCP works well for complex problems, but it is definitely not a panacea for all
of nonconvex optimization. SCP can fail to find a solution, but usually a
slight change to the parameters recovers convergence. This part of the article
provides the reader with all the necessary insights to get started with
SCP. The numerical examples in Part III provide a practical and open\dash
source implementation of the algorithms herein.

\subsection{Historical Development of SCP}

Tracing the origins of what we refer to as sequential convex programming is not
a simple task. Since the field of nonlinear programming gained traction as a
popular discipline in the 1960s and 70s, many researchers have explored the
solution of nonconvex optimization problems via convex approximations. This
section attempts to catalogue some of the key developments, with a focus on
providing insight into how the field moved toward the present day version of
sequential convex programming for trajectory generation.

The idea to solve a general (nonconvex) optimization problem by iteratively
approximating it as a convex program was perhaps first developed using
branch\dash and\dash bound
techniques~\cite{Falk1969,Soland1971,Horst1976,Horst1984}. Early results were of
mainly academic interest, and computationally tractable methods remained
elusive. One of the most important ideas that emerged from these early
investigations appears to be that of McCormick
relaxations~\cite{McCormick1976}. These are a set of atomic rules for
constructing convex/concave relaxations of a specific class of functions that
everywhere under-/over-estimate the original functions. These rules result in a
class of SCP methods, and algorithms based on McCormick relaxations continue to
be developed with increasing computational
capabilities~\cite{Mitsos2009,Tsoukalas2014,Singer2004,Singer2006}.

\begin{blurb}
  Because SCP is agnostic to the choice of subproblem optimizer, well\dash
  tested algorithms can be used, which is attractive for safety\dash critical
  applications.
\end{blurb}

Difference\dash of\dash convex programming is a related class of SCP
methods~\cite{Horst1999,Lipp2016}. These types of algorithms rely on the
formulation of nonconvex constraints as the difference between two convex
functions, say $f = f_1 - f_2$, where both $f_1$ and $f_2$ are convex
functions. The advantage of this decomposition is that only the function $f_2$
needs to be linearized in order to approximate the nonconvex function $f$. The
convex-concave procedure presented in~\cite{Yuille2003} is one example of a
successful implementation of this idea, and it has been applied, among other
places, in the field of machine learning to support vector machines and
principal component analysis~\cite{Lanckriet2009}.

Perhaps the earliest and simplest class of SCP methods whose structure
resembles that shown in \figref{scp_loop} is, unsurprisingly, sequential
\alert{linear} programming (SLP). These algorithms linearize all nonlinear
functions about a current reference solution so that each subproblem is a
linear program. These linear programs are then solved with a trust region to
obtain a new reference, and the process is repeated. Early developments came
from the petroleum industry and were intended to solve large\dash scale
problems~\cite{Palacios-Gomez1982,Byrd2003}. From a computational perspective,
SLP was initially attractive due to the maturity of the simplex algorithm. Over
time, however, solvers for more general classes of convex optimization problems
have advanced to the point that restricting oneself to linear programs to save
computational resources at \alglocation{\solveloc} in \figref{scp_loop} has become
unnecessary, except perhaps for very large\dash scale problems.

Another important class of SCP methods is that of sequential quadratic
programming (SQP). The works of Han~\cite{Han1977,Han1979},
Powell~\cite{Powell1978,Powell1978a,Powell1986}, Boggs and
Tolle~\cite{Boggs1982,Boggs1984,Boggs1989}, and Fukushima~\cite{Fukushima1986}
appear to have exerted significant influence on the early developments of
SQP-type algorithms, and their impact remains evident today. An excellent
survey was written by Boggs and Tolle~\cite{Boggs1995} and an exhaustive
monograph is available by Conn, Gould, and Toint \cite{ConnTrustRegionBook}. SQP
methods approximate a nonconvex program with a quadratic program using some
reference solution, and then use the solution to this quadratic program to
update the approximation. Byrd, Schnabel, and Schultz provide a general theory
for inexact SQP methods \cite{Byrd1988}. The proliferation of SQP-type
algorithms can be attributed to three main aspects: 1) their similarity with
the familiar class of Newton methods, 2) the fact that the initial reference
need not be feasible, and 3) the existence of algorithms to quickly and
reliably solve quadratic programs. In fact, the iterates obtained by SQP
algorithms can be interpreted either as solutions to quadratic programs or as
the application of Newton's method to the optimality conditions of the original
problem~\cite{NocedalBook}. SQP algorithms are arguably the most mature class of
SCP methods~\cite{Gill2005}, and modern developments continue to address both
theoretical and applied aspects~\cite{Betts1993,Lawrence2001}.


Their long history of successful deployment in NLP solvers
notwithstanding~\cite{Gill2005}, SQP methods do come with several
drawbacks. Gill and Wong nicely summarize the difficulties that can arise when
using SQP methods~\cite{Gill2012}, and we will only outline the basic ideas
here. Most importantly (and this goes for any ``second-order'' method), it is
difficult to accurately and reliably estimate the Hessian of the nonconvex
program's Lagrangian. Even if this is done, say, analytically, there is no
guarantee that it will be positive semidefinite, and an indefinite Hessian
results in an NP-hard nonconvex quadratic program. Hessian approximation
techniques must therefore be used, such as keeping only the positive
semidefinite part or using the BFGS update~\cite{Fletcher2013}. In the latter
case, additional conditions must be met to ensure that the Hessian remains
positive-definite. These impose both theoretical and computational challenges
which, if unaddressed, can both impede convergence and curtail the real-time
applicability of an SQP-type algorithm. Fortunately, a great deal of effort has
gone into making SQP algorithms highly practical, resulting in mature algorithm
packages like SNOPT~\cite{Gill2005}.

One truly insurmountable drawback of SQP methods for trajectory generation in
particular is that quadratic programs require all constraints to be affine in
the solution variable. Alas, many motion planning problems are naturally
subject to non\dash affine convex constraints. We have already seen an example
of a second-order cone constraint that arises from a spacecraft glideslope
requirement in \sbref{sidebar_affinestate}, shown
in~\figref{landing_glideslope_cyclic_shift}. 
For problems with non\dash affine convex constraints, the use of an SQP
algorithm may require more iterations to converge compared to a more general
SCP algorithm, leading to a reduction in computational efficiency. Moreover,
each SQP iterate is not guaranteed to be feasible with respect to the original
convex constraints, whereas the SCP iterates will be.

There are several classes of SCP algorithms that generalize the idea of SQP in
order to deal with exactly this limitation. Semidefinite programs are the most
general class of convex optimization problems for which efficient off-the-shelf
solvers are available. Fares et al. introduced sequential semidefinite
programming~\cite{Fares2002}, which uses matrix variables that are subject to
definiteness constraints. Such algorithms find application most commonly in
robust control, where problems are formulated as (nonconvex) semidefinite
programs with linear and bilinear matrix inequalities~\cite{BoydLMI}. Recent
examples have appeared for robust planetary rocket landing
\cite{Reynolds2021Funnel, ReynoldsThesis}. We can view sequential semidefinite
programming as the furthest possible generalization of SLP to the idea of
exploiting existing convexity in the subproblems.

This article focuses on the class of SCP methods that solve a general convex
program at each iteration, without a priori restriction to one of the
previously mentioned classes of convex programs (e.g., LPs, QPs, SOCPs, and
SDPs). This class of SCP methods has been developed largely over the last
decade, and represents the most active area of current development, with
successful applications in robot and spacecraft trajectory
optimization~\cite{Schulman2014,Liu2014,Liu2015,Liu2016,Lee2017,%
  SzmukReynolds2018,Reynolds2020,Reynolds2019b,Sagliano2017,%
  Simplicio2019,MalyutaARC}. We focus, in particular, on two specific
algorithms within this class of SCP methods: \scvx and \gusto. These two
algorithms are complementary in a number of ways, and enjoy favorable
theoretical guarantees. On the one hand, the theoretical analysis of \scvx
works with the temporally discretized problem and provides guarantees in terms
of the Karush-Kuhn-Tucker (KKT) optimality
conditions~\cite{Mao2016,Mao2017,Mao2018}. On the other hand, \gusto is analyzed
for the continuous\dash time problem and provides theoretical guarantees in
terms of the Pontryagin maximum principle
\cite{Bonalli2019a,Bonalli2019b,PontryaginBook,BerkovitzBook}. The numerical
examples in Part III of this article are solved using both \scvx and \gusto
exactly as they are presented here. These examples illustrate that the methods
are, to some degree, interchangeable.

\begin{csmfigure}[%
  caption={%
    Sequential convex programming can be placed atop a hierarchy of classical
    optimization algorithms. In this illustration, the ``width'' of each layer
    is representative of the corresponding algorithm's implementation and
    runtime complexity (to be used only as an intuitive guide). Each layer
    embeds within itself the algorithms from the layers below it.
  },%
  label={opti_alg_hierarchy}]%
  \centering%
  \ifmaketwocolcsm
  \includegraphics{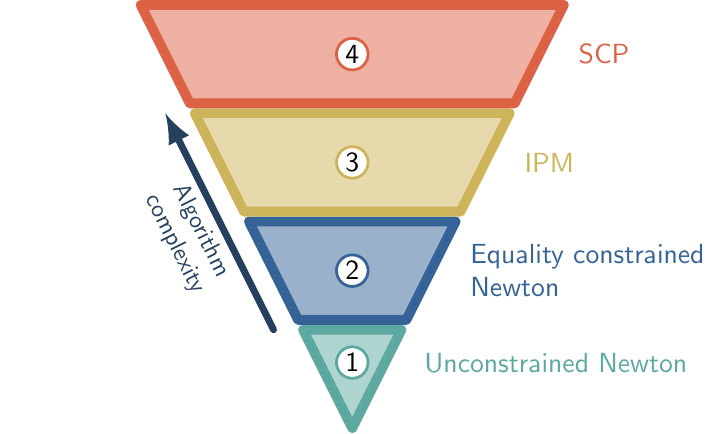}
  \else
  \includegraphics[scale=\csmpreprintfigscale]{opti_alg_hierarchy}
  \fi
\end{csmfigure}

Typically, although not necessarily, the convex solver used at
\alglocation{\solveloc} in \figref{scp_loop} is based on an interior point method
\cite{Nesterov1994,WrightIPMBook}. This leads to a nice interpretation of SCP as
the ``next layer up'' in a hierarchy of optimization algorithms described in
\cite[Chapter~11]{BoydConvexBook}, and which we illustrate in
\figref{opti_alg_hierarchy}. In the bottommost layer, we have the unconstrained
Newton's method, which solves a sequence of unconstrained QPs. The next layer
solves linear equality constrained convex problems. This again uses Newton's
method, but with a more complicated step
computation. 
The third layer is the IPM family of methods, which solve a convex problem with
linear equality and convex inequality constraints as a sequence of linear
equality constrained problems. Thus, we may think of IPMs as iteratively
calling the algorithm in layer \alglayer{beamerBlue}{2} of
\figref{opti_alg_hierarchy}. Analogously, SCP solves a nonconvex problem as a
sequence of convex problems with linear equality and convex inequality
constraints. Thus, SCP iteratively calls an IPM algorithm from layer
\alglayer{beamerYellow!85!black}{3}. Numerical experience has shown that for
most problems, IPMs require on the order of tens of iterations (i.e., calls to
layer \alglayer{beamerBlue}{2}) to converge \cite{BoydConvexBook}. Similarly,
our experience has been that SCP requires on the order of tens of iterations
(i.e., calls to layer \alglayer{beamerYellow!85!black}{3}) to converge.

The rest of this part of the article is organized as follows. We first state
the general continuous-time optimal control problem that we wish to solve, and
discuss the common algorithmic underpinnings of the SCP framework. We then
describe the \scvx and \gusto algorithms in full detail. At the end of Part II,
we compare \scvx and \gusto and give some advice on using SCP in the real
world. Part III will present two numerical experiments that provide practical
insight and highlight the capabilities of each algorithm.


\subsection{Problem Formulation}

The goal of SCP methods is to solve continuous-time optimal control problems of
the following form:
\begin{optimization}[%
  label={scp_gen_cont},%
  variables={u,p},%
  objective={\ocost(x,u,p)}]%
  \optilabel{dynamics}%
  & \dot{x}(t) = f\big(t,x(t),u(t),p\big), \\
  \optilabel{convex_path_constraints_X}%
  & \pare[big]{x(t),p} \in \set{X}(t), \\
  \optilabel{convex_path_constraints_U}%
  &\pare[big]{u(t),p} \in \set{U}(t), \\
  \optilabel{nonconvex_constraints}%
  & s\big(t,x(t),u(t),p\big) \leq 0, \\
  \optilabel{initial_conditions}%
  & \gic\big(x(0),p\big) = 0, \\
  \optilabel{final_conditions}%
  & \gtf\big(x(1),p\big) = 0,
\end{optimization}
where $x(\cdot)\in\reals^n$ is the state trajectory, $u(\cdot)\in\reals^m$ is
the control trajectory, and $\pk\in\real^d$ is a vector of parameters. The
function
$f : \real\times\real^n \times \real^m \times \real^d \rightarrow \real^n$
represents the (nonlinear) dynamics, which are assumed to be at least once
continuously differentiable. Initial and terminal boundary conditions are
enforced by using the continuously differentiable functions
$\gic:\real^n \times \real^d \rightarrow \real^{\dimgic}$ and
$\gtf : \real^n \times \real^d \rightarrow \real^{\dimgtf}$. We separate convex
and nonconvex path (i.e., state and control) constraints by using the convex
sets $\set{X}(t)$ and $\set{U}(t)$ to represent convex path constraints, and
the continuously differentiable function
$s : \real \times \real^n \times \real^m \times \real^d \rightarrow \real^{\dimss}$
to represent nonconvex path constraints. It is assumed that the sets
$\set{X}(t)$ and $\set{U}(t)$ are compact (i.e., closed and bounded). This
amounts to saying that the vehicle cannot escape to infinity or apply infinite
control action, which is obviously reasonable for all practical
applications. Finally, note that \pref{scp_gen_cont} is defined on the $[0,1]$ time
interval, and the constraints
\optieqref{scp_gen_cont}{dynamics}-\optieqref{scp_gen_cont}{nonconvex_constraints} have
to hold at each time instant.

We highlight that the parameter vector $\pk$ can be used, among other things,
to capture free initial and/or free final time problems by making $\tio$ and
$\tf$ elements of $\pk$. In particular, an appropriate scaling of time can
transform the $[0,1]$ time interval in \pref{scp_gen_cont} into a $[\tio,\tf]$
interval. This allows us to restrict the problem statement to the $[0,1]$ time
interval without loss of generality~\cite{Betts1998}. We make use of this
transformation in the numerical examples at the end of the article.

Hybrid systems like bouncing balls, colliding objects, and bipedal robots
require integer variables in their optimization models. The integer variable
type, however, is missing from \pref{scp_gen_cont}. Nevertheless, methods exist to
embed integer variables into the continuous\dash variable formulation. Among
these methods are state-triggered constraints
\cite{SzmukThesis,SzmukReynolds2018,Reynolds2019b,szmuk2019successive,
  szmuk2019real,malyuta2020fast,MalyutaJGCD}, and homotopy techniques such as
the relaxed autonomously switched hybrid system and composite smooth control
\cite{saranathan2018relaxed,taheri2020novel,MalyutaARC}. We shall therefore move
forward using \pref{scp_gen_cont} ``without loss of generality'', keeping in mind
that there are methods to embed integer solution variables exactly or as an
arbitrarily accurate approximation \cite{MalyutaARC}.

We also take this opportunity to note that \pref{scp_gen_cont} is not the most
general optimal control problem that SCP methods can solve. However, it is
general enough for the introductory purpose of this article, and can already
cover the vast majority of trajectory optimization problems
\cite{MalyutaARC}. The numerical implementation attached to this article (see
\figref{github_qr}) was applied to solve problems ranging from quadrotor trajectory
generation to spacecraft rendezvous and docking
\cite{malyuta2020fast,MalyutaJGCD}.

The cost function in~\optiobjref{scp_gen_cont} is assumed to be of the Bolza form
\cite{LiberzonBook}:
\begin{equation}
  \label{eq:ocost_nlin}
  \ocost(x,u,p) = \term(x(1),p) + \int_0^{1} \runn(x(t),u(t),p)\sdt,
\end{equation}
where the terminal cost $\term:\reals^n\times\reals^d\to\reals$ is assumed to
be a convex function and the running cost
$\runn:\reals^n\times\reals^m\times\reals^d\to\reals$ can be in general a
nonconvex function. Note that convexity assumptions on $\term$ are without loss
of generality. For example, a nonconvex $\term$ can be replaced by a linear
terminal cost $\tau_{\mathrm{f}}$ (where $\tau_{\mathrm{f}}$ becomes an element
of $p$), and a nonconvex terminal boundary condition is added to the definition
of $\gtf$ in \optieqref{scp_gen_cont}{final_conditions}:
\begin{equation}
  \label{eq:terminal_cost_ncvx_fix}
  \term(x(1),p)=\tau_{\mathrm{f}}.
\end{equation}

\subsection{SCP Algorithm Foundations}

All SCP methods work by solving a sequence of local convex approximations
to~\pref{scp_gen_cont}, which we call subproblems. As shown in \figref{scp_loop}, this
requires having access to an existing reference trajectory at location
\alglocation{\iterstartloc} of the figure. We will call this a \alert{reference
  solution}, with the understanding that this trajectory need not be a feasible
solution to the problem (neither for the dynamics nor for the constraints). SCP
methods update this reference solution after each passage around the loop of
\figref{scp_loop}, with the solution obtained at \alglocation{\solveloc} becoming
the reference for the next iteration. This begs the question: where does the
reference solution for the first iteration come from? 

\subsubsection{Initial Trajectory Guess}

A user-supplied initial trajectory guess is responsible for providing the first
SCP iteration with a reference solution. Henceforth, the notation $\trajohone$
shall denote a reference trajectory on the time interval $[0,1]$.

We will see in the following sections that the SCP algorithms that we discuss,
\scvx and \gusto, are guaranteed to converge almost regardless of the initial
trajectory guess. In particular, this guess can be grossly infeasible with
respect to both the dynamics \optieqref{scp_gen_cont}{dynamics} and the nonconvex
constraints \optieqref{scp_gen_cont}{nonconvex_constraints}\dash
\optieqref{scp_gen_cont}{final_conditions}. However, the algorithms do require the
guess to be feasible with respect to the convex path constraints
\CTNLconvexpath. Assuring this is almost always an easy task, either by
manually constructing a simplistic solution that respects the convex
constraints, or by projecting an infeasible guess onto the $\set{X}(t)$ and
$\set{U}(t)$ sets. For reference, both strategies are implemented in our
open\dash source code linked in \figref{github_qr}.

Numerical experience has shown that both \scvx and \gusto are extremely adept
at morphing coarse initial guesses into feasible and locally optimal
trajectories. This represents a significant algorithmic benefit, since most
traditional methods, like SQP and NLP, require good (or even feasible) initial
guesses, which can be very hard to come by \cite{Betts1998,Kelly2017}.

To give the reader a sense for what kind of initial guess can be provided, we
present an initialization method called \alert{straight\dash line
  interpolation}. We have observed that this technique works well for a wide
variety of problems, and we use it in the numerical examples at the end of this
article. However, we stress that this is merely a rule of thumb and not a
rigorously derived technique.

We begin by fixing the initial and final states $\xxic$ and $\xxfc$ that
represent either single-point boundary conditions or points in a desired
initial and terminal set defined by \optieqref{scp_gen_cont}{initial_conditions}
and \optieqref{scp_gen_cont}{final_conditions}. The state trajectory is then
defined as a linear interpolation between the two endpoints:
\begin{equation}
  \label{eq:state_initial_guess}
  \xb(t) = (1-t)\xxic+t\xxfc,~\textnormal{for}~t\in [0,1].
\end{equation}

If a component of the state is a non-additive quantity, such as a unit
quaternion, then linear interpolation is not the most astute choice. In such
cases, we opt for the simplest alternative to linear interpolation. For unit
quaternions, this would be spherical linear interpolation \cite{Shoemake1985}.

Whenever possible, we select the initial input trajectory based on insight from
the physics of the problem. For example, for an aerial vehicle we would choose
an input that opposes the pull of gravity. In the case of a rocket, the choice
can be $\uuic = -m_{\mathrm{wet}}\gI$ and $\uufc=-m_{\mathrm{dry}}\gI$, where
$m_{\mathrm{wet}}$ and $m_{\mathrm{dry}}$ are the initial and estimated final
masses of the vehicle, and $\gI$ is the inertial gravity vector. If the problem
structure does not readily admit a physics\dash based choice of control input,
our go\dash to approach is to set the input to the smallest feasible value that
is compliant with \optieqref{scp_gen_cont}{convex_path_constraints_U}. The
intuition is that small inputs are often associated with a small cost
\optiobjref{scp_gen_cont}. In any case, the initial control solution is
interpolated using a similar expression to~\eqref{eq:state_initial_guess}:
\begin{equation}
  \label{eq:input_initial_guess}
  \ub(t) = (1-t)\uuic+t\uufc,~\textnormal{for}~t\in [0,1].
\end{equation}

The initial guess for $\pb$ can have a significant impact on the number of SCP
iterations required to obtain a solution. For example, if $\pb$ represents the
final time of a free final time problem that evolves on the $[0,\tf]$ interval,
then it is best to guess a time dilation value that is reasonable for the
expected trajectory. Since parameters are inherently problem specific, however,
it is unlikely that any generic rule of thumb akin to
\eqref{eq:state_initial_guess} and \eqref{eq:input_initial_guess} will prove
reliable. Fortunately, since SCP runtime is usually on the order of a few
seconds or less, the user can experiment with different initial guesses for
$\pb$ and come up with a good initialization strategy relatively quickly.

For all but the simplest problems, the initial guess $\trajohone$ constructed
above is going to be (highly) infeasible with respect to the dynamics and
constraints of \pref{scp_gen_cont}. Nevertheless, SCP methods like \scvx and \gusto
can, and often do, converge to usable trajectories using such a coarse initial
guess. However, this does not relieve the user entirely from choosing an
initial guess that exploits the salient features of their particular problem. A
well-chosen initial guess will (likely) have the following three benefits for
the solution process:
\begin{itemize}
\item It will reduce the number of iterations and the time required to
  converge. This is almost always a driving objective in the design of a
  trajectory optimization algorithm since fast convergence is not only a
  welcome feature but also a hard requirement for onboard implementation in an
  autonomous system;
\item It will encourage the converged solution to be feasible for
  \pref{scp_gen_cont}. As mentioned, SCP methods like \scvx and \gusto will always
  converge to a trajectory, but without a guarantee that the trajectory will be
  feasible for the original problem. The fact that the solution often
  \textit{is} feasible with respect to \pref{scp_gen_cont} is a remarkable
  ``observation'' that researchers and engineers have made, and it is a driving
  reason for the modern interest in SCP methods. Nevertheless, an observation
  is not a proof, and there are limits to how bad an initial guess can be. The
  only rule of thumb that is always valid is that one should embed as much
  problem knowledge as possible in the initial guess;
\item A better initial guess may also improve the converged trajectory's
  optimality. However, the level of optimality is usually difficult to measure,
  because a globally optimal solution is rarely available for the kinds of
  difficult trajectory problems that we are concerned with using
  SCP. Nevertheless, some attempts to characterize the optimality level have
  been made in recent years \cite{Reynolds2020,Reynolds2020c}.
\end{itemize}

\subsubsection{Linearization}

Let us now place ourselves at location \alglocation{\iterstartloc} in
\figref{scp_loop}, and imagine that the algorithm is at some iteration during the
SCP solution process. The first task in the way of constructing a convex
subproblem is to remove the nonconvexities of \pref{scp_gen_cont}. For this
purpose, recall that the algorithm has access to the reference trajectory
$\trajohone$. If we replace every nonconvexity by its first-order approximation
around the reference trajectory, then we are guaranteed to generate convex
subproblems. Furthermore, these are computationally inexpensive to compute
relative to second-order approximations (i.e., those involving Hessian
matrices). As we alluded in the previous section on SCP history, linearization
of all nonconvex elements is not the only choice for SCP -- it is simply a very
common one, and we take it for this article. To formulate the linearized
nonconvex terms, the following Jacobians must be computed (the time argument is
omitted where necessary to keep the notation short):
\begin{subequations}
  \label{eq:scvx_lin_mats}
  \begin{align}
    \label{eq:scvx_lin_mats_d}
    A(t) &\definedas \diff{x}{f}(t, \xb(t),\ub(t),\pb), \\
    \label{eq:scvx_lin_mats_e}
    B(t) &\definedas \diff{u}{f}(t, \xb(t),\ub(t),\pb), \\
    \label{eq:scvx_lin_mats_f}
    F(t) &\definedas \diff{p}{f}(t, \xb(t),\ub(t),\pb), \\
    \label{eq:scvx_lin_mats_g}
    r(t) &\definedas f(t, \xb(t),\ub(t),\pb) - A \xb(t) - B \ub(t) - F\pb, \\
    \label{eq:scvx_lin_mats_h}
    C(t) &\definedas \diff{x}{s}(t, \xb(t),\ub(t),\pb), \\
    \label{eq:scvx_lin_mats_i}
    D(t) &\definedas \diff{u}{s}(t, \xb(t),\ub(t),\pb), \\
    \label{eq:scvx_lin_mats_j}
    G(t) &\definedas \diff{p}{s}(t, \xb(t),\ub(t),\pb), \\
    \label{eq:scvx_lin_mats_k}
    r\der(t) &\definedas s(t, \xb(t),\ub(t),\pb) - C \xb(t) - D \ub(t) - G\pb, \\
    \label{eq:scvx_lin_mats_l}
    H_0 &\definedas \diff{x}{\gic}(\xb(0),\pb), \\
    \label{eq:scvx_lin_mats_m}
    K_0 &\definedas \diff{p}{\gic}(\xb(0),\pb), \\
    \label{eq:scvx_lin_mats_n}
    \ell_0 &\definedas \gic(\xb(0),\pb) - H_0 \xb(0) - K_0 \pb, \\
    \label{eq:scvx_lin_mats_o}
    \Hf &\definedas \diff{x}{\gtf}(\xb(1),\pb), \\
    \label{eq:scvx_lin_mats_p}
    \Kf &\definedas \diff{p}{\gtf}(\xb(1),\pb), \\
    \label{eq:scvx_lin_mats_q}
    \ellf &\definedas \gtf(\xb(1),\pb) - \Hf \xb(1) - \Kf \pb.
  \end{align}
\end{subequations}

These matrices can be used to write down the first-order Taylor series
approximations for each of $f$, $s$, $\gic$, and $\gtf$. Note that we will not
linearize the cost function \eqref{eq:ocost_nlin} at this point, since \scvx and
\gusto make different assumptions about its particular form. Convexification of
the cost function will be tackled separately in later sections on \scvx and
\gusto.

Using the terms in \eqref{eq:scvx_lin_mats}, we obtain the following approximation
of~\pref{scp_gen_cont} about the reference trajectory:
\begin{optimization}[
  label={scp_gen_cvx},
  variables={u,p},
  objective={\ocost(x,u,p)}]%
  \optilabel{dynamics}
  & \dot{x}(t) = A(t)x(t)+B(t)u(t)+F(t)p+r(t),\hspace{-2mm} \\
  \optilabel{convex_path_constraints_X}%
  & \pare[big]{x(t),p} \in \set{X}(t), \\
  \optilabel{convex_path_constraints_U}%
  &\pare[big]{u(t),p} \in \set{U}(t), \\
  \optilabel{convexified_constraints}
  & C(t)x(t)+D(t)u(t)+G(t)p+r\der(t) \leq 0, \\
  \optilabel{initial_conditions}
  & H_0 x(0) + K_0 \pk + \ell_0 = 0, \\
  \optilabel{final_conditions}
  & \Hf x(1) + \Kf \pk + \ellf = 0.
\end{optimization}

\pref{scp_gen_cvx} is convex in the constraints and potentially nonconvex in the
cost. Note that the convex path constraints in \CTNLconvexpath are kept without
approximation. This is a key advantage of SCP over methods like SLP and SQP, as
was discussed in the previous section on SCP history.

Because the control trajectory $u(\cdot)$ belongs to an infinite\dash
dimensional vector space of continuous\dash time functions, \pref{scp_gen_cvx}
cannot be implemented and solved numerically on a digital computer. To do so,
we must consider a finite\dash dimensional representation of the control
function $u(t)$, which can be obtained via temporal discretization or direct
collocation~\cite{Betts1998,Malyuta2019b}. These representations turn the
original infinite\dash dimensional optimal control problem into a finite\dash
dimensional parameter optimization problem that can be solved on a digital
computer.

In general, and rather unsurprisingly, solutions to discretized problems are
only approximately optimal and feasible with respect to the original
problem. In particular, a discrete-time control signal has fewer degrees of
freedom than its continuous\dash time counterpart. Therefore, it may lack the
flexibility required to exactly match the true continuous\dash time optimal
control signal. By adding more temporal nodes, the approximation can become
arbitrarily accurate, albeit at the expense of problem size and computation
time.

Another problem with discretization is that the path constraints are usually
enforced only at the discrete temporal nodes, and not over the entire time
horizon. This can lead to (typically mild) constraint violation between the
discrete-time nodes, although some techniques exist to remedy this artifact
\cite{acikmese2008enhancements,Dueri2017Clipping}.

The bad news notwithstanding, there are well established discretization methods
that ensure exact satisfaction of the original continuous\dash time nonlinear
dynamics \optieqref{scp_gen_cont}{dynamics}. Thus, the discretized solution can
still produce strictly dynamically feasible continuous\dash time
trajectories. We refer the reader
to~\cite{SzmukReynolds2018,Reynolds2019b,Malyuta2019b,MalyutaARC} for detailed
explanations of discretization methods that ensure exact satisfaction of the
continuous\dash time nonlinear dynamics. An introduction to the technique that
we will use for the numerical examples at the end of this article is given in
\sbref{discretization}.

\declaresidebar[%
title={Discretizing Continuous-time Optimal Control Problems},
label={discretization}]{%
  sidebars/discretization.tex}

Our systematic linearization of all nonconvex elements has ensured that
\pref{scp_gen_cvx} is convex in the constraints, which is good news. However,
linearization unsurprisingly has a price. We have inadvertently introduced two
artifacts that must be addressed: artificial unboundedness and artificial
infeasibility.

\subsubsection{Artificial Unboundedness}

Linear approximations are only accurate in a neighborhood around the reference
solution $\trajohone$. Thus, for each $t\in [0,1]$, the subproblem solution
must be kept ``sufficiently close'' to the linearization point defined by the
reference solution. Another reason to not deviate too far from the reference is
that, in certain malicious cases, linearization can render the solution
unbounded below (i.e., the convex cost \optiobjref{scp_gen_cvx} can be driven to
negative infinity). We refer to this phenomenon as \alert{artificial
  unboundedness}. To mitigate this problem and to quantify the meaning of
``sufficiently close'', we add the following trust region constraint:
\begin{align}
  &\delta x(t) = x(t) - \xb(t), \nonumber \\
  &\delta u(t) = u(t) - \ub(t), \nonumber \\
  &\delta p = p - \pb, \nonumber \\
  \nonumber
  &\trx \norm[q]{\delta x(t)} +%
    \tru \norm[q]{\delta u(t)} + \\
  \label{eq:scp_trust_region}
  &\qquad\qquad
    \trp \norm[q]{\delta p} \leq \tr,~\textnormal{for}~t\in [0,1].
\end{align}
for some choice of $q\in\brac{1,2,2^{\scriptscriptstyle +},\infty}$ and
constants $\trx,\tru,\trp\in\{0,1\}$. We use $q=2^{\scriptscriptstyle +}$ to
denote the ability to impose the trust region as the quadratic two-norm
squared. The trust region radius $\tr$ is a fixed scalar that is updated
between SCP iterations (i.e., passages around the loop in \figref{scp_loop}). The
update rule associated with the trust region measures how well the
linearization approximates the original nonconvex elements at each
iterate. This informs the algorithm whether to shrink, grow, or maintain the
trust region radius. SCP methods can be differentiated by how they update the
trust region, and so the trust region update will be discussed separately for
\scvx and \gusto in the following sections.

\figref{scvx_infeas_tr} shows a two\dash dimensional toy problem that exemplifies a
single iteration of an SCP convergence process. In this example, the ``original
problem'' consists of one parabolic (nonconvex) equality constraint
({\color{tay_col_b}blue}), a convex equality constraint
({\color{tay_col_g}green}), and a convex halfspace inequality constraint
(feasible to the left of the vertical {\color{tay_col_r}red} dashed line). The
original problem is approximated about the reference solution $\bar{z}$,
resulting in the {\color{tay_col_b}blue} dash-dot equality constraint and the
same convex equality and inequality constraints. The trust region is shown as
the {\color{tay_col_r}red} circle, and represents the region in which the SCP
algorithm has deemed the convex approximation to be valid. Evidently, if the
new solution $z$ deviates too much from $\bar{z}$, the linear approximation of
the parabola becomes poor. Moreover, had the green equality constraint been
removed, removal of the trust region would lead to artificial unboundedness, as
the cost could be decreased indefinitely.

Clearly, there is another problem with the linearization in \figref{scvx_infeas_tr}
-- the resulting subproblem is infeasible. This is because the
{\color{tay_col_g}green} and {\color{tay_col_b}blue} dash-dot equality
constraints do not intersect inside the set defined by the trust region and the
convex inequality constraint halfspace. This issue is known as artificial
infeasibility.

\begin{csmfigure}[%
  caption={%
    A two\dash dimensional nonconvex toy problem that exemplifies a convex
    subproblem obtained during an SCP iteration. In this case, the cost
    function $L(z)=-z_2$ and level curves of the cost are shown as
    {\color{gray}gray} dashed lines. The {\color{tay_col_b}blue} curve
    represents a nonconvex equality constraint, and its linearization is shown
    as the {\color{tay_col_b}blue} dash-dot line. Another convex equality
    constraint is shown in {\color{tay_col_g}green}, and a convex inequality
    constraint is shown as the vertical {\color{tay_col_r}red} dashed line. The
    trust region is the {\color{tay_col_r}red} circle centered at the
    linearization point $\bar{z}$, and has radius~$\tr$. The optimal solution
    of the original (non-approximated) problem is shown as the black
    square. The convex subproblem is artificially infeasible. Without the trust
    region and the {\color{tay_col_g}green} constraint, it would also be
    artifically unbounded.%
  },%
  label={scvx_infeas_tr},position=t]%
  \centering
  \ifmaketwocolcsm
  \includegraphics[width=\textwidth]{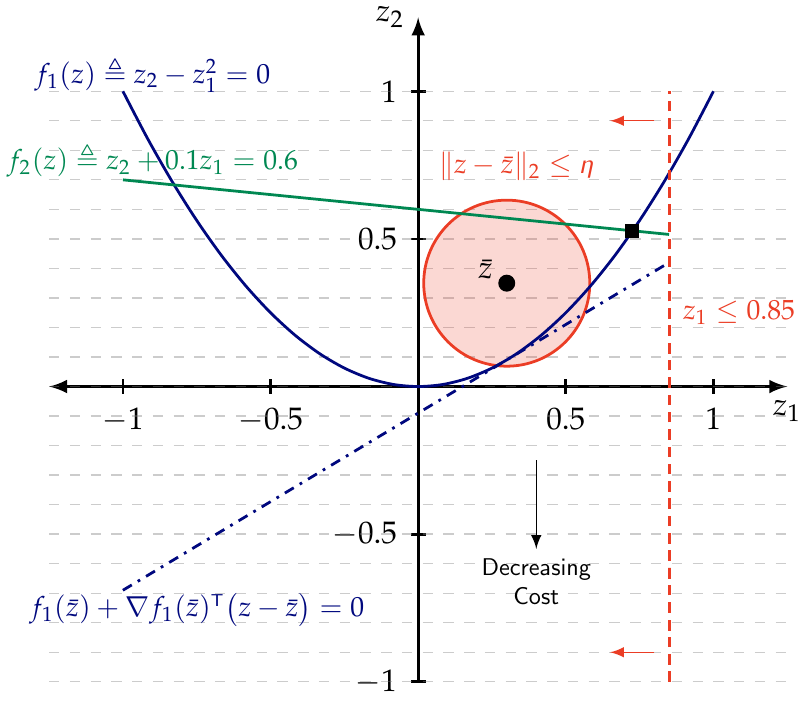}
  \else
  \includegraphics[width=0.8\textwidth]{tikz_scvx_infeasible_tr_depiction}
  \fi
\end{csmfigure}

\subsubsection{Artificial Infeasibility}

Linearization can make the resulting subproblem infeasible. Two independent
cases can arise wherein the constraints imposed in~\pref{scp_gen_cvx} are
inconsistent (i.e., no feasible solution exists), even though the original
constraints admit a non-empty feasible set:
\begin{itemize}
\item In the first case, the intersection of the convexified path constraints
  may be empty. This occurs in the example of~\figref{scvx_infeas_tr}, where no
  feasible solution exists because the linearized constraints
  ({\color{tay_col_g}green} and {\color{tay_col_b}blue} dash-dot) do not
  intersect to the left of the {\color{tay_col_r}red} inequality constraint;
\item In the second case, the trust region may be so small that it restricts
  the solution variables to a part of the solution space that is outside of the
  feasible set. In other words, the intersection of the trust region with the
  (non-empty) feasible region of the convexified constraints may itself be
  empty. This would have been the case in \figref{scvx_infeas_tr} if the
  {\color{tay_col_g}green} and {\color{tay_col_b}blue} dash-dot lines were to
  intersect outside of the {\color{tay_col_r}red} trust region circle, but to
  the left of the halfspace inequality (for example, if $\bar z$ was slightly
  further to the right).
\end{itemize}

The occurence of either case would prevent SCP from finding a new reference
solution at \alglocation{\solveloc} in \figref{scp_loop}. Thus, the solution loop
cannot continue, and the algorithm will fail. As a result, even if the original
problem admits a feasible solution (shown as the black square in
\figref{scvx_infeas_tr}), either of the two aforementioned scenarios would prevent
SCP from finding it. We refer to this phenomenon as \alert{artificial
  infeasibility}.

Artificial infeasibility in sequential convex programming was recognized early
in the development of SCP algorithms~\cite{Powell1978,Byrd2003}. 
Two equivalent strategies exist to counteract this issue. One approach adds an
unconstrained, but penalized, slack variable to each linearized
constraint. This variable is sometimes called a \alert{virtual control} when
applied to the dynamics constraint \optieqref{scp_gen_cont}{dynamics}, and a
virtual buffer when applied to other constraints \cite{Mao2018}. To keep our
language succinct, we will use virtual control in both cases.

The second approach penalizes constraint violations by augmenting the original
cost function with \alert{soft penalty} terms \optiobjref{scp_gen_cvx}. When the
same functions and weights are used as to penalize the virtual control terms,
this strategy results in the same optimality conditions. The ultimate result of
both strategies is that subproblem is guaranteed to be feasible.

Because the virtual control is a new term that we add to the problem, it
follows that the converged solution must have a zero virtual control in order
to be feasible and physically meaningful with respect to the original
problem. If, instead, the second strategy is used and constraint violations are
penalized in the cost, the converged solution must not violate the constraints
(i.e., the penalty terms should be zero). Intuitively, if the converged
solution uses a non-zero virtual control or has a non-zero constraint violation
penalty, then it is not a solution of the original optimal control problem.

The fact that trajectories can converge to an infeasible solution is one of the
salient limitations of SCP. However, it is not unlike the drawback of any NLP
optimization method, which may fail to find a solution entirely even if one
exists. When SCP converges to an infeasible solution, we call it a ``soft''
failure, since usually only a few virtual control terms are non-zero. A soft
failure carries important information, since the temporal nodes and constraints
with non-zero virtual control hint at how and where the solution is
infeasible. Usually, relatively mild tuning of the algorithm parameters or the
problem definition will recover convergence to a feasible solution. In relation
to the optimization literature at large, the soft failure exhibited by SCP is
related to one-norm regularization, lasso regression, and basis pursuit, used
to find sparse approximate solutions \cite[Section~11.4.1]{BoydConvexBook}.

The specific algorithmic choices made to address artificial unboundedness
(e.g., selection of the trust region radius update rule) and artificial
infeasibility (e.g., virtual control versus constraint penalization) lead to
SCP algorithms with different characteristics. The next two sections review two
such methods, \scvx and \gusto, and highlight their design choices and
algorithmic properties. To facilitate a concise presentation, we henceforth
suppress the time argument $t$ whenever possible.

\subsection{The \scvx Algorithm}

In light of the above discussion, the \scvx algorithm makes the following
algorithmic choices:
\begin{itemize}
\item The terminal and running costs in \eqref{eq:ocost_nlin} are both assumed to
  be convex functions. We already mentioned that this is without loss of
  generality for the terminal cost, and the same reasoning applies for the
  running cost. Any nonconvex term in the cost can be offloaded into the
  constraints, and an example was given in \eqref{eq:terminal_cost_ncvx_fix};
\item To handle artificial unboundedness, \scvx enforces
  \eqref{eq:scp_trust_region} as a hard constraint. While several choices are
  possible \cite{Mao2016, Mao2018}, this article uses $\trx=\tru=\trp=1$. The
  trust region radius $\eta$ is adjusted at each iteration by an update rule,
  which we discuss below;
\item To handle artificial infeasibility, \scvx uses virtual control terms.
\end{itemize}

Let us begin by describing how \scvx uses virtual control to handle artificial
infeasibility. This is done by augmenting the linear
approximations~\optieqref{scp_gen_cvx}{dynamics}
and~\optieqref{scp_gen_cvx}{convexified_constraints}-\optieqref{scp_gen_cvx}{final_conditions}
as follows:
\begin{subequations}
  \label{eq:scvx_lin_approxs_vc}
  \begin{align}
    \label{eq:scvx_lin_approxs_vc_f}
    \dot{x} &= A x + B u + F \pk + r + E \vc, \\
    \label{eq:scvx_lin_approxs_vc_s}
    \vcc &\geq C x + D u + G \pk + r\der, \\
    \label{eq:scvx_lin_approxs_vc_gic}
    0 &= H_0 x(0) + K_0 \pk + \ell_0 + \vccic, \\
    \label{eq:scvx_lin_approxs_vc_gtc}
    0 &= \Hf x(1) + \Kf \pk + \ellf + \vcctf.
  \end{align}
\end{subequations}
where $\vc(\cdot)\in\reals^{\dimvc}$, $\vcc(\cdot)\in\reals^{\dimss}$,
$\vccic\in\reals^{\dimgic}$, and $\vcctf\in\reals^{\dimgtf}$ are the virtual
control terms. To keep notation manageable, we will use the symbol $\vcany$ as
a shorthand for the argument list $(\vc,\vcc,\vccic,\vcctf)$.

The virtual control $\vc$ in \eqref{eq:scvx_lin_approxs_vc_f} can be viewed simply
as another control variable that can be used to influence the state
trajectory. Like the other virtual control terms, we would like $\vc$ to be
zero for any converged trajectory, because it is a synthetic input that cannot
be used in reality. Note that it is required that the pair $\pare{A,E}$ in
\eqref{eq:scvx_lin_approxs_vc_f} is controllable, which is easy to verify using any
of several available controllability tests \cite{Antsaklis2006}. A common choice
is $E=I_n$, in which case $\pare{A,E}$ is unconditionally controllable.

The use of non-zero virtual control in the subproblem solution is discouraged
by augmenting the cost function with a virtual control penalty
term. Intuitively, this means that virtual control is used only when it is
necessary to avoid subproblem infeasibility. To this end, we define a positive
definite penalty function $P : \real^{n} \times \real^p \rightarrow \nonneg$
where $p$ is any appropriate integer. The following choice is typical in
practice:
\begin{equation}
  \label{eq:P_penalty_def}
  P(y, z) = \norm[1]{y}+\norm[1]{z},
\end{equation}
where $y$ and $z$ are placeholder arguments. The cost
function~\optiobjref{scp_gen_cvx} is augmented using the penalty function as
follows:
\begin{align}
  \nonumber
  \ocostw(x,u,p,\vcany)
  &\definedas \termw(x(1),p,\vccic,\vcctf)+\\
  \label{eq:scvx_Lpen}
  &\pushright{\int_0^1\runnw(x,u,p,E\vc,\vcc)\sdt,\qquad\quad} \\
  \nonumber
  \termw(x(1),p,\vccic,\vcctf) &= \term(x(1),p)+\Jw P(0,\vccic)+
                                 \Jw P(0,\vcctf), \\
  \label{eq:scvx_Lrunnw}
  \runnw(x,u,p,E\vc,\vcc) &= \runn(x,u,p)+\lambda P\pare[big]{E \vc, \vcc}.
\end{align}

The positive weight $\lambda\in\pos$ is selected by the user, and must be
sufficiently large. We shall make this statement more precise later in the
section on \scvx convergence guarantees. For now, we note that in practice it
is quite easy to find an appropriately large $\lambda$ value by selecting a
power of ten. In general, $\lambda$ can be a function of time, but we rarely do
this in practice.

We also point to an important notational feature of \eqref{eq:scvx_Lrunnw}, where
we used $E\vc$ in place of $\vc$ for the argument list of $\runnw$. This will
help later on to highlight that the continuous\dash time matrix $E$ is
substituted with its discrete\dash time version after temporal discretization
(see ahead in \eqref{eq:scvx_costs_runnw_trapz}).

The continuous-time convex subproblem that is solved at each iteration of the
\scvx algorithm can then be stated formally as:
\begin{optimization}[
  label={subproblem_scvx_ct},
  variables={u,p,\vcany},
  objective={\ocostw(x,u,p,\vcany)}]%
  \optilabel{dynamics}
  & \dot{x} = A x + B u + F \pk + r + E \vc, \\
  \optilabel{xu_constraints_1}
  & (x, p) \in \set{X}, \quad (u, p) \in\set{U}, \\
  \optilabel{xu_constraints_2}
  & C x + D u + G \pk + r\der \leq \vcc, \\
  & H_0 x(0) + K_0 \pk + \ell_0 + \vccic = 0, \\
  & \Hf x(1) + \Kf \pk + \ellf + \vcctf = 0, \\
  \optilabel{xu_constraints_3}
  & \norm[q]{\delta x}+\norm[q]{\delta u}+\norm[q]{\delta p} \leq \tr.
\end{optimization}

It was mentioned in the previous section that \pref{subproblem_scvx_ct} is not
readily implementable on a computer because it is a continuous-time, and hence
infinite-dimensional, optimization problem. To solve the problem numerically, a
temporal discretization is applied, such as the one discussed in
\sbref{discretization}. In particular, we select a set of temporal nodes
$t_k\in[0,1]$ for $k=1,\ldots,N$ and recast the subproblem as a parameter
optimization problem in the (overloaded) variables $x = \{\xk\}_{k=1}^{N}$,
$u = \{\uk\}_{k=1}^{N}$, $p$, $\vc=\brac{\vck}_{k=1}^{N-1}$,
$\vcc=\brac{\vcck}_{k=1}^N$, $\vccic$, and $\vcctf$.

Depending on the details of the discretization scheme, there may be fewer
decision variables than there are temporal nodes. In particular, for simplicity
we will use a zeroth-order hold (ZOH) assumption (i.e., a piecewise constant
function) to discretize the dynamics virtual control $\vc(\cdot)$. This means
that the virtual control takes the value $\vc(t)=\vck$ inside each time
interval $[t_k,t_{k+1})$, and the discretization process works like any other
interpolating polynomial method from \sbref{discretization}. Because this leaves
$\vc_N$ undefined, we take $\vc_N=0$ for notational convenience whenever it
appears in future equations.

The continuous-time cost function~\eqref{eq:scvx_Lpen} can be discretized using any
appropriate method. Pseudospectral methods, for example, specify a numerical
quadrature that must be used to discretize the integral \cite{Kelly2017}. For
simplicity, we shall assume that the time grid is uniform (i.e.,
$t_{k+1}-t_k=\timeintvl$ for all $k=1,\dots,N-1$) and that trapezoidal
numerical integration is used. This allows us to write the discrete-time
version of \eqref{eq:scvx_Lpen} as:
\begin{align}
  \label{eq:scvx_costs_L}
  \mflw(x,u,p,\vcany)
  &= \termw\pare[big]{x(1),p,\vccic,\vcctf}+\mathtt{trapz}(\runnw^N), \\
  \label{eq:scvx_costs_runnw_trapz}
  \runnw[,k]^N &= \runnw(\xk,\uk,p,E_k\vck,\vcck).
\end{align}
where trapezoidal integration is implemented by the function
$\mathtt{trapz}:\reals^N\to\reals$, defined as follows:
\begin{equation}
  \label{eq:trapz}
  \mathtt{trapz}(z) \definedas \frac{\timeintvl}{2}
  \sum_{k=1}^{N-1}z_k+z_{k+1}.
\end{equation}

We call \eqref{eq:scvx_costs_L} the \alert{linear augmented cost function}. This
name is a slight misnomer, because \eqref{eq:scvx_costs_L} is in fact a general
nonlinear convex function. However, we use the ``linear'' qualifier to
emphasize that the cost relates to the convex subproblem for which all
nonconvexities have been linearized. In particular, the virtual control terms
can be viewed as linear measurements of dynamic and nonconvex path and boundary
constraint infeasibility.

Lastly, the constraints~\optieqref{subproblem_scvx_ct}{xu_constraints_1},
\optieqref{subproblem_scvx_ct}{xu_constraints_2},
and~\optieqref{subproblem_scvx_ct}{xu_constraints_3} are enforced at the discrete
temporal nodes $t_k$ for each $k=1,\ldots,N$. In summary, the following
discrete-time convex subproblem is solved at each \scvx iteration (i.e.,
location \alglocation{\solveloc} in \figref{scp_loop}):%
\begin{optimization}[
  label={subproblem_scvx_dt},
  variables={x,u,p,\vcany},
  objective={\mflw(x,u,p,\vcany)}]%
  \optilabel{dynamics}
  & \xkp = A_k \xk + B_k \uk + F_k \pk + r_k + E_k \vck, \\
  \optilabel{convex_path}
  & (\xk, \pk) \in \set{X}_k, \quad (\uk, \pk) \in\set{U}_k, \\
  \optilabel{nonconvex_path}
  & C_k \xk + D_k \uk + G_k \pk + r_k\der \leq \vcck, \\
  & H_0 x_1 + K_0 \pk + \ell_0 + \vccic = 0, \\
  \optilabel{final_conditions}
  & \Hf x_N + \Kf \pk + \ellf + \vcctf = 0, \\
  \optilabel{trust_region}
  & \norm[q]{\delta \xk}+\norm[q]{\delta \uk}+\norm[q]{\delta p} \leq \tr.
\end{optimization}

We want to clarify that \optieqref{subproblem_scvx_dt}{dynamics} is written as a
shorthand convenience for discretized dynamics, and is not representative of
every possible discretization choice. For example,
\optieqref{subproblem_scvx_dt}{dynamics} is correct if ZOH discretization is
used. However, as specified in \eqref{eq:sidebar_dynamics_dt}, FOH discretization
would lead to the following constraint that replaces
\optieqref{subproblem_scvx_dt}{dynamics}:
\begin{equation*}
  \xkp = A_k \xk + B_k^{\scriptscriptstyle -} \uk +
  B_k^{\scriptscriptstyle +} \ukp + F_k p + r_k + E_k \vck.
\end{equation*}

The most general interpolating polynomial discretization fits into the
following discrete\dash time dynamics constraint:
\begin{equation*}
  \xkp = A_k \xk + \sum_{j=1}^N B_{k}^{j} u_j + F_k p + r_k + E_k \vck,
\end{equation*}
where the $j$ superscript indexes different input coefficient matrices. Other
discretization methods lead to yet other affine equality constraints, all of
which simply replace \optieqref{subproblem_scvx_dt}{dynamics}
\cite{Malyuta2019b}. With this in mind, we will continue to write
\optieqref{subproblem_scvx_dt}{dynamics} for simplicity. Furthermore, it is
implicitly understood that the constraints
\optieqref{subproblem_scvx_dt}{dynamics}-\optieqref{subproblem_scvx_dt}{nonconvex_path}
and \optieqref{subproblem_scvx_dt}{trust_region} hold at each temporal node.

Because \pref{subproblem_scvx_dt} is a finite-dimensional convex optimization
problem, it can be solved to global optimality using an off-the-shelf convex
optimization solver \cite{BoydConvexBook}. We shall denote the optimal solution
by $x^* = \brac{\xk^*}_{k=1}^{N}$, $u^* = \brac{\uk^*}_{k=1}^{N}$, $p^*$,
$\vc^*=\brac{\vck^*}_{k=1}^{N-1}$, $\vcc^*=\brac{\vcck^*}_{k=1}^N$, $\vccic^*$,
and $\vcctf^*$.

\subsubsection{\scvx Update Rule}

At this point, we know how to take a nonconvex optimal control problem like
\pref{scp_gen_cont} and: 1) convexify it to \pref{scp_gen_cvx}, 2) add a trust region
\eqref{eq:scp_trust_region} to avoid artificial unboundedness, 3) add virtual
control terms \eqref{eq:scvx_lin_approxs_vc} to avoid artificial infeasibility, 4)
penalize virtual control usage in the cost \eqref{eq:scvx_Lpen}, and 5) apply
discretization to obtain a finite-dimensional convex \pref{subproblem_scvx_dt}. In
fact, this gives us all of the necessary ingredients to go around the loop in
\figref{scp_loop}, except for one thing: how to update the trust region radius
$\eta$ in \eqref{eq:scp_trust_region}. In general, the trust region changes after
each pass around the loop. In this section, we discuss this missing ingredient.

To begin, we define a linearization accuracy metric called the \alert{defect}:
\begin{equation}
  \label{eq:scvx_defect}
  \defectk \definedas \xkp - \flow(t_k,t_{k+1},x_{k},u,p)
\end{equation}
for $k=1,\dots,N-1$. The function $\flow$ is called the \alert{flow map} and
its role is to ``propagate'' the control input $u$ through the continuous\dash
time nonlinear dynamics \optieqref{scp_gen_cont}{dynamics}, starting at state
$x_k$ at time $t_k$ and evolving until the next temporal grid node $t_{k+1}$
\cite{HirschDynamicsBook}. It is important that the flow map is implemented in a
way that is consistent with the chosen discretization scheme, as defined below.

\begin{definition}[scvx_flow_consistency]
  The flow map $\flow$ in \eqref{eq:scvx_defect} is \alert{consistent} with the
  discretization used for \pref{subproblem_scvx_dt}, if the following equation
  holds for all $k=1,\dots,N-1$:
  \begin{equation}
    \label{eq:consistency}
    \flow(t_k,t_{k+1},\xbk,\ub,\pb) =
    A_k \xbk + B_k \ubk + F_k \pbk + r_k.
  \end{equation}
\end{definition}

\begin{csmfigure}[%
  caption={Illustration of the flow map consistency property in
    \dref{scvx_flow_consistency}. When the flow map is consistent with the
    discretization scheme, the state $\tilde x_{k+1}$ propagated through the
    flow map ({\color{beamerYellow}yellow} circle) and the state $\hat x_{k+1}$
    propagated through the discrete-time linearized update equation (dashed
    {\color{beamerGreen}green} circle) match. When the reference trajectory is
    not dynamically feasible, it generally deviates from the flow map
    trajectory, hence $\tilde x_{k+1}\ne \bar x_{k+1}$.%
  },%
  label={scvx_consistency}]%
  \centering%
  \ifmaketwocolcsm
  \includegraphics{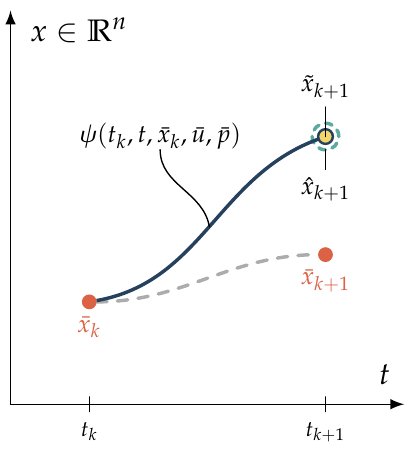}%
  \else
  \includegraphics[scale=1.4]{scvx_consistency}%
  \fi
\end{csmfigure}

There is an intuitive way to think about the consistency property of
$\flow$. The reader may follow along using the illustration in
\figref{scvx_consistency}. On the one hand, $\flow(t_k,t_{k+1},\xbk,\ub,\pb)$ maps
an initial state $\xbk$ through the continuous\dash time nonlinear dynamics
\optieqref{scp_gen_cont}{dynamics} to a new state $\tilde x_{k+1}$. On the other
hand, the right-hand side of \eqref{eq:consistency} does the same, except that it
uses the linearized and discretized dynamics and outputs a new state
$\hat x_{k+1}$. Because the linearization is being evaluated at the reference
trajectory (i.e., at the linearization point), the linearized continuous-time
dynamics will yield the exact same trajectory. Thus, the only difference
between the left- and right-hand sides of \eqref{eq:consistency} is that the
right-hand side works in discrete-time. Consistency, then, simply means that
propagating the continuous-time dynamics yields the same point as performing
the discrete-time update (i.e., $\tilde x_{k+1}=\hat x_{k+1}$). For every
discretization method that is used to construct \pref{subproblem_scvx_dt}, there
exists a consistent flow map.

The reader is likely already familiar with flow maps, even if the term sounds
new. Consider the following concrete examples. When using forward Euler
discretization, the corresponding consistent flow map is simply:
\begin{equation}
  \label{eq:consistent_flow_map_forward_euler}
  \flow(t_k,t_{k+1},x_{k},u,p) = x_k+\timeintvl f(t_k,x_k,u_k,p).
\end{equation}

When using an interpolating polynomial discretization scheme like the one
described in \sbref{discretization}, the corresponding consistent flow map is the
solution to the dynamics \optieqref{scp_gen_cont}{dynamics} obtained through
numerical integration. In other words, the flow map satisfies the following
conditions at each time instant $t\in[t_k,t_{k+1}]$:
\begin{subequations}
  \label{eq:consistent_flow_map_interpolating_poly}
  \begin{align}
    \flow(t_k,t_k,x_{k},u,p)
    &= x_k, \\
    \dot\flow(t_k,t,x_{k},u,p)
    &= f\pare[big]{t,\flow(t_k,t,x_{k},u,p),u,p}.
  \end{align}
\end{subequations}

\begin{csmfigure}[%
  caption={Illustration of defect calculation according to
    \eqref{eq:scvx_defect}. The flow map computation restarts at each discrete-time
    node. At each temporal node, the defect is computed as the difference
    between the discrete solution output by \pref{subproblem_scvx_dt} and the
    corresponding flow map value.%
  },%
  label={scvx_defects}]%
  \centering%
  \ifmaketwocolcsm
  \includegraphics{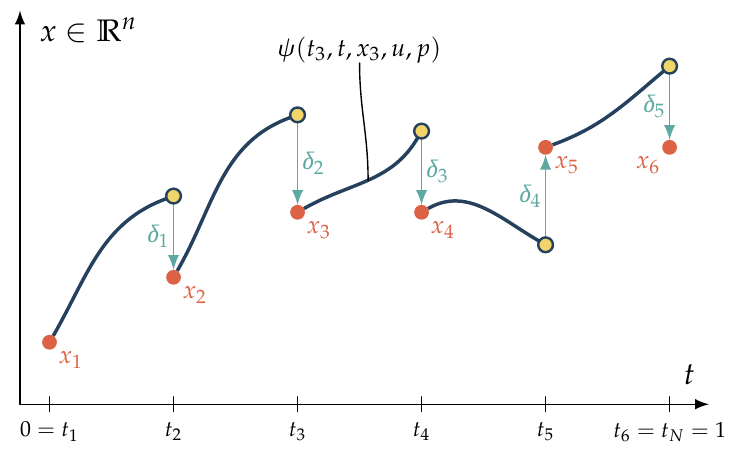}%
  \else
  \includegraphics[scale=1.4]{scvx_defects}%
  \fi
\end{csmfigure}

As illustrated in \figref{scvx_defects}, the defect \eqref{eq:scvx_defect} captures the
discrepancy between the next discrete-time state $x_{k+1}$ and the state
obtained by using the flow map starting at time $t_k$. The defect has the
following interpretation: a non-zero defect indicates that the solution to the
subproblem is dynamically infeasible with respect to the original nonconvex
dynamics. For a dynamically feasible subproblem solution, the flow map
trajectories in \figref{scvx_defects} coincide with the discrete-time states at the
discrete\dash time nodes. This is a direct consequence of the consistency
property from \dref{scvx_flow_consistency}. We shall see this happen for the
converged solutions of the numerical examples presented in Part III of this
article (e.g., see~\figref{ex_quad_pos,ex_ff_pos}).

We now know how to compute a consistent flow map and how to use it to calculate
defects using \eqref{eq:scvx_defect}. We will now leverage defects to update the
trust region radius in \scvx. First, define a nonlinear version of
\eqref{eq:scvx_costs_L} as follows:
\begin{align}
  \notag
  \mfnlw(x,u,p,\vcany)
  &= \termw\pare[big]{x(1),p,\gic(x(0),p),\gtf(x(1),p)}+ \\
  \label{eq:scvx_costs_J}
  &\pushright{\mathtt{trapz}(\runnw^N),\qquad} \\
  \notag
  \runnw[,k]^N
  &= \runnw\big(\xk,\uk,p,\defectk, \\
  \label{eq:scvx_costs_J_trap}
  &\pushright{\brak[big]{s(t_k,\xk,\uk,p)}^+\big).\qquad\qquad}
\end{align}
where the positive\dash part function $\brak{\cdot}^+$ returns zero when its
argument is negative, and otherwise just returns the argument. A salient
feature of \eqref{eq:scvx_costs_J} is that it evaluates the penalty function
\eqref{eq:P_penalty_def} based on how well the actual nonconvex constraints are
satisfied. To do so, when compared to \eqref{eq:scvx_costs_L}, $E_k\vck$ is
replaced with the defect $\defectk$ measuring dynamic infeasibility, $\vcck$ is
replaced with the actual nonconvex path constraints
\optieqref{scp_gen_cont}{nonconvex_constraints}, while $\vccic$ and $\vcctf$ are
replaced by the actual boundary conditions
\optieqref{scp_gen_cont}{initial_conditions} and
\optieqref{scp_gen_cont}{final_conditions}. Because evaluation of the defect and
the nonconvex path and boundary constraints is a nonlinear operation, we call
\eqref{eq:scvx_costs_J} the nonlinear augmented cost function.

\begin{csmfigure}[%
  caption={%
    The \scvx trust region update rule. The accuracy metric $\rho$ is defined
    in~\eqref{eq:scvx_ratio} and provides a measure of how accurately the convex
    subproblem given by \pref{subproblem_scvx_dt} describes the original
    \pref{scp_gen_cont}. Note that Case 1 actually rejects the solution to
    \pref{subproblem_scvx_dt} and shrinks the trust region before proceeding to the
    next iteration. In this case, the convex approximation is deemed so poor
    that it is unusable.%
  },
  label={scvx_updates},
  columns=2,position=t]%
  \centering
  \includegraphics{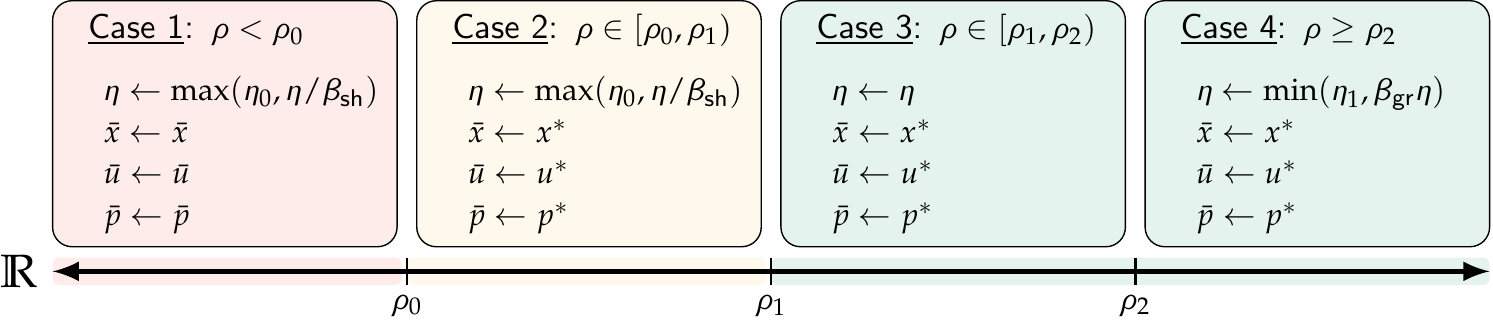}
\end{csmfigure}

The \scvx algorithm uses the linear augmented cost \eqref{eq:scvx_costs_L} and the
nonlinear augmented cost \eqref{eq:scvx_costs_J} to, roughly speaking, measure the
accuracy of the convex approximation of \pref{scp_gen_cont} by
\pref{scp_gen_cvx}. Using the reference solution and the optimal solution to
\pref{subproblem_scvx_dt}, \scvx defines the following scalar metric to measure
convexification accuracy:
\begin{equation}
  \label{eq:scvx_ratio}
  \rat \definedas
  \frac{\mfnlw(\xb,\ub,\pb) - \mfnlw(x^*,u^*,p^*)}
  {\mfnlw(\xb,\ub,\pb) - \mflw(x^*,u^*,p^*,\vc^*,{\vcc}^*)}.
\end{equation}

Let us carefully analyze the elements of \eqref{eq:scvx_ratio}. First of all, the
denominator is always nonnegative because the following relation holds
\cite[Theorem~3.10]{Mao2018}:
\begin{equation}
  \label{eq:scvx_ratio_denominator_nonneg}
  \mfnlw(\xb,\ub,\pb) \geq%
  \mflw(x^*,u^*,p^*,\vc^*,{\vcc}^*).
\end{equation}

The proof of \eqref{eq:scvx_ratio_denominator_nonneg} works by constructing a
feasible subproblem solution that matches the cost $\mfnlw(\xb,\ub,\pb)$. Begin
by setting the state, input, and parameter values to the reference solution
$\trajohone$. We now have to choose the virtual controls that make this
solution feasible for \pref{subproblem_scvx_dt} and that yield a matching cost
\optiobjref{subproblem_scvx_dt}. The consistency property from
\dref{scvx_flow_consistency} ensures that this is always possible to do. In
particular, choose the dynamics virtual control to match the defects, and the
other virtual controls to match the constraint values at the reference
solution. This represents a feasible solution of \pref{subproblem_scvx_dt} whose
cost equals $\mfnlw(\xb,\ub,\pb)$. The optimal cost for the subproblem cannot
be worse, so the inequality \eqref{eq:scvx_ratio_denominator_nonneg} follows.

If the denominator of \eqref{eq:scvx_ratio} is zero, it follows from the above
discussion that the reference trajectory is an optimal solution of the convex
subproblem. This signals that it is appropriate to terminate \scvx and to exit
the loop in \figref{scp_loop}. Hence, we can use the denominator of
\eqref{eq:scvx_ratio} as part of a stopping criterion that avoids a division by
zero. The stopping criterion is discussed further in the next section.

Looking at \eqref{eq:scvx_ratio} holistically, it is essentially a ratio between
the actual cost improvement (the numerator) and the predicted cost improvement
(the denominator), achieved during one \scvx iteration
\cite{Kochenderfer2019}:
\begin{equation}
  \label{eq:scvx_rat_holistic}
  \rat = \frac{\textnormal{actual improvement}}{%
    \textnormal{predicted improvement}}.
\end{equation}

In fact, the denominator can be loosely interpreted as a lower bound
prediction: the algorithm ``thinks'' that the actual cost improves by at least
that much. We can thus adopt the following intuition based on seeing
\pref{subproblem_scvx_dt} as a local model of \pref{scp_gen_cont}. A small
$\rho$ value indicates that the model is inaccurate because the actual cost
decrease is much smaller than predicted. If $\rho$ is close to unity, the model
is accurate because the actual cost decrease is similar to the prediction. If
the $\rho$ value is greater than unity, the model is ``conservative'' because
it underestimates the cost reduction. As a result, large $\rho$ values
incentivize growing the trust region because the model is trustworthy and we
want to utilize more of it. On the other hand, small $\rho$ values incentivize
shrinking the trust region in order to not ``overstep'' an inaccurate model
\cite{Mao2018}.

The \scvx update rule for the trust region radius $\tr$ formalizes the above
intuition. Using three user-defined scalars $\rat_0,\,\rat_1,\,\rat_2\in (0,1)$
that split the real number line into four parts, the trust region radius and
reference trajectory are updated at the end of each \scvx iteration according
to~\figref{scvx_updates}. The user-defined constants $\trshrink,\,\trgrow > 1$ are
the trust region shrink and growth rates, respectively. Practical
implemenations of \scvx also let the user define a minimum and a maximum trust
region radius by using $\tr_0,\,\tr_1> 0$.

The choice of the user-defined scalars $\rat_0$, $\rat_1$, and $\rat_2$ greatly
influences the algorithm runtime. Indeed, as shown in~\figref{scvx_updates}, when
$\rat<\rat_0$ the algorithm will actually outright reject the subproblem
solution, and will resolve the same problem with a smaller trust region. This
begs the question: can rejection go on indefinitely? If this occurs, the
algorithm will be stuck in an infinite loop. The answer is no, the
metric~\eqref{eq:scvx_ratio} must eventually rise above $\rat_0$
\cite[Lemma~3.11]{Mao2018}.

\subsubsection{\scvx Stopping Criterion}

The previous sections provide a complete description of the
\textsf{\color{beamerBlue}Starting} and \textsf{\color{beamerRed}Iteration}
regions in \figref{scp_loop}. A crucial remaining element is how and when to exit
from the ``SCP loop'' of the \textsf{\color{beamerRed}Iteration} region. This
is defined by a stopping criterion (also called an exit criterion), which is
implemented at location \alglocation{\testloc} in \figref{scp_loop}. When the
stopping criterion triggers, we say that the algorithm has converged, and the
final iteration's trajectory $\{x^*(t),\,u^*(t),\,p^*\}_0^1$ is output.

The basic idea of the stopping criterion is to measure how different the new
trajectory $\{x^*(t),\,u^*(t),\,p^*\}_0^1$ is from the reference trajectory
$\trajohone$. Intuitively, if the difference is small, then the algorithm
considers the trajectory not worthy of further improvement and it is
appropriate to exit the SCP loop. The formal \scvx exit criterion uses the
denominator of \eqref{eq:scvx_ratio} (i.e, the predicted cost improvement) as the
stopping criterion \cite{Mao2018,Mao2016,Mao2017,Dueri2017Clipping}:
\begin{equation}
  \label{eq:scvx_stopping_criterion_official}
  \mfnlw(\xb,\ub,\pb) - \mflw(x^*,u^*,p^*,\vc^*,{\vcc}^*)\le\varepsilon,
\end{equation}
where $\varepsilon\in\pos$ is a user-defined (small) threshold. For notational
simplicity, we will write \eqref{eq:scvx_stopping_criterion_official} as:
\begin{equation*}
  \bar{\mfnlw} - \mflw^*\le\varepsilon.
\end{equation*}

Numerical experience has shown a control-dependent stopping criterion to
sometimes lead to an unnecessarily conservative definition of convergence
\cite{Reynolds2019b}. For example, in an application like rocket landing, common
vehicle characteristics (e.g., inertia and engine specific impulse) imply that
relatively large changes in control can have little effect on the state
trajectory and optimality. When the cost is control dependent (which it often
is), \eqref{eq:scvx_stopping_criterion_official} may be a conservative choice that
will result in more iterations without providing a more optimal solution. We
may thus opt for the following simpler and less conservative stopping
criterion:
\begin{equation}
  \label{eq:scvx_stopping_criterion}
  \norm[\hat q]{p^* - \pb} +%
  \max_{k\in\{1,\dots,N\}} \norm[\hat q]{\xk^* - \xbk}%
  \leq \varepsilon,
\end{equation}
where $\hat q\in\brac{1,2,2^{\scriptscriptstyle +},\infty}$ defines a norm
similarly to \eqref{eq:scp_trust_region}. Importantly, the following correspondence
holds between \eqref{eq:scvx_stopping_criterion_official} and
\eqref{eq:scvx_stopping_criterion}. For any $\varepsilon$ choice in
\eqref{eq:scvx_stopping_criterion}, there is a (generally different)
$\varepsilon$ choice in \eqref{eq:scvx_stopping_criterion_official} such that: if
\eqref{eq:scvx_stopping_criterion_official} holds, then
\eqref{eq:scvx_stopping_criterion} holds. We call this an ``implication
correspondence'' between stopping criteria.

\scvx guarantees that there will be an iteration for which
\eqref{eq:scvx_stopping_criterion_official} holds. By implication correspondence,
this guarantee extends to \eqref{eq:scvx_stopping_criterion}. In general, the user
can define yet another stopping criterion that is tailored to the specific
trajectory problem, as long as implication correspondence holds. In fact, we do
this for the numerical examples at the end of the article, where we use the
following stopping criterion that combines
\eqref{eq:scvx_stopping_criterion_official} and \eqref{eq:scvx_stopping_criterion}:
\begin{equation}
  \label{eq:scvx_numerical_example_stopping_criterion}
  \textnormal{\eqref{eq:scvx_stopping_criterion} holds or}~%
  \bar{\mfnlw} - \mflw^*%
  \le \varepsilon_{\mathrm{r}}\abso[big]{\bar{\mfnlw}},
\end{equation}
where $\varepsilon_{\mathrm{r}}\in\pos$ is a user-defined (small) threshold on
the relative cost improvement. The second term in
\eqref{eq:scvx_numerical_example_stopping_criterion} allows the algorithm to
terminate when relatively little progress is being made in decreasing the cost,
a signal that it has achieved local optimality.

\subsubsection{\scvx Convergence Guarantee}

The \scvx trust region update rule in \figref{scvx_updates} is designed to permit a
rigorous convergence analysis of the algorithm. A detailed discussion is given
in \cite{Mao2018}. To arrive at the convergence result, the proof requires the
following (mild) technical condition that is common to most if not all
optimization algorithms. We really do mean that the condition is ``mild''
because we almost never check it in practice and \scvx just works.

\begin{condition}[scvx_licq]
  The gradients of the equality and active inequality constraints of
  \pref{scp_gen_cont} must be linearly independent for the final converged
  trajectory output by \scvx. This is known as a linear independence constraint
  qualification (LICQ) \cite[Chapter~12]{NocedalBook}.
\end{condition}

It is also required that the weight $\lambda$ in \eqref{eq:scvx_Lpen} is
sufficiently large. This ensures that the integral penalty term in
\eqref{eq:scvx_Lpen} is a so-called exact penalty function. A precise condition for
``large enough'' is provided in \cite[Theorem~3.9]{Mao2018}, and a possible
strategy is outlined following the theorem in that paper. However, in practice,
we simply iterate over a few powers of ten until \scvx converges with zero
virtual control. An approximate range that works for most problems is $10^2$ to
$10^4$. Once a magnitude is identified, it usually works well across a wide
range of parameter variations for the problem.

The \scvx convergence guarantee is stated below, and is proved in
\cite[Theorem~3.21]{Mao2018}. Beside \conref{scvx_licq}, the theorem requires a few
more mild assumptions on the inequality constraints in \pref{scp_gen_cont} and on
the penalized cost \eqref{eq:scvx_costs_L}. We do not state them here due to their
technical nature, and refer the reader to the paper. It is enough to say that,
like \conref{scvx_licq}, these assumptions are mild enough that we rarely check
them in practice.

\begin{theorem}[scvx_convergence]
  Suppose that \conref{scvx_licq} holds and the weight $\lambda$ in
  \eqref{eq:scvx_Lpen} is large enough. Regardless of the initial reference
  trajectory provided in \figref{scp_loop}, the \scvx algorithm will always
  converge at a superlinear rate by iteratively
  solving~\pref{subproblem_scvx_dt}. Furthermore, if the virtual controls are zero
  for the converged solution, then it is a stationary point of \pref{scp_gen_cont}
  in the sense of having satisfied the KKT conditions.
\end{theorem}

\tref{scvx_convergence} confirms that \scvx solves the KKT conditions of the
original \pref{scp_gen_cont}. Because these are first-order necessary conditions of
optimality, they are also satisfied by local maxima and saddle
points. Nevertheless, because each convex subproblem is minimized by \scvx, the
event of converging to a stationary point that is not a local minimum is very
small.

\tref{scvx_convergence} also states that the convergence rate is superlinear
\cite{NocedalBook}, which is to say that the distance away from the converged
solution decreases superlinearly \cite[Theorem~4.7]{Mao2018}. This is better
than general NLP methods, and on par with SQP methods that usually also attain
at most superlinear convergence \cite{Boggs1995}. 



In conclusion, note that \tref{scvx_convergence} is quite intuitive and confirms
our basic intuition. If we are ``lucky'' to get a solution with zero virtual
control, then it is a local minimum of the original nonconvex
\pref{scp_gen_cont}. The reader will be surprised, however, at just how easy it is
to be ``lucky''. In most cases, it really is a simple matter of ensuring that
the penalty weight $\lambda$ in \eqref{eq:scvx_Lpen} is large enough. If the
problem is feasible but \scvx is not converging or is converging with non-zero
virtual control, the first thing to try is to increase $\lambda$. In more
difficult cases, some nonconvex constraints may need to be reformulated as
equivalent versions that play nicer with the linearization process (the
free-flyer example in Part III discusses this in detail). However, let us be
clear that there is no a priori hard guarantee that the converged solution will
satisfy $\vcany=0$. In fact, since \scvx always converges, even an infeasible
optimal control problem can be solved, and it will return a solution for which
$\vcany$ is non-zero.

\subsection{The \gusto Algorithm}

The \gusto algorithm is another SCP method that can be used for trajectory
generation. The reader will find that \gusto has a familiar feel to that of
\scvx. Nevertheless, the algorithm is subtly different from both computational
and theoretical standpoints. For example, while \scvx works directly with the
temporally discretized \pref{subproblem_scvx_dt} and derives its convergence
guarantees from the KKT conditions, \gusto performs all of its analysis in
continuous-time using Pontryagin's maximum principle
\cite{Mao2018,Bonalli2019a,PontryaginBook,BerkovitzBook}. Temporal
discretization is only introduced at the very end to enable numerical solutions
to the problem. \gusto applies to versions of \pref{scp_gen_cont} where the running
cost in \eqref{eq:ocost_nlin} is quadratic in the control variable:
\begin{equation}
  \label{eq:gusto_running_cost}
  \runn(x,u,p) = u^\transp \Jq(p) u + u^\transp \Jl (x,p) + \Jc(x,p),
\end{equation}
where the parameter $p$, as before, can be used to capture free final time
problems. The functions $\Jq$, $\Jl$, and $\Jc$ must all be continuously
differentiable. Furthermore, the function $\Jq$ must be positive
semidefinite. The dynamics must be affine in the control variable:
\begin{equation}
  \label{eq:control_affine_gusto}
  f(t,x,u,p) = f_0(t,x,p) + \sum_{i=1}^{m} u_i f_i(t,x,p),
\end{equation}
where the $f_i:\reals\times\reals^n\times\reals^d\to\reals^n$ are nonlinear
functions representing a decomposition of the dynamics into terms that are
control-independent and terms that linearly depend on each of the control
variables $u_i$. Note that any Lagrangian mechanical system can be expressed in
the control affine form, and so \eqref{eq:control_affine_gusto} is applicable to
the vast majority of real\dash world vehicle trajectory generation applications
\cite{Bullo2005}. Finally, the nonconvex path constraints
\optieqref{scp_gen_cont}{nonconvex_constraints} are independent of the control
terms, i.e., $s(t,x,u,p)=s(t,x,p)$. Altogether, these assumptions specialize
\pref{scp_gen_cont} to problems that are convex in the control variable.

At its core, \gusto is an SCP trust region method just like \scvx. Thus, it has
to deal with the same issues of artificial infeasibility and artificial
unboundedness. To this end, the \gusto algorithm makes the following
algorithmic choices:
\begin{itemize}
\item To handle artificial unboundedness, \gusto augments the cost function
  with a soft penalty on the violation of \eqref{eq:scp_trust_region}. Because the
  original \pref{scp_gen_cont} is convex in the control variables, there is no need
  for a trust region with respect to the control and we use $\tru=0$. In its
  standard form, \gusto works with default values
  $\alpha_x=\alpha_p=1$. However, one can choose different values
  $\alpha_x, \alpha_p>0$ without affecting the convergence guarantees;
\item To handle artificial infeasibility, \gusto augments the cost function
  with a soft penalty on nonconvex path constraint violation. As the algorithm
  progresses, the penalty weight increases.
\end{itemize}

From these choices, we can already deduce that a salient feature of the \gusto
algorithm is its use of soft penalties to enforce nonconvex constraints. Recall
that in \scvx, the trust region \eqref{eq:scp_trust_region} is imposed exactly, and
the linearized constraints are (equivalently) relaxed by using virtual control
\eqref{eq:scvx_lin_approxs_vc}. By penalizing constraints, \gusto can be analyzed
in continuous\dash time via the classical Pontryagin maximum principle
\cite{Bonalli2019a,Bonalli2019b,PontryaginBook,BerkovitzBook}. An additional
benefit of having moved the state constraints into the cost is that the virtual
control term employed in the dynamics \eqref{eq:scvx_lin_approxs_vc_f} can be
safely removed, since linearized dynamics are almost always controllable
\cite{BonalliLewTAC2021}.

Let us begin by formulating the augmented cost function used by the \gusto
convex subproblem at \alglocation{2} in \figref{scp_loop}. This involves three
elements: the original cost \eqref{eq:ocost_nlin}, soft penalties on path
constraints that involve the state, and a soft penalty on trust region
violation. We will now dicuss the latter two elements.

To formulate the soft penalties, consider a continuously differentiable,
convex, and nondecreasing penalty function $\hpen:\reals\to\nonneg$ that
depends on a scalar weight $\Jw\ge 1$ \cite{Bonalli2019a,BonalliLewTAC2021}. The
goal of $\hpen$ is to penalize any positive value and to be agnostic to
nonpositive values. Thus, a simple example is a quadratic rectifier:
\begin{equation}
  \label{eq:gusto_hpen}
  \hpen(z) = \Jw\pare[big]{\brak{z}^+}^2,
\end{equation}
where higher $\lambda$ values indicate higher penalization. Another example is
the softplus function \cite{SoftplusNIPS}:
\begin{equation}
  \label{eq:gusto_hpen_softplus}
  \hpen(z) = \Jw \sigma\inv\log\pare[big]{1+\exp{\sigma z}},
\end{equation}
where $\sigma\in\pos$ is a sharpness parameter. As $\sigma$ grows, the softplus
function becomes an increasingly accurate approximation of
$\Jw\max\{0,z\}$. 

We use $\hpen$ to enforce soft penalties for violating the constraints
\optieqref{scp_gen_cont}{convex_path_constraints_X} and
\optieqref{scp_gen_cont}{nonconvex_constraints} that involve the state and
parameter. To this end, let $\indi:\reals^n\times\reals^d\to\reals^{\dimindi}$
be a convex continuously differentiable indicator function that is nonpositive
if and only if the convex state constraints in
\optieqref{scp_gen_cont}{convex_path_constraints_X} are satisfied:
\begin{equation*}
  \indi(x,p)\le 0~\iff~(x, p)\in\set X.
\end{equation*}

Note that because $\set X$ is a convex set, such a function always
exists. Using $\indi$, $s$ from \optieqref{scp_gen_cont}{nonconvex_constraints},
and $\hpen$, we can lump all of the state constraints into a single soft
penalty function:
\begin{equation}
  \label{eq:gusto_soft_penalty_state}
  \Jpen(t,x,p)\definedas\sum_{i=1}^{\dimindi}\hpen\pare[big]{\indi_i(t,x)}+
  \sum_{i=1}^{\dimss}\hpen\pare[big]{s_i(t,x,p)}.
\end{equation}


Another term is added to the augmented cost function to handle artificial
unboundedness. This term penalizes violations of the trust region constraint,
and is defined in a similar fashion as \eqref{eq:gusto_soft_penalty_state}:
\begin{equation}
  \label{eq:gusto_soft_penalty_tr}
  \Jtr(x,p)\definedas \hpen\pare[big]{\norm[q]{\delta x}+
    \norm[q]{\delta p}-\tr},
\end{equation}
although we note that hard enforced versions of the trust region constraint can
also be used \cite{BonalliLewTAC2021}.

The overall augmented cost function is obtained by combining the original cost
\optiobjref{scp_gen_cont} with \eqref{eq:gusto_soft_penalty_state} and
\eqref{eq:gusto_soft_penalty_tr}. Because the resulting function is generally
nonconvex, we take this opportunity to decompose it into its convex and
nonconvex parts:
\begin{subequations}
  \begin{align}
    \label{eq:gusto_Jpen}
    \ocostw(x,u,p)
    &\definedas \ocostwCvx(x,p)+\ocostwNCvx(x,u,p), \\
    \nonumber
    \ocostwCvx(x,p)
    &= \term(x(1),p)+\int_0^1 \Jtr(x,p)+ \\
    \label{eq:gusto_Jpen_cvx}
    &\pushright{\sum_{i=1}^{\dimindi}\hpen\pare{\indi_i(t,x)}\sdt,\quad} \\
    \label{eq:gusto_Jpen_ncvx}
    \ocostwNCvx(x,u,p)
    &= \int_0^1 \runn(x,u,p)+
      \sum_{i=1}^{\dimss}\hpen\pare{s_i(t,x,p)}\sdt.
  \end{align}
\end{subequations}

The terms $\Jtr(x, p)$ and $\hpen(\indi_i(x))$ in \eqref{eq:gusto_Jpen_cvx} are
convex functions, since they are formed by composing a convex nondecreasing
function $\hpen$ with a convex function \cite{BoydConvexBook}. Thus,
\eqref{eq:gusto_Jpen_cvx} is the convex part of the cost, while
\eqref{eq:gusto_Jpen_ncvx} is the nonconvex part.

Our ultimate goal is to construct a convex subproblem that can be solved at
location \alglocation{\solveloc} in \figref{scp_loop}. Thus, the nonconvex part of
the cost \eqref{eq:gusto_Jpen_ncvx} must be convexified around the reference
trajectory $\trajohone$. This requires the following Jacobians in addition to
the initial list \eqref{eq:scvx_lin_mats}:
\begin{subequations}
  \label{eq:gusto_lin_mats}
  \begin{align}
    \label{eq:scvx_lin_mats_a_lambda}
    \diffLx &\definedas \diff{x}\runn(\xb,\ub,\pb), \\
    \label{eq:scvx_lin_mats_b_lambda}
    \diffLu &\definedas \diff{u}\runn(\xb,\ub,\pb), \\
    \label{eq:scvx_lin_mats_c_lambda}
    \diffLp &\definedas \diff{p}\runn(\xb,\ub,\pb).
  \end{align}
\end{subequations}

Using \eqref{eq:scvx_lin_mats_h}, \eqref{eq:scvx_lin_mats_j}, and
\eqref{eq:gusto_lin_mats}, we can write the convex approximation of
\eqref{eq:gusto_Jpen_ncvx}:
\begin{align}
  \nonumber
  \ocostLwNCvx(x,u,p)
  &= \int_0^1
    \runn(\xb,\ub,\pb)+\diffLx\delta x+\diffLu\delta u+\diffLp\delta p
    + \\
  \label{eq:gusto_Jpen_ncvx_lin}
  &\pushright{
    \sum_{i=1}^{\dimss}\hpen\pare[big]{
    s_i(t,\xb,\pb)+C_i\delta x+G_i\delta p
    }\sdt,\quad}
\end{align}
where $C_i$ and $G_i$ are the $i$-th rows of the Jacobians
\eqref{eq:scvx_lin_mats_h} and \eqref{eq:scvx_lin_mats_j}, respectively. Note that,
strictly speaking, $\ocostLwNCvx$ is not a linearized version of $\ocostwNCvx$
because the second term in \eqref{eq:gusto_Jpen_ncvx} is only linearized inside the
$\hpen(\cdot)$ function.

Replacing $\ocostwNCvx$ in \eqref{eq:gusto_Jpen} with $\ocostLwNCvx$, we obtain a
convexified augmented cost function:
\begin{equation}
  \label{eq:gusto_Lpen}
  \ocostLw(x,u,p) = \ocostwCvx(x,p)+\ocostLwNCvx(x,u,p).
\end{equation}

To summarize the above discussion, the continuous-time convex subproblem that
is solved at each iteration of the \gusto algorithm can be stated formally as:
\begin{optimization}[
  label={subproblem_gusto_ct},
  variables={u,p},
  objective={\ocostLw(x,u,p)}]%
  \optilabel{dynamics}
  & \dot{x} = A x + B u + F \pk + r,  \\
  \optilabel{u_constraints}
  & (u,p) \in\set{U}, \\
  \optilabel{initial_condition}
  & H_0 x(0) + K_0 \pk + \ell_0 = 0, \\
  \optilabel{terminal_condition}
  & \Hf x(1) + \Kf \pk + \ellf = 0.
\end{optimization}

\pref{subproblem_gusto_ct} can be compared to the \scvx continuous\dash time convex
subproblem, given by \pref{subproblem_scvx_ct}. Broadly speaking, the problems are
quite similar. Their main differences stem from how the \scvx and \gusto
algorithms handle artificial infeasibility and artificial unboundedness. In the
case of \scvx, virtual control terms are introduced and a hard trust region
\optieqref{subproblem_scvx_ct}{xu_constraints_3} is imposed. In the case of \gusto,
everything is handled via soft penalties in the augmented cost. The result is
that penalty terms in the cost \optiobjref{subproblem_gusto_ct} replace the
constraints \optieqref{subproblem_scvx_ct}{xu_constraints_1},
\optieqref{subproblem_scvx_ct}{xu_constraints_2}, and
\optieqref{subproblem_scvx_ct}{xu_constraints_3}.

There is another subtle difference between \pref{subproblem_gusto_ct} and
\pref{subproblem_scvx_ct}, which is that \gusto does not use virtual control terms
for the linearized boundary conditions
\optieqref{subproblem_gusto_ct}{initial_condition} and
\optieqref{subproblem_gusto_ct}{terminal_condition}. In \scvx, these virtual
control terms maintain subproblem feasibility when the linearized boundary
conditions define hyperplanes that do not intersect with the hard-enforced
convex state path constraints in
\optieqref{subproblem_scvx_ct}{xu_constraints_1}. This is not a problem in \gusto,
since the convex state path constraints are only penalized in the cost
\eqref{eq:gusto_soft_penalty_state}, and violating them for the sake of retaining
feasibility is allowed. Furthermore, the linearized dynamics
\optieqref{subproblem_gusto_ct}{dynamics} are theoretically guaranteed to be almost
always controllable \cite{BonalliLewTAC2021}. This implies that the dynamics
\optieqref{subproblem_gusto_ct}{dynamics} can always be steered between the
linearized boundary conditions \optieqref{subproblem_gusto_ct}{initial_condition}
and \optieqref{subproblem_gusto_ct}{terminal_condition} \cite{Antsaklis2006}.

Similar to how we treated \pref{subproblem_scvx_ct}, a temporal discretization
scheme must be applied in order to numerically solve the subproblem.
Discretization proceeds in the same way as for \scvx: we select a set of
temporal points $t_k\in[0,1]$ for $k=1,\ldots,N$, and recast
\pref{subproblem_gusto_ct} as a parameter optimization problem in the (overloaded)
variables $x = \{\xk\}_{k=1}^{N}$, $u = \{\uk\}_{k=1}^{N}$, and $p$. The same
discretization choices are available as for \scvx, such as described in
\sbref{discretization}.


The integrals of the continuous-time cost function \optiobjref{subproblem_gusto_ct}
also need to be discretized. This is done in a similar fashion to the way
\eqref{eq:scvx_costs_L} was obtained from \eqref{eq:scvx_Lpen} for \scvx. For
simplicity, we will again assume that the time grid is uniform with step size
$\timeintvl$ and that trapezoidal numerical integration \eqref{eq:trapz} is
used. For notational convenience, by combining the integral terms in
\eqref{eq:gusto_Jpen_cvx} and \eqref{eq:gusto_Jpen_ncvx_lin}, we can write
\eqref{eq:gusto_Lpen} compactly as:
\begin{equation}
  \label{eq:gusto_Lpen_concat}
  \ocostLw(x,u,p) = \term(x(1),p)+\int_0^1 \ocostLwrunn(x,u,p)\sdt,
\end{equation}
where the convex function $\ocostLwrunn$ is formed by summing the integrands of
\eqref{eq:gusto_Jpen_cvx} and \eqref{eq:gusto_Jpen_ncvx_lin}. We can then compute the
discrete-time version of \eqref{eq:gusto_Lpen_concat}, which is the \gusto
equivalent of its \scvx counterpart \eqref{eq:scvx_costs_L}:
\begin{align}
  \label{eq:gusto_costs_L}
  \mflw(x,u,p) &=\term(x(1),p)+\mathtt{trapz}(\ocostLwrunnN), \\
  \ocostLwrunnN[,k] &= \ocostLwrunn(\xk,\uk,p).
\end{align}

Lastly, like in \scvx, the constraint
\optieqref{subproblem_gusto_ct}{u_constraints} is enforced only at the discrete
temporal nodes. In summary, the following discrete-time convex subproblem is
solved at each \gusto iteration (i.e., location \alglocation{\solveloc} in
\figref{scp_loop}):
\begin{optimization}[
  label={subproblem_gusto_dt},
  variables={x,u,p},
  objective={\mflw(x,u,p)}]%
  \optilabel{dynamics}
  & \xkp = A_k \xk + B_k \uk + F_k \pk + r_k, \\
  \optilabel{input_constraints}
  & (\uk, \pk) \in\set{U}_k, \\
  & H_0 x_1 + K_0 \pk + \ell_0 = 0, \\
  & \Hf x_N + \Kf \pk + \ellf = 0,
\end{optimization}
where it is implicitly understood that the constraints
\optieqref{subproblem_gusto_dt}{dynamics} and
\optieqref{subproblem_gusto_dt}{input_constraints} hold at each temporal node.

\subsubsection{\gusto Update Rule}
\label{subsubsec:gusto_update}

We are now at the same point as we were in the \scvx section: by using
\pref{subproblem_gusto_dt} as the discrete-time convex subproblem, all of the
necessary elements are available to go around the loop in \figref{scp_loop}, except
for how to update the trust region radius $\eta$ and the penalty weight
$\Jw$. Both values are generally different at each \gusto iteration, and we
shall now describe the method by which they are updated. The reader will find
the concept to be similar to how things worked for \scvx.

First, recall that \gusto imposes the trust region \eqref{eq:scp_trust_region} as a
soft constraint via the penalty function \eqref{eq:gusto_soft_penalty_tr}. This
means that the trust region constraint can possible be violated. If the
solution to \pref{subproblem_gusto_dt} violates \eqref{eq:scp_trust_region}, \gusto
rejects the solution and increases the penalty weight $\lambda$ by a
user-defined factor $\Jwgrow>1$. Otherwise, if \eqref{eq:scp_trust_region} holds,
the algorithm proceeds by computing the following convexification accuracy
metric that is analogous to \eqref{eq:scvx_ratio}:
\begin{equation}
  \label{eq:gusto_ratio}
  \rat \definedas
  \frac{\abso[big]{\ocostw(x^*,u^*,p^*) -
      \ocostLw(x^*,u^*,p^*)} + \varTheta^*}
  {\abso[big]{\ocostLw (x^*,u^*,p^*)} +
    \int_0^1 \norm[2]{\dot{x}^*}\sdt},
\end{equation}
where, recalling \optieqref{scp_gen_cont}{dynamics} and
\optieqref{subproblem_gusto_ct}{dynamics}, we have defined:
\begin{subequations}
  \label{eq:gusto_ratio_defs}
  \begin{align}
    \label{eq:gusto_integrate_x_star}
    \dot{x}^*   &\definedas A x^* + Bu^* + Fp^* + r, \\
    \label{eq:gusto_integrate_Theta_star}
    \varTheta^* &\definedas \int_0^1 \norm[2]{f(t,x^*,u^*,p^*) - \dot{x}^*}\sdt.
  \end{align}
\end{subequations}

Equation \eqref{eq:gusto_integrate_x_star} integrates to yield the continuous-time
state solution trajectory. This is done in accordance with the temporal
discretization scheme, such as the one described in \sbref{discretization}. As a
result, the value of $\varTheta^*$ is nothing but a measure of the total
accumulated error that results from linearizing the dynamics along the
subproblem's optimal solution. In a way, $\varTheta^*$ is the counterpart of
\scvx defects defined in \eqref{eq:scvx_defect}, with the subtle difference that
while defects measure the discrepancy between the linearized and nonlinear
state trajectories, $\varTheta^*$ measures the discrepancy between the
linearized and nonlinear state dynamics.

The continuous-time integrals in \eqref{eq:gusto_ratio} can, in principle, be
evaluated exactly (i.e., to within numerical precision) based on the
discretization scheme used. In practice, we approximate them by directly
numerically integrating the solution of \pref{subproblem_gusto_dt} using (for
example) trapezoidal integration \eqref{eq:trapz}. This means that we evaluate the
following approximation of \eqref{eq:gusto_ratio}:
\begin{align}
  \label{eq:gusto_ratio_approx}
  \rat
  &\approx
    \frac{\abso[big]{\mfnlw(x^*,u^*,p^*) -
    \mflw(x^*,u^*,p^*)} + \varTheta^*}
    {\abso[big]{\mflw (x^*,u^*,p^*)} +
    \mathtt{trapz}(\mathsf{x}^*)}, \\
  \nonumber
  \varTheta^*
  &\approx
    \mathtt{trapz}(\Delta f^*), \\
  \nonumber
  \Delta f^*_k
  &= \norm[2]{f(t_k,\xk^*, \uk^*, p^*)-\dot{x}_k^*},\quad%
    \mathsf{x}^*_k = \norm[2]{\dot{x}_k^*}.
\end{align}

The nonlinear augmented cost $\mfnlw$ in the numerator of
\eqref{eq:gusto_ratio_approx} is a temporally discretized version of
\eqref{eq:gusto_Jpen}. In particular, combine the integrals of
\eqref{eq:gusto_Jpen_cvx} and \eqref{eq:gusto_Jpen_ncvx} to express \eqref{eq:gusto_Jpen}
as:
\begin{equation}
  \label{eq:gusto_pen_concat}
  \ocostw(x,u,p) = \term(x(1),p)+\int_0^1 \ocostwrunn(x,u,p)\sdt,
\end{equation}
which allows us to compute $\mfnlw$ as follows:
\begin{align}
  \label{eq:gusto_costs_J}
  \mfnlw(x,u,p) &=\term(x(1),p)+\mathtt{trapz}(\ocostwrunnN), \\
  \label{eq:gusto_costs_J_trapz}
  \ocostwrunnN[,k] &= \ocostwrunn(\xk,\uk,p).
\end{align}

Looking at \eqref{eq:gusto_ratio} holistically, it can be seen as a normalized
error that results from linearizing the cost and the dynamics:
\begin{equation}
  \label{eq:gusto_rat_holistic}
  \rat = \frac{\textnormal{cost error}+\textnormal{dynamics error}}{%
    \textnormal{normalization term}}.
\end{equation}

Note that as long as the solution of \pref{subproblem_gusto_dt} is nontrivial
(i.e., $x^*(t)=0$ for all $t\in [0,1]$ does not hold), the normalization term
is guaranteed to be strictly positive. Thus, there is no danger of dividing by
zero.

For comparison, \scvx evaluates convexification accuracy through
\eqref{eq:scvx_rat_holistic}, which measures accuracy as a relative error in the
cost improvement prediction. This prediction is a ``higher order'' effect:
linearization error indirectly influences cost prediction error through the
optimization of \pref{subproblem_scvx_dt}. \gusto takes a more direct route with
\eqref{eq:gusto_rat_holistic}, and measures convexification accuracy directly as a
normalized error that results from linearizing both the cost and the dynamics.

Looking at \eqref{eq:gusto_rat_holistic}, we can adopt the following intuition
about the size of $\rat$. As for \scvx, let us view \pref{subproblem_gusto_dt} as a
local model of \pref{scp_gen_cont}. A large $\rat$ value indicates an inaccurate
model, since the linearization error is large. A small $\rat$ value indicates
an accurate model, since the linearization error is relatively small compared
to the normalization term. Hence, large $\rat$ values incentivize shrinking the
trust region in order to not overstep the model, while small $\rat$ values
incentivize growing it in order to exploit a larger region of an accurate
model. Note that the opposite intuition holds for \scvx, where small $\rat$
values are associated with shrinking the trust region.

\begin{csmfigure}[%
  caption={%
    The \gusto trust region update rule. The accuracy metric $\rat$ is defined
    in~\eqref{eq:gusto_ratio} and provides a measure of how accurately the convex
    subproblem given by \pref{subproblem_gusto_dt} describes the original
    \pref{scp_gen_cont}. Note that Cases 3 and 4 reject the solution to
    \pref{subproblem_gusto_dt}. In Case 3, this is due to the convex approximation
    being deemed so inaccurate that it is unusable, and the trust region is
    shrunk accordingly. In Case 4, this is due to the trust region constraint
    \eqref{eq:scp_trust_region} being violated.%
  },%
  label={gusto_updates},%
  columns=2]%
  \centering%
  \includegraphics{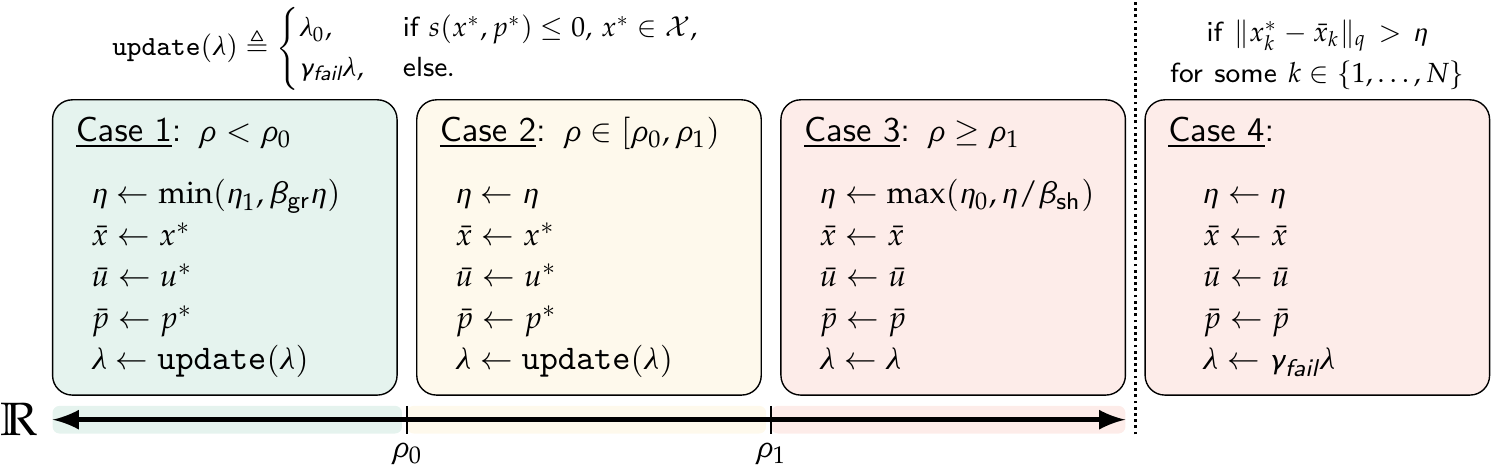}
\end{csmfigure}

The \gusto update rule formalizes the above intuition. Using two user-defined
constants $\rat_0,\,\rat_1\in (0,1)$ that split the real number line into three
parts, the trust region radius $\tr$ and the penalty weight $\Jw$ are updated
at the end of each \gusto iteration according to \figref{gusto_updates}. Just like
in \scvx, the user-defined constants $\trshrink,\,\trgrow > 1$ are the trust
region shrink and growth rates, respectively. The sizes of the trust region and
the penalty weight are restricted by user-defined constants: $\tr_0$ is the
minimum trust region radius, $\tr_1$ is the maximum trust region radius, and
$\Jw_0\ge 1$ is the minimum penalty weight. Importantly, for cases 1 and 2 in
\figref{gusto_updates}, whenever the solution of \pref{subproblem_gusto_dt} is
accepted, the penalty weight $\Jw$ is increased by a factor $\Jwgrow$ if any of
the nonconvex state constraints in \optieqref{scp_gen_cont}{nonconvex_constraints}
are violated. This incentivizes subsequent iterates to become feasible with
respect to \optieqref{scp_gen_cont}{nonconvex_constraints}.

In addition, \gusto requires that $\eta$ eventually shrinks to zero as the SCP
iterations progress (we cover this in more detail in a later section on \gusto
convergence) \cite{BonalliLewTAC2021}. However, if the trajectory is converging,
then $\delta x,\,\delta u,\,\delta p\to 0$ and the linearizations in
\pref{subproblem_gusto_ct} are bound to become more accurate. Hence, we expect
near-zero $\rat$ values when the trajectory is converging, which means that
\figref{gusto_updates} will grow the trust region instead of shrinking it. A simple
remedy is to apply the following exponential shrink factor following the update
in \figref{gusto_updates}:
\begin{equation}
  \label{eq:gusto_exp_shrink}
  \tr \gets \shrinkrate^{\brak{1+k-\shrinkk}^+}\tr,
\end{equation}
where $\shrinkrate\in (0,1)$ is an exponential shrink rate and $\shrinkk\ge 1$
is the first SCP iteration where shrinking is applied. The user sets both
$\shrinkrate$ and $\shrinkk$, and in this manner can regulate how fast the
trust region shrinks to zero. We view low $\shrinkrate$ values and high
$\shrinkk$ values as setting the algorithm up for a longer ``exploration''
phase prior to tightening the trust region.

\subsubsection{\gusto Stopping Criterion}

To complete the description of the \gusto algorithm, it remains to define the
stopping criterion used at location \alglocation{\testloc} of
\figref{scp_loop}. The formal \gusto exit criterion directly checks if the control
and parameters are all almost unchanged \cite{Bonalli2019a,LewBonalliECC2020}:
\begin{align}
  &\biggg(\norm[\hat q]{\pb-p^*}\le\varepsilon
  ~\textnormal{and} \notag \\
  \label{eq:gusto_stopping_criterion_official}
  &\qquad\qquad\int_0^1\norm[\hat q]{u^*(t)-\ub(t)}\le\varepsilon\biggg)
  ~\textnormal{or}~\Jw>\Jw_{\max},
\end{align}
where $\varepsilon\in\pos$ and
$\hat q\in\brac{2,2^{\scriptscriptstyle +},\infty}$ have the same meaning as
before in \eqref{eq:scvx_stopping_criterion}, and the new parameter
$\Jw_{\max}\in\pos$ is a (large) maximum penalty weight.

When \eqref{eq:gusto_stopping_criterion_official} triggers due to $\Jw_{\max}$,
\gusto exits with an unconverged trajectory that violates the state and/or
trust region constraints. This is equivalent to \scvx exiting with non-zero
virtual control, and indicates that the algorithm failed to solve the problem
due to inherent infeasibility or numerical issues.

Computing \eqref{eq:gusto_stopping_criterion_official} can be computationally
expensive due to numerical integration of the control deviation. We can
simplify the calculation by directly using the discrete\dash time
solution. This leads to the following stopping criterion:
\begin{align}
  \label{eq:gusto_stopping_criterion}
  &\norm[\hat q]{p^* - \pb} +
    \mathtt{trapz}(\Delta u^*) \leq \varepsilon~\textnormal{or}~\Jw>\Jw_{\max}, \\
  \notag
  &\Delta u_k^* = \norm[\hat q]{\uk^* - \ubk}.
\end{align}

Just as discussed for \scvx, an implication correspondence holds between the
stopping criteria \eqref{eq:gusto_stopping_criterion_official} and
\eqref{eq:gusto_stopping_criterion}. In practice, we follow the example of
\eqref{eq:scvx_numerical_example_stopping_criterion} and add the option for exiting
when relatively little progress is being made in decreasing the cost, which can
often signal local optimality sooner than \eqref{eq:gusto_stopping_criterion} is
satisfied:
\begin{equation}
  \label{eq:gusto_numerical_example_stopping_criterion}
  \textnormal{\eqref{eq:gusto_stopping_criterion} holds or}~%
  \abso[big]{\bar{\mfnlw} - \mfnlw^*}%
  \le \varepsilon_{\mathrm{r}}\abso[big]{\bar{\mfnlw}},
\end{equation}
where, as before, $\varepsilon_{\mathrm{r}}\in\pos$ is a user-defined (small)
threshold on the relative cost improvement. The numerical examples at the end
of the article implement
\eqref{eq:gusto_numerical_example_stopping_criterion}.

\subsubsection{\gusto Convergence Guarantee}

Each iteration of the \gusto numerical implementation can be seen as composed
of three stages. Looking at \figref{scp_loop}, first the convex
\pref{subproblem_gusto_dt} is constructed and solved with a convex optimizer at
location \alglocation{\solveloc}. Using the solution, the stopping criterion
\eqref{eq:gusto_numerical_example_stopping_criterion} is checked at
\alglocation{\testloc}. If the test fails, the third and final stage updates
the trust region radius and soft penalty weight according to \figref{gusto_updates}
and \eqref{eq:gusto_exp_shrink}. In this context, a convergence guarantee ensures
that the stopping criterion at \alglocation{\testloc} in \figref{scp_loop}
eventually triggers.

 \gusto is an SCP algorithm that was designed and analyzed in
continuous-time using the Pontryagin maximum principle
\cite{Bonalli2019a,BonalliLewTAC2021}. Thus, the first convergence guarantee
that we are able to provide assumes that the \gusto algorithm solves the
continuous\dash time subproblem (i.e., \pref{subproblem_gusto_ct}) at each
iteration \cite[Corollary~III.1]{Bonalli2019a},
\cite[Theorem~3.2]{BonalliLewTAC2021}.

\begin{theorem}[gusto_convergence]
  Regardless of the initial reference trajectory provided in \figref{scp_loop}, the
  \gusto algorithm will always converge by iteratively
  solving~\pref{subproblem_scvx_ct}. Furthermore, if $\Jw\le\Jw_{\max}$ (in
  particular, the state constraints are exactly satisfied), then the solution
  is a stationary point of \pref{scp_gen_cont} in the sense of having satisfied the
  necessary optimality conditions of the Pontryagin maximum principle.
\end{theorem}

 We will emphasize once again the following duality between the
\gusto and \scvx convergence guarantees. Both algorithms can converge to
infeasible points of the original \pref{scp_gen_cont}. For \scvx, this corresponds
to a non\dash zero virtual control (recall \tref{scvx_convergence}). For \gusto,
this corresponds to $\Jw>\Jw_{\max}$. The situations are completely
equivalent, and correspond simply to the different choices made by the
algorithms to impose nonconvex constraints either with virtual control (as in
\scvx) or as soft cost penalties (as in \gusto).

\tref{gusto_convergence} should be viewed as the \gusto counterpart of
\tref{scvx_convergence} for \scvx. Despite the similarities between the two
statements, there are three important nuances about \tref{gusto_convergence} that
we now discuss.

The first difference concerns how the \gusto update rule in \figref{gusto_updates}
and \eqref{eq:gusto_exp_shrink} play into the convergence proof. In \scvx, the
update rule from \figref{scvx_updates} plays a critical role in proving
convergence. Thus, \scvx is an SCP algorithm that is intimately linked to its
trust region update rule. This is not the case for \gusto, whose convergence
proof does not rely on the update rule definition at all. The only thing
assumed by \tref{gusto_convergence} is that the trust region radius $\eta$
eventually shrinks to zero \cite{BonalliLewTAC2021}. Ensuring this is the reason
for \eqref{eq:gusto_exp_shrink}. Thus, \gusto is an SCP algorithm that accepts any
update rule that eventually shrinks to zero. Clearly, some update rules may
work better than others, and the one presented in this article is simply a
choice that we have implemented in practice. Another simple update rule is
given in \cite{BonalliLewTAC2021}. Overall, the \gusto update rule can be viewed
as a mechanism by which to accept or reject solutions based on their
convexification accuracy, and not as a function designed to facilitate formal
convergence proofs.

The reason why it is necessary to shrink the trust region to zero is the second
nuanced detail of \tref{gusto_convergence}. In nonlinear optimization literature,
most claims of convergence fall into the so-called weak convergence category
\cite{Absil2005,NocedalBook}. This is a technical term which means that a
\textit{subsequence} among the trajectory iterates $\{x^i(t),u^i(t),p^i\}_0^1$
converges. For example, the subsequence may be $i=1,3,5,7,\dots$, while the
iterates $i=2,4,6,\dots$ behave differently. Both \scvx and \gusto provide, at
the baseline, the weak convergence guarantee. The next step is to provide a
strong convergence guarantee, which ensures that the full sequence of iterates
converges. For \scvx, this is possible by leveraging properties of the update
rule in \figref{scvx_updates}. As we mentioned in the above paragraph, \gusto is
more flexible in its choice of update rule. To get a converging iterate
sequence, the additional element \eqref{eq:gusto_exp_shrink} is introduced to force
all subsequences to converge to the same value. Because our analysis for \scvx
has shown the strong convergence guarantee to be tied exclusively to the choice
of update rule, \gusto can likely achieve a similar guarantee with an
appropriately designed update rule.

The third and final nuanced detail of \tref{gusto_convergence} is that it assumes
continuous\dash time subproblems (i.e., \pref{subproblem_gusto_ct}) are solved. In
reality, the discretized \pref{subproblem_gusto_dt} is implemented on the
computer. This difference between proof and reality should ring a familiar tone
to Part I on \lcvx, where proofs are also given using the continuous-time
Pontryagin's maximum principle, even though the implementation is in
discrete-time. If the discretization error is small, the numerically
implemented \gusto algorithm will remain in the vicinity of a convergent
sequence of continuous\dash time subproblem solutions. Thus, given an accurate
temporal discretization, \tref{gusto_convergence} can be reasonably assumed to
apply for the numerical \gusto algorithm \cite{BonalliLewTAC2021}.

Our default choice in research has been to use an interpolating polynomial
method such as described in \sbref{discretization}. This approach has three
advantages \cite{Malyuta2019b}: 1) the continuous\dash time dynamics
\optieqref{subproblem_gusto_ct}{dynamics} are satisfied exactly, 2) there is a
cheap way to convert the discrete\dash time numerical solution into a
continuous-time control signal, and 3) the discretized dynamics
\optieqref{subproblem_gusto_dt}{dynamics} result in a more sparse problem than
alternative formulations (e.g., pseudospectral), which benefits real\dash time
solution. With this choice, discretization introduces only two artifacts: the
control signal has fewer degrees of freedom, and the objective function
\optiobjref{subproblem_gusto_dt} is off from \optiobjref{subproblem_gusto_ct} by a
discretization error. Thus, by using an interpolating polynomial discretization
we can rigorously say that \pref{subproblem_gusto_dt} finds a ``local'' optimum for
a problem that is ``almost'' \pref{subproblem_gusto_ct}. We say ``local'' because
the control function has fewer DoFs, and ``almost'' due to discretization error
in the objective function. At convergence, this means that the \gusto solution
satisfies the Pontryagin maximum principle for a slightly different problem
than \pref{scp_gen_cont}. In practice, this technical discrepancy makes little to
no difference.



\subsection{Implementation Details}

We have now seen a complete description of the \scvx and \gusto algorithms. In
particular, we have all the elements that are needed to go around, and to
eventually exit, the SCP loop in \figref{scp_loop}. In this section, we discuss two
implementation details that significantly improve the performance of both
algorithms. The first detail concerns the temporal discretization procedure,
and the second detail is about variable scaling.

\subsubsection{Temporal Discretization}

The core task of discretization is to convert a continuous\dash time LTV
dynamics constraint into a discrete\dash time constraint. For \scvx, this means
converting \optieqref{subproblem_scvx_ct}{dynamics} to
\optieqref{subproblem_scvx_dt}{dynamics}. For \gusto, this means converting
\optieqref{subproblem_gusto_ct}{dynamics} to
\optieqref{subproblem_gusto_dt}{dynamics}. In all cases, our approach is to find
the equivalent discrete-time representation of the consistent flow map $\flow$
from \dref{scvx_flow_consistency}. The details of this conversion depend entirely
on the type of discretization. One example was given in detail for FOH, which
is an interpolating polynomial method, in \sbref{discretization}. This example
encapsulates a core issue with many other discretization schemes, so we will
work with it for concreteness.

In FOH discretization, the integrals in \eqref{eq:sidebar_dynamics_integrals}
require the continuous\dash time reference trajectory $\trajohone$ in order to
evaluate the corresponding Jacobians in \eqref{eq:scvx_lin_mats}. However, the
convex subproblem solution only yields a discrete\dash time trajectory. To
obtain the continuous\dash time reference, we use the continuous\dash time
input obtained directly from \eqref{eq:sidebar_control_foh} and integrate
\eqref{eq:consistent_flow_map_interpolating_poly} in tandem with the integrals in
\eqref{eq:sidebar_dynamics_integrals}. This operation is implemented over a time
interval $[t_k,t_{k+1}]$ as one big integration of a concatenated time
derivative composed of \eqref{eq:consistent_flow_map_interpolating_poly} and all
the integrands in \eqref{eq:sidebar_dynamics_integrals} \cite{Reynolds2020}. This
saves computational resources by not repeating multiple integrals, and has the
numerical advantage of letting an adaptive step integrator automatically
regulate the temporal resolution of the continuous\dash time reference
trajectory. Because the integration is reset at the end of each time interval,
the continuous\dash time reference state trajectory is discontinuous, as
illustrated in \figref{scvx_defects}. Further details are provided in
\cite{Malyuta2019b,Reynolds2020} and the source code of our numerical examples,
which is linked in \figref{github_qr}.

In theory, discretization occurs in the forward path of the SCP loop, just
before subproblem solution, as shown in \figref{scp_loop}. However, by integrating
\eqref{eq:consistent_flow_map_interpolating_poly} as described above, we can
actually obtain the \scvx defects that are needed at stage
\alglocation{\testloc} in \figref{scp_loop}. Integrating the flow map twice, once
for discretization and another time for the defects, is clearly wasteful.
Hence, in practice, we implement discretization at stage \alglocation{\testloc}
in \figref{scp_loop}, and store the result in memory to be used at the next
iteration. Note that there is no computational benefit when the flow map does
not need to be integrated, such as in
\eqref{eq:consistent_flow_map_forward_euler}. This is the case for many
discretization schemes like forward Euler, Runge-Kutta, and pseudospectral
methods. The unifying theme of these methods is that the next discrete\dash
time state is obtained algebraically as a function of (a subset of) the other
discrete\dash time states, rather than through a numerical integration process.



\subsubsection{Variable Scaling}


The convex SCP subproblems consist of the following optimization variables: the
states, the inputs, the parameter vector, and possibly a number of virtual
controls. The \alert{variable scaling} operation uses an invertible function to
transform the optimization variables into a set of new ``scaled''
variables. The resulting subproblems are completely equivalent, but the
magnitudes of the optimization variables are different
\cite{BoydConvexBook,NocedalBook}. 

\scvx and \gusto are not scale\dash invariant algorithms. This means that good
variable scaling can significantly impact not only how quickly a locally
optimal solution is found, but also the solution quality (i.e., the level of
optimality it achieves). Hence, it is imperative that the reader applies good
variable scaling when using SCP to solve trajectory problems. We will now
present a standard variable scaling technique that we have used successfully in
practice.

To motivate our scaling approach, let us review the two major effects of
variable magnitude on SCP algorithms. The first effect is common throughout
scientific computing and arises from the finite precision arithmetic used by
modern computers
\cite{SchorghoferBook,NocedalBook,BettsBook,GillPracticalBook,Ross2018}. When
variables have very different magnitudes, a lot of numerical error can
accumulate over the iterations of a numerical optimization algorithm
\cite{NocedalBook}. Managing variables of very different magnitudes is a common
requirement for numerical optimal control, where optimization problems describe
physical processes. For example, the state may include energy (measured in
Joules) and angle (measured in radians). The former might be on the order of
$10^6$ while the latter is on the order of $10^0$ \cite{MalyutaAMZ}. Most
algorithms will struggle to navigate the resulting decision space, which is
extremely elongated along the energy axis. 

The second effect of different variable magnitudes occurs in the formulation of
the trust region constraint \eqref{eq:scp_trust_region}. This constraint mixes all
of the states, inputs, and parameters into a single sum on its left\dash hand
side. We have long been taught not to compare apples and oranges, and yet
(without scaling), this is exactly what we are doing in
\eqref{eq:scp_trust_region}. Reusing the previous example, a trust region radius
$\eta=1$ means different things for an angle state than for an energy state. It
would effectively bias progress to the angle state, while allowing almost no
progress in the energy state. However, if we scale the variables to
nondimensionalize their values, then the components in the left-hand side of
\eqref{eq:scp_trust_region} become comparable and the sum is valid. Thus, variable
scaling plays an important role in ensuring that the trust region radius $\eta$
is ``meaningful'' across states, inputs, and parameters. Without variable
scaling, SCP algorithms usually have a very hard time converging to feasible
solutions.

We now have an understanding that variable scaling should seek to
nondimensionalize the state, input, and parameter vector in order to make them
comparable. Furthermore, it should bound the values to a region where finite
precision arithmetic is accurate. To this end, we have used the following
affine transformation with success across a broad spectrum of our trajectory
optimization research:
\begin{subequations}
  \label{eq:variable_scaling}
  \begin{align}
    \label{eq:x_scaling}
    x &= S_x\hat x+c_x, \\
    \label{eq:u_scaling}
    u &= S_u\hat u+c_u, \\
    \label{eq:p_scaling}
    p &= S_p\hat p+c_p,
  \end{align}
\end{subequations}
where $\hat x\in\reals^n$, $\hat u\in\reals^m$, and $\hat p\in\reals^d$ are the
new scaled variables. The user-defined matrices and offset vectors in
\eqref{eq:variable_scaling} are chosen so that the state, input, and parameter
vector are roughly bounded by a unit hypercube: $x\in [0,1]^n$, $u\in [0,1]^m$,
and $\pk\in [0,1]^d$. Another advantage to using a $[0,1]$ interval is that box
constraint lower bounds on the variables can be enforced ``for free'' if the
low\dash level convex solver operates on nonnegative variables
\cite{Reynolds2020}.

To give a concrete example of \eqref{eq:variable_scaling}, suppose that the state
is composed of two quantities: a position that takes values in
$[100,1000]~\si{\meter}$, and a velocity that takes values in
$[-10,10]~\si{\meter\per\second}$. We would then choose
$S_x=\diag\pare{900,20}\in\reals^{2\times 2}$ and
$c_x=(100,-10)\in\reals^2$. Most trajectory problems have enough
problem-specific information to find an appropriate scaling. When exact bounds
on possible variable values are unknown, an engineering approximation is
sufficient.

\subsection{Discussion}

The previous two sections make it clear that \scvx and \gusto are two instances
of SCP algorithms that share many practical and theoretical properties. The
algorithms even share similar parameters sets, which are listed
in \tabref{algo_params}.

From a practical point of view, \scvx can be shown to have superlinear
convergence rates under some mild additional assumptions. On the other hand,
the theory of \gusto allows one to leverage additional information from dual
variables to accelerate the algorithm, providing quadratic convergence
rates. Finally, both methods can be readily extended to account for
manifold-type constraints, which can be thought of as implicit, nonlinear
equality constraints, and stochastic problem settings, e.g.,
\cite{Bonalli2019a,Bonalli2019b,BonalliLewSICON2021,LewBonalliECC2020}.

Although theoretical convergence proofs of \scvx and \gusto do rely on a set of
assumptions, these assumptions are not strictly necessary for the algorithms to
work well in practice. Much of our numerical experience suggests that SCP
methods can be deployed to solve diverse and challenging problem formulations
that do not necessarily satisfy the theory of the past two sections. It is
often the case that what works best in practice, we cannot (yet) prove
rigorously in theory \cite{SzmukThesis}. As Renegar writes regarding convex
optimization \cite[p.~51]{RenegarBook}, ``It is one of the ironies of the IPM
literature that algorithms which are more efficient in practice often have
somewhat worse complexity bounds.'' The same applies for SCP methods, where
algorithms that work better in practice may admit weaker theoretical guarantees
in the general case \cite{Reynolds2020b}.

The reader should thus feel free to modify and adapt the methods based on the
requirements of their particular problem. In general, we suggest adopting a
modify\dash first\dash prove\dash later approach. For example, when using
\scvx, we are often able to significantly speed up convergence by entirely
replacing the trust region update step with a soft trust region penalty in the
cost function \cite{Reynolds2020,Reynolds2019b,SzmukReynolds2018}.

To conclude Part II, we note that \scvx, \gusto, and related SCP methods have
been used to solve a plethora of trajectory generation problems, including:
reusable launch vehicles~\cite{Goddard1920,Lu1997,Bonnard2003}, robotic
manipulators~\cite{Ratliff2009,Kalakrishnan2011,Schulman2014}, robot motion
planning~\cite{Kavraki1996,Lavalle2001,Majumdar2017}, and other examples
mentioned in the article's introduction and at the beginning of Part II. All of
these SCP variants and applications should inspire the reader to come up with
their own SCP method that works best for their particular application.

\begin{csmtable}[%
  caption={%
    A summary of the user-selected parameters for the \scvx and \gusto
    algorithms.%
  },%
  label={algo_params}]%
  \centering%
  \ifmaketwocolcsm
  \includegraphics[width=\textwidth]{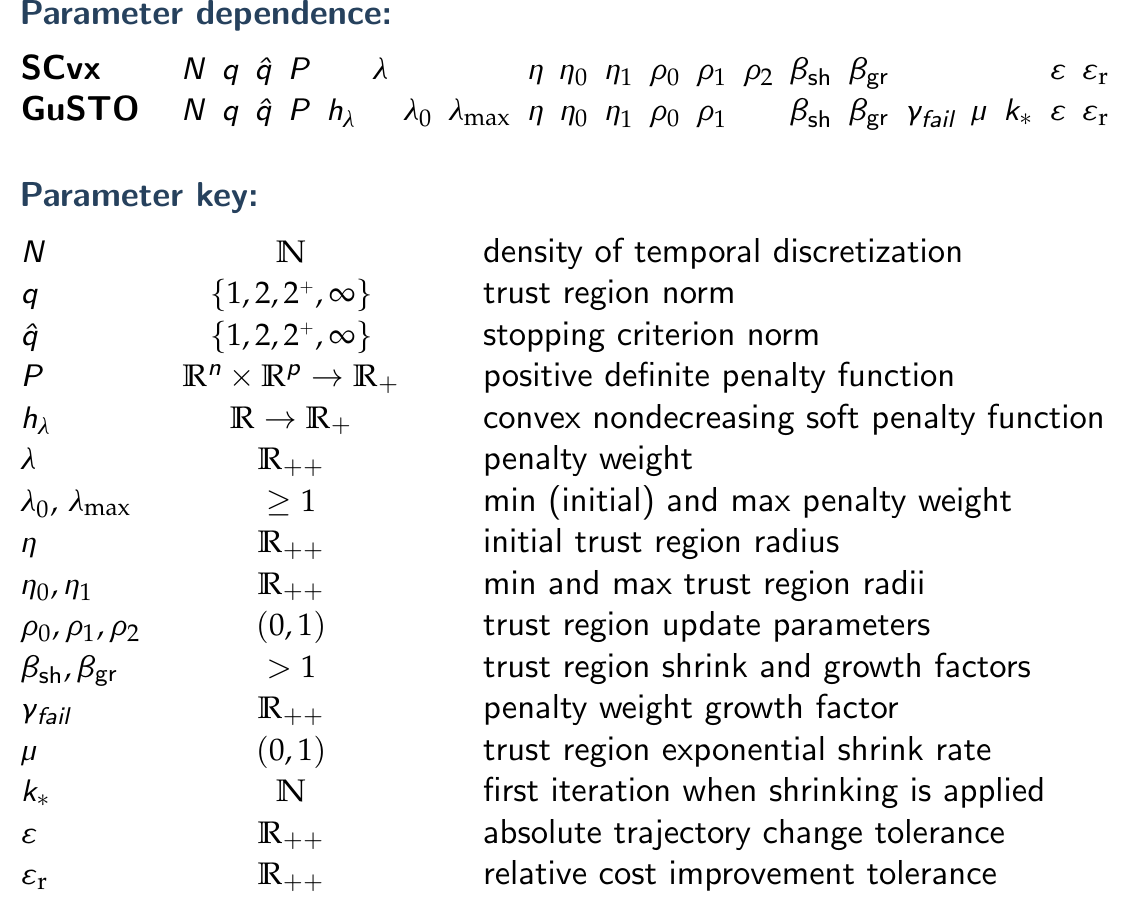}
  \else
  \includegraphics[scale=\csmpreprintfigscale]{tikz_scp_parameters}
  \fi
\end{csmtable}


\section{Part III: Application Examples}

Parts I and II of this article provided the theoretical background necessary to
start solving nonconvex trajectory generation problems using \lcvx, \scvx, and
\gusto. This final part of the article is dedicated to doing just that:
providing examples of how trajectory generation problems can be solved using
the three algorithms. We cover rocket landing first, followed by quadrotor
flight with obtacle avoidance, and finally six degree\dash of\dash freedom
flight of a robot inside a space station. The implementation code that produces
the exact plots seen in this part of the article is entirely available at the
link in \figref{github_qr}. We provide a CVX-like parser interface to \scvx and
\gusto, so the reader can leverage the code to begin solving their own
problems.

\subsection{\lcvx: 3-DoF Fuel-Optimal Rocket Landing}

Rocket\dash powered planetary landing guidance, also known as
\defintext{powered descent guidance} (PDG), was the original motivation for the
development of lossless convexification. It makes use of several \lcvx results
that we presented, making it a fitting first example.

Work on PDG began in the late 1960s~\cite{Meditch1964}, and has since been
extensively studied \cite{Ploen2006}. The objective is to find a sequence of
thrust commands that transfer the vehicle from a specified initial state to a
desired landing site without consuming more propellant than what is
available. Using optimization to compute this sequence can greatly increase the
range of feasible landing sites and the precision of the final landing location
\cite{Carson2011b,Wolf2006,Wolf2012}.

Lossless convexification for PDG was first introduced in
\cite{Acikmese2005,Acikmese2007}, where minimizing fuel usage was the
objective. It was the first time that convex optimization was shown to be
applicable to PDG. This discovery unlocked a polynomial time algorithm with
guaranteed convergence properties for generating optimal landing
trajectories. The work was later expanded to handle state constraints
\cite{Acikmese2011,Harris2013b,Harris2013a}, to solve a mininimum landing\dash
error problem \cite{Blackmore2010}, to include nonconvex pointing constraints
\cite{Carson2011,acikmese2013flight}, and to handle nonlinear terms in the
dynamics such as aerodynamic drag and nonlinear gravity \cite{Blackmore2012}.

Today, we know that SpaceX uses convex optimization for the Falcon 9 rocket
landing algorithm \cite{lars2016autonomous}. Flight tests have also been
performed using lossless convexification in a collaboration between NASA and
Masten Space Systems \cite{Scharf2017}, as shown in \figref{lcvx_masten}.

\begin{csmfigure}[%
  caption={%
    The Masten Xombie rocket near the end of a 750~m divert maneuver. Figure
    reproduced with permission from \cite[Fig.~1]{Scharf2017}.%
  },
  label={lcvx_masten}]%
  \centering%
  \ifmaketwocolcsm
  \includegraphics[width=0.9\columnwidth]{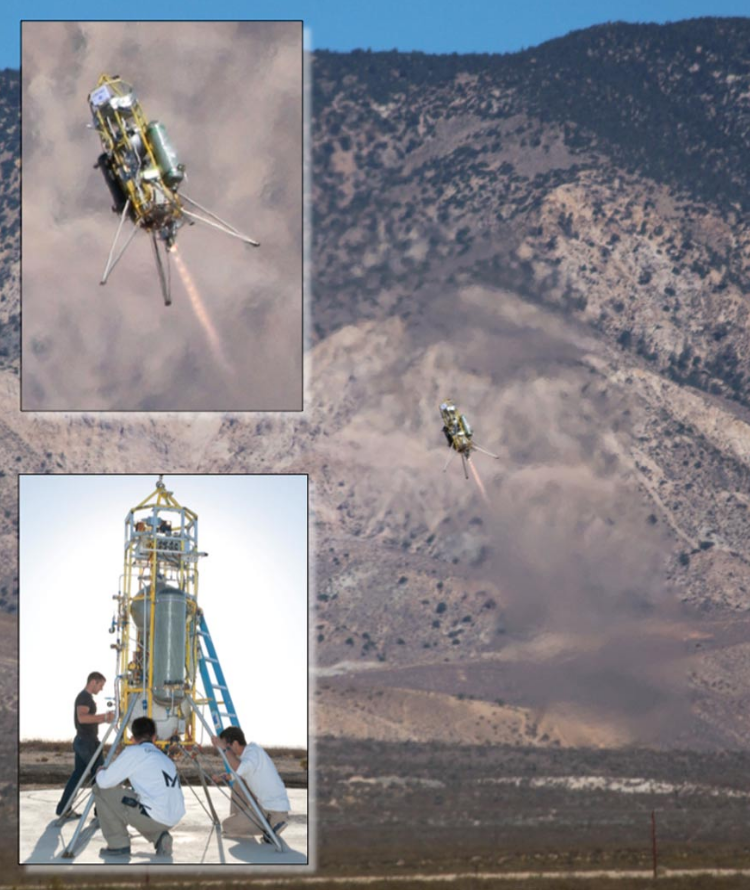}
  \else
  \includegraphics[width=0.7\columnwidth]{masten_adapt}
  \fi
\end{csmfigure}

We will now present a lossless convexification PDG example based on a mixture
of original ideas from \cite{Acikmese2007,Carson2011}. Note that \lcvx considers
the 3\dash degree\dash of\dash freedom (DoF) PDG problem, where the vehicle is
modeled as a point mass. This model is accurate as long as attitude can be
controlled in an inner loop faster than the outer translation control
loop. This is a valid assumption for many vehicles, including rockets and
aircraft. The thrust vector is taken as the control input, where its direction
serves as a proxy for the vehicle attitude. We begin by stating the raw minimum
fuel 3\dash DoF PDG problem.

\begin{csmfigure}[%
  caption={%
    Illustration of the 3-DoF PDG problem, showing some of the relevant
    constraints on the rocket-powered lander's trajectory. The thrust direction
    $T_c(t)/\norm[2]{T_c(t)}$ servers as a proxy for vehicle attitude.%
  },
  label={rocket_landing_setup},
  columns=2]%
  \centering%
  \ifmaketwocolcsm
  \includegraphics[width=0.85\textwidth]{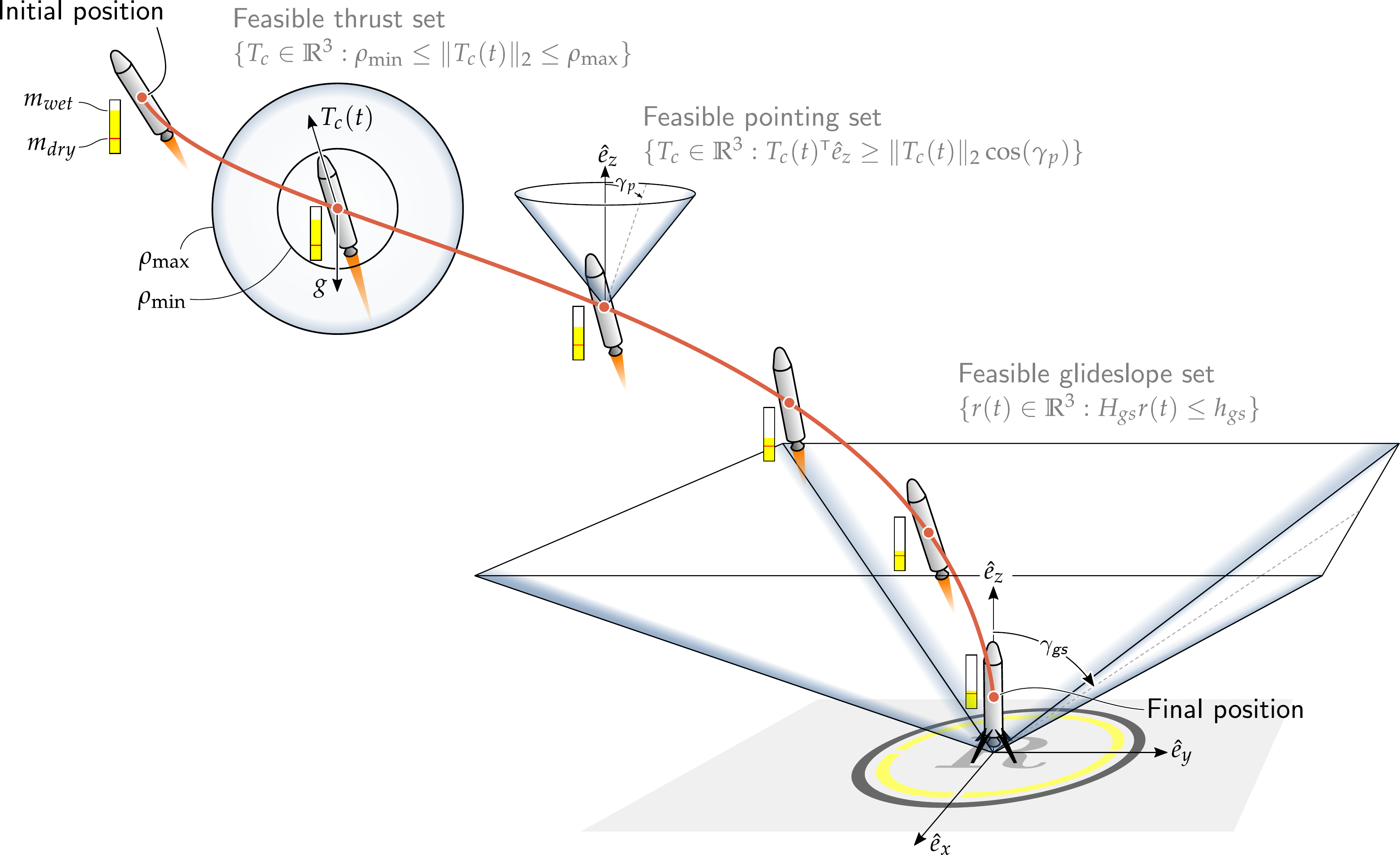}
  \else
  \includegraphics[width=\textwidth]{rocket_landing_setup}
  \fi
\end{csmfigure}

\begin{optimization}[
  label={lcvx_o_pdg},
  variables={T_c,t_f},
  objective={\int_0^{t_f}\norm[2]{T_c(t)}\dt}]
  \optilabel{rdot}
  & \dot r(t) = v(t), \\
  \optilabel{vdot}
  & \dot v(t) = g+\frac{T_c(t)}{m(t)}-\omega^\times\omega^\times
  r(t)-2\omega^\times v(t), \\
  \optilabel{mdot}
  & \dot m(t) = -\alpha\norm[2]{T_c(t)}, \\
  \optilabel{thrust_bounds}
  & \rho_{\min}\le\norm[2]{T_c(t)}\le\rho_{\max}, \\
  \optilabel{pointing}
  & T_c(t)\T\hat e_z\ge \norm[2]{T_c(t)}\cos(\gamma_{p}), \\
  \optilabel{glideslope}
  & H_{gs}r(t)\le h_{gs}, \\
  \optilabel{speed}
  & \norm[2]{v(t)} \le v_{\max}, \\
  \optilabel{dry}
  & m_{dry}\le m(t_f), \\
  \optilabel{bdry1}
  & r(0)=r_0,~v(0)=v_0,~m(0)=m_{wet}, \\
  \optilabel{bdry2}
  & r(t_f)=v(t_f)=0.
\end{optimization}

The vehicle translational dynamics correspond to a double integrator with
variable mass, moving in a constant gravitational field, viewed in the planet's
rotating frame. In particular, $r\in\real^3$ is the position, $v\in\real^3$ is
the velocity, $g\in\real^3$ is the gravitational acceleration, and
$\omega\in\real^3$ is the planet's constant angular velocity. The notation
$\omega^\times$ denotes the skew-symmetric matrix representation of the cross
product $\omega\times(\cdot)$. The mass $m\in\real$ is depleted by the rocket
engine according to \optieqref{lcvx_o_pdg}{mdot} with the fuel consumption rate
\begin{equation}
  \alpha\definedas\frac{1}{I_{sp}g_e},
\end{equation}
where $g_e\approx 9.807~\si{\meter\per\second\squared}$ is the Earth's standard
gravitational acceleration and $I_{sp}~\si{\second}$ is the rocket engine's
specific impulse. As illustrated in \figref{rocket_landing_setup},
$T_c\in\reals^3$ is the rocket engine thrust vector, which is upper and lower
bounded via \optieqref{lcvx_o_pdg}{thrust_bounds}. The lower bound was motivated in
\sbref{inputrelax}. The thrust vector, and therefore the vehicle attitude, also has
a tilt angle constraint \optieqref{lcvx_o_pdg}{pointing} which prevents the vehicle
from deviating by more than angle $\gamma_p$ away from the vertical. Following
the discussion in \sbref{sidebar_affinestate}, we also impose an affine glideslope
constraint via \optieqref{lcvx_o_pdg}{glideslope} with a maximum glideslope angle
$\gamma_{gs}$. The velocity is constrained to a maximum magnitude $v_{\max}$ by
\optieqref{lcvx_o_pdg}{speed}. The final mass must be greater than $m_{dry}$
\optieqref{lcvx_o_pdg}{dry}, which ensures that no more fuel is consumed than what
is available. Constraints \optieqref{lcvx_o_pdg}{bdry1} and
\optieqref{lcvx_o_pdg}{bdry2} impose fixed boundary conditions on the rocket's
state. In particular, the starting mass is $m_{wet}>m_{dry}$.

To begin the lossless convexification process, the standard input slack
variable $\sigma\in\real$ from \sbref{inputrelax} is introduced in order to remove
the nonconvex lower bound in \optieqref{lcvx_o_pdg}{thrust_bounds}. As a
consequence, the familiar \lcvx equality constraint appears in
\optieqref{lcvx_r1_pdg}{lcvx_equality}. Following \sbref{pointingrelax}, this also
replaces $\norm[2]{T_c(t)}$ in \optieqref{lcvx_o_pdg}{pointing} with
$\sigma(t)$ in \optieqref{lcvx_r1_pdg}{pointing}.

\begin{optimization}[
  label={lcvx_r1_pdg},
  variables={\sigma,T_c,t_f},
  objective={\int_0^{t_f}\sigma(t)\dt}
  ]
  & \dot r(t) = v(t), \\
  \optilabel{vdot}
  & \dot v(t) = g+\frac{T_c(t)}{m(t)}-\omega^\times\omega^\times
  r(t)-2\omega^\times v(t), \\
  & \dot m(t) = -\alpha\sigma(t), \\
  \optilabel{sigma_bounds}
  & \rho_{\min}\le\sigma(t)\le\rho_{\max}, \\
  \optilabel{lcvx_equality}
  & {\color{lcvxColor}\norm[2]{T_c(t)}\le\sigma(t)}, \\
  \optilabel{pointing}
  & T_c(t)\T\hat e_z\ge \sigma(t)\cos(\gamma_{p}), \\
  & H_{gs}r(t)\le h_{gs}, \\
  & \norm[2]{v(t)} \le v_{\max}, \\
  & m_{dry}\le m(t_f), \\
  & r(0)=r_0,~v(0)=v_0,~m(0)=m_{wet}, \\
  & r(t_f)=v(t_f)=0.
\end{optimization}

Next, we approximate the nonlinearity $T_c/m$ in
\optieqref{lcvx_r1_pdg}{vdot}. To this end, \cite{Acikmese2007} showed that the
following change of variables can be made:
\begin{equation}
  \label{eq:lcvx_pdg_cov}
  \xi\definedas\frac{\sigma}{m},~u\definedas\frac{T_c}{m},~
  z\definedas\ln(m).
\end{equation}

Using the new variables, note that
\begin{equation}
  \frac{\dot m(t)}{m(t)} = -\alpha\xi(t)~\implies~z(t)=z(0)-\alpha\int_0^{t}\xi(\tau)\dd\tau.
\end{equation}

Since the cost in \optiobjref{lcvx_r1_pdg} maximizes $m(t_f)$, and since
$\alpha>0$, an equivalent cost is to minimize $\int_0^{t_f}\xi(t)\dt$. It turns
out that the new variables linearize all constraints except for the upper bound
part of \optieqref{lcvx_r1_pdg}{sigma_bounds}. In the new variables,
\optieqref{lcvx_r1_pdg}{sigma_bounds} becomes:
\begin{equation}
  \label{eq:lcvx_pdg_nonlinear_exp}
  \rho_{\min}\exp{-z(t)}\le\xi(t)\le\rho_{\max}\exp{-z(t)}.
\end{equation}

To keep the optimization problem an SOCP, it is desirable to also do something
about the lower bound in \eqref{eq:lcvx_pdg_nonlinear_exp}, which is a convex
exponential cone. In \cite{Acikmese2007}, it was shown that the following Taylor
series approximation is accurate to within a few percent, and is therefore
acceptable:
\begin{subequations}
  \label{eq:exp_approx}
  \begin{align}
    \label{eq:1}
    \mu_{\min}(t)\left[1-\delta z(t)+\frac 12 \delta z(t)^2\right]
    &\le \xi(t), \\
    \mu_{\max}(t)\bigl[1-\delta z(t)\bigr]
    &\ge \xi(t),
  \end{align}
\end{subequations}
where:
\begin{subequations}
  \begin{align}
    \mu_{\min}(t) &= \rho_{\min}\exp{-z_0(t)}, \\
    \mu_{\max}(t) &= \rho_{\max}\exp{-z_0(t)}, \\
    \delta z(t) &= z(t)-z_0(t), \\
    z_0(t) &= \ln(m_{wet}-\alpha\rho_{\max}t).
  \end{align}
\end{subequations}

The reference profile $z_0(\cdot)$ corresponds to the maximum fuel rate, thus
$z_0(t)$ lower bounds $z(t)$. To ensure that physical bounds on $z(t)$ are not
violated, the following extra constraint is imposed:
\begin{equation}
  \label{eq:z_phys_bnd}
  z_0(t)\le z(t)\le\ln(m_{wet}-\alpha\rho_{\min}t).
\end{equation}

The constraints \eqref{eq:exp_approx} and \eqref{eq:z_phys_bnd} together approximate
\eqref{eq:lcvx_pdg_nonlinear_exp} and have the important property of being
conservative with respect to the original constraint
\cite[Lemma~2]{Acikmese2007}. Thus, using this approximation \emphasize{will not
  generate solutions that are infeasible for the original problem}. We can now
write the finalized convex relaxation of \pref{lcvx_o_pdg}.

\begin{optimization}[
  label={lcvx_r2_pdg},
  variables={\xi,u,t_f},
  objective={\int_0^{t_f}\xi(t)\dt}]
  \optilabel{rdot}
  & \dot r(t) = v(t), \\
  \optilabel{vdot}
  & \dot v(t) = g+u(t)-\omega^\times\omega^\times r(t)-2\omega^\times v(t), \\
  \optilabel{zdot}
  & \dot z(t) = -\alpha\xi(t), \\
  \optilabel{xi_lb}
  & \mu_{\min}(t)\left[1-\delta z(t)+\frac 12 \delta z(t)^2\right]\le\xi(t), \\
  \optilabel{xi_ub}
  & \mu_{\max}(t)\bigl[1-\delta z(t)\bigr]\ge\xi(t), \\
  \optilabel{lcvx_equality}
  & {\color{lcvxColor}\norm[2]{u(t)}\le\xi(t)}, \\
  \optilabel{pointing}
  & u(t)\T\hat e_z\ge\xi(t)\cos(\gamma_p), \\
  \optilabel{glideslope}
  & H_{gs}r(t)\le h_{gs}, \\
  \optilabel{speed}
  & \norm[2]{v(t)} \le v_{\max}, \\
  \optilabel{z_dry}
  & \ln(m_{dry})\le z(t_f), \\
  \optilabel{z_bounds}
  & z_0(t)\le z(t)\le\ln(m_{wet}-\alpha\rho_{\min}t), \\
  & r(0)=r_0,~v(0)=v_0,~z(0)=\ln(m_{wet}), \\
  & r(t_f)=v(t_f)=0.
\end{optimization}

Several conditions must now be checked to ascertain lossless convexification,
i.e., that the solution of \pref{lcvx_r2_pdg} is globally optimal for
\pref{lcvx_o_pdg}. To begin, let us view $z$ in \optieqref{lcvx_r2_pdg}{zdot} as a
``fictitious'' state that is used for concise notation. In practice, we know
explicitly that $z(t)=\ln(m_{wet})-\alpha\int_0^{t}\xi(t)\dt$ and thus we can
replace every instance of $z(t)$ in
\optieqref{lcvx_r2_pdg}{xi_lb}-\optieqref{lcvx_r2_pdg}{z_bounds} with this
expression. Thus, we do \emphasize{not} consider $z$ as part of the state and,
furthermore, \optieqref{lcvx_r2_pdg}{xi_lb}, \optieqref{lcvx_r2_pdg}{xi_ub},
\optieqref{lcvx_r2_pdg}{z_dry}, and \optieqref{lcvx_r2_pdg}{z_bounds} are
\emphasize{input} and not state constraints. Defining the state as
$x=(r,v)\in\reals^6$, the state\dash space matrices are
\begin{equation}
  A =
  \Matrix{
    0 & I_3 \\
    -\omega^\times\omega^\times & -2\omega^\times
  },~
  B =
  \Matrix{
    0 \\ I_3
  },~E=B.
\end{equation}

\begin{csmfigure}[%
  caption={%
    A force balance in the normal direction $\hat n_{gs}$ can be used to
    guarantee that the glideslope constraint can only be activated
    instantaneously.%
  },
  label={rocket_glideslope_balance}]%
  \centering%
  \ifmaketwocolcsm
  \includegraphics[width=0.5\columnwidth]{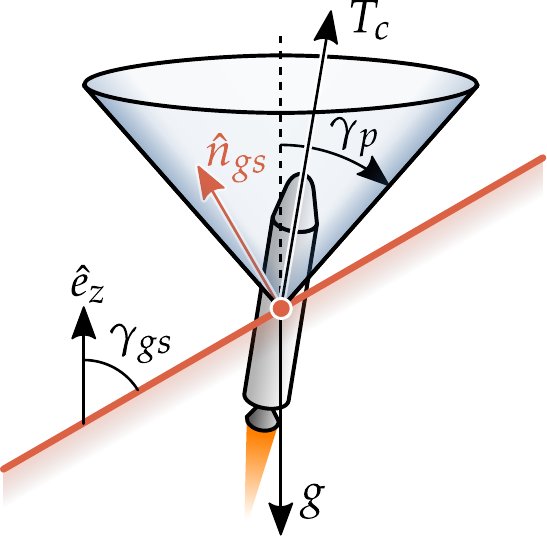}
  \else
  \includegraphics[width=0.4\columnwidth]{rocket_glideslope_balance}
  \fi
\end{csmfigure}

The pair $\{A,B\}$ is unconditionally controllable and therefore
\conref{lcvx_nostate_controllability} is satisfied. We also need to verify
\conref{lcvx_nostate_controllability_pointing} due to the presence of the pointing
constraint \optieqref{lcvx_r2_pdg}{pointing}. In this case $\hat n_u=\hat e_z$ and
$N=\Matrix{\hat e_x & \hat e_y}$. \conref{lcvx_nostate_controllability_pointing}
holds as long as $\omega^\times \hat e_z\ne 0$, which means that the planet
does not rotate about the local vertical of the landing frame.

The glideslope constraint \optieqref{lcvx_r2_pdg}{glideslope} can be treated
either via \conref{lcvx_linstate_controllability} or by checking that it can
only be instantaneously active. In this case, we can prove the latter by
considering a force balance in the glideslope plane's normal direction
$\hat n_{gs}$, as illustrated in \figref{rocket_glideslope_balance}. We
also invoke the fact that the optimal thrust profile for the rocket landing
problem is bang-bang \cite{Acikmese2007}. This allows to distill the
verification down to a conservative condition that the following inequalities
hold for all
$\theta\in [\frac\pi 2-\gamma_p-\gamma_{gs},\frac\pi 2+\gamma_p-\gamma_{gs}]$:
\begin{subequations}
  \begin{align}
    \rho_{\min}\cos(\theta) &< m_{dry}\norm[2]{g}\sin(\gamma_{gs}), \\
    \rho_{\max}\cos(\theta) &> m_{wet}\norm[2]{g}\sin(\gamma_{gs}).
  \end{align}
\end{subequations}

\begin{csmfigure}[%
  caption={%
    Various views of the globally optimal rocket landing trajectory obtained
    via lossless convexification for \pref{lcvx_o_pdg}. Circular markers show the
    discrete\dash time trajectory directly returned from the convex
    solver. Lines show the continuous\dash time trajectory, which is obtained
    by numerically integrating the ZOH discretized control signal through the
    dynamics \optieqref{lcvx_o_pdg}{rdot}-\optieqref{lcvx_o_pdg}{mdot}. The exact match
    at every discrete\dash time node demonstrates that the solution is
    dynamically feasible for the actual continuous\dash time vehicle.
  },%
  label={lcvx_pdg_results},%
  columns=2,
  position={!b}]%
  \centering
  \includegraphics[width=1\textwidth]{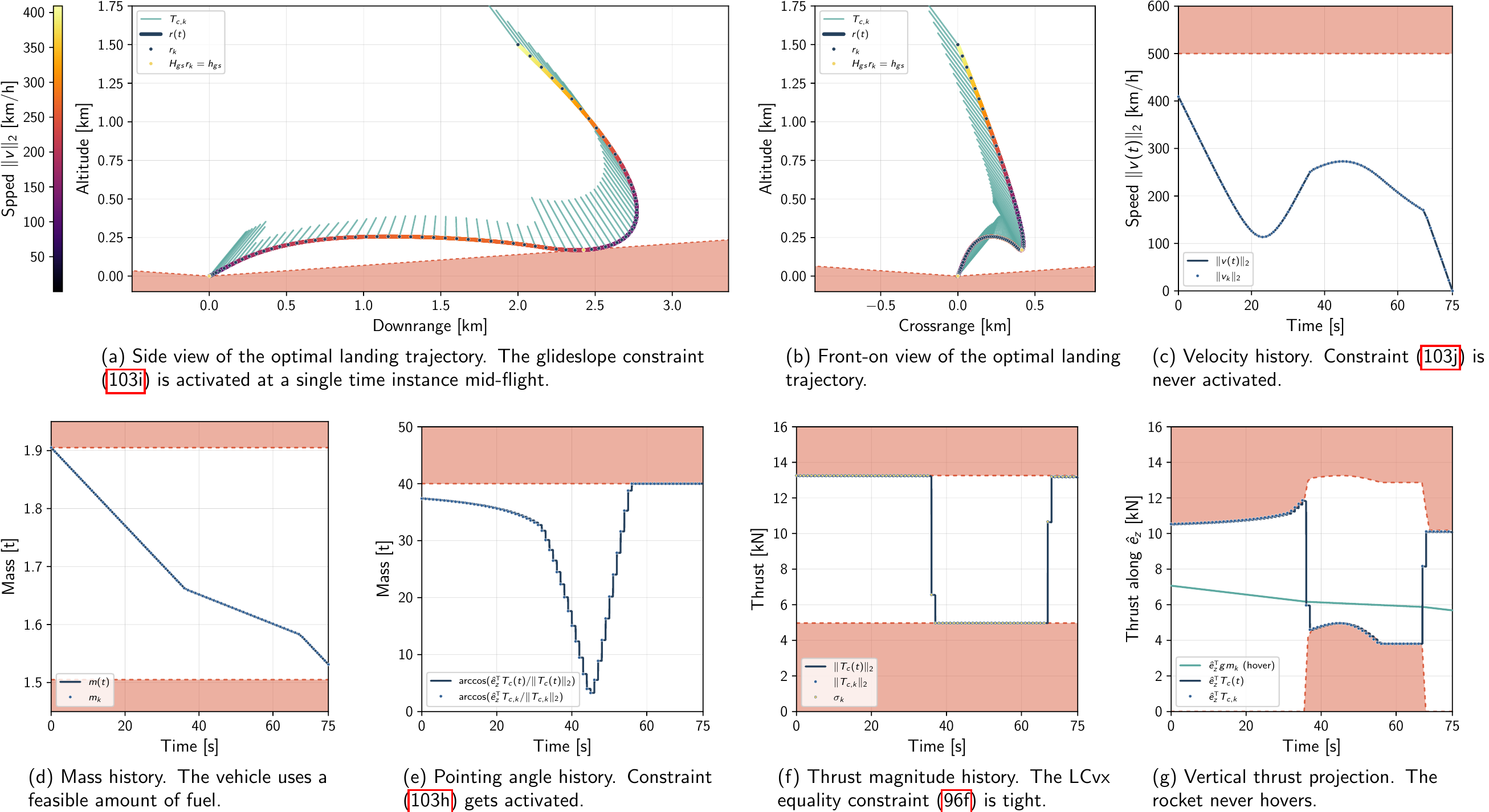}
\end{csmfigure}

It turns out that for the problem parameters \eqref{eq:lcvx_3dof_pdg_parameters}
of this example, the above inequalities hold. Therefore
\optieqref{lcvx_r2_pdg}{glideslope} can only be instantaneously active and does
not pose a threat to lossless convexification.

We also need to check \conref{lcvx_nostate_linindep}, which relates to the
transversality condition of the maximum principle
\cite{PontryaginBook,BerkovitzBook,LiberzonBook}. In the case of
\pref{lcvx_r2_pdg}, we have $m[t_f]=0$ and $b[t_f]=x(t_f)$, hence
\begin{equation}
  m_{\text{\lcvx}} =
  \Matrix{
    0\in\real^{6} \\
    \xi(t_f)
  },~
  B_{\text{\lcvx}} =
  \Matrix{
    I_6 \\
    0
  }.
\end{equation}

Hence, as long as $\xi(t_f)>0$, \conref{lcvx_nostate_linindep} holds. Since
$\xi(t_f)>0$ is guaranteed by the fact that the lower bound in
\eqref{eq:exp_approx} is greater than the exact lower bound
$\rho_{\min}\exp{-z(t_f)}>0$ \cite[Lemma~2]{Acikmese2007},
\conref{lcvx_nostate_linindep} holds.

The remaining constraint to be checked is the maximum velocity bound
\optieqref{lcvx_r2_pdg}{speed}. Although it can be written as the quadratic
constraint $v_{\max}^{-2}v(t)\T v(t)\le 1$, the dynamics
\optieqref{lcvx_r2_pdg}{rdot}-\optieqref{lcvx_r2_pdg}{vdot} do not match the form
\optieqref{lcvx_o_quadstate}{dynamics_1}-\optieqref{lcvx_o_quadstate}{dynamics_2}
required by the \lcvx result for quadratic state constraints. Thus, we must
resort to the more restricted statement for general state constraints in
\tref{genstate}. According to \tref{genstate}, we conclude that lossless
convexification will hold as long as the maximum velocity bound
\optieqref{lcvx_r2_pdg}{speed} is activated at most a discrete number of times. In
summary, we can conclude that the solution of \pref{lcvx_r2_pdg} is guaranteed to
be globally optimal for \pref{lcvx_o_pdg} as long as \optieqref{lcvx_r2_pdg}{speed} is
never persistently active.

We now present a trajectory result obtained by solving
\pref{lcvx_r2_pdg} using a ZOH discretized input signal with a
$\Delta t=1~\si{\second}$ time step (see \sbref{discretization}). We use the
following numerical parameters:
\begin{subequations}
  \label{eq:lcvx_3dof_pdg_parameters}
  \begin{gather}
    \label{eqref:pdg_pars}
    g = -3.71\hat e_z~\si{\meter\per\second\squared},~
    m_{dry} = 1505~\si{\kilogram}, \\
    m_{wet} = 1905~\si{\kilogram},~
    I_{sp} = 225~\si{\second}, \\
    \omega = (3.5\hat e_x+2\hat e_z)\cdot 10^{-3}%
    ~\si{\degree\per\second}, \\
    \rho_{\min} = 4.971~\si{\kilo\newton},~
    \rho_{\max} = 13.258~\si{\kilo\newton}, \\
    \gamma_{gs} = 86~\si{\degree},~%
    \gamma_{p} = 40~\si{\degree},~
    v_{\max} = 500~\si{\kilo\meter\per\hour}, \\
    r_0 = 2\hat e_x+1.5\hat e_z~\si{\kilo\meter}, \\
    v_0 = 288\hat e_x+108\hat e_y-270\hat e_z~\si{\kilo\meter\per\hour}.
  \end{gather}
\end{subequations}

Golden search \cite{Kochenderfer2019} is used to find the optimal time of flight
$t_f$. This is a valid choice because the cost function is unimodal with respect
to $t_f$ \cite{Blackmore2010}. For the problem parameters in
\eqref{eq:lcvx_3dof_pdg_parameters}, an optimal rocket landing trajectory is
found with a minimum fuel time of flight $\optimal{t_f}=75~\si{\second}$.

\figref{lcvx_pdg_results} visualizes the computed optimal landing trajectory. From
\fakesubfigref{lcvx_pdg_results}{a}, we can clearly see that the glide slope constraint is
not only satisfied, it is also activated only twice (once mid\dash flight and
another time at touchdown), as required by the \lcvx guarantee. From
\fakesubfigref{lcvx_pdg_results}{c}, we can see that the velocity constraint is never
activated. Hence, in this case, it does not pose a threat to lossless
convexification. Several other interesting views of the optimal trajectory are
plotted in \fakesubfigref[g]{lcvx_pdg_results}{d}. In all cases, the state and input
constraints of \pref{lcvx_o_pdg} hold. The fact that \pref{lcvx_r2_pdg} is a lossless
convexification of \pref{lcvx_o_pdg} is most evident during the minimum-thrust
segment from about 40 to 70 seconds in \fakesubfigref{lcvx_pdg_results}{f}. Here, it is
important realize that $\norm[2]{T_c(t)}< \rho_{\min}$ is feasible for the
relaxed problem. The fact that this never occurs is a consequence of the
lossless convexification guarantee that \optieqref{lcvx_r2_pdg}{lcvx_equality}
holds with equality.

In conclusion, we emphasize that the rocket landing trajectory presented in
\figref{lcvx_pdg_results} is not just a feasible solution, nor even a locally
optimal one, but it is a \textit{globally} optimal solution to this rocket
landing problem. This means that one can do no better given these parameters
and problem description.


\subsection{SCP: Quadrotor Obstacle Avoidance}

\begin{csmfigure}[%
  caption={%
    A quadrotor at the University of Washington's Autonomous Controls
    Laboratory, executing a collision-free trajectory computed by \scvx
    \cite{Szmuk2017,mceowen_tablet,ACLYoutubeObstacles}.
  },
  label={scvx_quadrotor_obstacle},
  position={!t}]%
  \centering%
  \includegraphics[width=0.9\columnwidth]{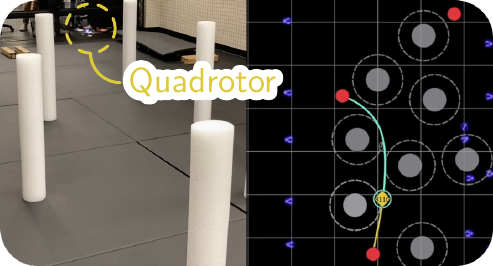}
\end{csmfigure}

We now move on to trajectory generation for problems which cannot be handled by
\lcvx. The objective of this first example is to compute a trajectory for a
quadrotor that flies from one point to another through an obstacle-filled
flight space. This example is sourced primarily from~\cite{Szmuk2017}, and a
practical demonstration is shown in \figref{scvx_quadrotor_obstacle}.

The quadrotor is modeled as a point mass, which is a reasonable approximation
for small and agile quadrotors whose rotational states evolve over much shorter
time scales than the translational states \cite{Acikmese2007}. We express the
equations of motion in an East-North-Up (ENU) inertial coordinate system, and
take the state vector to be the position and velocity. Using a simple
double-integrator model, the continuous\dash time equations of motion are
expressed in the ENU frame as:
\begin{equation}
  \label{eq:ex_quad_dynamics}
  \ddot{r}(\tabs) = a(\tabs) - g\hat{n},
\end{equation}
where $r\in\real^3$ denotes the position, $g\in\reals$ is the (constant)
acceleration due to gravity, $\hat{n} = (0,0,1)\in\real^3$ is the ``up''
direction, and $a(\tabs)\in\real^3$ is the commanded acceleration, which is the
control input. The time $\tabs$ spans the interval $[0,\tf]$. Because the
previous section on SCP used a normalized time $t\in [0,1]$, we call $\tabs$
the absolute time and use the special font to denote it. We allow the
trajectory duration to be optimized, and impose the following constraint to
keep the final time bounded:
\begin{equation}
  \label{eq:ex_quad_tf_bounds}
  \tf[,\min] \le \tf \le \tf[,\max],
\end{equation}
where $\tf[,\min]\in\reals$ and $\tf[,\max]\in\reals$ are user-defined
parameters. Boundary conditions on the position and velocity are imposed to
ensure that the vehicle begins and ends at the desired states:
\begin{subequations}
  \label{eq:ex_quad_bcs}
  \begin{alignat}{3}
    r(\ti) &= r_0, \quad \dot{r}(\ti) &&= v_0, \\
    r(\tf) &= r_f, \quad \dot{r}(\tf) &&= v_f.
  \end{alignat}
\end{subequations}

The magnitude of the commanded acceleration is limited from above and below by
the electric motor and propeller configuration, and its direction is
constrained to model a tilt angle constraint on the quadrotor. In effect, the
acceleration direction is used as a proxy for the vehicle attitude. This is an
accurate approximation for a ``flat'' quadrotor configuration where the
propellers are not canted with respect to the plane of the quadrotor
body. Specifically, we enforce the following control constraints:
\begin{subequations}
  \label{eq:ex_quad_cc}
  \begin{align}
    \label{eq:ex_quad_nrm_bnd}
    a_{\min} \leq \norm[2]{a(\tabs)} &\leq a_{\max}, \\
    \label{eq:ex_quad_tilt_bnd}
    \norm[2]{a(\tabs)} \cos \tiltmax &\leq \hat{n}^\transp a(\tabs),
  \end{align}
\end{subequations}
where $0<a_{\min}<a_{\max}$ are the acceleration bounds and
$\tiltmax\in(0\si{\degree},180\si{\degree}]$ is the maximum angle by which the
acceleration vector is allowed to deviate from the ``up'' direction.

Finally, obstacles are modeled as three-dimensional ellipsoidal keep-out
zones, described by the following nonconvex constraints:
\begin{equation}
  \label{eq:ex_quad_sc_obs}
  \norm[2]{H_j \big( r(\tabs) - c_j \big)} \geq 1, \quad
  j=1,\ldots,\Nobs,
\end{equation}
where $c_j\in\real^3$ denotes the center, and $H_j\in\pd^3$ defines the size
and shape of the $j$-th obstacle. Note that the formulation allows the
obstacles to intersect. 

Using the above equations, we wish to solve the following free final time
optimal control problem that minimizes control energy:
\begin{optimization}[%
  label={ex_quad_ocp},
  variables={u,\tf},
  objective={\int_{0}^{\tf}\norm[2]{a(\tabs)}^2\sdd\tabs}]%
  &%
  \textnormal{\eqref{eq:ex_quad_dynamics}-\eqref{eq:ex_quad_sc_obs}.}
\end{optimization}

Due to the presence of nonconvex state constraints \eqref{eq:ex_quad_sc_obs}, only
embedded \lcvx applies to the above problem. In particular, \lcvx can be used
to handle the nonconvex input lower bound in \eqref{eq:ex_quad_nrm_bnd} along with
the tilt constraint in \eqref{eq:ex_quad_tilt_bnd}. This removes some of the
nonconvexity. However, it is only a partial convexification which still leaves
behind a nonconvex optimal control problem. Thus, we must resort to SCP
techniques for the solution.

\subsubsection{\scvx Formulation}

\begin{csmtable}[
  caption={Algorithm parameters for the quadrotor obstacle avoidance example.},
  label={ex_quad_params},
  position={!t}]%
  \centering%
  \ifmaketwocolcsm
  \includegraphics[width=\textwidth]{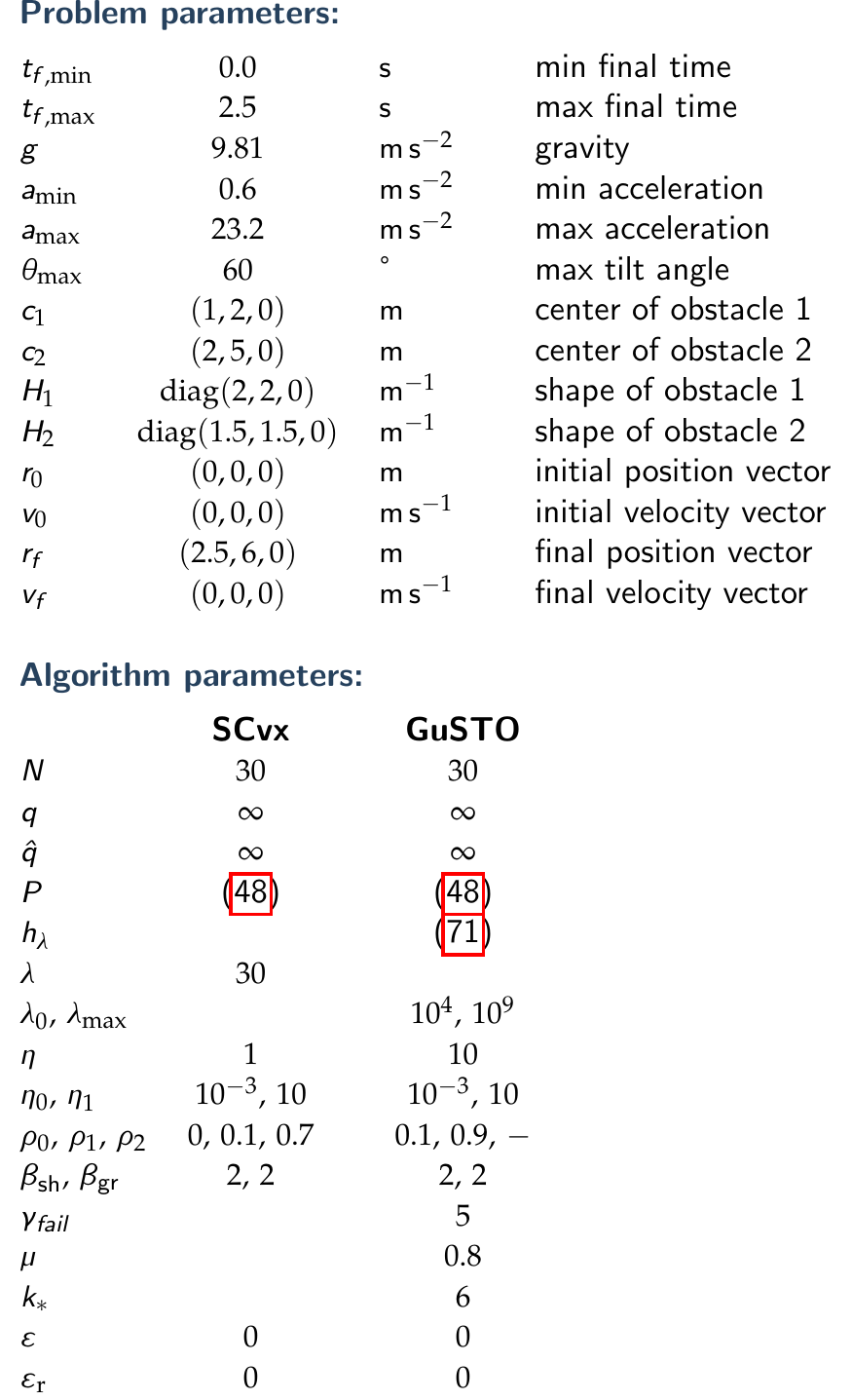}
  \else
  \includegraphics[scale=\csmpreprintfigscale]{tikz_ex_quad_parameters}
  \fi
\end{csmtable}

We begin by demonstrating how \scvx can be used to solve \pref{ex_quad_ocp}. To
this end, we describe how the problem can be cast into the template of
\pref{scp_gen_cont}. Once this is done, the rest of the solution process is
completely automated by the mechanics of the \scvx algorithm as described in
the previous section. To make the notation lighter, we will omit the argument
of time whenever possible.

Let us begin by defining the state and input vectors. For the input vector,
note that the nonconvex input constraints in \eqref{eq:ex_quad_cc} can be
convexified via embedded \lcvx. In particular, the relaxation used for
\pref{lcvx_r_nostate_pointing} can losslessly convexify both input constraints by
introducing a slack input $\sigma\in\real$ and rewriting~\eqref{eq:ex_quad_cc} as:
\begin{subequations}
  \label{eq:ex_quad_cc_cvx}
  \begin{align}
    a_{\min} \leq \sigma &\leq a_{\max}, \\
    \label{eq:ex_quad_cc_cvx_lcvx_equality}
    \color{red}\norm[2]{a} &\color{red}\leq \sigma{\color{black},} \\
    \sigma \cos \tiltmax &\leq \hat{n}^\transp a,
  \end{align}
\end{subequations}
where \eqref{eq:ex_quad_cc_cvx_lcvx_equality} is the familiar \lcvx equality
constraint. Thus, we can define the following state and ``augmented'' input
vectors:
\begin{equation}
  \label{eq:ex_quad_scvx_state_input}
  x=\Matrix{r \\ \dot r}\in\reals^6,\quad
  u=\Matrix{a \\ \sigma}\in\reals^4.
\end{equation}

Next, we have to deal with the fact that \pref{scp_gen_cont} uses normalized time
$t\in [0,1]$, while \pref{ex_quad_ocp} uses absolute time $\tabs\in [0,\tf]$. To
reconcile the two quantities, we use a one\dash dimensional parameter vector
$p\in\reals$. The parameter defines a \alert{time dilation} such that the
following relation holds:
\begin{equation}
  \label{eq:ex_quad_time_dilation}
  \tabs=p t,
\end{equation}
from which it follows that $p\equiv\tf$. In absolute time, the dynamics are
given directly by writing \eqref{eq:ex_quad_dynamics} in terms of
\eqref{eq:ex_quad_scvx_state_input}, which gives a set of time-invariant
first-order ordinary differential equations:
\begin{equation}
  \label{eq:ex_quad_dynamics_abs_time}
  f\pare{x, u} = \Matrix{
    \dot r \\
    a-g\hat n}.
\end{equation}

For \pref{scp_gen_cont}, we need to convert the dynamics to normalized time. We do
so by applying the chain rule:
\begin{equation}
  \label{eq:ex_quad_dynamics_chain_rule}
  \frac{\dd x}{\dt} = \frac{\dd x}{\dd\tabs}\frac{\dd\tabs}{\dt} =
  p f\pare[big]{x, u},
\end{equation}
which yields the dynamics \optieqref{scp_gen_cont}{dynamics} in normalized time:
\begin{equation}
  \label{eq:ex_quad_dynamics_normalized_time}
  f\pare{x, u, p} = p f\pare[big]{x, u}.
\end{equation}

The convex path constraints \CTNLconvexpath are easy to write. Although there
are no convex state constraints, there are convex final time bounds
\eqref{eq:ex_quad_tf_bounds}. These can be included as a convex state path
constraint \optieqref{scp_gen_cont}{convex_path_constraints_X}, which is mixed in
the state and parameter. Using \eqref{eq:ex_quad_time_dilation}, we define the
convex state path constraints as:
\begin{equation}
  \label{eq:ex_quad_state_path}
  \set{X}=\brac[big]{(x,p)\in\reals^6\times\reals\where
    \tf[,\min] \le p \le \tf[,\max]}.
\end{equation}

On the other hand, the convex input constraint set $\set{U}$ is given by all
the input vectors that satisfy \eqref{eq:ex_quad_cc_cvx}.

The nonconvex path constraints \optieqref{scp_gen_cont}{nonconvex_constraints} are
given by the vector function $s:\reals^3\to\reals^{\Nobs}$, whose elements
encode the obstacle avoidance constraints:
\begin{equation}
  \label{eq:ex_quad_s_obs_avoid}
  s_j\pare{r} = 1-\norm[2]{H_j ( r - c_j )},\quad j=1,\dots,\Nobs.
\end{equation}

We will briefly mention how to evaluate the Jacobian \eqref{eq:scvx_lin_mats_h} for
\eqref{eq:ex_quad_s_obs_avoid}. Suppose that a reference position trajectory
$\{\bar r(t)\}_0^1$ is available, and consider the $j$-th obstacle constraint
in \eqref{eq:ex_quad_s_obs_avoid}. The following gradient then allows to evaluate
\eqref{eq:scvx_lin_mats_h}:
\begin{equation}
  \grad s_j(r) = -\frac{H_j^\transp H_j \big( r - c_j \big)}%
  {\norm[2,big]{H_j \big( r - c_j \big)}}.
\end{equation}

The boundary conditions \optieqref{scp_gen_cont}{initial_conditions} and
\optieqref{scp_gen_cont}{final_conditions} are obtained from \eqref{eq:ex_quad_bcs}:
\begin{subequations}
  \label{eq:ex_quad_gic_gtc}
  \begin{align}
    \gic\pare[big]{x(0),p} &= \Matrix{r(0)-r_0 \\ \dot{r}(0)-v_0}, \\
    \gtf\pare[big]{x(1),p} &= \Matrix{r(1)-r_f \\ \dot{r}(1)-v_f}.
  \end{align}
\end{subequations}

Lastly, we have to convert the cost \optiobjref{ex_quad_ocp} into the Bolza form
\eqref{eq:ocost_nlin}. There is no terminal cost, hence $\term\equiv 0$. On the
other hand, the direct transcription of the running cost would be
$\runn(x, u, p) = p\sigma^2$. However, we will simplify this by omitting the
time dilation parameter. This simplification makes the problem easier by
removing nonconvexity from the running cost, and numerical results still show
good resulting trajectories. Furthermore, since \scvx augments the cost with
penalty terms in \eqref{eq:scvx_Lpen}, we will normalize the running cost by its
nominal value in order to make it roughly unity for a ``mild'' trajectory. This
greatly facilitates the selection of the penalty weight $\lambda$. By taking
the nominal value of $\sigma$ as the hover condition, we define the following
running cost:
\begin{equation}
  \label{eq:ex_quad_running_cost}
  \runn(x, u, p) = \pare[biggg]{\frac{\sigma}{g}}^2.
\end{equation}

We now have a complete definition of \pref{scp_gen_cont} for the quadrotor obstacle
avoidance problem. The only remaining task is to choose the \scvx algorithm
parameters listed in \tabref{algo_params}. The rest of the solution process is
completely automated by the general \scvx algorithm description in Part II.

\subsubsection{\gusto Formulation}

The \gusto algorithm can also be used to solve \pref{ex_quad_ocp}. The formulation
is almost identical to \scvx, which underscores the fact that the two
algorithms can be used interchangeably to solve many of the same problems. The
quadratic running cost \eqref{eq:gusto_running_cost} encodes
\eqref{eq:ex_quad_running_cost} as follows:
\begin{subequations}
  \label{eq:ex_quad_gusto_running_cost}
  \begin{align}
    \Jq(p) &= \diag\pare[big]{0,0,0,g^{-2}}, \\
    \Jl (x,p) &= 0, \\
    \Jc(x,p) &= 0.
  \end{align}
\end{subequations}

We can also cast the dynamics \eqref{eq:ex_quad_dynamics_abs_time} into the control
affine form \eqref{eq:control_affine_gusto}:
\begin{subequations}
  \label{eq:ex_quad_gusto_dynamics}
  \begin{align}
    f_0(x) &= \Matrix{\dot r \\ -g\hat n}, \\
    f_i(x) &= \Matrix{0 \\ e_i},\quad i=1,2,3, \\
    f_4(x) &= 0,
  \end{align}
\end{subequations}
where $e_i\in\reals^3$ is the $i$-th standard basis vector. The reader may be
surprised that these are the only changes required to convert the \scvx
formulation from the previous section into a form that can be ingested by
\gusto. The only remaining task is to choose the \gusto algorithm parameters
listed in \tabref{algo_params}. Just like for \scvx, the rest of the solution
process is completely automated by the general \gusto algorithm description in
Part II.

\subsubsection{Initial Trajectory Guess}

The initial state reference trajectory is obtained by a simple straight-line
interpolation as provided by \eqref{eq:state_initial_guess}:
\begin{equation}
  \label{eq:ex_quad_initial_state}
  x(t)=(1-t)\Matrix{r_0\\ v_0}+t\Matrix{r_f\\ v_f},%
  ~\textnormal{for}~t\in [0,1].
\end{equation}

The initial parameter vector, which is just the time dilation, is chosen to be
in the middle of the allowed trajectory durations:
\begin{equation}
  \label{eq:ex_quad_initial_parameter}
  p = \frac{\tf[,\min] + \tf[,\max]}{2}.
\end{equation}

The initial input trajectory is guessed based on insight about the quadrotor
problem. Ignoring any particular trajectory task, we know that the quadrotor
will generally have to support its own weight under the influence of
gravity. Thus, we choose a constant initial input guess that would make a
static quadrotor hover:
\begin{equation}
  \label{eq:ex_quad_initial_input}
  a(t)=g\hat n,~\sigma(t)=g,~\textnormal{for}~t\in [0,1].
\end{equation}

This initial guess is infeasible with respect to both the dynamics and the
obstacle constraints, and it is extremely cheap to compute. The fact that it
works well in practice highlights two facts: SCP methods are relatively easy to
initialize, and they readily accept infeasible initial guesses. Note that in
the particular case of this problem, an initial trajectory could also be
computed using a convex version of the problem obtained by removing the
obstacle constraints \eqref{eq:ex_quad_sc_obs} and applying the \lcvx relaxation
\eqref{eq:ex_quad_cc_cvx}.

\subsubsection{Numerical Results}

\begin{csmtable}[
  caption={Breakdown of subproblem size for the quadrotor obstacle avoidance
    example.},
  label={ex_quad_subproblem_size},
  position={!t}]%
  \centering%
  \ifmaketwocolcsm
  \includegraphics[width=0.7\columnwidth]{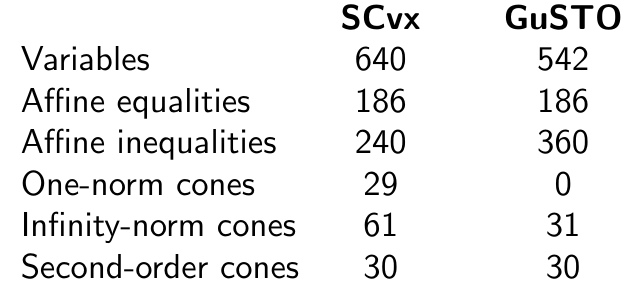}
  \else
  \includegraphics[scale=\csmpreprintfigscale]{scp_quad_sizes}
  \fi
\end{csmtable}

We now have a specialized instance of \pref{scp_gen_cont} and an initialization
strategy for the quadrotor obstacle avoidance problem. The solution is obtained
via SCP according to the general algorithm descriptions for \scvx and
\gusto. For temporal discretization, we use the FOH interpolating polynomial
method from \sbref{discretization}. The algorithm parameters are provided in
\tabref{ex_quad_params}, and the full implementation is available in the code
repository linked in \figref{github_qr}. ECOS is used as the numerical convex
optimizer at \alglocation{\solveloc} in \figref{scp_loop}
\cite{domahidi2013ecos}. The timing results correspond to a Dell XPS 13 9260
laptop powered by an Intel Core i5-7200U CPU clocked at
2.5~\si{\giga\hertz}. The computer has 8~\si{\gibi\byte} LPDD3 RAM and
128~\si{\kibi\byte} L1, 512~\si{\kibi\byte} L2, and 3~\si{\mebi\byte} L3 cache.

\begin{csmfigure}[%
  caption={%
    The position trajectory evolution (left) and the final converged trajectory
    (right) for the quadrotor obstacle avoidance problem. The continuous\dash
    time trajectory in the right plots is obtained by numerically integrating
    the dynamics \eqref{eq:ex_ff_dynamics}. The fact that this trajectory passes
    through the discrete\dash time subproblem solution confirms dynamic
    feasibility. The {\color{beamerRed}red} lines in the right plots show the
    acceleration vector as seen from above.
  },
  label={ex_quad_pos},
  position={!t}]%
  \centering
  \begin{subfigure}[t]{\columnwidth}
    \centering
    \includegraphics[width=\columnwidth,page=1]{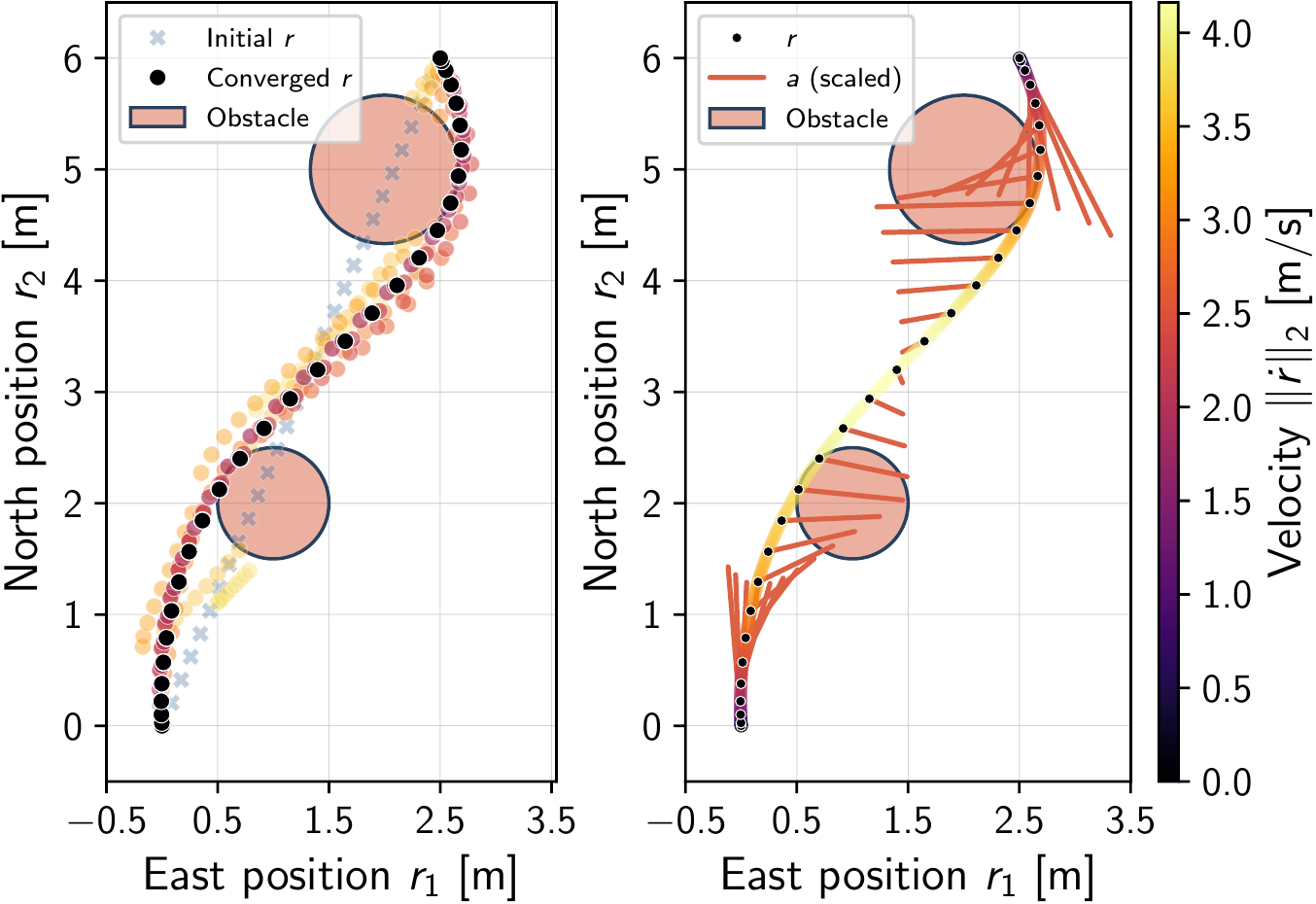}%
    \caption{\scvx.}
    \label{fig:ex_quad_pos_scvx}
  \end{subfigure}%
  \vspace{2mm}%

  \begin{subfigure}[t]{\columnwidth}
    \centering
    \includegraphics[width=\columnwidth,page=2]{scp_quad_pos}%
    \caption{\gusto.}
    \label{fig:ex_quad_pos_gusto}
  \end{subfigure}%
\end{csmfigure}

\begin{csmfigure}[%
  caption={%
    Convergence and runtime performance for the quadrotor obstacle avoidance
    problem. Both algorithms take a similar amount of time and number of
    iterations to converge to a given tolerance. The runtime subplots in the
    bottom row show statistics on algorithm performance over 50 executions. The
    solution times per subproblem are roughly equal, which shows that the
    subproblem difficulty stays constant over the
    iterations. {\color{beamerDarkBlue}Formulate} measures the time taken to
    parse the subproblem into the input format of the convex optimizer;
    {\color{beamerYellow}Discretize} measures the time taken to temporally
    discretize the linearized dynamics \optieqref{scp_gen_cvx}{dynamics}; and
    {\color{beamerRed}Solve} measures the time taken by the core convex
    numerical optimizer. Each bar shows the median time. The
    {\color{beamerBlue}blue} trace across the diagonal shows the cumulative
    time obtained by summing the runtimes of all preceding iterations. Its
    markers are placed at the median time, and the error bars show the $10\%$
    (bottom) and $90\%$ (top) quantiles. Because small runtime differences
    accumulate over iterations, the error bars grow with iteration count.
  },
  label={ex_quad_convergence},
  columns=2,
  position={!t}]%
  \centering
  \ifmaketwocolcsm
  \def\subfigwidth{0.45\textwidth}
  \else
  \def\subfigwidth{0.48\textwidth}
  \fi
  \def\figwidth{\columnwidth}
  \begin{subfigure}[t]{\subfigwidth}
    \centering
    \includegraphics[width=\figwidth]{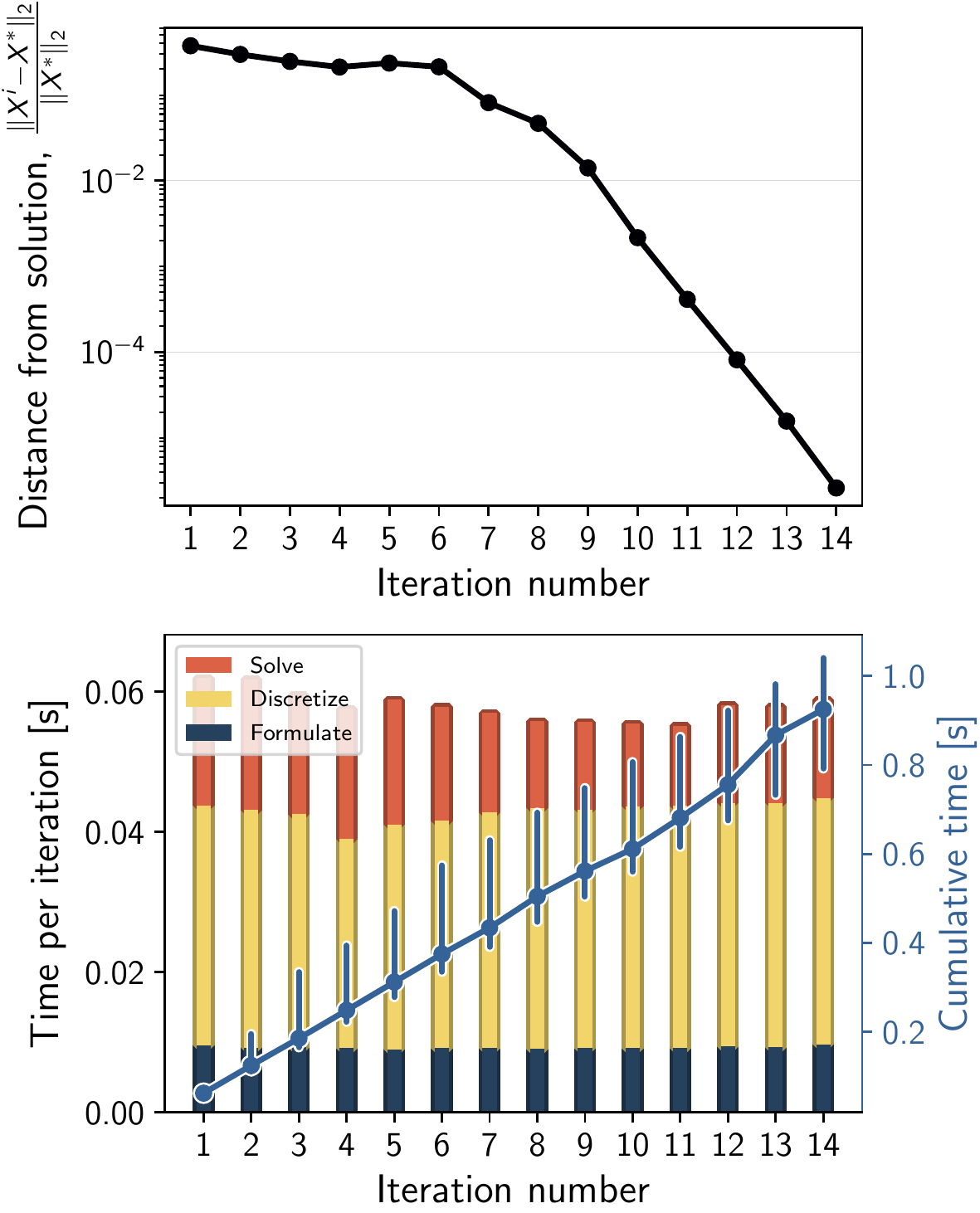}%
    \caption{\scvx.}
    \label{fig:ex_quad_convergence_scvx}
  \end{subfigure}%
  \ifmaketwocolcsm
  \hspace{1cm}%
  \else
  \hfill%
  \fi
  \begin{subfigure}[t]{\subfigwidth}
    \centering
    \includegraphics[width=\figwidth]{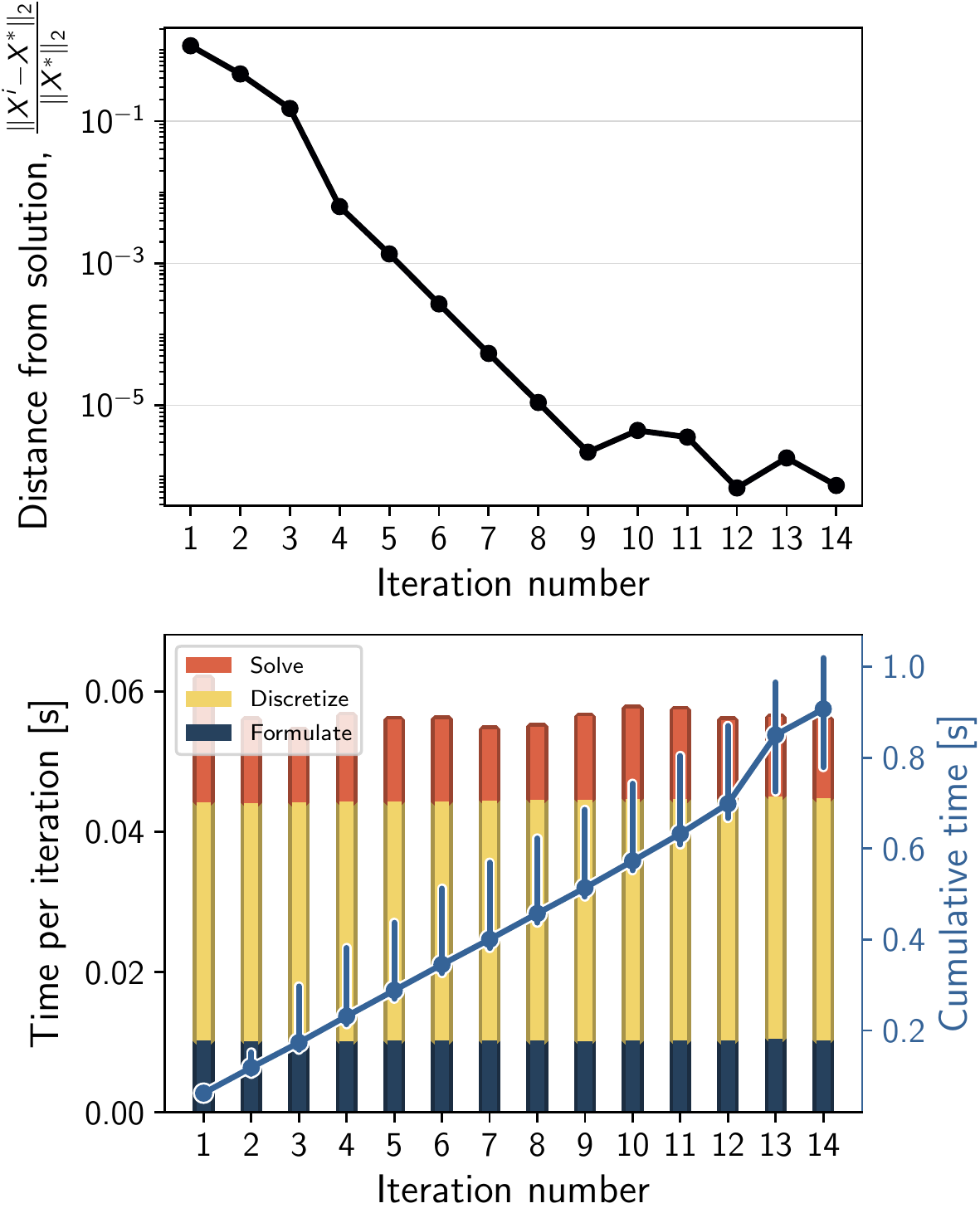}%
    \caption{\gusto.}
    \label{fig:ex_quad_convergence_gusto}
  \end{subfigure}%
\end{csmfigure}

The convergence process for both algorithms is illustrated in
\figref{ex_quad_convergence}. We have set the convergence tolerances
$\varepsilon=\varepsilon_{\mrm{r}}=0$ and we terminate both algorithms after 15
iterations. At each iteration, the algorithms have to solve subproblems of
roughly equal sizes, as shown in \tabref{ex_quad_subproblem_size}. Differences in
the sizes arise from the slightly different subproblem formulations of each
algorithm. Among the primary contributors are how artificial infeasibility and
unboundedness are treated, and differences in the cost penalty terms. Notably,
\gusto has no one-norm cones because it does not have a dynamics virtual
control term like \scvx (compare \optieqref{subproblem_scvx_dt}{dynamics} and
\optieqref{subproblem_gusto_dt}{dynamics}). Despite their differences,
\figref{ex_quad_convergence} shows that \gusto and \scvx are equally fast.

The trajectory solutions for \scvx and \gusto are shown in
\figref{ex_quad_pos,ex_quad_timeseries}. Recall that this is a free final time
problem, and both algorithms are able to increase the initial guess
\eqref{eq:ex_quad_initial_parameter} until the maximum allowed flight time
$\tf[,\max]$. Note that this is optimal, since a quadrotor minimizing control
energy \optiobjref{ex_quad_ocp} will opt for a slow trajectory with the lowest
acceleration.

The \scvx and \gusto solutions are practically identical. The fact that they
were produced by two different algorithms can only be spotted from the
different convergence histories on the left side in \figref{ex_quad_pos}. It is not
always intuitive how the initial guess morphs into the final trajectory. It is
remarkable that the infeasible trajectories of the early iterations morph into
a smooth and feasible trajectory. Yet, this is guaranteed by the SCP
convergence theory in Part II. For this and many other trajectory problems, we
have observed time and again how SCP is adept at morphing rough initial guesses
into fine-tuned feasible and locally optimal trajectories.

Finally, we will make a minor note of that temporal discretization results in
some clipping of the obstacle keep-out zones in \figref{ex_quad_pos}. This is the
direct result of imposing constraints only at the discrete\dash time
nodes. Various strategies exist to mitigate the clipping effect, such as
increasing the radius of the keep-out zones, increasing the number of
discretization points, imposing a sufficiently low maximum velocity constraint,
or numerically minimizing a function related to the state transition matrix
\cite{acikmese2008enhancements,Dueri2017Clipping}.

\begin{csmfigure}[%
  caption={%
    The acceleration norm and tilt angle time histories for the converged
    trajectory of the quadrotor obstacle avoidance problem. These are visually
    identical for \scvx and \gusto, so we only show a single plot. The
    continuous\dash time acceleration norm is obtained from the FOH assumption
    \eqref{eq:sidebar_control_foh}, while the continuous\dash time tilt angle is
    obtained by integrating the solution through the nonlinear
    dynamics. Similar to \figref{lcvx_pdg_results}, the acceleration time history
    plot confirms that lossless convexification holds (i.e., the constraint
    \eqref{eq:ex_quad_cc_cvx_lcvx_equality} holds with equality).
    %
  },
  label={ex_quad_timeseries},
  position={!t}]%
  \centering
  \ifmaketwocolcsm
  \def\figwidth{0.9\columnwidth}
  \else
  \def\figwidth{0.8\columnwidth}
  \fi
  \includegraphics[width=\figwidth,page=1]{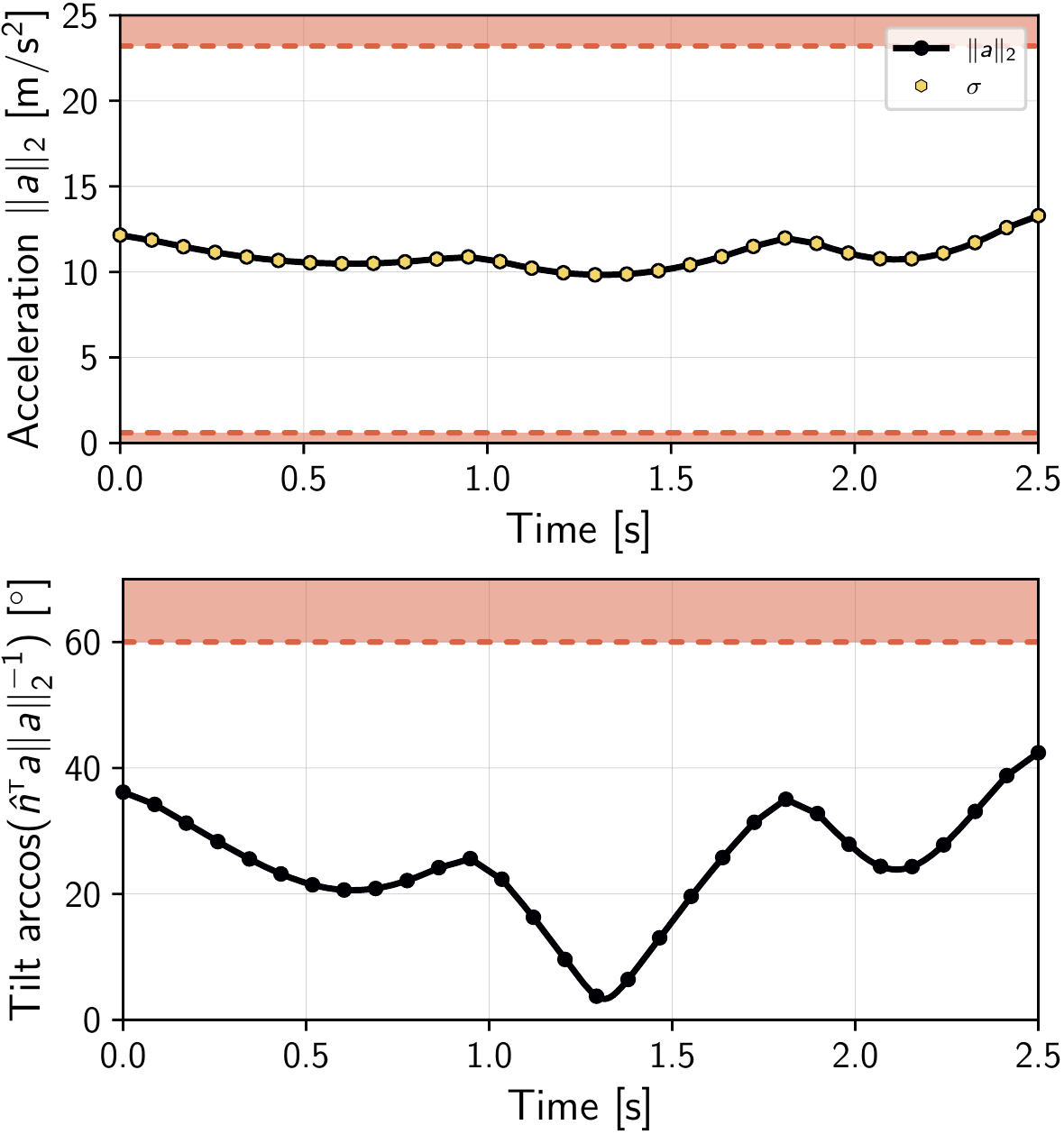}%
\end{csmfigure}

\subsection{SCP: 6-DoF Free-flyer}

\begin{csmfigure}[%
  caption={%
    Two examples of free-flyer robots at the International Space Station, the
    JAXA Int-Ball (left) and the Naval Postgraduate School/NASA Astrobee
    (right) \cite{AstrobeeImage,IntBallImage}. These robots provide a helping
    hand in space station maintenance tasks.
    %
  },
  label={ff_examples}]%
  \centering%
  \includegraphics[width=0.9\columnwidth]{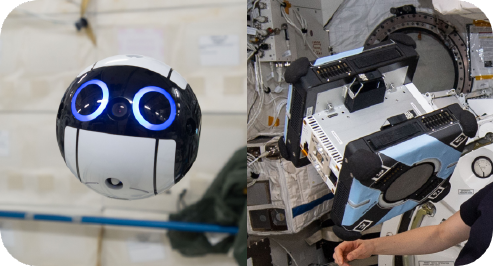}
\end{csmfigure}

Having demonstrated the use of sequential convex programming on a relatively
simple quadrotor trajectory, we now present a substantially more challenging
example involving nonlinear 6-DoF dynamics and a more complex set of obstacle
avoidance constraints.

The objective is to compute a trajectory for a 6-DoF free-flying robotic
vehicle that must navigate through an environment akin to the International
Space Station (ISS). Free-flyers are robots that provide assistance to human
operators in micro\dash gravity environments~\cite{Estrada2016,Mote2020}. As
shown in \figref{ff_examples}, Astrobee and the JAXA Internal Ball Camera
(Int-Ball) are two recent successful deployments of such robots. Their goals
include filming the ISS and assisting with maintenance
tasks~\cite{Smith2016}. The particulars of this SCP example are taken primarily
from~\cite{Bonalli2019a}.

The quadrotor in the previous section was modeled as a point mass whose
attitude is approximated by the direction of the acceleration vector. The
free-flyer, on the other hand, is a more complex vehicle that must generally
perform coupled translational and rotational motion using a multi-thruster
assembly. Maneuvers may require to point a camera at a fixed target, or to
emulate nonholonomic behavior for the sake of predictability and operator
comfort \cite{Barlow2016Astrobee,Szafir2015}. This calls for whole-body motion
planning, for which we model the free-flyer as a full 6-DoF rigid body with
both translational and rotational dynamics.

To describe the equations of motion, we need to introduce two reference
frames. First, let $\Finertial$ be an inertial reference frame with a
conveniently positioned, but otherwise arbitrary, origin. Second, let $\Fbody$
be a rotating reference frame affixed to the robot's center of mass, and whose
unit vectors are aligned with the robot's principal axes of inertia. We call
$\Finertial$ the inertial frame and $\Fbody$ the body frame. Correspondingly,
vectors expressed in $\Finertial$ are inertial vectors and carry an $\inertial$
subscript (e.g., $x_{\inertial}$), while those expressed in $\Fbody$ are body
vectors and carry a $\body$ subscript (e.g., $x_{\body}$). For the purpose of
trajectory generation, we encode the orientation of $\Fbody$ with respect to
$\Finertial$ using a vector representation of a unit quaternion,
$\qIB\in\real^4$.

Our convention is to represent the translational dynamics in $\Finertial$ and
the attitude dynamics in $\Fbody$.  This yields the following Newton\dash Euler
equations that govern the free-flyer's motion \cite{SzmukReynolds2018}:
\begin{subequations}
  \label{eq:ex_ff_dynamics}
  \begin{align}%
    \drI(\tabs)
    &= \vI(\tabs), \\
    \dvI(\tabs)
    &= m\inv T_{\inertial}(\tabs), \\
    \label{eq:ex_ff_dynamics_quaternion}
    \dqIB(\tabs)
    &= \frac{1}{2} \qIB(\tabs) \otimes \wB(\tabs), \\
    \dwB(\tabs)
    &= J^{-1} \pare[big]{ M_{\body}(\tabs) -
      \wB(\tabs)\skew J \wB(\tabs) }.
  \end{align}
\end{subequations}

The state variables in the above equations are the inertial position
$\rI\in\reals^3$, the inertial velocity $\vI\in\reals^3$, the aforementioned
unit quaternion attitude $\qIB\in\reals^4$, and the body angular velocity
$\wB\in\reals^3$. The latter variable represents the rate at which $\Fbody$
rotates with respect to $\Finertial$. The free-flyer's motion is controlled by
an inertial thrust vector $T_{\inertial}\in\real^3$ and a body torque vector
$M_{\body}\in\real^3$. The physical parameters of the free-flyer, the mass
$m>0$ and the principal moment of inertia matrix $J\in\reals^{3\times 3}$, are
fixed. The dynamics are written in absolute time $\tabs$ that spans the
interval $[0,\tf]$. We allow the final time $\tf$ to be optimized, and bound it
using the previous constraint \eqref{eq:ex_quad_tf_bounds}.

The initial and final conditions for each state are specified in this example
to be fixed constants:
\begin{subequations}%
  \label{eq:ex_ff_bcs}
  \begin{alignat}{4}
    \label{eq:ex_ff_bcs_pos}
    \rI(0)  &= \ric, \quad &\rI(\tf)  &= \rfc, \\
    \vI(0)  &= \vic, \quad &\vI(\tf)  &= \vfc, \\
    \label{eq:ex_ff_bcs_quat}
    \qIB(0) &= \qic, \quad &\qIB(\tf) &= \qfc, \\
    \wB(0)  &= 0,    \quad &\wB(\tf)  &= 0.
  \end{alignat}
\end{subequations}

The free-flyer robot implements a 6-DoF holonomic actuation system based on a
centrifugal impeller that pressurizes air, which can then be vented from a set
of nozzles distributed around the body \cite{Bualat2015,JAXAIntBall}. Holonomic
actuation means that the thrust and torque vectors are independent from the
vehicle's attitude \cite{Roque2016}. The capability of this system can be
modeled by the following control input constraints:
\begin{equation}
  \label{eq:ex_ff_cc}
  \norm[2]{T_{\inertial}(\tabs)} \leq T_{\max}, \quad
  \norm[2]{M_{\body}(\tabs)} \leq M_{\max},
\end{equation}
where $T_{\max} > 0$ and $M_{\max} > 0$ are user-defined constants representing
the maximum thrust and torque.

This problem involves both convex and nonconvex state constraints. Convex
constraints are used to bound the velocity and angular velocity magnitudes to
user\dash defined constants:
\begin{equation}
  \label{eq:ex_ff_cvx_sc}
  \norm[2]{\vI(\tabs)} \leq v_{\max}, \quad
  \norm[2]{\wB(\tabs)} \leq \omega_{\max}.
\end{equation}

Nonconvex state constraints are used to model the (fictional) ISS flight space
and to avoid floating obstacles. The latter are modeled exactly as in the
quadrotor example using the constraints in \eqref{eq:ex_quad_sc_obs}. The flight
space, on the other hand, is represented by a union of rectangular rooms. This
is a difficult nonconvex constraint and its efficient modeling requires some
explanation.

\begin{csmfigure}[%
  caption={%
    A heatmap visualization of the exact SDF \eqref{eq:sdf_iss} and the approximate
    SDF \eqref{eq:sdf_iss_approx} for several values of the sharpness parameter
    $\sigma$. Each plot also shows the SDF zero-level set boundary as a dashed
    line. This boundary encloses the feasible flight space, which corresponds
    to nonnegative values of the SDF.
  },%
  label={sdf_illustration},%
  columns=2]%
  \centering%
  \ifmaketwocolcsm
  \includegraphics[width=0.9\textwidth]{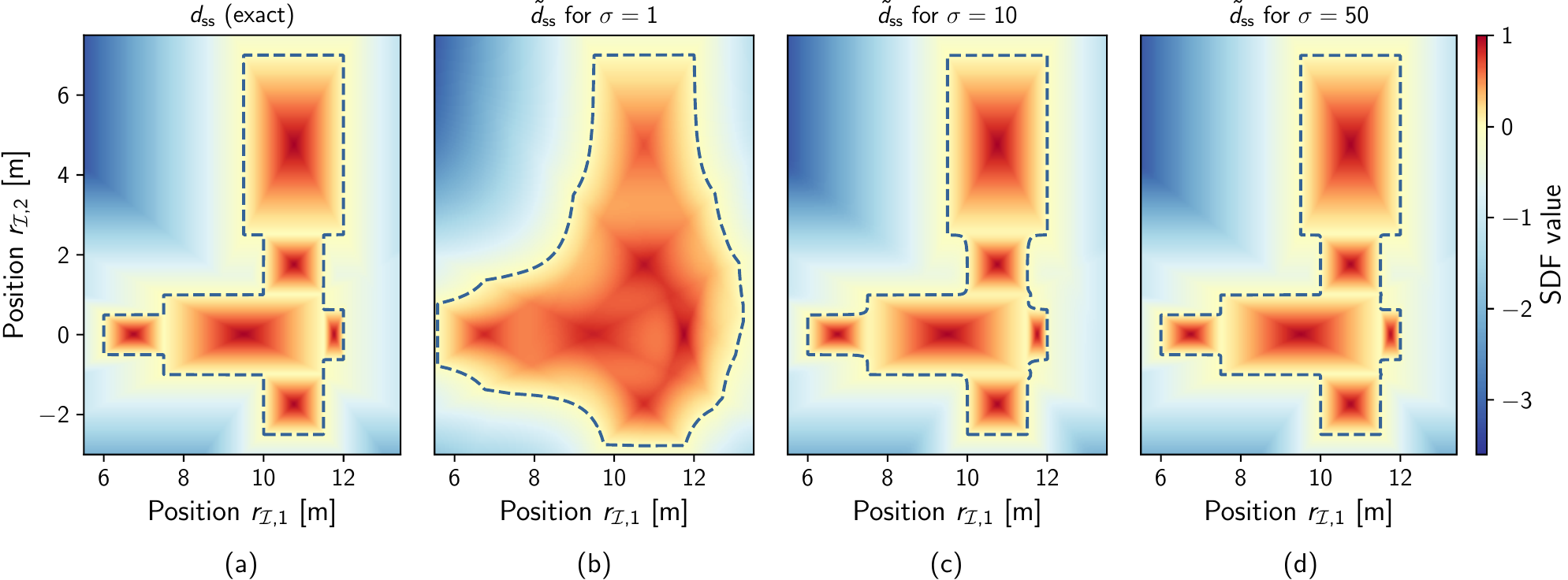}
  \else
  \includegraphics[width=\textwidth]{sdf_illustration}
  \fi
\end{csmfigure}

At the conceptual level, the space station flight space is represented by a
function $\diss:\reals^3\to\reals$ that maps inertial position to a scalar
number. This is commonly referred to as a signed distance function (SDF)
\cite{AkenineMoller2018}. Let us denote the set of positions that are within the
flight space by $\ObsISS\subset\reals^3$. A valid SDF is given by any function
that satisfies the following property:
\begin{equation}
  \label{eq:sdf_property}
  \rI\in\ObsISS~\iff~\diss(\rI)\ge 0.
\end{equation}

If we can formulate a continuously differentiable $\diss$, then we can model
the flight space using \eqref{eq:sdf_property} as the following nonconvex path
constraint \optieqref{scp_gen_cont}{nonconvex_constraints}:
\begin{equation}
  \label{eq:sdf_constraint_exact}
  \diss(\rI)\ge 0.
\end{equation}

Open-source libraries such as Bullet~\cite{Coumans2020} are available to compute
the SDF for obstacles of arbitrary shape. In this work, we will use a simpler
custom implementation. To begin, let us model the space station as an assembly
of several rooms. This is expressed as a set union:
\begin{equation}
  \label{eq:iss_volume}
  \ObsISS \definedas \Union_{i=1}^{\niss}\Obsi,
\end{equation}
where each room $\Obsi$ is taken to be a rectangular box:
\begin{equation}
  \label{eq:iss_room}
  \Obsi \definedas \brac[big]{
    \rI\in\reals^3 : l_i^{\ISS} \le \rI \le u_i^{\ISS}}.
\end{equation}

The coordinates $l_i^{\ISS}$ and $u_i^{\ISS}$ represent the ``rear bottom
right'' and the ``ahead top left'' corners of the $i$-th room, when looking
along the positive axis directions. Equivalently, but more advantageously for
implementing the SDF, the set \eqref{eq:iss_room} can be represented as:
\begin{equation}
  \label{eq:iss_room_inf_norm}
  \Obsi
  \definedas \brac[biggg]{
    \rI\in\reals^3 : \norm[\infty,biggg]{%
    \frac{\rI-c_i^{\ISS}}{s_i^{\ISS}}}\le 1
    },
\end{equation}
where vector division is entrywise, and the new terms $c_i^{\ISS}$ and
$s_i^{\ISS}$ are the room's centroid and diagonal:
\begin{equation}
  c_i^{\ISS} = \frac{u_i^{\ISS}+l_i^{\ISS}}{2},\quad%
  s_i^{\ISS} = \frac{u_i^{\ISS}-l_i^{\ISS}}{2}.
\end{equation}

To find the SDF for the overall flight space, we begin with the simpler task of
writing an SDF for a single room. This is straightforward using
\eqref{eq:iss_room_inf_norm}:
\begin{equation}
  \label{eq:sdf_i}
  \diss[i](\rI) \definedas
  1-\norm[\infty,biggg]{\frac{\rI-c_i^{\ISS}}{s_i^{\ISS}}},
\end{equation}
which is a concave function that satisfies a similar property to
\eqref{eq:sdf_property}:
\begin{equation}
  \label{eq:sdf_i_property}
  \rI\in\Obsi~\iff~\diss[i](\rI)\ge 0.
\end{equation}

Because $\diss[i]$ is concave, the constraint on the right side of
\eqref{eq:sdf_i_property} is convex. This means that constraining the robot to be
inside room $\Obsi$ is a convex operation, which makes sense since $\Obsi$ is a
convex set.

As the free-flyer traverses the flight space, one can imagine the room SDFs to
evolve based on the robot's position. When the robot enters the $i$-th room,
$\diss[i]$ becomes positive and grows up to a maximum value of one as the robot
approaches the room center. As the robot exits the room, $\diss[i]$ becomes
negative and decreases in value as the robot flies further away. To keep the
robot inside the space station flight space, the room SDFs have to evolve such
that there is always at least one nonnegative room SDF. This requirement
precisely describes the overall SDF, which can be encoded mathematically as a
maximization:
\begin{equation}
  \label{eq:sdf_iss}
  \diss(\rI) \definedas \max_{i=1,\dots,\niss} \diss[i](\rI).
\end{equation}

We now have an SDF definition which satisfies the required property
\eqref{eq:sdf_property}. \fakesubfigref{sdf_illustration}{a} shows an example scalar field
generated by \eqref{eq:sdf_iss} for a typical space station layout. Visually, when
restricted to a plane with a fixed altitude $\rI[3]$, the individual room SDFs
form four\dash sided ``pyramids'' above their corresponding room.

The room SDFs $\diss[i]$ are concave, however the maximization in
\eqref{eq:sdf_iss} generates a nonconvex function. Two possible strategies to
encode \eqref{eq:sdf_iss} are by introducing integer variables or by a smooth
approximation \cite{Blackmore2012}. The former strategy generates a
mixed-integer convex subproblem, which is possible to solve but does not fit
the SCP algorithm mold of this article. As mentioned before for
\pref{scp_gen_cont}, our subproblems do not involve integer variables. We thus
pursue the smooth approximation strategy, which yields an arbitrarily accurate
approximation of the feasible flight space $\ObsISS$ and carries the
significant computational benefit of avoiding mixed-integer programming.

The crux of our strategy is to replace the maximization in \eqref{eq:sdf_iss} with
the softmax function. Given a general vector $v\in\reals^n$, this function is
defined by:
\begin{equation}
  \label{eq:softmax}
  \softmax[\sigma](v) = \sigma\inv \log\sum_{i=1}^n \exp{\sigma v_i},
\end{equation}
where $\sigma>0$ is a sharpness parameter such that $\softmax[\sigma]$ upper
bounds the exact $\max$ with an additive error of at most $\log(n)/\sigma$. To
develop intuition, consider the SDF \eqref{eq:sdf_iss} for two adjacent rooms and
restricted along the $j$-th axis of the inertial frame. We can then write the
SDF simply as:
\begin{equation}
  \label{eq:sdf_iss_restricted}
  \diss(\rI[j]) = \max\brac[Big]{
    1-\abso[Big]{\frac{\rI[j]-c_{1j}^{\ISS}}{s_{1j}^{\ISS}}},
    1-\abso[Big]{\frac{\rI[j]-c_{2j}^{\ISS}}{s_{2j}^{\ISS}}}
  }.
\end{equation}

\begin{csmfigure}[%
  caption={%
    Illustration of the correspondence between between the exact SDF
    \eqref{eq:sdf_iss} and its approximation $\sdiss$ using the softmax function
    \eqref{eq:softmax}. As the sharpness parameter $\sigma$ increases, the
    approximation quickly converges to the exact SDF. This figure illustrates a
    sweep for $\sigma\in [1, 5]$, where lower values are associated with darker
    colors.
  },%
  label={softmax_1d_intuition},%
  position={!b}]%
  \centering%
  \ifmaketwocolcsm
  \includegraphics[width=0.9\columnwidth]{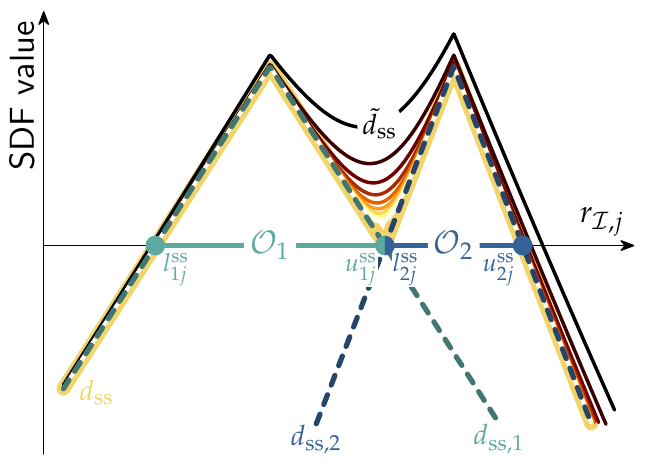}
  \else
  \includegraphics[width=0.8\columnwidth]{softmax_1d_intuition}
  \fi
\end{csmfigure}

\figref{softmax_1d_intuition} illustrates the relationship between the exact SDF
\eqref{eq:sdf_iss_restricted} and its approximation, which is obtained by replacing
the $\max$ operator with $\softmax[\sigma]$. We readily observe that $\diss$ is
indeed highly nonconvex, and that the approximation quickly converges to the
exact SDF as $\sigma$ increases. We now generalize this one-dimensional example
and replace \eqref{eq:sdf_iss} with the following approximation:
\begin{subequations}
  \label{eq:sdf_iss_approx}
  \begin{align}
    \label{eq:sdf_iss_smooth}
    \sdiss(\rI) &\definedas \softmax[\sigma]\pare[big]{\disslb(\rI)}, \\
    \label{eq:sdf_room_convex}
    \disslb[i](\rI) &\le \diss[i](\rI),\quad i=1,\dots,\niss,
  \end{align}
\end{subequations}
where the new functions $\disslb[i]$ as so-called ``slack'' room SDFs. The
model \eqref{eq:sdf_iss_approx} admits several favorable properties. First,
\eqref{eq:sdf_iss_smooth} is smooth in the new slack SDFs and can be included
directly in \optieqref{scp_gen_cont}{nonconvex_constraints} as the following
nonconvex path constraint:
\begin{equation}
  \label{eq:ex_ff_iss}
  \sdiss(\rI)\ge 0.
\end{equation}

Second, the constraints in \eqref{eq:sdf_room_convex} are convex and can be
included directly in \optieqref{scp_gen_cont}{convex_path_constraints_X}. Overall,
the approximate SDF \eqref{eq:sdf_iss_smooth} satisfies the following property:
\begin{equation}
  \label{eq:sdf_iss_approx_property}
  \rI\in\sObsISS~\iff~\exists\disslb(\rI)~\textnormal{such that
  \eqref{eq:ex_ff_iss} holds}.
\end{equation}
where $\sObsISS\subset\reals^3$ is an approximation of $\ObsISS$ that becomes
arbitrarily more accurate as the sharpness parameter $\sigma$ increases. The
geometry of this convergence process is illustrated in \fakesubfigref[c]{sdf_illustration}{b}
for a typical space station layout. Crucially, the fact that $\softmax[\sigma]$
is an upper bound of the $\max$ means that the approximate SDF $\sdiss$ is
nonnegative at the interfaces of adjacent rooms. In other words, the passage
between adjacent rooms is not artificially blocked by our approximation.

Summarizing the above discussion, we wish to solve the following free final
time optimal control problem that minimizes control energy:
\begin{optimization}[
  label={ex_ff_ocp},
  variables={T_{\inertial},M_{\body},\tf},
  objective={\int_{0}^{\tf} \norm[2]{T_{\inertial}(\tabs)}^2 +%
    \norm[2]{M_{\body}(\tabs)}^2\sdd\tabs}]
  & 
  \textnormal{\eqref{eq:ex_ff_dynamics}-%
    \eqref{eq:ex_ff_cvx_sc},~\eqref{eq:ex_ff_iss},~%
    \eqref{eq:ex_quad_tf_bounds},~%
    \eqref{eq:ex_quad_sc_obs}.}
\end{optimization}

\subsubsection{\scvx Formulation}

\begin{csmtable}[%
  caption={%
    Algorithm parameters for the 6-DoF free\dash flyer example.
  },%
  label={ex_ff_params},
  position={!t}]%
  \centering%
  \ifmaketwocolcsm
  \includegraphics[width=\textwidth]{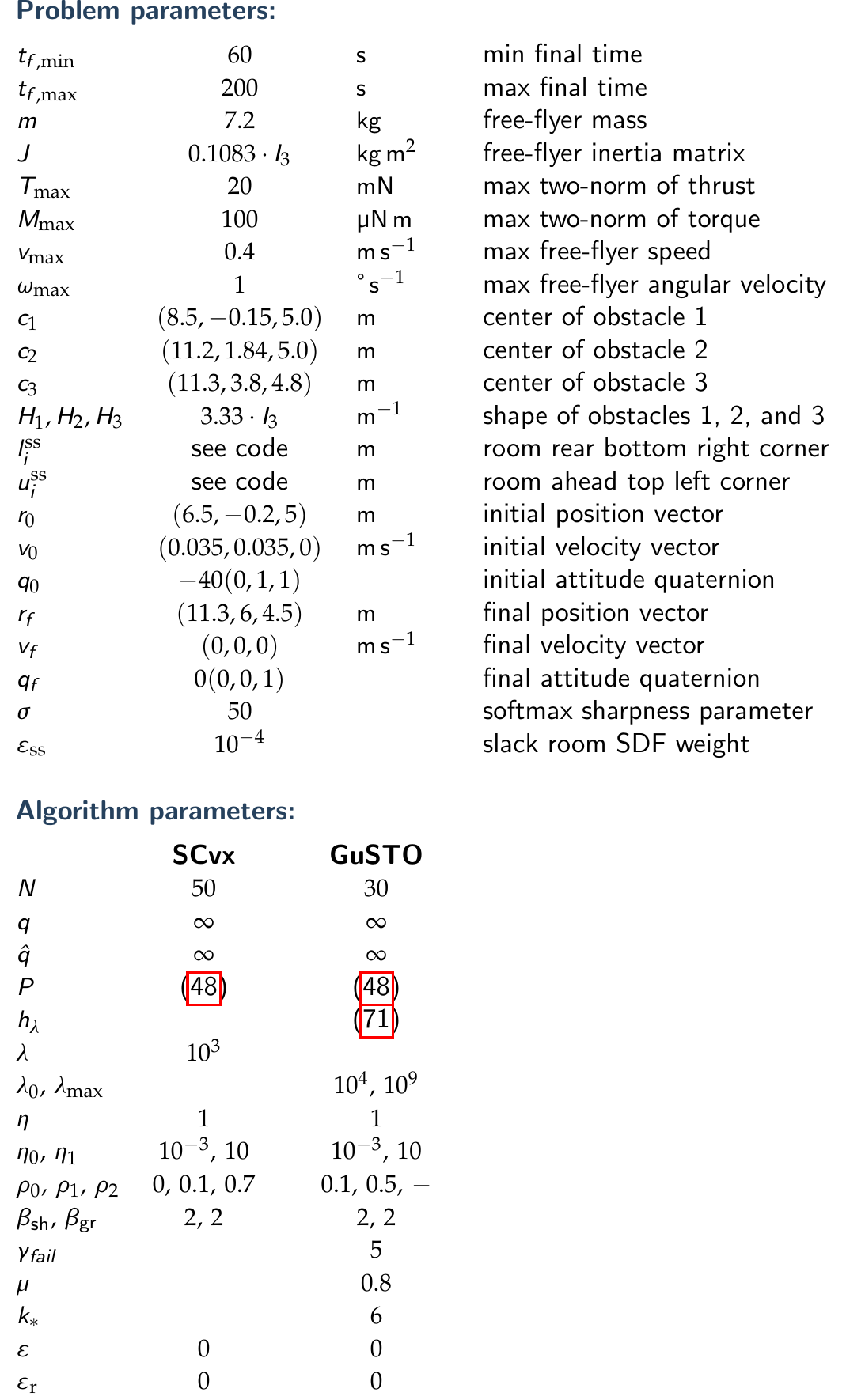}
  \else
  \includegraphics[scale=\csmpreprintfigscale]{tikz_ex_ff_parameters}
  \fi
\end{csmtable}

Like for the quadrotor example, we begin by demonstrating how to cast
\pref{ex_ff_ocp} into the standard template of \pref{scp_gen_cont}. While this process
is mostly similar to that of the quadrotor, the particularities of the flight
space constraint \eqref{eq:ex_ff_iss} will reveal a salient feature of efficient
modeling for SCP. Once the modeling step is done and an initial guess
trajectory is defined, the solution process is completely automated by the
general \scvx algorithm description in Part II. To keep the notation light, we
omit the argument of time where possible.

Looking at the dynamics \eqref{eq:ex_ff_dynamics}, we define the following state
and control vectors:
\begin{subequations}
  \label{eq:ex_ff_scvx_state_input}
  \begin{align}
    x &= \pare[big]{\rI,~\vI,~\qIB,~\wB}\in\reals^{13},\quad \\
    u &= \pare[big]{T_{\inertial},~M_{\body}}\in\reals^6.
  \end{align}
\end{subequations}

Next, we define the parameter vector to serve two purposes. First, as for the
quadrotor, we define a time dilation $\tdil$ such that
\eqref{eq:ex_quad_time_dilation} holds, yielding $\tabs=\tdil t$. Second, we take
advantage of the fact that \optieqref{scp_gen_cont}{convex_path_constraints_X} is
mixed in the state and parameter in order to place the slack room SDFs in
\eqref{eq:sdf_room_convex} into the parameter vector. In particular, we recognize
that according to \optieqref{subproblem_scvx_dt}{convex_path}, the constraint
\eqref{eq:sdf_room_convex} is imposed only at the discrete\dash time grid
nodes. Thus, for a grid of $N$ nodes, there are only $N\niss$ instances of
\eqref{eq:sdf_room_convex} to be included. We can therefore define the following
vector:
\begin{equation}
  \label{eq:ex_ff_vectorized_diss_lb}
  \disslbvec \definedas
  \pare[bigg]{
    \disslb^1, \dots, \disslb^N
  } \in \reals^{N\niss},
\end{equation}
where $\disslb^k\equiv\disslb\pare[big]{\rI(t_k)}$ and
$\disslbvec[i+(k-1)\niss]$ denotes the slack SDF value for the $i$-th room at
time $t_k$. To keep the notation concise, we will use the shorthand
$\disslbvec[ik]\equiv \disslbvec[i+(k-1)\niss]$. The overall parameter vector
is then given by:
\begin{equation}
  \label{eq:ex_ff_parameter}
  p = \Matrix{\tdil \\ \disslbvec}\in\reals^{1+N\niss}.
\end{equation}

In absolute time, the dynamics \optieqref{scp_gen_cont}{dynamics} are given
directly by \eqref{eq:ex_ff_dynamics}. As for the quadrotor, this forms a set of
time\dash invariant first\dash order ordinary differential equations:
\begin{equation}
  \label{eq:ex_ff_dynamics_abs_time}
  f\pare{x, u} = \Matrix{
    \vI \\
    m\inv T_{\inertial} \\
    \frac{1}{2} \qIB \otimes \wB \\
    J^{-1} \pare[big]{ M_{\body} - \wB\skew J \wB }
  }.
\end{equation}

The boundary conditions \optieqref{scp_gen_cont}{initial_conditions} and
\optieqref{scp_gen_cont}{final_conditions} are obtained from \eqref{eq:ex_ff_bcs}:
\begin{subequations}
  \label{eq:ex_ff_gic_gtc}
  \begin{align}
    \gic\pare[big]{x(0),p} &= \Matrix{\rI(0)-\ric \\
    \vI(0)-\vic \\ \qIB(0) - \qic \\ \wB(0)}, \\
    \gtf\pare[big]{x(1),p} &= \Matrix{\rI(1)-\rfc \\
    \vI(1)-\vfc \\ \qIB(1) - \qfc \\ \wB(1)}.
  \end{align}
\end{subequations}

The dynamics are converted to normalized time in the same way as
\eqref{eq:ex_quad_dynamics_normalized_time}:
\begin{equation}
  \label{eq:ex_ff_dynamics_normalized_time}
  f\pare{x, u, p} = \tdil f\pare{x, u}.
\end{equation}

The convex state and input path constraints \CTNLconvexpath are
straightforward. As for the quadrotor example, we leverage the mixed state\dash
parameter nature of \optieqref{scp_gen_cont}{convex_path_constraints_X} to include
all of the convex state and parameter constraints. In particular, these are
\eqref{eq:ex_quad_tf_bounds}, \eqref{eq:ex_ff_cvx_sc}, and
\eqref{eq:sdf_room_convex}. Using the definition of time dilation, we translate
\eqref{eq:ex_quad_tf_bounds} into the constraint:
\begin{equation}
  \label{eq:ex_ff_tdil_constraint}
  \tf[,\min] \le \tdil \le \tf[,\max].
\end{equation}

Using the definition of the concatenated slack SDF vector
\eqref{eq:ex_ff_vectorized_diss_lb}, we translate \eqref{eq:sdf_room_convex} into the
following constraints:
\begin{align}
  \label{eq:ex_ff_slack_sdf_constraints}
  \disslbvec[ik] \le \diss[i]\pare[big]{\rI(t_k)},\quad
  &i=1,\dots,\niss, \\
  \notag
  &k=1,\dots,N.
\end{align}

Consequently, the convex path constraint set $\set X$ in
\optieqref{scp_gen_cont}{convex_path_constraints_X} is given by:
\begin{align}
  \notag
  \set{X}
  = \big\{(x,p)\in\reals^{13}\times\reals^{1+N\niss}\where~%
  &\textnormal{\eqref{eq:ex_ff_cvx_sc}, \eqref{eq:ex_ff_tdil_constraint},} \\
  \label{eq:ex_ff_state_path}
  \textnormal{and~}
  &\textnormal{\eqref{eq:ex_ff_slack_sdf_constraints} hold}\big\}.
\end{align}

The convex input constraint set $\set{U}$ in
\optieqref{scp_gen_cont}{convex_path_constraints_U} is given simply by all the
input vectors that satisfy \eqref{eq:ex_ff_cc}.

The nonconvex path constraints \optieqref{scp_gen_cont}{nonconvex_constraints} for
the free-flyer problem involve the ellipsoidal floating obstacles
\eqref{eq:ex_quad_sc_obs} and the approximate flight space constraint
\eqref{eq:ex_ff_iss}. The floating obstacle constraints are modeled exactly like
for the quadrotor using \eqref{eq:ex_quad_s_obs_avoid}. For the flight space
constraint, we leverage the concatenated slack SDF vector
\eqref{eq:ex_ff_vectorized_diss_lb} and the eventual temporal discretization of the
problem in order to impose \eqref{eq:ex_ff_iss} at each temporal grid node as
follows:
\begin{equation}
  \label{eq:ex_ff_scvx_iss}
  \softmax[\sigma]\pare[big]{\disslb^k}\ge 0,~
  \quad k=1,\dots,N.
\end{equation}

Hence, the nonconvex path constraint function in
\optieqref{scp_gen_cont}{nonconvex_constraints} can be written as
$s:\reals^3\times\reals^{N\niss}\to\reals^{\Nobs+1}$. The first $\Nobs$
components are given by \eqref{eq:ex_quad_s_obs_avoid} and the last component is
given by the negative left\dash hand side of \eqref{eq:ex_ff_scvx_iss}.

It remains to define the running cost of the Bolza cost function
\eqref{eq:ocost_nlin}. Like for the quadrotor, we scale the integrand in
\optiobjref{ex_ff_ocp} to be mindful of the penalty terms which the \scvx algorithm
will add. Furthermore, we simplify the cost by neglecting time dilation and
directly associating the absolute\dash time integral of \optiobjref{ex_ff_ocp} with
the normalized\dash time integral of \eqref{eq:ocost_nlin}. This yields the
following convex running cost definition:
\begin{equation}
  \label{eq:ex_ff_running_cost}
  \runn(x, u, p) =
  \pare[Big]{\frac{\norm[2]{T_{\inertial}}}{T_{\max}}}^2+
  \pare[Big]{\frac{\norm[2]{M_{\body}}}{M_{\max}}}^2.
\end{equation}

At this point, it may seem as though we are finished with formulating
\pref{scp_gen_cont} for \scvx. However, the seemingly innocuous flight space
constraint \eqref{eq:ex_ff_scvx_iss} actually hides an important difficulty that we
will now address. The importance of the following discussion cannot be
overstated, as it can mean the difference between successful trajectory
generation, and convergence to an infeasible trajectory (i.e., one with
non\dash zero virtual control). In the case of the 6-DoF free\dash flyer,
omission of the following discussion incurs a $40\si{\percent}$ optimality
penalty for the value of \optiobjref{ex_ff_ocp}.

We begin investigating \eqref{eq:ex_ff_scvx_iss} by writing down its Jacobians,
which \scvx will use for linearizing the constraint. Note that
\eqref{eq:ex_ff_scvx_iss} is a function of only $\disslb^k\in\reals^{\niss}$, which
is part of the concatenated slack SDF \eqref{eq:ex_ff_vectorized_diss_lb}, and thus
resides in the parameter vector. Hence, only the Jacobian
\eqref{eq:scvx_lin_mats_j} is non\dash zero. Using the general softmax definition
\eqref{eq:softmax}, the $i$-th element of
$\grad\softmax[\sigma]\pare[big]{\disslb^k}$ is given by:
\begin{equation}
  \label{eq:ex_ff_softmax_iss}
  \frac{\partial\softmax[\sigma]\pare[big]{\disslb^k}}{\partial\disslb[i]^k}
  = \pare[Big]{\sum_{j=1}^{\niss}
    \exp{\sigma\disslb[j]^k}}\inv
  \exp{\sigma\disslb[i]^k}.
\end{equation}

\begin{csmfigure}[%
  caption={%
    Visualization of the effect of slackness in the SDF lower\dash bound
    constraint \eqref{eq:sdf_room_convex} on the gradient of the approximate SDF
    \eqref{eq:sdf_iss_smooth}. This plot is obtained by setting $\niss=6$ and
    $\disslb[j]=-1$ for all $j\ne i$. The individual curves are obtained by
    gradually reducing $\disslb[i]$ from its maximum value of $\diss[i]$. As
    the sharpness parameter $\sigma$ increases, a ``cutoff'' value appears,
    below which the approximate SDF becomes insensitive to changes in
    $\disslb[i]$.
  },%
  label={softmax_gradient_problem},%
  position={!t}]%
  \centering%
  \includegraphics[width=\columnwidth]{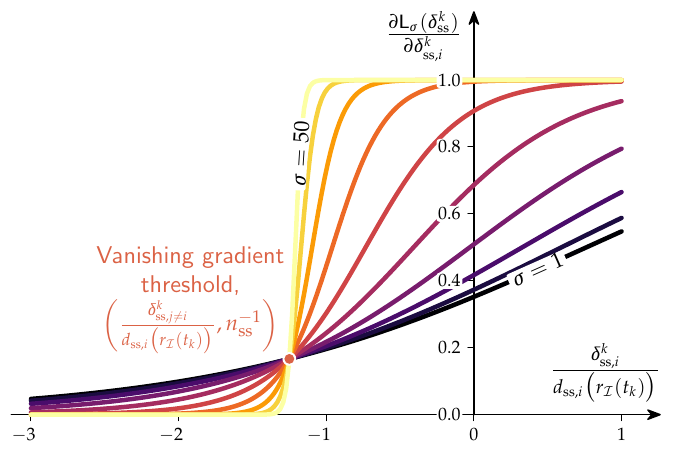}
\end{csmfigure}

When the slack SDF satisfies the lower\dash bound
\eqref{eq:ex_ff_slack_sdf_constraints} with equality, the Jacobian
\eqref{eq:ex_ff_softmax_iss} is an accurate representation of how the overall SDF
\eqref{eq:sdf_iss_smooth} changes due to small perturbations in the robot's
position. The problematic case occurs when this bound is loose. To illustrate
the idea, suppose that the robot is located near the center of $\Obsi$, such
that $\diss[i]\pare[big]{\rI(t_k)}=0.8$. We assume that the rooms do not
overlap and that the slack SDF values of the other rooms satisfy
$\disslb[j]^k=\diss[j]\pare[big]{\rI(t_k)}=-1$ for all $j\ne i$. Since the
robot is uniquely inside $\Obsi$, the exact SDF \eqref{eq:sdf_iss} is locally a
linear function of $\diss[i]$ and has a gradient
$\partial\diss/\partial\diss[i]=1$. Since we want the SCP subproblem to be an
accurate local approximation of the nonconvex problem, we expect the same
behavior for the approximate SDF \eqref{eq:sdf_iss_smooth} for high values of
$\sigma$. However, this may not be the case.

\figref{softmax_gradient_problem} illustrates what happens to the approximate SDF
gradient \eqref{eq:ex_ff_softmax_iss} as the slackness in
\eqref{eq:ex_ff_slack_sdf_constraints} increases. First, we note that when there is
no slackness, increasing the $\sigma$ parameter does indeed make the
approximate gradient converge to the exact value of one. However, as slackness
grows, there is a distinct cutoff value below which \eqref{eq:ex_ff_softmax_iss}
becomes zero. This is known as a \textbf{vanishing gradient} problem, and has
been studied extensively for machine learning \cite{GoodfellowDeepLearning}. The
core issue is that SCP relies heavily on gradient information in order to
determine how to improve the feasibility and optimality of the subproblem
solution. As an analogy, the gradient acts like a torchlight that illuminates
the local surroundings in a dark room and allows one to take a step closer to a
light switch. When the gradient vanishes, so does the torchlight, and SCP no
longer has information about which direction is best to take. Unless the
solution is already locally optimal, a vanished gradient most often either
blocks SCP from finding more optimal solutions, or forces it to use non\dash
zero virtual control. The result is that the converged trajectory is either
(heavily) suboptimal or even infeasible.

\begin{blurb}
  Favorable gradient behavior is instrumental for good performance.
\end{blurb}

Looking at \figref{softmax_gradient_problem}, one may ask, in order to recover
gradient information, why does SCP not simply increase $\disslb[i]^k$ above the
vanishing threshold? But remember, it is the gradient that indicates that
increasing $\disslb[i]^k$ is a good approach in the first place. The situation
is much like focusing your eyes on the flat region of the $\sigma=50$ curve on
the very left in \figref{softmax_gradient_problem}. If you only saw the part of the
curve for $\disslb[i]^k/\diss[i]\pare[big]{\rI(t_k)}\in [-3, -2]$, you would
also not know whether it is best to increase or decrease $\disslb[i]^k$.

Fortunately, the remedy is quite simple. Because \eqref{eq:sdf_property} is a
necessary and sufficient condition, we know that slackness in
\eqref{eq:ex_ff_slack_sdf_constraints} cannot be used to make the trajectory more
optimal. In other words, a trajectory with non\dash zero slackness will not
achieve a lower cost \optiobjref{ex_ff_ocp}. Hence, we simply need some way to
incentivize the convex subproblem optimizer to make
\eqref{eq:ex_ff_slack_sdf_constraints} hold with equality. Our approach is to
introduce a terminal cost that maximizes the concatenated slack SDF:
\begin{equation}
  \label{eq:ex_ff_scvx_terminal_cost}
  \term\pare{\disslbvec} = -\sdissweight
  \sum_{k=1}^{N}\sum_{i=1}^{\niss}\disslbvec[ik],
\end{equation}
where $\sdissweight\in\pos$ is any user\dash chosen positive number. To make
sure that \eqref{eq:ex_ff_scvx_terminal_cost} does not interfere with the extra
penalty terms introduced by \scvx, we set $\sdissweight$ to a very small but
numerically tolerable value as shown in \tabref{ex_ff_params}.

In summary, our investigation into \eqref{eq:ex_ff_scvx_iss} allowed us to identify
a vanishing gradient issue. This resulted in a simple yet effective remedy in
the form of a terminal cost \eqref{eq:ex_ff_scvx_terminal_cost}. The discussion
hopefully highlights three salient features of good modeling for SCP\dash based
trajectory generation. First, SCP does not have equal performance for
mathematically equivalent problem formulations (such as the free\dash flyer
problem with and without \eqref{eq:ex_ff_scvx_terminal_cost}). Second, favorable
gradient behavior is instrumental for good performance. Third, remedies to
recover good performance for difficult problems are often surprisingly simple.



\subsubsection{\gusto Formulation}

Like for the quadrotor example, the \gusto formulation is very similar to
\scvx. We can express \eqref{eq:ex_ff_running_cost} as the quadratic running cost
\eqref{eq:gusto_running_cost} as follows:
\begin{subequations}
  \label{eq:ex_ff_gusto_running_cost}
  \begin{align}
    \Jq(p) &= \diag\pare[bigg]{T_{\max}\inv[2] I_3, M_{\max}\inv[2] I_3}, \\
    \Jl (x,p) &= 0, \\
    \Jc(x,p) &= 0.
  \end{align}
\end{subequations}

The dynamics \eqref{eq:ex_ff_dynamics_abs_time} are also cast into the control
affine form \eqref{eq:control_affine_gusto}:
\begin{subequations}
  \begin{align}
    \label{eq:ex_ff_gusto_f0}
    f_0\pare{x, p}
    &= \Matrix{
    \vI \\
    0 \\
    \frac{1}{2} \qIB \otimes \wB \\
    -J^{-1} \pare[big]{ \wB\skew J \wB }
    }, \\
    f_i\pare{x, p}
    &= \Matrix{
    0 \\
    m\inv e_i \\
    0 \\
    0
    },\quad i=1,2,3, \\
    f_i\pare{x, p}
    &= \Matrix{
    0 \\
    0 \\
    0 \\
    J^{-1} e_i
    },\quad i=4,5,6,
  \end{align}
\end{subequations}
where $e_i\in\reals^3$ is the $i$-th standard basis vector. Just like for the
quadrotor, we are ``done'' at this point and the rest of the optimization model
is exactly the same as for \scvx in the last section.

\subsubsection{Initial Trajectory Guess}

The initial trajectory guess is based on some simple intuition about what a
feasible free-flyer trajectory might look like. Although this guess is more
complicated than the straight\dash line initialization used for the quadrotor,
it is based on purely kinematic considerations. This makes the guess quick to
compute but also means that it does not satisfy the dynamics and obstacle
constraints. The fact that SCP readily morphs this coarse guess into a feasible
and locally optimal trajectory corroborates the effectiveness of SCP methods
for high\dash dimensional nonconvex trajectory generation tasks.

To begin, the time dilation $\tdil$ is obtained by averaging the final time
bounds as in \eqref{eq:ex_quad_initial_parameter}. An ``L\dash shape'' path is then
used for the position trajectory guess. In particular, recall that according to
\eqref{eq:ex_ff_bcs_pos}, the free\dash flyer has to go from $\ric$ to $\rfc$. We
define a constant velocity trajectory which closes the gap between $\ric$ and
$\rfc$ along the first axis, then the second, and finally the third, in a total
time of $\tdil$ seconds. The trajectory thus consists of three straight legs
with sharp 90 degree turns at the transition points, which is akin to the
Manhattan or taxicab geometry of the one\dash norm
\cite{LaVallePlanningBook}. The velocity is readily derived from the position
trajectory, and is a constant\dash norm vector whose direction changes twice to
align with the appropriate axis in each trajectory leg. Furthermore, we
initialize the concatenated slack SDF parameter vector
\eqref{eq:ex_ff_vectorized_diss_lb} by evaluating \eqref{eq:sdf_i} along the position
trajectory guess for each room and discrete\dash time grid node.

The attitude trajectory guess is only slightly more involved, and it is a
general procedure that we can recommend for attitude trajectories. According to
\eqref{eq:ex_ff_bcs_quat}, the free\dash flyer has to rotate between the attitudes
encoded by $\qic$ and $\qfc$. Since quaternions are not additive and must
maintain a unit norm to represent rotation, straight\dash line interpolation
from $\qic$ to $\qfc$ is not an option. Instead, we use spherical linear
interpolation (SLERP) \cite{Shoemake1985,Sola2017}. This operation performs a
continuous rotation from $\qic$ to $\qfc$ at a constant angular velocity around
a fixed axis. To define the operation, we introduce the exponential and
logarithmic maps for unit quaternions:
\begin{subequations}
  \label{eq:q_exp_log_maps}
  \begin{align}
    \label{eq:q_exp_map}
    \Expquat(q)
    &\definedas \alpha u,~\textnormal{where}~\alpha\in\reals,~u\in\reals^3, \\
    \label{eq:q_log_map}
    \Logquat\pare[big]{\alpha u}
    &\definedas \Matrix{u \sin\pare{\alpha/2} \\ \cos\pare{\alpha/2}}.
  \end{align}
\end{subequations}

The exponential map converts a unit quaternion to its equivalent angle\dash
axis representation. The logarithmic map converts an angle\dash axis rotation
back to a quaternion, which we write here in the vectorized form used to
implement $\qIB$ in \eqref{eq:ex_ff_dynamics_quaternion}. SLERP for the attitude
quaternion $\qIB$ can then be defined by leveraging \eqref{eq:q_exp_log_maps}:
\begin{subequations}
  \label{eq:slerp}
  \begin{align}
    \label{eq:slerp_error_quaternion}
    q_e &= \qic\qconj\otimes \qfc, \\
    \qIB(t) &= q_0\otimes \Expquat\pare[big]{t\Logquat\pare{q_e}},
  \end{align}
\end{subequations}
where $q_e$ is the error quaternion between $\qfc$ and $\qic$, and $t\in [0,1]$
is an interpolation parameter such that $\qIB(0)=\qic$ and $\qIB(1)=\qfc$. The
angular velocity trajectory guess is simple to derive, since SLERP performs a
constant velocity rotation around a fixed axis. Using
\eqref{eq:slerp_error_quaternion}:
\begin{equation}
  \label{eq:ex_ff_ang_vel_guess}
  \wB(t) = \tdil\inv \Logquat\pare{q_e}.
\end{equation}

The free\dash flyer is a very low thrust vehicle to begin with. By using
\optiobjref{ex_ff_ocp}, we are in some sense searching for the lowest of low thrust
trajectories. Hence, we expect the control inputs $T_{\inertial}$ and
$M_{\body}$ to be small. Without any further insight, it is hard to guess what
the thrust and torque would look like for a 6-DoF vehicle in a micro\dash
gravity environment. Hence, we simply set the initial control guess to zero.

\subsubsection{Numerical Results}

\begin{csmtable}[
  caption={Breakdown of subproblem size for the 6-DoF free-flyer
    example.},
  label={ex_ff_subproblem_size},
  position={!b}]%
  \centering%
  \ifmaketwocolcsm
  \includegraphics[width=0.7\columnwidth]{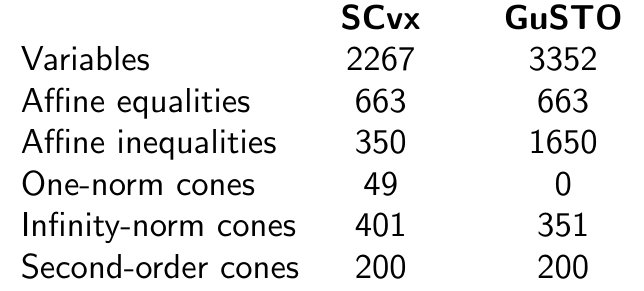}
  \else
  \includegraphics[scale=\csmpreprintfigscale]{scp_ff_sizes}
  \fi
\end{csmtable}

We now have a specialized instance of \pref{scp_gen_cont} and an initialization
strategy for the 6-DoF free\dash flyer problem. The trajectory solution is
generated using \scvx and \gusto with temporal discretization performed using
the FOH interpolating polynomial method in \sbref{discretization}. The algorithm
parameters are provided in \tabref{ex_ff_params}, where the initial and final
quaternion vectors are expressed in degrees using the angle\dash axis
representation of \eqref{eq:q_exp_map}. ECOS is used as the numerical convex
optimizer, and the full implementation is available in the code repository
linked in \figref{github_qr}.

\begin{csmfigure}[%
  caption={%
    Convergence and runtime performance for the 6-DoF free-flyer problem. Both
    algorithms take a similar amount of time to converge. The runtime subplots
    in the bottom row show statistics on algorithm performance over 50
    executions. \gusto converges slightly faster for this example and, although
    both algorithms reach numerical precision for practical purposes, \gusto
    converges all the way down to a $10^{-14}$ tolerance. The plots are
    generated according to the same description as for
    \figref{ex_quad_convergence}.
  },
  label={ex_ff_convergence},
  columns=2,
  position={!t}]%
  \centering
  \ifmaketwocolcsm
  \def\subfigwidth{0.45\textwidth}
  \else
  \def\subfigwidth{0.48\textwidth}
  \fi
  \def\figwidth{\columnwidth}
  \begin{subfigure}[t]{\subfigwidth}
    \centering
    \includegraphics[width=\figwidth]{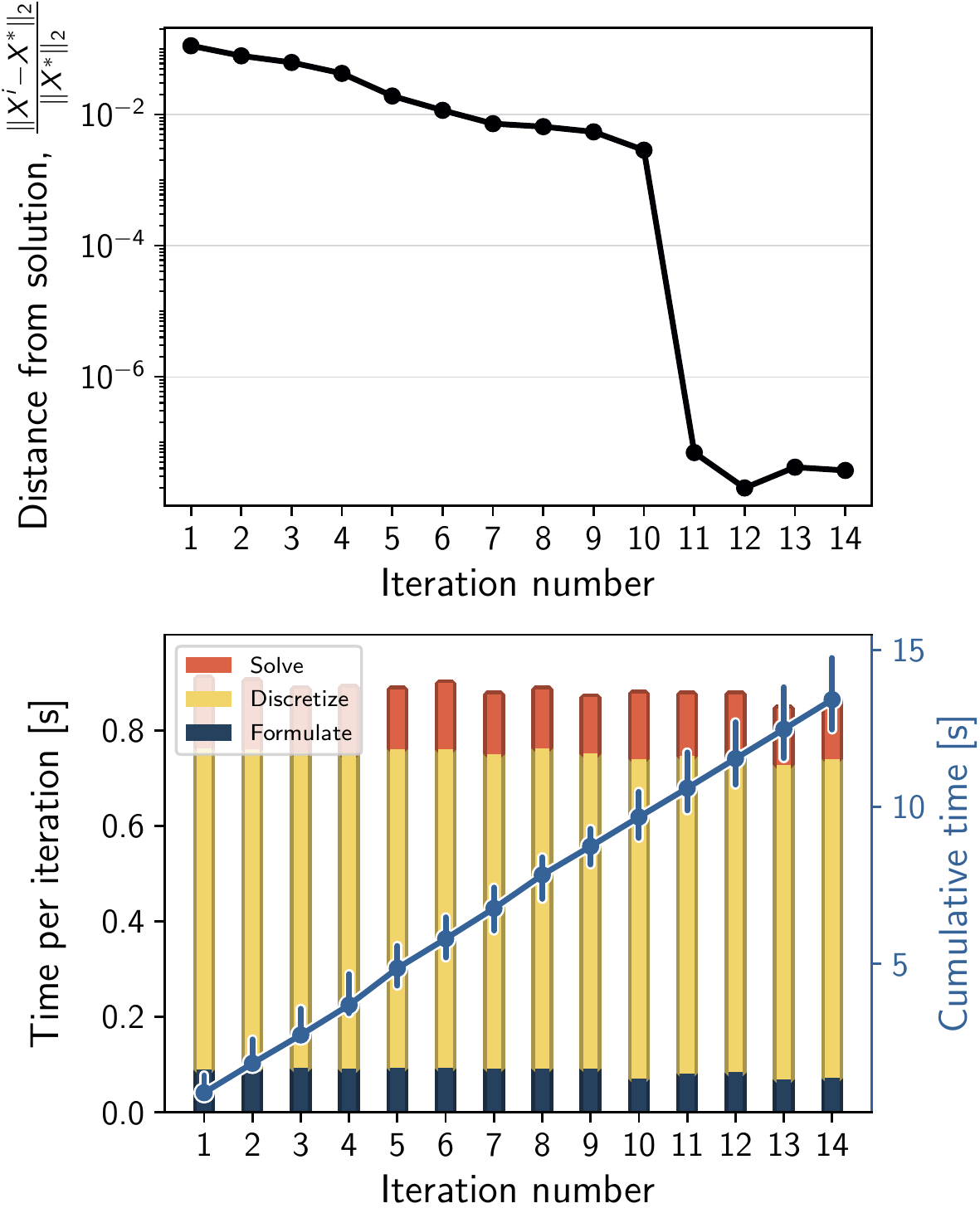}%
    \caption{\scvx.}
    \label{fig:ex_ff_convergence_scvx}
  \end{subfigure}%
  \ifmaketwocolcsm
  \hspace{1cm}%
  \else
  \hfill%
  \fi
  \begin{subfigure}[t]{\subfigwidth}
    \centering
    \includegraphics[width=\figwidth]{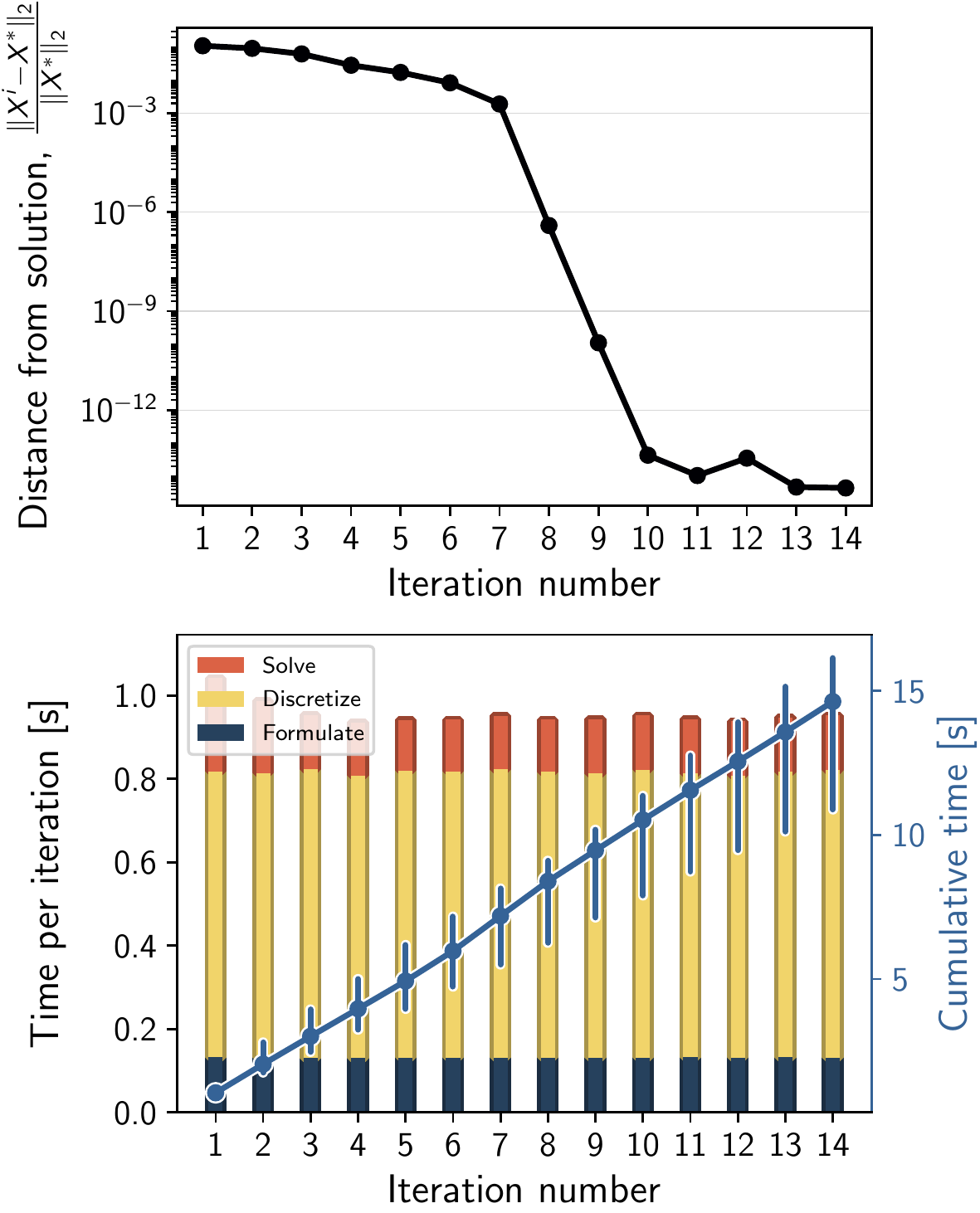}%
    \caption{\gusto.}
    \label{fig:ex_ff_convergence_gusto}
  \end{subfigure}%
\end{csmfigure}

The convergence processes for \scvx and \gusto are shown in
\figref{ex_ff_convergence}. We have again set
$\varepsilon=\varepsilon_{\mrm{r}}=0$ so that we can observe the convergence
process for exactly 15 iterations. At each iteration, the algorithms solve a
convex subproblem whose size is documented in \tabref{ex_ff_subproblem_size}. Note
that the subproblems of both algorithms are substantially larger than for the
quadrotor example, and represent a formidable increase in dimensionality for
the numerical problem. However, modern IPMs easily handle problems of this size
and we will see that the increased variable and constraint count is of little
concern. We further note that the larger number of variables and affine
inequalities for \gusto is due to how our implementation uses extra slack
variables to encode the soft penalty function \eqref{eq:gusto_hpen}. Because \gusto
does not use a dynamics virtual control, it has no one\dash norm cones, while
\scvx has several such cones to model the virtual control penalty
\eqref{eq:P_penalty_def}. Due to its larger subproblem size and slightly more
complicated code for including constraints as soft penalties, the ``solve'' and
``formulate'' times per subproblem are slightly larger for \gusto in this
example. Nevertheless, both algorithms have roughly equal runtimes, and \gusto
has the advantage of converging to a given tolerance in slightly fewer
iterations.


\begin{csmfigure}[%
  caption={%
    The position trajectory evolution (left) and the final converged trajectory
    (right) for the 6-DoF free-flyer problem. In the right plots for each
    algorithm, the continuous\dash time trajectory is obtained by numerically
    integrating the dynamics \eqref{eq:ex_ff_dynamics}. The fact that this
    trajectory passes through the discrete\dash time subproblem solution
    confirms dynamic feasibility.
  },
  label={ex_ff_pos},
  columns=2]%
  \ifmaketwocolcsm
  \def\subfigwidth{0.5\textwidth}
  \def\subfigfigwidth{0.98\textwidth}
  \else
  \def\subfigwidth{\textwidth}
  \def\subfigfigwidth{0.7\textwidth}
  \fi
  \begin{subfigure}[t]{\subfigwidth}
    \centering
    \includegraphics[width=\subfigfigwidth,page=1]{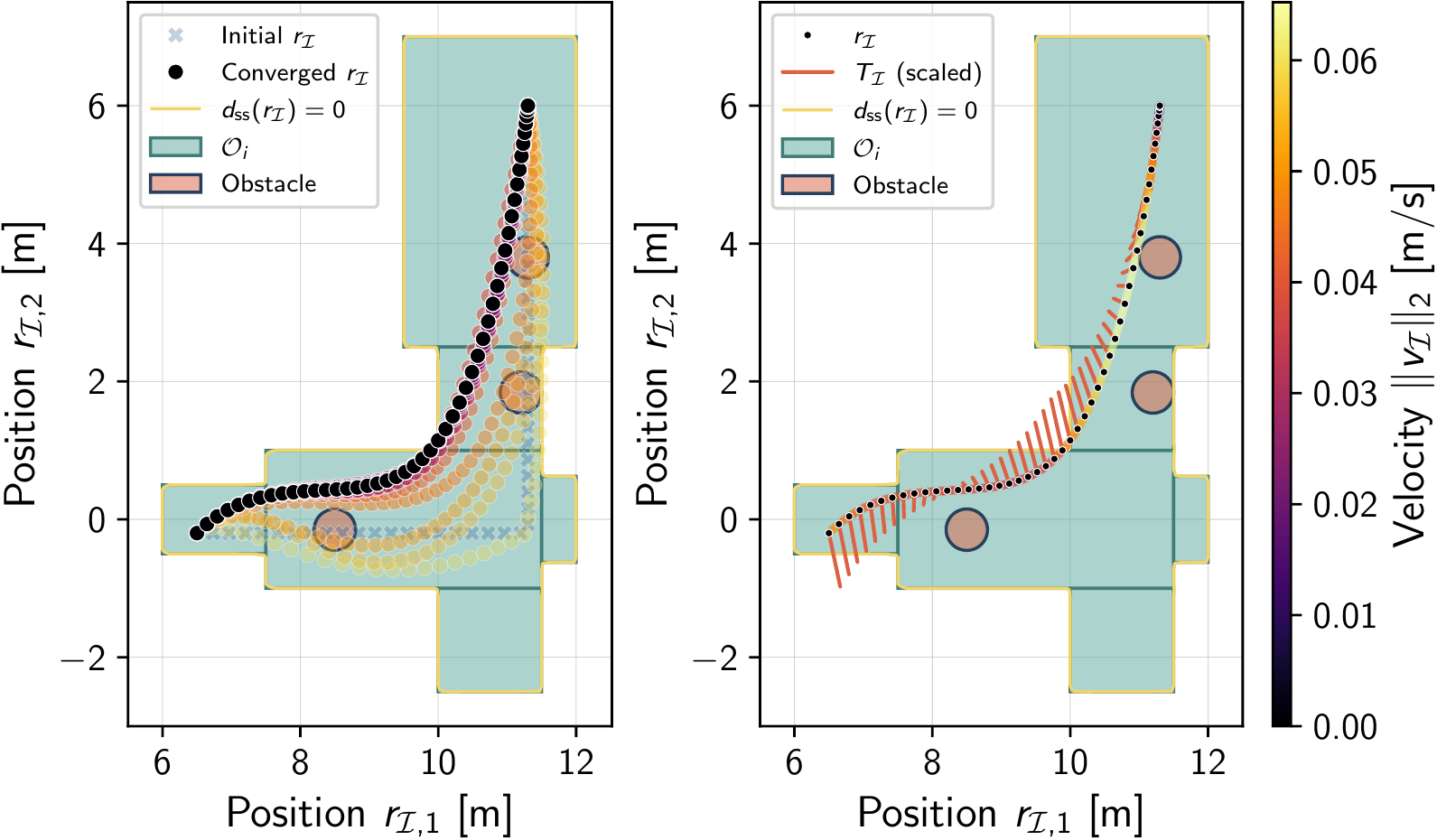}%
    \caption{\scvx.}
    \label{fig:ex_ff_pos_scvx}
  \end{subfigure}%
  \ifmaketwocolcsm
  \hfill%
  \else
  \vspace{3mm}

  \fi
  \begin{subfigure}[t]{\subfigwidth}
    \centering
    \includegraphics[width=\subfigfigwidth,page=2]{scp_ff_pos}%
    \caption{\gusto.}
    \label{fig:ex_ff_pos_gusto}
  \end{subfigure}%
\end{csmfigure}

The converged trajectories are plotted in \figref{ex_ff_pos,ex_ff_timeseries}. The
left subplots in \figref{ex_ff_pos_scvx,ex_ff_pos_gusto} show a remarkably similar
evolution of the initial guess into the converged trajectory. The final
trajectories are visually identical, and both algorithms discover that the
maximum allowed flight time of $\tf[,\max]$ is control energy\dash optimal, as
expected.

Lastly, \figref{ex_ff_obstacles} plots the evolution of the nonconvex flight space
and obstacle avoidance inequalities \eqref{eq:ex_ff_iss} and
\eqref{eq:ex_quad_sc_obs}. Our first observation is that the constraints hold at
the discrete\dash time nodes, and that the free\dash flyer approaches the
ellipsoidal obstacles quite closely. This is similar to how the quadrotor
brushes against the obstacles in \figref{ex_quad_pos}, and is a common feature of
time- or energy\dash optimal trajectories in a cluttered environment. Our
second observation concerns the sawtooth\dash like nonsmooth nature of the SDF
time history. Around $\tabs=25~\si{\second}$ and $\tabs=125~\si{\second}$, the
approximate SDF comes close to zero even though the position trajectory in
\figref{ex_ff_pos} is not near a wall at those times. This is a consequence of our
modeling, since the SDF is near\dash zero at the room interfaces (see
\figref{sdf_illustration,softmax_1d_intuition}), even though these are not physical
``walls''. However, around $\tabs=100~\si{\second}$, the flight space
constraint \eqref{eq:ex_ff_iss} is actually activated as the free\dash flyer rounds
a corner. Roughly speaking, this is intuitively the optimal thing to do. Like a
Formula One driver rounding a corner by following the racing line, the
free\dash flyer spends less control effort by following the shortest path. An
unfortunate consequence is that this results in minor clipping of the
continuous\dash time flight space constraint. The issue can be mitigated by the
same strategies as proposed in the last section for the quadrotor example
\cite{acikmese2008enhancements,Dueri2017Clipping}.

\begin{csmfigure}[%
  caption={%
    State and control time histories for the converged trajectory of the 6-DoF
    free-flyer problem. These are visually identical for \scvx and \gusto, so
    we only show a single plot. Euler angles using the intrinsic Tait-Bryan
    convention are shown in place of the quaternion attitude. Like in
    \figref{ex_quad_timeseries}, the dots represent the discrete\dash time solution
    while the continuous lines are obtained by propagating the solution through
    the actual continuous\dash time dynamics \optieqref{scp_gen_cont}{dynamics}.
  },
  label={ex_ff_timeseries},
  position={!t}]%
  \centering
  \includegraphics[width=\columnwidth]{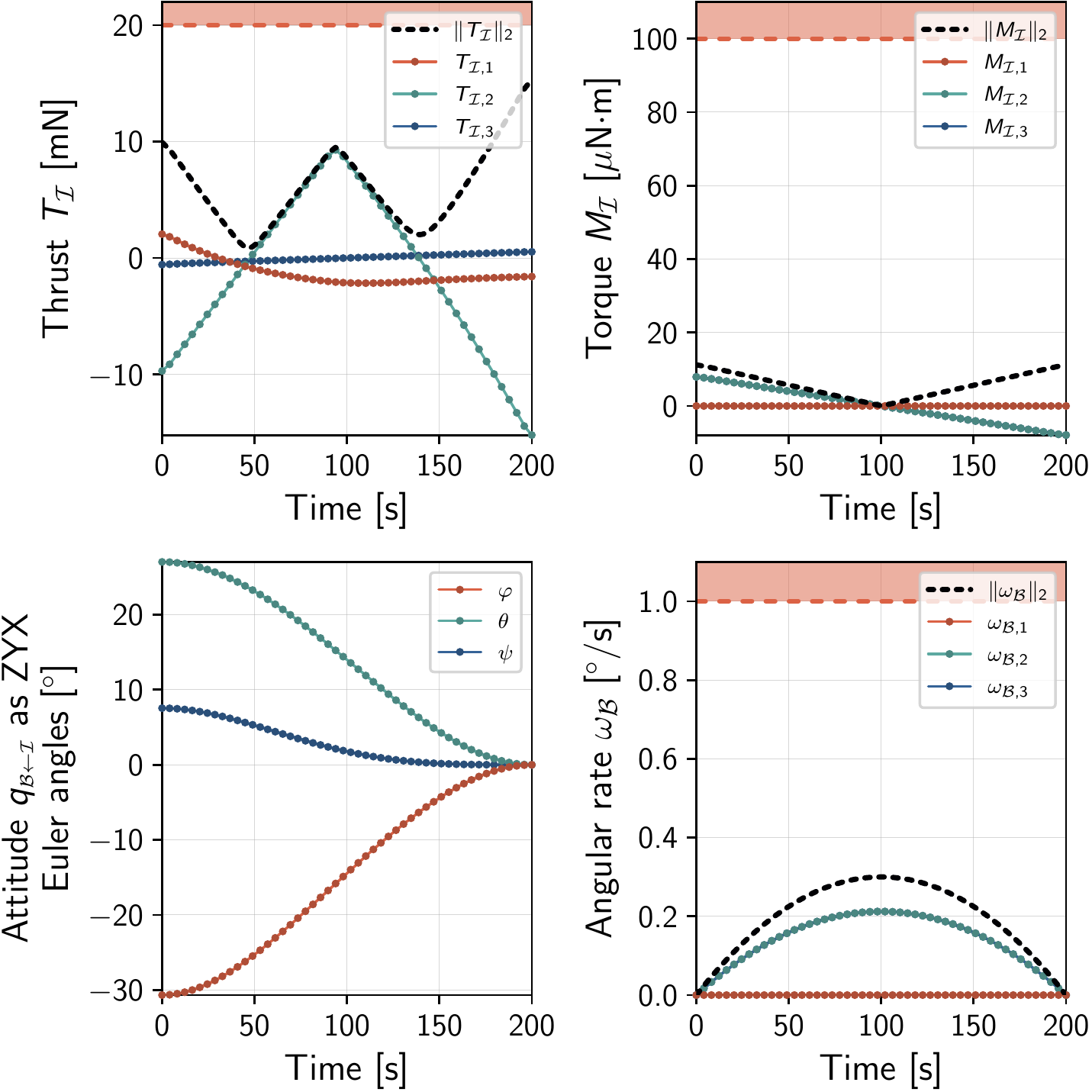}%
\end{csmfigure}

\begin{csmfigure}[%
  caption={%
    Signed distance function and obstacle avoidance time histories for the
    converged trajectory of the 6-DoF free-flyer problem. Like for
    \figref{ex_ff_timeseries}, these are visually identical for \scvx and \gusto,
    so we only show a single plot. Note the highly nonlinear nature of the SDF,
    whose time history exhibits sharp corners as the robot traverses the
    feasible flight space. Although the SDF constraint \eqref{eq:ex_ff_scvx_iss} is
    satisfied at the discrete\dash time nodes, minor inter-sample constraint
    clipping occurs around 100 seconds as the robot rounds a turn in the middle
    of its trajectory (see \figref{ex_ff_pos}).
  },
  label={ex_ff_obstacles},
  position={!t}]%
  \centering
  \includegraphics[width=\columnwidth]{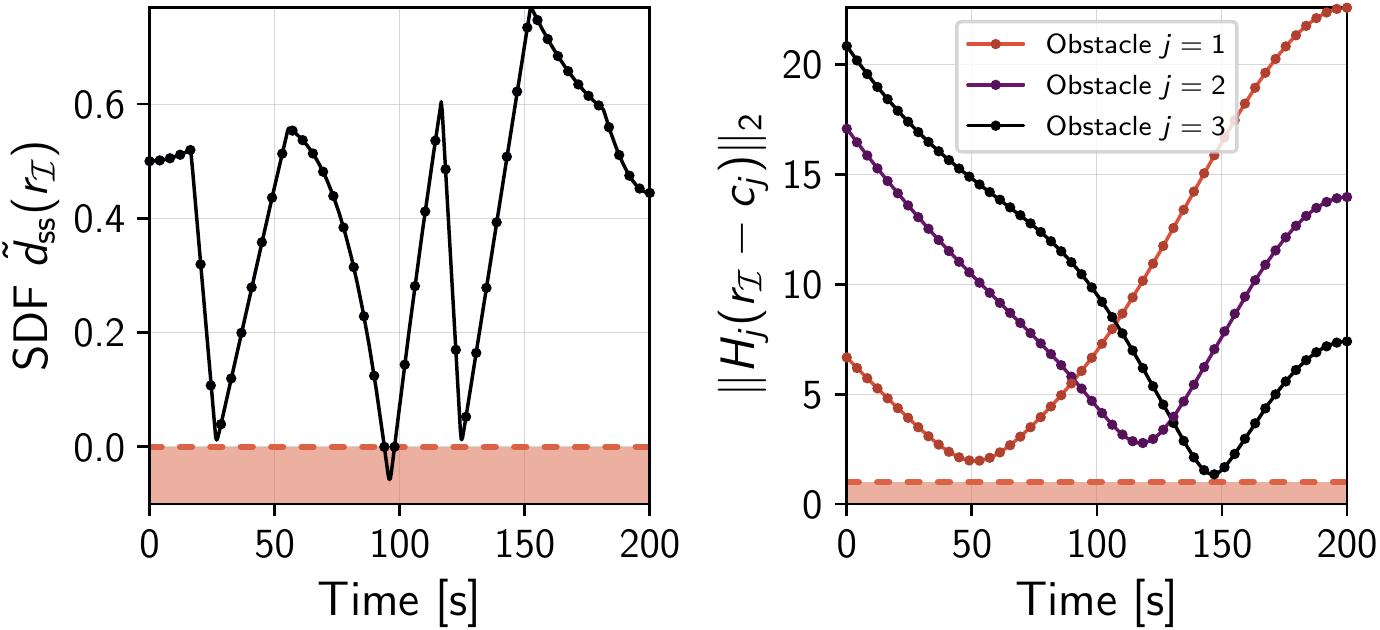}%
\end{csmfigure}


\section{Conclusions}

Modern vehicle engineering is moving in the direction of increased
autonomy. This includes aerospace, automotive, and marine transport, as well as
robots on land, in the air, and in our homes. No matter the application, a
common feature across autonomous systems is the basic requirement to generate
trajectories. In a general sense, trajectories serve like plans to be executed
in order for the system to complete its task. Due to the large scale of
deployment and/or the safety\dash critical nature of the system, reliable
real\dash time onboard trajectory generation has never been more important.

This article takes the stance that convex optimization is a prime contender for
the job, thanks to 40 years of optimization research having produced a
remarkable suite of numerical methods for quickly and reliably solving convex
problems \cite{NocedalBook,BoydConvexBook,wright2005interior}. Many of these
methods are now packaged as either commercial or open\dash source off\dash
the\dash shelf codes \cite{domahidi2013ecos,Zanelli2017}. This makes the
injection of convex optimization into an autonomous system easier than ever
before, provided that the right high\dash level algorithms exist to leverage
it.

To leverage convex optimization for the difficult task of nonconvex trajectory
generation, this article provides an expansive tutorial of three
algorithms. First, the lossless convexification (\lcvx) algorithm is introduced
to remove acute nonconvexities in the control input constraints. This provides
an optimal control theory\dash backed way to transform certain families of
nonconvex trajectory generation tasks into ones that can be solved in one shot
by a convex optimizer. A variable\dash mass rocket landing example at the end
of the article illustrates a real\dash world application of the \lcvx method.

Not stopping there, the article then motivates an entire family of optimization
methods called sequential convex programming (SCP). These methods use a
linearize\dash solve loop whereby a convex optimizer is called several times
until a locally optimal trajectory is obtained. SCP strikes a compelling
middleground between ``what is possible'' and ``what is acceptable'' for
safety\dash critical real\dash time trajectory generation. In particular, SCP
inherits much from trust region methods in numerical optimization, and its
performance is amenable to theoretical analysis using standard tools of
analysis, constrained optimization, and optimal control theory
\cite{NocedalBook,ConnTrustRegionBook,PontryaginBook}. This articles provides a
detailed overview of two specific and closely related SCP algorithms called
\scvx and \gusto \cite{Mao2018,Bonalli2019a}. To corroborate their effectiveness
for difficult trajectory generation tasks, two numerical examples are presented
based on a quadrotor and a space\dash station free\dash flyer maintenance
robot.

The theory behind \lcvx, \scvx, and \gusto is relatively new and under active
research, with the oldest method in this article (i.e., classical \lcvx
\cite{Acikmese2005}) being just 15 years old. We firmly believe that neither one
of the methods has attained the limits of its capabilities, and this presents
the reader with an exciting opportunity to contribute to the effort. It is
clear to us that convex optimization has a role to play in the present and
future of advanced trajectory generation. With the help of this article and the
associated source code, we hope that the reader now has the knowledge and tools
to join in the adventure.

\section{Acknowledgements}


This work was supported in part by the National Science Foundation,
Cyber-Physical Systems (CPS) program (award 1931815), by the Office of Naval
Research, ONR YIP Program (contract N00014-17-1-2433), and by the King
Abdulaziz City for Science and Technology (KACST). The authors would like to
extend their gratitude to Yuanqi Mao for his invaluable inputs on sequential
convex programming algorithms, to Abhinav Kamath for his meticulous review of
every detail, and to Jonathan P. How for the initial encouragement to write
this article.


\section{Author Information}

\begin{authorbio}[Danylo Malyuta]%
  received the B.Sc. degree in Mechanical Engineering from EPFL and the
  M.Sc. degree in Robotics, Systems and Control from ETH Z{\"u}rich. He is
  currently a Ph.D. candidate in the Autonomous Controls Lab at the Department
  of Aeronautics and Astronautics of the University of Washington. His research
  is primarily focused on computationally efficient optimization-based control
  of dynamical systems. Danylo has held internship positions at the NASA Jet
  Propulsion Laboratory, NASA Johnson Space Center, and Amazon Prime Air.
\end{authorbio}

\begin{authorbio}[Taylor P. Reynolds]%
  received the B.S. degree in Mathematics \& Engineering from Queen's
  University in 2016. He received the Ph.D. degree from the Department of
  Aeronautics \& Astronautics at the University of Washington in 2020 under the
  supervision of Mehran Mesbahi. During his Ph.D., Taylor worked with NASA JSC
  and Draper Laboratories to develop advanced guidance algorithms for planetary
  landing on the SPLICE project, and also co-founded the Aeronautics \&
  Astronautics CubeSat Team at the University of Washington. He now works as a
  Research Scientist at Amazon Prime Air.
\end{authorbio}

\begin{authorbio}[Michael Szmuk]%
  received the B.S. and M.S. degrees in Aerospace Engineering from the
  University of Texas at Austin. In 2019, he received the Ph.D. degree while
  working in the Autonomous Controls Lab at the Department of Aeronautics and
  Astronautics of the University of Washington, under the supervision of
  \Behcet{} \Acikmese{}. During his academic career, he completed internships
  at NASA, AFRL, Emergent Space, Blue Origin, and Amazon Prime Air. He now
  works as a Research Scientist at Amazon Prime Air, specializing in the design
  of flight control algorithms for autonomous air delivery vehicles.
\end{authorbio}

\begin{authorbio}[Thomas Lew]%
  is a Ph.D. candidate in Aeronautics and Astronautics at Stanford
  University. He received the B.Sc. degree in Microengineering from {\'E}cole
  Polytechnique F{\'e}d{\'e}rale de Lausanne in 2017 and the M.Sc. degree in
  Robotics from ETH Z{\"u}rich in 2019. His research focuses on the
  intersection between optimal control and machine learning techniques for
  robotics and aerospace applications.
\end{authorbio}

\begin{authorbio}[Riccardo Bonalli]%
  obtained the M.Sc. degree in Mathematical Engineering from Politecnico di
  Milano, Italy in 2014 and the Ph.D. degree in applied mathematics from
  Sorbonne Universit{\'e}, France in 2018 in collaboration with ONERA -- The
  French Aerospace Lab. He is recipient of the ONERA DTIS Best Ph.D. Student
  Award 2018. He is now a postdoctoral researcher at the Department of
  Aeronautics and Astronautics at Stanford University. His main research
  interests concern theoretical and numerical robust optimal control with
  applications in aerospace systems and robotics.
\end{authorbio}

\begin{authorbio}[Marco Pavone]%
  is an Associate Professor of Aeronautics and Astronautics at Stanford
  University, where he is the Director of the Autonomous Systems
  Laboratory. Before joining Stanford, he was a Research Technologist within
  the Robotics Section at the NASA Jet Propulsion Laboratory.  He received the
  Ph.D. degree in Aeronautics and Astronautics from the Massachusetts Institute
  of Technology in 2010.  His main research interests are in the development of
  methodologies for the analysis, design, and control of autonomous systems,
  with an emphasis on self-driving cars, autonomous aerospace vehicles, and
  future mobility systems.  He is a recipient of a number of awards, including
  a Presidential Early Career Award for Scientists and Engineers, an ONR YIP
  Award, an NSF CAREER Award, and a NASA Early Career Faculty Award.  He was
  identified by the American Society for Engineering Education (ASEE) as one of
  America’s 20 most highly promising investigators under the age of 40.  He is
  currently serving as an Associate Editor for the IEEE Control Systems
  Magazine.
\end{authorbio}

\begin{authorbio}[\Behcet{} \Acikmese{}]%
  is a Professor at the University of Washington, Seattle. He received the
  Ph.D. degree in Aerospace Engineering from Purdue University.  He was a
  senior technologist at the NASA Jet Propulsion Laboratory (JPL) and a
  lecturer at the California Institute of Technology. At JPL, he developed
  control algorithms for planetary landing, spacecraft formation flying, and
  asteroid and comet sample return missions. He developed the ``flyaway''
  control algorithms used successfully in NASA’s Mars Science Laboratory (MSL)
  and Mars 2020 missions during the landings of Curiosity and Perseverance
  rovers on Mars. He is a recipient of the NSF CAREER Award, the IEEE Award for
  Technical Excellence in Aerospace Control, and numerous NASA Achievement
  awards for his contributions to NASA missions and technology development. His
  research interests include optimization\dash based control, nonlinear and
  robust control, and stochastic control.
\end{authorbio}


\csmprintmainbib

\end{document}
